\documentclass[11pt,a4paper]{article}
\usepackage{amsfonts}

\usepackage{indentfirst}
\usepackage{cite}
\usepackage{amsmath,amsthm}
\allowdisplaybreaks[2]
\usepackage{amssymb}
\usepackage{mathrsfs}

\usepackage[hidelinks,bookmarksnumbered,hypertexnames=false]{hyperref}
\hypersetup{pdftitle={Two-dimensional full compressible flow with near-linear bulk viscosity and temperature-dependent heat conduction}}
\usepackage[
  hmargin=2.2cm,
  top=2.4cm,
  bottom=2.6cm,
  footskip=0.9cm
]{geometry}

\newtheorem{theorem}{Theorem}[section]
\newtheorem{remark}{Remark}[section]
\newtheorem{definition}{Definition}[section]
\newtheorem{lemma}[theorem]{Lemma}

\newtheorem{proposition}[theorem]{Proposition}

\newcommand{\n}{\rho}
\newcommand{\mr}{\mathbb{R}}
\newcommand{\lm}{\lambda}
\renewcommand{\div}{ {\rm div }  }
\newcommand{\na}{\nabla }
\newcommand{\pa}{\partial}
\newcommand{\ga}{\gamma}
\newcommand{\OM}{\Omega}
\newcommand{\curl}{{\rm curl} }
\newcommand{\la}{\label}
\newcommand{\ol}{\overline}
\newcommand{\T}{\mathbb{T}}
\newcommand{\ba}{\begin{aligned}}
\newcommand{\ea}{\end{aligned}}
\newcommand{\be}{\begin{equation}}
\newcommand{\ee}{\end{equation}}
\def\p{\partial}
\def\lam{\lambda}

\def\o{\omega}
\def\rr{\mathbb{R}^2}

\makeatletter
\@addtoreset{equation}{section}
\makeatother

\title{Global Existence and Large-Time Behavior of Strong Solutions to the Two-Dimensional Full Compressible Navier-Stokes Equations with Degenerate Heat Conductivity and Large Initial Data}

\author{Qinghao L{\small EI}$^{a,b}$, Xin S{\small{ONG}}$^{c}$ \thanks{Email addresses: leiqinghao22@mails.ucas.ac.cn (Q. H. Lei), songxin2357@email.ncu.edu.cn (X. Song)} \\
{\normalsize a. School of Mathematical Sciences,}\\
{\normalsize University of Chinese Academy of Sciences, Beijing 100049, P. R. China;} \\
{\normalsize b. Institute of Applied Mathematics,} \\
{\normalsize Academy of Mathematics and Systems Science,} \\
{\normalsize Chinese Academy of Sciences, Beijing 100190, P. R. China;} \\
{\normalsize c. School of Mathematics and Computer Sciences,}\\
{\normalsize  Nanchang University, Nanchang 330031, P. R. China.}
}

\date{}
\begin{document}
\maketitle
\begin{abstract}
In this paper, we investigate the two-dimensional full compressible Navier--Stokes system with degenerate heat conductivity in three settings: the periodic domain, a general bounded simply connected domain subject to Navier-slip boundary conditions, and the whole space with non-vacuum far-field density.
We assume that the shear viscosity $\mu$ is a positive constant, the bulk viscosity $\lambda(\rho)\ge0$ depends on the density, and the heat conductivity takes the form $\kappa(\theta)=\kappa_0 \theta^\alpha$ with $\kappa_0>0$.
For a broad class of density-dependent bulk viscosities, we establish the global existence and large-time behavior of strong solutions for arbitrarily large initial data allowing vacuum.
In particular, our results include the sublinear growth law
$$
\lambda(\rho)=\frac{\rho}{(1+\rho)^\varepsilon \log\log(e^e+\rho)},
$$
where $\varepsilon>0$ is a sufficiently small constant depending only on the domain, the initial data, and the physical parameters.
\end{abstract}


\section{Introduction and main results}
We consider a compressible, viscous, heat-conducting polytropic fluid in a domain $\Omega \subset \mathbb{R}^2$.
Its motion is governed by the full compressible Navier--Stokes system:
\be\ba\la{ffns}
\begin{cases}
\rho_t + \div(\rho u) = 0, \\
(\n u)_t + \div(\n u\otimes u) + \na P = \div \mathbb{S}, \\
(\n E)_t + \div( \n E u + P u) - \div(\kappa \na \theta) = \div(\mathbb{S}u),
\end{cases}
\ea\ee
where the viscous stress tensor $\mathbb{S}$ and the specific total energy $E$ are given by
\be\nonumber\ba
\mathbb{S} = 2 \mu D(u) + \lambda \div u\,\mathbb{I}_2, \quad E = e + \frac{1}{2} |u|^2.
\ea\ee
Here $D(u) = \frac{1}{2} \left[ \na u + (\na u)^{T} \right]$ is the deformation tensor, and $\mathbb{I}_2$ denotes the $2 \times 2$ identity matrix.
Throughout, $t \ge 0$ and $x \in \OM \subset \rr$ denote the time and spatial variables, respectively, while
$\n$, $u=(u^1,u^2)$, $e$, $P$, and $\theta$ represent the fluid density, velocity, specific internal energy, pressure, and absolute temperature.

We assume constant shear viscosity and power-law heat conductivity,
\begin{equation}\label{ext:heatlaw}
 \mu=\text{constant}>0,\qquad
 \kappa(\theta)=\kappa_0\theta^\alpha,\qquad \kappa_0>0,\quad \alpha>0 .
\end{equation}
The bulk viscosity satisfies, for a constant $C_\lambda>0$,
\begin{equation}\label{visc:regularity}
 \lambda\in C^1([0,\infty))\cap C^2((0,\infty)),\quad
 \lambda\ge0,\quad r|\lambda'(r)|\le C_\lambda(2\mu+\lambda(r))
 \quad(r>0).
\end{equation}
Define
\begin{equation}\label{visc:growth1}
 \Lambda(r)=\int_1^r\frac{2\mu+\lambda(z)}{z}\,dz,\qquad r>0 .
\end{equation}

For an ideal polytropic gas, the equations of state are given by
\be\la{sd2}\ba
P(\n,e) = (\ga-1)\n e = R \n \theta, \quad e = \frac{R \theta}{\ga-1},
\ea\ee
where $\ga>1$ is the adiabatic exponent and $R>0$ is the gas constant.
Without loss of generality, we normalize $R=1$ and set
\be\la{sd3}\ba
c_v \triangleq \frac{R}{\ga-1}, \quad K(\theta) \triangleq \int_0^\theta \kappa(s)\,ds.
\ea\ee
For strong solutions, the system \eqref{ffns} can then be rewritten as
\be\ba\la{fns}
\begin{cases}
\rho_t + \div(\rho u) = 0, \\
\n(u_t+u \cdot \na u) = \mu \Delta u + \na\big((\mu+\lm)\div u\big) - \na P, \\
c_v\n(\theta_t + u \cdot \na \theta) - \Delta K(\theta) = 2\mu |D(u)|^2 + \lambda(\div u)^2 - P \div u,
\end{cases}
\ea\ee
supplemented with the initial data
\be\la{cztj}
\rho(0)=\rho_0, \quad (\rho u)(0)=\rho_0 u_0, \quad (\rho\theta)(0)=\rho_0\theta_0.
\ee
In this paper, we consider the following three settings:

\textit{Case 1. Periodic boundary conditions.}
In this case,
\be\la{zqbjtj}\ba
\OM=\T^2 \triangleq \mr^2/\mathbb{Z}^2, \ \text{and} \  (\n,u,\theta) \ \text{is $1$-periodic in each spatial direction}.
\ea\ee

\textit{Case 2. Navier-slip boundary conditions.}
In this case, $\OM$ is a simply connected bounded $C^{3,1}$ domain in $\rr$, and
\be\la{tri:bd-boundary}\ba
 u\cdot n=0,\qquad \omega=-A u\cdot\tau,\qquad \nabla\theta\cdot n=0
 \quad\text{on }\partial\Omega\times(0,\infty),
\ea\ee
where $\omega \triangleq \pa_2 u^1 - \pa_1 u^2$ is the vorticity, $n$ is the unit outer normal to $\partial\Omega$, $\tau=(n_2,-n_1)$ is the unit tangent vector, $\chi=(\partial_\tau n)\cdot\tau$ denotes the boundary curvature, and $A$ is a prescribed smooth nonnegative function on $\partial\Omega$.

\textit{Case 3. Cauchy problem.}
In this case, $\OM=\rr$, and $(\n,u,\theta)$ satisfies the far-field condition
\be\la{tri:farfield}\ba
 (\rho,u,\theta)(x,t)\longrightarrow(\rho_\infty,0,\theta_\infty)
 \quad\text{as }|x|\to\infty,
\ea\ee
where $\rho_\infty$ and $\theta_\infty$ are fixed positive constants.

There is an extensive literature on the existence and large-time behavior of solutions to the full compressible Navier--Stokes equations.
In one space dimension, Kazhikhov and Shelukhin~\cite{KazhikhovShelukhin1977} first proved the global existence and uniqueness of strong solutions to the initial-boundary value problem with constant viscosity and heat conductivity in the absence of vacuum; Kazhikhov~\cite{Kazhikhov1982} subsequently extended this theory to the Cauchy problem for large, regular perturbations of a constant state with positive density and temperature.
For density-dependent viscosity and temperature-dependent heat conductivity, Duan, Guo, and Zhu~\cite{DuanGuoZhu2017} obtained global well-posedness under stress-free, thermally insulated boundary conditions, assuming the initial density is bounded away from zero.
For the multi-dimensional problem, early results include Serrin's uniqueness theorem for smooth non-vacuum solutions~\cite{Serrin1959} and Nash's local existence result for the three-dimensional Cauchy problem with positive density~\cite{Nash1962}.
When vacuum is allowed, Cho and Kim~\cite{ChoKim2006} constructed local strong solutions under compatibility conditions on the initial momentum and thermal energy equations.
In the small-perturbation regime, Matsumura and Nishida~\cite{MatsumuraNishida1980} established the existence and uniqueness of global classical solutions in $\mathbb{R}^3$ for $H^3$ perturbations of a constant state with positive density and temperature, whereas Hoff~\cite{Hoff1997} constructed global weak solutions admitting discontinuities in two and three dimensions.
As for global weak solutions with large data, Feireisl~\cite{Feireisl2004} established the existence of global variational solutions for a class of pressure laws excluding the ideal polytropic gas, while Bresch and Desjardins~\cite{BreschDesjardins2007} proved the stability of weak solutions under a structural relation linking density-dependent shear and bulk viscosities.
Allowing vacuum and large oscillations, Huang and Li~\cite{HuangLi2018} established global classical solutions in $\mathbb{R}^3$ for smooth initial data with small initial energy and positive far-field density and temperature; this result was subsequently extended to bounded domains by Li, L\"u, and Wang~\cite{LLW} and to exterior domains by Li, Li, and L\"u~\cite{LiLiLu2025}.
Under different smallness mechanisms, Wen and Zhu~\cite{WenZhu2017} obtained global strong solutions for sufficiently small initial mass, whereas Li~\cite{Li2020} proved global well-posedness under a scaling-invariant smallness condition assuming $2\mu>\lambda$, a restriction later removed by Wen~\cite{Wen2025}.
In three-dimensional bounded domains with Dirichlet boundary conditions, Shen, Xu, and Zhang~\cite{ShenXuZhang2025} also obtained global strong solutions with density-dependent viscosity and initial vacuum for sufficiently small initial energy.

In the isentropic case with $\lambda(\rho)=\rho^\beta$, the system (\ref{fns}) corresponds to the model introduced by Vaigant-Kazhikhov \cite{VK}, who established the existence and uniqueness of global strong solutions for large initial data with density away from vacuum in rectangular domains under the assumption that $\beta>3$.
Subsequently, Jiu, Wang, and Xin~\cite{JWX1} extended this result to allow initial vacuum on the periodic torus.
Huang and Li~\cite{HL2,HL3} (see also Jiu, Wang, and Xin~\cite{JWX2}) relaxed the crucial condition from $\beta>3$ to $\beta>\frac{4}{3}$ by developing new techniques based on commutator theory and blow-up criteria in periodic domains and the whole space, allowing the density to vanish.
Later, Fan, Li, and Li~\cite{FLL} established global strong solutions for $\beta>\frac{4}{3}$ in general bounded simply connected domains with Navier-slip boundary conditions, and Fan, Li, and Wang~\cite{FLW} derived time-independent density upper bounds and exponential decay in both periodic and bounded simply connected domains.
For non-isentropic flows, Fan and Huang~\cite{FanHuang} recently established the global existence and large-time behavior of strong solutions to the two-dimensional compressible Navier--Stokes equations with transported entropy and $\lambda(\rho)=\rho^\beta$ ($\beta>\frac{4}{3}$).
Their result holds for arbitrarily large initial data with vacuum on both the torus and bounded simply connected domains.
In the bounded-domain case, they exploited the Navier-slip boundary condition to derive a new commutator estimate.
More recently, for a broad class of density-dependent bulk viscosities, including all superlinear power laws and near-linear laws with logarithmic gains, Fan, Jiang, Li, and L\"u~\cite{FJLL1,FJLL2} established the global existence and large-time behavior of strong solutions in periodic domains, the whole space, and smooth bounded domains with Navier-slip boundary conditions.
Subsequently, for the linear bulk viscosity $\lambda(\rho)=\rho$, Li, Liang, and Ye~\cite{LLY} and Fan, Huang, and Lei~\cite{FHL} proved global existence and uniqueness of strong solutions for large initial data in $\mathbb{R}^2$ with vacuum far-field density and on $\mathbb{T}^2$, respectively.

Extending this large-data theory to the full heat-conducting system \eqref{fns} with degenerate thermal conductivity introduces substantial new difficulties: the pressure work $P\div u$ and viscous dissipation nonlinearly couple the momentum and temperature equations, while the degeneracy $\kappa(0)=0$ leads to a loss of uniform parabolicity as the temperature approaches zero.
Addressing these difficulties requires delicate estimates to handle both the effects of vacuum and high-temperature growth without imposing any smallness assumptions on the initial data.
The aim of this paper is to establish the global existence and large-time behavior of strong solutions to the two-dimensional full compressible Navier--Stokes equations \eqref{fns} with degenerate heat conductivity for arbitrarily large initial data allowing vacuum.

\begin{definition}[Coefficient and data classes]\label{coeff:admissible}
Set
$$h(r)\triangleq 2\mu+\lambda(r), \quad L(r)\triangleq \log(e+r).$$
Fix the domain, the physical parameters, $C_\lambda$ in \eqref{visc:regularity}, and $H_0\ge\sup_{0\le r\le 2}h(r)$. Assume \eqref{ext:heatlaw} and \eqref{visc:regularity}, together with one of the following four conditions:

\smallskip\noindent\textup{(A) Pointwise lower bound ($\alpha>1$).}
Let $r_*\ge2$, $c_*>0$, and let
$D:[2,\infty)\to[1,\infty)$ be a nondecreasing function such that,
for some $A_D>0$ and $0\le\beta_D<1/2$,
\begin{equation}\label{new:den}
 D(r)\le A_D L(r)^{\beta_D},
 \qquad
 \int_2^\infty\frac{dr}{(e+r)L(r)D(r)}=\infty.
\end{equation}
We consider the following two alternatives.

For a fixed viscosity law and arbitrary compatible initial data, assume
\begin{equation}\label{ext:iter}
 h(r)\ge\frac{c_*r}{D(r)}
 \qquad (r\ge r_*).
\end{equation}

For the data-class alternative, fix $D,c_*,r_*$ and the constants
$A_D,\beta_D$ in \eqref{new:den}, in addition to the domain, physical
parameters, $C_\lambda$, and $H_0$ fixed above. Given $M>0$, consider
all admissible viscosity laws satisfying
\begin{equation}\label{ext:powerfamily}
 h(r)\ge\frac{c_*r}{(1+r)^\varepsilon D(r)}
 \qquad (r\ge r_*),\qquad
 0\le\varepsilon\le\varepsilon_*(M),
\end{equation}
and all corresponding compatible initial data satisfying
\begin{equation}\label{ext:dataM}
\begin{gathered}
 \|\rho_0\|_{L^\infty}+\|u_0\|_{H^1}+\|\theta_0\|_{L^\infty}
 +(\operatorname*{ess\,inf}\theta_0)^{-1}
 +\|g_1\|_{L^2}+\|\nabla\theta_0\|_{L^2}\le M,\\
 M_0+M_0^{-1}\le M\quad (|\Omega|<\infty),\qquad
 \|\rho_0-\rho_\infty\|_{L^2}
 +\|\theta_0-\theta_\infty\|_{L^2}\le M
 \quad (\Omega=\mathbb R^2).
\end{gathered}
\end{equation}
For each $M>0$, the threshold $\varepsilon_*(M)>0$ can be chosen to depend only on $M$ and the quantities fixed above.

\smallskip\noindent\textup{(B) Integral criterion ($\alpha>1$).}
Set $D_h(r)\triangleq 1+\sup_{0\le z\le r}\frac{z}{\mu+\lambda(z)}$ and assume that
\begin{equation}\label{new:fixed}
 \lim_{r\to\infty}\frac{\log D_h(r)}{\log L(r)}=0,\qquad
 \int_2^\infty\frac{\Lambda(r)}{r^2 L(r)}\,dr=\infty.
\end{equation}

\smallskip\noindent\textup{(C) Logarithmic viscosity gain.}
Let $0\le d<1$, $c_*>0$, and $r_*\ge 2$. We consider the following two alternatives.

For a fixed viscosity law and arbitrary compatible initial data, assume
\begin{equation}\label{new:gain}
 \mu+\lambda(r)\ge c_* r L(r)^d\qquad (r\ge r_*),
\end{equation}
together with
\[
 \alpha\ge \frac{2(1-d)}{2+d} \ \textnormal{(global existence)},\quad
 \alpha>1-d \ \textnormal{(time-uniform bounds)}.
\]

For the data-class alternative, fix $0<d<1$, $\alpha>1-d$, $c_*>0$, and $r_*\ge2$, in addition to the domain, physical parameters, $C_\lambda$, and $H_0$ fixed above. Given $M>0$, consider all admissible viscosity laws satisfying
\begin{equation}\label{ext:gainfamily}
 \mu+\lambda(r)\ge\frac{c_*rL(r)^d}{(1+r)^\varepsilon}
 \qquad(r\ge r_*),\qquad
 0\le\varepsilon\le\varepsilon_{*,C}(M),
\end{equation}
and all corresponding compatible initial data satisfying \eqref{ext:dataM}. For each $M>0$, the threshold $\varepsilon_{*,C}(M)>0$ can be chosen to depend only on $M$ and the quantities fixed above.

\smallskip\noindent\textup{(D) Linear heat conductivity ($\alpha=1$).}
Set $J(r)\triangleq \log\log(e^e+r)$, and assume that there exist $c>0$, $r_*>1$, $p\ge 0$, and $\Psi\in C^1([1,\infty))$ with $\Psi\ge 1$ such that
\begin{equation}\label{crit:viscosity}
 \mu+\lambda(r)\ge \frac{c r}{\Psi(J(r))}\qquad (r\ge r_*),
\end{equation}
and
\begin{equation}\label{crit:profile}
 0\le \frac{x\Psi'(x)}{\Psi(x)}\le p\quad (x\ge 1),\qquad
 \int_2^\infty \frac{dx}{\Psi(x)\log(e+x)}=\infty.
\end{equation}
\end{definition}

\begin{remark}[Admissible exponents and further examples]\label{rem:alpha-ranges}
\normalfont
Let $\lambda_0\ge 0$ and $\kappa(\theta)=\kappa_0\theta^\alpha$ with $\kappa_0>0$. For arbitrary compatible initial data, Definition~\ref{coeff:admissible} yields the following representative ranges of $\alpha$:
\begin{center}
\renewcommand{\arraystretch}{1.35}
\begin{tabular}{@{}l@{\hspace{1em}}c@{\hspace{1em}}c@{\hspace{1em}}c@{}}
\hline
Bulk viscosity $\lambda(r)$ & Global existence & Time-uniform bounds & Case \\
\hline\noalign{\vskip 4pt}
$\lambda_0+\dfrac{r}{[\log\log(e^e+r)]^\eta}$, $0<\eta\le 1$ & $\alpha\ge 1$ & $\alpha>1$ & \textup{(A)}, \textup{(D)} \\[7pt]
$\lambda_0+r$ & $\alpha\ge 1$ & $\alpha>1$ & \textup{(A)}, \textup{(D)} \\[3pt]
$\lambda_0+r[\log(e+r)]^d$, $0<d<1$ & $\alpha\ge\dfrac{2(1-d)}{2+d}$ & $\alpha>1-d$ & \textup{(C)} \\[6pt]
$\lambda_0+r[\log(e+r)]^d$, $d\ge 1$ & $\alpha>0$ & $\alpha>0$ & \textup{(C)} \\[2pt]
$\lambda_0+r\log(e+r)[\log\log(e^e+r)]^\delta$, $\delta\in\mathbb R$ & $\alpha>0$ & $\alpha>0$ & \textup{(C)} \\[2pt]
$\lambda_0+r^m[\log(e+r)]^s$, $m>1$, $s\in\mathbb R$ & $\alpha>0$ & $\alpha>0$ & \textup{(C)} \\
\hline
\end{tabular}
\end{center}
The first row follows by taking $D(r)=J(r)^\eta$ in \textup{(A)} and $\Psi(x)=x^\eta$ in \textup{(D)}; in particular, it also covers $\lambda_0+r/\log\log(e+r)$ (with $\lambda(0)=\lambda_0+e$ by continuity).

To construct further examples involving iterated logarithms, define $L_1(r)=\log(e+r)$ and $L_{k+1}(r)=1+\log L_k(r)$ for $k\ge 1$, with the convention that an empty product equals $1$. For any integer $m\ge 2$, one may take
\[
 \begin{aligned}
  \textup{(A)}\ (\alpha>1):\quad & D(r)=\prod_{k=2}^{m} L_k(r)
   \qquad \Bigl(L_{m+1}'(r)=\frac{1}{(e+r)L_1(r)D(r)}\Bigr),\\
  \textup{(D)}\ (\alpha=1):\quad & \Psi(x)=x\prod_{k=2}^{m} L_k(x),\\
  &\lambda(r)=\lambda_0+\frac{r}{L_2(r)^\beta\prod_{k=3}^{m}L_k(r)^{s_k}}
   \quad (0\le \beta<1,\ s_3,\ldots,s_m\in\mathbb R).
 \end{aligned}
\]
In particular, setting $m=2$ in \textup{(D)} yields the denominator $J(r)L_2(J(r))$ in \eqref{crit:viscosity}.
\end{remark}

Before stating our main results, we introduce the notation and conventions used throughout the paper. We define
\be\ba\nonumber
\na^\bot \triangleq (\p_2,-\p_1),\quad
\int f\,dx \triangleq \int_{\OM} f\,dx,\quad
\ol{f} \triangleq \frac{1}{|\OM|}\int f\,dx \quad (\text{when } |\OM|<\infty).
\ea\ee
For a positive integer $s$ and $1\leq r\leq \infty$, we denote the standard Lebesgue and Sobolev spaces by
\be\ba\nonumber
L^r =L^r(\OM),\quad W^{s,r} =W^{s,r}(\OM),\quad H^s =W^{s,2}.
\ea\ee
The material derivative is defined by
\be\ba\nonumber
\frac{D}{Dt}f=\dot{f} \triangleq f_t + u\cdot\na f,
\ea\ee
and the effective viscous flux $G$ and vorticity $\o$ are given by
\be\ba\la{gw}
G \triangleq (2\mu+\lam)\div u-P, \quad \o \triangleq \na^\bot \cdot u = \pa_2 u^1 - \pa_1 u^2.
\ea\ee

We now give the definition of strong solutions to \eqref{fns}.
\begin{definition}
If all derivatives involved in \eqref{fns} for $(\n,u,\theta)$ are regular distributions and the equations in \eqref{fns} hold almost everywhere in $\OM \times(0,T)$, then $(\n,u,\theta)$ is called a strong solution.
\end{definition}

Our first result concerns the global existence and large-time behavior of strong solutions on the torus.

\begin{theorem}\la{th1}
Let $\OM=\T^2$. Assume one of the conditions \textup{(A)--(D)} in Definition~\ref{coeff:admissible}.
Let the initial data $(\n_0,u_0,\theta_0)$ satisfy, for some $q>2$ and a positive constant $\underline\theta_0$,
\be\la{th12}\ba
& \n_0 \in W^{1,q}(\T^2), \quad \n_0 \ge 0, \quad M_0 \triangleq \int_{\T^2} \n_0\,dx > 0, \\
& u_0 \in H^2(\T^2), \quad \theta_0 \in H^1(\T^2) \cap L^\infty(\T^2), \quad \inf_{x \in \T^2} \theta_0(x) \ge \underline\theta_0>0,
\ea\ee
and the momentum compatibility condition
\be\la{th13}\ba
-\mu \Delta u_0 - \na\big((\mu+\lambda(\n_0)) \div u_0\big) + \na(\n_0 \theta_0) = \sqrt{\n_0}\,g_1,
\ea\ee
with $g_1\in L^2(\T^2;\rr)$.
The problem \eqref{fns}--\eqref{zqbjtj} admits a unique global strong solution $(\n,u,\theta)$ in $\T^2 \times (0,\infty)$ satisfying, for any $0<\tau<T<\infty$,
\be\la{th15}\ba
\begin{cases}
\rho\in C([0,T];W^{1,q}), \quad \n_t \in L^\infty(0,T;L^q), \\
u \in C([0,T];H^1) \cap L^\infty(0,T;H^2) \cap L^2(0,T;W^{2,q}), \\
\theta\in L^\infty(0,T;H^1)\cap L^2(0,T;H^2),\quad \rho\theta\in C([0,T];L^2), \\
\theta\in C([\tau,T];H^2)\cap L^2(\tau,T;H^3),\quad \theta_t\in L^2(\tau,T;H^1), \\
\sqrt\rho\,u_t\in L^\infty(0,T;L^2),\quad u_t\in L^2(0,T;H^1), \\
\sqrt\rho\,\theta_t\in L^2(0,T;L^2)\cap L^\infty(\tau,T;L^2).
\end{cases}
\ea\ee
For every $T>0$ there are constants $c_T,C_T>0$ such that
\begin{equation}\label{main:finite-bounds}
 0\le\rho(x,t)\le C_T,\qquad c_T\le\theta(x,t)\le C_T,
 \qquad (x,t)\in\Omega\times[0,T].
\end{equation}
Suppose now that \textup{(A)}, \textup{(B)}, or \textup{(C)} with $\alpha>1-d$ holds. Set
\be\la{ltp-equilibrium}\ba
u_\infty = \frac{1}{M_0}\int_{\T^2}\rho_0 u_0\,dx,\quad
\theta_\infty = \frac{1}{c_v M_0}\int_{\T^2}
     \left(c_v\rho_0\theta_0+\frac12\rho_0|u_0-u_\infty|^2\right)dx .
\ea\ee
Then there exist constants $C>0$ and $\theta_*>0$, independent of time, such that
\be\la{ltp-main-uniform}\ba
0\le\rho(x,t)\le C,\qquad
\theta_*\le\theta(x,t)\le C,\qquad (x,t)\in\T^2\times[0,\infty).
\ea\ee
Moreover, for each $1\le p<\infty$, there exist positive constants $C_p$ and $\eta_p$ such that, for all $t\ge1$,
\be\la{ltp-main-decay}\ba
&\|\rho(t)-\ol{\rho}_0\|_{L^p}
 +\|u(t)-u_\infty\|_{W^{1,p}}
 +\|\theta(t)-\theta_\infty\|_{W^{1,p}}
 \le C_p e^{-\eta_p t}.
\ea\ee
In the data-class assertions of \textup{(A)} and \textup{(C)}, there exists a constant $C_\rho(M)>0$ such that $\sup_{t\ge0}\|\rho(t)\|_{L^\infty}\le C_\rho(M)$ uniformly over the corresponding admissible coefficients and initial data.
\end{theorem}

The second result establishes the global existence and exponential decay of strong solutions in bounded domains.

\begin{theorem}\la{tri:main}
Let $\Omega\subset\mathbb R^2$ be a bounded, simply connected $C^{3,1}$ domain. Assume that $A+2\chi\ge0$ on $\partial\Omega$ and that one of the conditions \textup{(A)--(D)} in Definition~\ref{coeff:admissible} holds.
Let the initial data $(\n_0,u_0,\theta_0)$ satisfy, for some $q>2$ and a positive constant $\underline\theta_0$,
\be\la{tri:bd-data}\ba
 &\rho_0\in W^{1,q}(\Omega),\quad\rho_0\ge0,\quad M_0\triangleq\int_\Omega\rho_0\,dx>0, \quad u_0\in H^2(\Omega),\\
 &\theta_0 \in H^1(\Omega) \cap L^\infty(\Omega),\quad \theta_0\ge\underline\theta_0,\quad
 u_0\cdot n=0,\quad\omega_0=-A u_0\cdot\tau,
\ea\ee
together with \eqref{th13} for some $g_1\in L^2(\Omega;\mathbb R^2)$.
The thermal initial space is the $H^1$ form domain of the Neumann Laplacian. Smooth approximations satisfy $\partial_n\theta_{0,j}=0$; for data with a classical compatible normal trace take $\partial_n\theta_0=0$. The problem \eqref{fns}, \eqref{cztj}, and \eqref{tri:bd-boundary} admits a unique global strong solution with the regularity in \eqref{th15} and the bounds \eqref{main:finite-bounds} on every finite interval $[0,T]$.
Assume in addition that \textup{(A)}, \textup{(B)}, or \textup{(C)} with $\alpha>1-d$ holds. Then there exist constants $C>0$ and $\theta_*>0$, independent of time, such that
\be\la{tri:bd-uniform}\ba
 0\le\rho(x,t)\le C,\qquad 0<\theta_*\le\theta(x,t)\le C,
 \qquad(x,t)\in\Omega\times[0,\infty).
\ea\ee
For every $1\le p<\infty$, there exist positive constants $C_p$ and $\eta_p$ such that, for all $t\ge1$,
\be\la{tri:bd-decay}\ba
 \|\rho-\rho_*\|_{L^p}+\|u\|_{W^{1,p}}+\|\theta-\theta_{\rm eq}\|_{W^{1,p}}\le C_p e^{-\eta_p t},
\ea\ee
where $\rho_*=M_0/|\Omega|$ and
\be\la{tri:bd-equilibrium}\ba
 \theta_{\rm eq}=\frac{1}{c_vM_0}\left[\int_\Omega\left(\frac12\rho_0|u_0|^2+c_v\rho_0\theta_0\right)dx-\mu\int_0^\infty\!\int_{\partial\Omega}(A+2\chi)|u\cdot\tau|^2\,dSdt\right]>0.
\ea\ee
In the data-class assertions of \textup{(A)} and \textup{(C)}, there exists a constant $C_\rho(M)>0$ such that $\sup_{t\ge0}\|\rho(t)\|_{L^\infty}\le C_\rho(M)$ uniformly over the corresponding admissible coefficients and initial data.
\end{theorem}

\begin{remark}\label{rem:bd-equilibrium-positive}
\normalfont
Under \textup{(A)}, \textup{(B)}, or \textup{(C)} with $\alpha>1-d$, conservation of mass and \eqref{tri:bd-uniform} give
\[
 E(t)=\int_\Omega\left(\frac12\rho|u|^2+c_v\rho\theta\right)dx
 \ge c_v\int_\Omega\rho\theta\,dx\ge c_v M_0\theta_*,\qquad t\ge0.
\]
By the total energy balance under the thermally insulated Navier-slip boundary condition, we have
$E'(t)=-\mu\int_{\partial\Omega}(A+2\chi)|u\cdot\tau|^2\,dS\le0$.
Consequently, $\lim_{t\to\infty}E(t)$ exists, and \eqref{tri:bd-equilibrium} yields
\[
 \theta_{\rm eq}=\frac{\lim_{t\to\infty}E(t)}{c_vM_0}\ge\theta_*>0 .
\]
\end{remark}

Finally, for the Cauchy problem \eqref{fns}, \eqref{cztj}, \eqref{tri:farfield}, we have the following global existence result.

\begin{theorem}\la{tri:whole}
Let $\Omega=\mathbb R^2$.
Assume one of the conditions \textup{(A)--(D)} in Definition~\ref{coeff:admissible}.
Let the initial data $(\n_0,u_0,\theta_0)$ satisfy, for some $q>2$ and a positive constant $\underline\theta_0$,
\be\la{tri:whole-data}\ba
 &\rho_0\ge0,\quad\rho_0-\rho_\infty\in L^2(\mathbb R^2)\cap W^{1,q}(\mathbb R^2), \quad u_0\in H^2(\mathbb R^2), \\
 &\theta_0-\theta_\infty \in H^1(\mathbb{R}^2) \cap L^\infty(\mathbb{R}^2), \quad\theta_0\ge\underline\theta_0,
\ea\ee
together with \eqref{th13} for some $g_1\in L^2(\mathbb R^2;\mathbb R^2)$.
The problem \eqref{fns}, \eqref{cztj}, and \eqref{tri:farfield} admits a unique global strong solution $(\n,u,\theta)$ in $\rr \times (0,\infty)$ satisfying, for any $0<\tau<T<\infty$,
\be\la{th35}\ba
\begin{cases}
\rho - \rho_\infty \in C([0,T];W^{1,q}), \quad \n_t \in L^\infty(0,T;L^q), \\
u \in C([0,T];H^1) \cap L^\infty(0,T;H^2) \cap L^2(0,T;W^{2,q}), \\
\theta-\theta_\infty\in L^\infty(0,T;H^1)\cap L^2(0,T;H^2),\quad \rho\theta-\rho_\infty\theta_\infty\in C([0,T];L^2), \\
\theta-\theta_\infty\in C([\tau,T];H^2)\cap L^2(\tau,T;H^3),\quad \theta_t\in L^2(\tau,T;H^1), \\
\sqrt\rho\,u_t\in L^\infty(0,T;L^2),\quad u_t\in L^2(0,T;H^1), \\
\sqrt\rho\,\theta_t\in L^2(0,T;L^2)\cap L^\infty(\tau,T;L^2).
\end{cases}
\ea\ee
The bounds \eqref{main:finite-bounds} hold on every finite interval $[0,T]$. Assume in addition that \textup{(A)}, \textup{(B)}, or \textup{(C)} with $\alpha>1-d$ holds. Then there exist constants $C>0$ and $\theta_*>0$, independent of time, such that
\be\la{tri:whole-uniform}\ba
 0\le\rho(x,t)\le C,\qquad 0<\theta_*\le\theta(x,t)\le C,
 \quad(x,t)\in\mathbb R^2\times[0,\infty).
\ea\ee
Moreover, for every $2\le p<\infty$, it holds that
\be\la{tri:whole-decay}\ba
 \lim_{t\to\infty}\left(\|\rho-\rho_\infty\|_{L^p}+\|u\|_{W^{1,p}}+\|\theta-\theta_\infty\|_{W^{1,p}}\right)=0.
\ea\ee
In the data-class assertions of \textup{(A)} and \textup{(C)}, there exists a constant $C_\rho(M)>0$ such that $\sup_{t\ge0}\|\rho(t)\|_{L^\infty}\le C_\rho(M)$ uniformly over the corresponding admissible coefficients and initial data.
\end{theorem}

\section{Preliminaries}

In this section, we recall some known facts and elementary inequalities that will be used frequently later.

We begin with the local existence of strong solutions, which can be established by arguments similar to those in \cite{WX,LXZthermal}.

\begin{lemma}\label{lct}
Under the assumptions on the coefficients and initial data in Theorems~\ref{th1}--\ref{tri:whole}, there exists a small time $T>0$ such that the problem \eqref{fns}--\eqref{cztj}, supplemented with the corresponding boundary or far-field condition \eqref{zqbjtj}, \eqref{tri:bd-boundary}, or \eqref{tri:farfield}, admits a unique strong solution $(\rho,u,\theta)$ on $\Omega\times(0,T]$ enjoying the regularity in \eqref{th15} (or \eqref{th35} when $\Omega=\mathbb{R}^2$).
\end{lemma}

Next, we derive elementary growth and comparison estimates relating the bulk viscosity $\lambda$ to its primitive $\Lambda$ under assumption \eqref{visc:regularity}.

\begin{lemma}\label{lem:visc-envelope}
Suppose that $\lambda$ satisfies \eqref{visc:regularity}, and let $\Lambda$ be defined by \eqref{visc:growth1}. Then there exist constants $B_h\ge 2$ and $C>0$, depending only on $C_\lambda$ and $\sup_{0\le s\le 1}(2\mu+\lambda(s))$, such that
\begin{equation}\label{visc:envelope}
\begin{aligned}
 2\mu+\lambda(r)&\le C(1+r)^{B_h}
     && (r\ge0),\\
 2\mu+\lambda(r)&\le C\Lambda(r),\quad
 2\mu+\lambda(r)\le Cr\Lambda(r)^{1-1/B_h}
     && (r\ge2).
\end{aligned}
\end{equation}
\end{lemma}

\begin{proof}
Set $h(r)\triangleq 2\mu+\lambda(r)\ge 2\mu>0$. By \eqref{visc:regularity}, we have
\[
 \left|\frac{d}{dr}\log h(r)\right|
 =\frac{|\lambda'(r)|}{h(r)}
 \le \frac{C_\lambda}{r}
 \qquad (r>0).
\]
For any $r\ge 1$, integrating this bound over $[1,r]$ yields
\[
 \log\frac{h(r)}{h(1)}
 =\int_1^r \frac{h'(z)}{h(z)}\,dz
 \le C_\lambda \int_1^r \frac{dz}{z}
 = C_\lambda \log r,
\]
and hence
\begin{equation*}
 h(r)\le h(1)\,r^{C_\lambda}\qquad (r\ge 1).
\end{equation*}
Setting $B_h\triangleq \max\{2,C_\lambda\}$ and noting that $h(r)\le \sup_{0\le s\le 1}h(s)$ for $r\in[0,1]$, we immediately obtain the first inequality in \eqref{visc:envelope}:
\[
 h(r)\le C(1+r)^{B_h}\qquad (r\ge 0).
\]
Next, for any $r\ge 2$ and $z\in [r/2,r]$, integrating over $[z,r]$ gives
\[
 \log\frac{h(r)}{h(z)}
 =\int_z^r \frac{h'(s)}{h(s)}\,ds
 \le C_\lambda \int_z^r \frac{ds}{s}
 = C_\lambda \log\frac{r}{z}
 \le C_\lambda \log 2,
\]
which implies
\[
 h(z)\ge 2^{-C_\lambda} h(r)\qquad \text{for all } z\in[r/2,r].
\]
Since $[r/2,r]\subset [1,r]$ for $r\ge 2$ and $h>0$, it follows from \eqref{visc:growth1} that
\[
 \Lambda(r)
 =\int_1^r \frac{h(z)}{z}\,dz
 \ge \int_{r/2}^r \frac{h(z)}{z}\,dz
 \ge 2^{-C_\lambda}h(r)\int_{r/2}^r \frac{dz}{z}
 =(2^{-C_\lambda}\log 2)\,h(r).
\]
This yields the second inequality in \eqref{visc:envelope}:
\[
 h(r)\le \frac{2^{C_\lambda}}{\log 2}\,\Lambda(r)\le C\Lambda(r)\qquad (r\ge 2).
\]

Finally, for $r\ge 2$, since $B_h\ge C_\lambda$ and $r\ge 1$, the bound $h(r)\le h(1)r^{C_\lambda}\le h(1)r^{B_h}$ implies
\[
 h(r)^{1/B_h}\le h(1)^{1/B_h}\,r,
\]
whereas the second bound in \eqref{visc:envelope} gives
\[
 h(r)^{1-1/B_h}\le C^{1-1/B_h}\,\Lambda(r)^{1-1/B_h}.
\]
Multiplying these two inequalities yields
\[
 h(r)=h(r)^{1/B_h}\,h(r)^{1-1/B_h}\le Cr\Lambda(r)^{1-1/B_h}\qquad (r\ge 2),
\]
which completes the proof of \eqref{visc:envelope}.
\end{proof}

We next recall the two-dimensional Gagliardo--Nirenberg inequality \cite{NI,TG}.
\begin{lemma}\la{gn1}
Assume that $\OM=\T^2$, $\OM=\rr$, or $\OM$ is a bounded smooth domain in $\rr$.
For $2<p<\infty$, there exists a positive constant $C$ depending only on $\OM$ such that for any $f \in H^1(\OM)$,
\be\ba\la{gn11}
\| f \|_{L^p} \le Cp^{1/2}\| f \|^{2/p}_{L^2} \| f \|^{1-2/p}_{H^1}.
\ea\ee
Moreover, the term $\| f \|_{H^1}$ on the right-hand side of \eqref{gn11} can be replaced by $\| \na f \|_{L^2}$ whenever $\OM=\rr$, or $\ol{f}=0$, or $f\in H^1(\OM;\rr)$ satisfies $f \cdot n|_{\p \OM}=0$.
\end{lemma}

The following Poincar\'e-type inequality controls the $H^1$-norm of a nonnegative function bounded on a subset of positive measure.

\begin{lemma}\la{1pi}
Let $\OM=\T^2$ or let $\OM \subset \rr$ be a bounded Lipschitz domain.
Assume that $f \in H^1(\OM)$ is nonnegative and that there exist a constant $M\ge 0$ and a measurable subset $\OM_0 \subset \OM$ with $|\OM_0|>0$ such that
\be\la{1pi01}\ba
f \le M \quad \textnormal{a.e. in } \OM_0.
\ea\ee
Then there exists a positive constant $C$ depending only on $|\OM_0|$, $\OM$, and $M$ such that
\be\la{1pi02}\ba
\| f \|_{H^1} \le C + C \| \na f \|_{L^2}.
\ea\ee
\end{lemma}
\begin{proof}
Since $f \ge 0$, we deduce from (\ref{1pi01}) that
\be\la{1pi11}\ba
\ol{f} - M \le \ol{ f } - f \ \text{ in } \OM_0.
\ea\ee
By (\ref{1pi11}) and Poincar\'e's inequality, we have
\be\la{1pi12}\ba
\left( \ol{ f } - M \right)_+^2
\le \frac{1}{|\OM_0|} \int_{\OM_0} \left| \ol{ f } - f \right|^2 dx
\le C \int_{\OM} |\na f|^2 dx,
\ea\ee
which yields
\be\la{1pi13}\ba
\ol{ f } \le M + C \| \na f \|_{L^2}.
\ea\ee
Using (\ref{1pi13}) and Poincar\'e's inequality, we arrive at
\be\la{1pi14}\ba
\| f \|_{H^1} \le \| f - \ol{ f } \|_{L^2} + |\OM|^{\frac{1}{2}} \ol{ f } + \| \na f \|_{L^2}
\le C + C \| \na f \|_{L^2}.
\ea\ee
This completes the proof.
\end{proof}

For a periodic vector field $v$, the Fourier representation of its divergence and curl gives, for $1<p<\infty$,
\be\la{dc1}\ba
 \|\nabla v\|_{L^p}\le C_p\big(\|\div v\|_{L^p}+\|\curl v\|_{L^p}\big).
\ea\ee
Here and below $\curl v=\partial_2v^1-\partial_1v^2$.
    
Let $\mathcal{H}^1(\T^2)$ and $\mathcal{BMO}(\T^2)$ denote the usual Hardy and BMO spaces.
Given a function $b$, define the commutator
\be\ba\la{p4}
[b,R_iR_j](f) \triangleq bR_i\circ R_j(f) - R_i \circ R_j(bf),\quad i,j=1,2,
\ea\ee
where $R_i = (-\Delta)^{-\frac{1}{2}}\pa_i$ is the usual Riesz transform on $\T^2$.

Now, we present the following well-known commutator estimates established by Coifman-Rochberg-Weiss \cite{CRW} and Coifman-Meyer \cite{CM}.
\begin{lemma}\la{jhzyl}
Let $b,f\in C^\infty(\T^2)$.
Then for $p\in (1,\infty)$, there is a positive constant $C$ depending only on $p$ such that
\be\ba\la{jhz1}
\|[b,R_iR_j](f)\|_{L^p}\leq C \| b\|_{\mathcal{BMO}}\|f\|_{L^p}.
\ea\ee
Moreover, for $p,q,r \in (1,\infty)$ with $\frac{1}{r} = \frac{1}{p}+ \frac{1}{q}$, there exists a positive constant $C$ depending only on $p$, $q$, and $r$ such that
\be\ba\la{jhz2}
\|\na[b,R_iR_j](f)\|_{L^r}\leq C \|\na b\|_{L^p}\|f\|_{L^q}.
\ea\ee
\end{lemma}

We also recall the following Beale--Kato--Majda inequality.
\begin{lemma}\la{bkm}
Let $2<q<\infty$.
Assume that $\OM$ is a simply connected bounded $C^{3,1}$ domain in $\rr$ and $u\in  \left\{ W^{2, q}(\Omega ) \big| u \cdot n=0, \curl u = -A u \cdot n^\bot \,\,\,\text{on}\,\,\, \partial\Omega \right\}$,
or $u \in W^{2,q}(\T^2)$, or $\na u \in L^2(\rr) \cap D^{1,q}(\rr)$.
Then there exists a positive constant $C$ depending only on $\OM$ and $q$ such that
\be\ba\la{bkm1}
\|\na u\|_{L^\infty} \le C \left( \|\div u \|_{L^\infty}
+ \|\curl u\|_{L^\infty} \right) \log \left(e+ \|\na^2 u\|_{L^q} \right)+ C\|\na u\|_{L^2}+C.
\ea\ee
\end{lemma}

\section{A priori estimates on the torus}\label{per-apriori}

In this section, we assume the coefficient hypotheses of Theorem~\ref{th1} and that the initial data $(\n_0,u_0,\theta_0)$ satisfy (\ref{th12}) with $\inf_{x \in \T^2} \n_0(x)>0$.
Let $(\n,u,\theta)$ be the strong solution to (\ref{fns})--(\ref{zqbjtj}) on $\T^2 \times (0,T]$, provided by Lemma \ref{lct}.

Moreover, since the total momentum is conserved for strong solutions of (\ref{fns}) and the system is Galilean invariant, we may assume without loss of generality that
\be\ba\nonumber
\int_{\T^2} \n u dx = \int_{\T^2} \n_0 u_0 dx = 0 \quad \text{ for all } t \in [0,T].
\ea\ee

We set
\be\ba\la{a1a2}
A_1^2(t) \triangleq 1 + \int \left( \frac{|G|^2}{2\mu+\lam} + \mu \o^2 \right) dx, \quad
A_2^2(t) \triangleq \int \n |\dot{u}|^2 dx,
\ea\ee
\be\ba\la{fp52}
 Q \triangleq 2\mu |D(u)|^2+\lambda(\div u)^2,\quad
 m(t) \triangleq \frac1{M_0}\int\rho\theta\,dx,\quad
 \mathcal D_\theta(t) \triangleq \int_\Omega
       \frac{\kappa(\theta)}{\theta^2}|\nabla\theta|^2\,dx,
\ea\ee
and
\be\ba\la{mdsj}
R_T \triangleq 2 + \sup_{0 \le t \le T} \| \n(t) \|_{L^\infty}.
\ea\ee

\subsection{Low-order estimates}

We first establish the basic energy, entropy, and temperature lower bound estimates.

\begin{lemma}\label{fpl1}
There exists a positive constant $C$, independent of $T$, such that
\be\la{fp01}\ba
 \sup_{0\le t\le T}\int\rho\bigl(1+|u|^2+\theta+(-\log\theta)_+
       +\log(1+\rho)\bigr)\,dx
 +\int_0^T\left(\int \frac{Q}{\theta}\,dx+\mathcal D_\theta(t)\right)dt\le C .
\ea\ee
If $\rho(x,t)\le r$ on $\T^2\times[0,T]$, then
\be\la{fp117}\ba
 \theta(x,t)\ge y(t)\triangleq(\underline\theta_0^{-1}+a_0t)^{-1},
 \quad a_0\triangleq\frac{1}{4c_v}\sup_{0\le z\le r}\frac{z}{\mu+\lambda(z)} .
\ea\ee
Moreover, there exist positive constants $M_1,M_2,h_0$ depending only on the initial energy, entropy, and mass such that, with $L_0\triangleq M_1$ and $\Gamma_t\triangleq\{x\in\T^2:\theta(x,t)\le M_1,\ \rho(x,t)\le M_2\}$,
\be\la{fp27}\ba
 |\Gamma_t|\ge \frac{M_0}{2M_2},\quad |\{\theta\ge h_0\}|\ge \frac{M_0}{2M_2},
\ea\ee
and for every $a\in(0,h_0)$ and every $n\ge 1$,
\be\la{extlow:cold}\ba
 |\{\theta\le a\}|\le C_{a,n}\,\mathcal D_\theta(t)^n .
\ea\ee
In addition, there exists a positive constant $C$, independent of $T$, such that
\be\la{ltp-low-u}\ba
 \|u\|_{L^j}\le C j^{1/2}\|\nabla u\|_{L^2},\quad 2\le j<\infty .
\ea\ee
\end{lemma}

\begin{proof}
First, integrating $(\ref{fns})_1$ over $\T^2$ gives the conservation of mass,
\be\nonumber\ba
 \int \rho\,dx = \int \rho_0\,dx = M_0.
\ea\ee
Multiplying $(\ref{fns})_2$ by $u$ and integrating by parts over $\T^2$, we derive
\be\nonumber\ba
 \frac{d}{dt}\int \frac{1}{2}\rho|u|^2\,dx
 +\int\bigl(\mu|\nabla u|^2+(\mu+\lambda)(\div u)^2\bigr)\,dx
 =\int P\div u\,dx.
\ea\ee
Integrating $(\ref{fns})_3$ over $\T^2$ yields
\be\nonumber\ba
 \frac{d}{dt}\int c_v\rho\theta\,dx
 =\int Q\,dx - \int P\div u\,dx.
\ea\ee
Noting that periodic integration by parts gives
\be\la{fp13}\ba
 \int Q\,dx = \int\bigl(\mu|\nabla u|^2+(\mu+\lambda)(\div u)^2\bigr)\,dx,
\ea\ee
we add the kinetic and internal energy identities to obtain
\be\la{fp15}\ba
E_* \triangleq \int\left(\frac{1}{2}\rho|u|^2+c_v\rho\theta\right) dx
= \int \left(\frac{1}{2}\rho_0|u_0|^2+c_v\rho_0\theta_0\right)dx.
\ea\ee

Next, to prove \eqref{fp117}, we use the two-dimensional pointwise bound $Q\ge(\mu+\lambda)(\div u)^2$ together with Young's inequality and $\rho\le r$ to deduce
\be\nonumber\ba
 P\div u = \rho\theta\div u
 \le (\mu+\lambda)(\div u)^2 + \frac{P^2}{4(\mu+\lambda)}
 \le Q + c_v a_0\rho\theta^2.
\ea\ee
Substituting this into $(\ref{fns})_3$ yields
\be\la{fp19}\ba
 c_v\rho\dot\theta - \Delta K(\theta)
 = Q - P\div u
 \ge -\frac{P^2}{4(\mu+\lambda)}
 \ge -c_v a_0\rho\theta^2.
\ea\ee
Observe that $y(t)=(\underline\theta_0^{-1}+a_0t)^{-1}$ solves $y_t+a_0 y^2=0$ with $\nabla y=0$ and $y(0)=\underline\theta_0$. Subtracting $c_v\rho(y_t+u\cdot\nabla y)-\Delta K(y)+c_v a_0\rho y^2=0$ from \eqref{fp19} gives
\be\nonumber\ba
 c_v\rho\bigl((\theta-y)_t+u\cdot\nabla(\theta-y)\bigr)
 -\div\bigl(\kappa(\theta)\nabla(\theta-y)\bigr)
 +c_v a_0\rho(\theta+y)(\theta-y)\ge 0.
\ea\ee
While $\theta>0$, set $z_-\triangleq\min\{\theta-y,0\}$. Multiplying the above inequality by $z_-$, integrating by parts over $\T^2$, and using $(\ref{fns})_1$, we arrive at
\be\la{fp115}\ba
 \frac{c_v}{2}\frac{d}{dt}\int\rho z_-^2\,dx
 +\int\kappa(\theta)|\nabla z_-|^2\,dx
 +c_v a_0\int\rho(\theta+y)z_-^2\,dx\le 0.
\ea\ee
Since $z_-(0)=0$ and $y(T)>0$, integrating \eqref{fp115} in time gives $z_-\equiv 0$.

To establish \eqref{fp01}, we multiply $(\ref{fns})_3$ by $1-\theta^{-1}$, integrate by parts over $\T^2$, add the kinetic energy identity, and use $(\ref{fns})_1$ to rewrite
\be\nonumber\ba
-\int\rho\div u\,dx=\frac{d}{dt}\int(\rho\log\rho-\rho+1)\,dx.
\ea\ee
This gives
\be\la{fp118}\ba
 \frac{d}{dt}\int\left[\frac{1}{2}\rho|u|^2
 +c_v\rho(\theta-\log\theta-1)+\rho\log\rho-\rho+1\right]dx
 +\int\frac{Q}{\theta}\,dx+\mathcal D_\theta(t)=0.
\ea\ee
Since $\theta-\log\theta-1\ge 0$, $(-\log\theta)_+\le(\theta-\log\theta-1)+1$, and $\rho\log(1+\rho)\le(\rho\log\rho-\rho+1)+\rho$, integrating \eqref{fp118} over $(0,T)$ and combining with \eqref{fp15} yields \eqref{fp01}.

By \eqref{fp01}, there exists a constant $C_0>0$ such that
\be\nonumber\ba
 \sup_{0\le t\le T}\int\rho\bigl(\theta+(-\log\theta)_++\log(1+\rho)\bigr)\,dx\le C_0.
\ea\ee
Taking $M_1=4C_0/M_0$, $M_2=e^{4C_0/M_0}-1$, and $h_0=e^{-4C_0/M_0}$, we infer from Chebyshev's inequality that each of the sets $\{\theta>M_1\}$, $\{\rho>M_2\}$, and $\{\theta<h_0\}$ carries at most $M_0/4$ of the total mass. Consequently,
\be\nonumber\ba
 \frac{M_0}{2}\le \int_{\Gamma_t}\rho\,dx\le M_2|\Gamma_t|,\quad
 \frac{M_0}{2}\le \int_{\{\theta\ge h_0,\ \rho\le M_2\}}\rho\,dx\le M_2|\{\theta\ge h_0\}|,
\ea\ee
which proves \eqref{fp27}. For any $a\in(0,h_0)$, consider the nonnegative truncation $w_a\triangleq(h_0-\max\{\theta,a\})_+$. Since $w_a=0$ on $\{\theta\ge h_0\}$, Lemma~\ref{1pi} (with $M=0$) and the Sobolev embedding yield, for every $n\ge 1$,
\be\nonumber\ba
 (h_0-a)^{2n}|\{\theta\le a\}|
 \le \int w_a^{2n}\,dx
 \le C_n\|w_a\|_{H^1}^{2n}
 \le C_n\|\nabla w_a\|_{L^2}^{2n}.
\ea\ee
Noting that $\kappa(\theta)/\theta^2=\kappa_0\theta^{\alpha-2}$ is bounded below by a positive constant depending only on $a$ and $h_0$ on $\{a<\theta<h_0\}$, we have
\be\nonumber\ba
 \|\nabla w_a\|_{L^2}^2
 = \int_{\{a<\theta<h_0\}}|\nabla\theta|^2\,dx
 \le C_a\int\frac{\kappa(\theta)}{\theta^2}|\nabla\theta|^2\,dx
 = C_a\,\mathcal D_\theta(t),
\ea\ee
which together with the previous inequality gives \eqref{extlow:cold}.

Finally, for any $v\in H^1(\T^2)$ with $\|\nabla v\|_{L^2}>0$, the Trudinger--Moser inequality on $\T^2$ \cite{Fontana93} provides a universal constant $c_0>0$ such that
\be\nonumber\ba
 \int\exp\left(\frac{c_0|v-\bar v|^2}{\|\nabla v\|_{L^2}^2}\right)dx \le C.
\ea\ee
Applying the elementary inequality $s_1 s_2\le s_1\log(1+s_1)+e^{s_2}$ ($s_1,s_2\ge 0$) and using \eqref{fp01}, we obtain
\be\nonumber\ba
 \frac{c_0}{\|\nabla v\|_{L^2}^2}\int\rho|v-\bar v|^2\,dx
 \le \int\rho\log(1+\rho)\,dx
 +\int\exp\left(\frac{c_0|v-\bar v|^2}{\|\nabla v\|_{L^2}^2}\right)dx \le C,
\ea\ee
which also holds trivially when $\nabla v\equiv 0$. Taking $v=u$ and using $\int\rho u\,dx=0$ yields
\be\nonumber\ba
 M_0|\bar u| = \left|\int\rho(\bar u-u)\,dx\right|
 \le M_0^{1/2}\left(\int\rho|u-\bar u|^2\,dx\right)^{1/2}
 \le C\|\nabla u\|_{L^2}.
\ea\ee
Therefore $\|u\|_{L^2}\le \|u-\bar u\|_{L^2}+|\bar u|\le C\|\nabla u\|_{L^2}$, and Lemma~\ref{gn1} gives \eqref{ltp-low-u}. This completes the proof of Lemma~\ref{fpl1}.
\end{proof}

Now we fix the coefficient profile $D$: take $D$ from \eqref{new:den} in \textup{(A)}, $D=D_h$ in \textup{(B)}, $D\equiv 1$ in \textup{(C)}, and $D=\Psi\circ J$ in \textup{(D)}, extended constantly to small densities.
In the data-class alternatives \eqref{ext:powerfamily} and \eqref{ext:gainfamily}, the estimates below are first derived on time intervals on which \eqref{ext:iter} and \eqref{new:gain}, respectively, hold with $c_*/2$ on the attained density range; these restrictions are removed by the first-attainment arguments in Proposition~\ref{order:periodic-cap}.
On the attained density range we have $h(r)\ge cr/D(r)$ for $r\ge r_*$, where $D$ is nondecreasing and satisfies $1\le D(r)\le C L(r)^{\beta_D}$ with $0\le\beta_D<\frac{1}{2}$.
Fix a threshold $r_1\ge\max\{2,2r_*\}$. Then for all $r\ge r_1$ in the attained range,
\be\la{extlow:struct}\ba
 \frac{r}{\mu+\lambda(r)}\le CD(r),\quad
 \Lambda(r)-\Lambda(r/2)\ge \frac{cr}{D(r)},\quad
 h(r)\le Cr\Lambda(r)^{1-1/B_h}.
\ea\ee
Indeed, $h\le 2(\mu+\lambda)$ gives the first inequality; integrating $h(z)/z\ge c/D(z)\ge c/D(r)$ over $[r/2,r]$ gives the second; and the third is directly \eqref{visc:envelope}.

For brevity, we set
\be\nonumber\ba
 \mathfrak d(t)\triangleq\|\nabla u(t)\|_{L^2}^2+\int Q\,dx+\mathcal D_\theta(t),\quad
 U\triangleq 1+\int_0^T\mathfrak d(t)\,dt,\quad
 N_p\triangleq 2+\sup_{0\le t\le T}\int_{\{\rho>r_1\}}\rho^p\,dx,
\ea\ee
where $p>2$.

\begin{lemma}\label{low:moment-result}
For each $p>2$, there exist positive constants $C_p$ and $c_p$, independent of $T$ and the provisional density ceiling, such that
\be\la{extlow:NvsU}\ba
 \log N_p\le C_p(1+U)^{c_p}.
\ea\ee
\end{lemma}
\begin{proof}
First, from (\ref{gw}) and $(\ref{fns})_2$, we deduce that
\be\la{fp31}\ba
 \rho\dot u=\nabla G+\mu\nabla^\perp\omega,\quad \Delta G=\div(\rho\dot u).
\ea\ee
Thus,
\be\la{fp33}\ba
G - \ol{G} + \frac{D}{Dt} \left( (- \Delta)^{-1} \div(\n u) \right)
= F,
\ea\ee
where
\be\la{fp34}\ba
F \triangleq \sum_{i,j=1}^{2} [u^i,R_iR_j](\n u^j).
\ea\ee
By $(\ref{fns})_1$, (\ref{visc:growth1}), and (\ref{fp33}), we have
\be\la{fp35}\ba
\frac{D}{Dt} (\Lambda(\n)-\psi) + P = - F - \ol{G},
\ea\ee
where
\be\la{fp36}\ba
\Lambda(\n) \triangleq \int_1^{\n}\frac{2\mu+\lambda(z)}{z}\,dz, \quad
\psi \triangleq (- \Delta)^{-1} \div(\n u).
\ea\ee

Since $\beta_D<\frac{1}{2}$, we may choose $m\in(\frac{1}{2}+\beta_D,1)$, $\xi\in(m,1)$, and $k_*>\max\{8,p,B_h\}$ large enough so that
\be\nonumber\ba
 m+\frac{m+\beta_D}{k_*}<\xi<1.
\ea\ee
In particular, $k>m(k+1)$ for all $k\in[k_*,k_*+1]$. All constants below are uniform for $k\in[k_*,k_*+1]$.
Let $K_0\ge r_1$ (to be fixed below) and $N\triangleq 2+\Lambda(K_0)$. Define
\be\nonumber\ba
 \phi\triangleq(\Lambda(\rho)-\psi-N)_+,\quad
 I_k\triangleq\int\rho\phi^k\,dx,\quad
 M_k\triangleq 1+I_k,\quad
 J_k\triangleq 1+\int_{\{\rho>K_0\}}\rho\Lambda(\rho)^k\,dx.
\ea\ee
Using Sobolev embedding, H\"older's inequality, and \eqref{fp01}, we arrive at
\be\nonumber\ba
 \|\psi\|_{L^{2k}}
 \le C\|\rho u\|_{L^{\frac{2k}{k+1}}}
 \le C\|\sqrt\rho\,u\|_{L^2}\|\rho\|_{L^k}^{1/2}
 \le C\|\rho\|_{L^k}^{1/2}.
\ea\ee
By \eqref{extlow:struct}, we have $\Lambda(\rho)\ge c\rho/D(\rho)$ on $\{\rho>K_0\}$, and hence
\be\nonumber\ba
 \|\rho\|_{L^k}^k + \int_{\{\rho>K_0\}}\frac{\rho^{k+1}}{D(\rho)^k}\,dx \le C J_k.
\ea\ee
Using $\Lambda(\rho)\le N+\phi+|\psi|$ on $\{\rho>K_0\}$, together with H\"older's inequality, we derive
\be\nonumber\ba
 J_k
 &\le C M_k + C\int_{\{\rho>K_0\}}\rho|\psi|^k\,dx
 \le C M_k + C\|\rho\mathbf 1_{\{\rho>K_0\}}\|_{L^2}\|\psi\|_{L^{2k}}^k \\
 &\le C M_k + C\|\rho\|_{L^1}^{\frac{k-2}{2(k-1)}}\|\rho\|_{L^k}^{\frac{k}{2(k-1)}}\|\rho\|_{L^k}^{\frac{k}{2}}
 \le C M_k + C J_k^{\frac{1}{2}+\frac{1}{2(k-1)}},
\ea\ee
which together with Young's inequality yields
\be\la{extlow:moment}\ba
 J_k+\|\rho\|_{L^k}^k+\|\psi\|_{L^{2k}}^{2k}\le C M_k .
\ea\ee

Next, we estimate the commutator term $F$ in \eqref{fp35}. Denote $T_u(g)\triangleq\sum_{i,j=1}^2[u^i,R_iR_j](g^j)$, so that $F=T_u(\rho u)$, and introduce the Orlicz function
\be\nonumber\ba
 B_k(s)\triangleq \int_0^s t^k L(t)^{-m(k+1)}\,dt,\quad s\ge 0.
\ea\ee
Since $k>m(k+1)$, $B_k'$ is strictly increasing and $B_k$ is strictly convex, with
\be\nonumber\ba
 B_k(s)\asymp s^{k+1}L(s)^{-m(k+1)},\quad
 B_k^*(s)\asymp s^{1+\frac{1}{k}}L(s)^{m+\frac{m}{k}},
\ea\ee
where $B_k^*$ denotes the Legendre conjugate of $B_k$. Fix $2<q_-<k_*+1<k_*+2<q_+$. By Lemma~\ref{jhzyl} and $\|u\|_{\mathcal{BMO}}\le C\|\nabla u\|_{L^2}$, we have
\be\nonumber\ba
 \|T_u(g)\|_{L^{q_\pm}}\le C\|\nabla u\|_{L^2}\|g\|_{L^{q_\pm}}.
\ea\ee
If $\|\nabla u\|_{L^2}=0$, then $F=0$. For $\|\nabla u\|_{L^2}>0$, splitting $g=g\mathbf 1_{\{|g|>t\}}+g\mathbf 1_{\{|g|\le t\}}$ at each level $t>0$ and applying Chebyshev's inequality to the $L^{q_-}$ and $L^{q_+}$ bounds yields
\be\nonumber\ba
 \left|\left\{|T_u(g)|>2\|\nabla u\|_{L^2}t\right\}\right|
 \le \frac{C}{t^{q_-}}\int_{\{|g|>t\}}|g|^{q_-}\,dx
 +\frac{C}{t^{q_+}}\int_{\{|g|\le t\}}|g|^{q_+}\,dx.
\ea\ee
Multiplying by $B_k'(t)$, integrating over $t\in(0,\infty)$, and applying Fubini's theorem give
\be\nonumber\ba
 \int B_k\left(\frac{|T_u(g)|}{2\|\nabla u\|_{L^2}}\right)dx
 &\le C\int |g|^{q_-}\!\left(\int_0^{|g|}t^{-q_-}B_k'(t)\,dt\right)dx
 +C\int |g|^{q_+}\!\left(\int_{|g|}^\infty t^{-q_+}B_k'(t)\,dt\right)dx \\
 &\le C\int B_k(|g|)\,dx,
\ea\ee
where the last step uses $q_-<k+1<q_+$, so that
$$\int_0^s t^{-q_-}B_k'(t)\,dt\le Cs^{-q_-}B_k(s),\quad \int_s^\infty t^{-q_+}B_k'(t)\,dt\le Cs^{-q_+}B_k(s).$$
By \eqref{ltp-low-u}, there exists $c_1>0$ such that
$$\int\exp\left( \frac{c_1|u|^2}{\|\nabla u\|_{L^2}^2}\right)\,dx\le C.$$
Choose a large constant $C_1>0$ so that $(k_*+3)/C_1^2\le c_1$.
On $\{|u|\le C_1\|\nabla u\|_{L^2}\sqrt{L(\rho)}\}$, using $m>\frac{1}{2}+\beta_D$ and $D(\rho)\le C L(\rho)^{\beta_D}$, we have
\be\nonumber\ba
 B_k\left(\frac{|\rho u|}{C_1\|\nabla u\|_{L^2}}\right)
 \le B_k\bigl(\rho\sqrt{L(\rho)}\bigr)
 \le C\rho^{k+1}L(\rho)^{\left(\frac{1}{2}-m\right)(k+1)}
 \le C\left(1+\frac{\rho^{k+1}}{D(\rho)^k}\right).
\ea\ee
On the complement $\{|u|>C_1\|\nabla u\|_{L^2}\sqrt{L(\rho)}\}$, we have $\rho\le\exp\bigl(|u|^2/(C_1\|\nabla u\|_{L^2})^2\bigr)$, and hence
\be\nonumber\ba
 B_k\left(\frac{|\rho u|}{C_1\|\nabla u\|_{L^2}}\right)
 \le C\left(\frac{|u|}{C_1\|\nabla u\|_{L^2}}\right)^{k+1}
 \exp\left(\frac{(k+1)|u|^2}{C_1^2\|\nabla u\|_{L^2}^2}\right)
 \le C\exp\left(\frac{c_1|u|^2}{\|\nabla u\|_{L^2}^2}\right).
\ea\ee
Applying the modular bound to $g=\rho u/(C_1\|\nabla u\|_{L^2})$ and using \eqref{extlow:moment}, we obtain
\be\nonumber\ba
 \int B_k\left(\frac{|F|}{2C_1\|\nabla u\|_{L^2}^2}\right)dx
 \le C\int B_k\left(\frac{|\rho u|}{C_1\|\nabla u\|_{L^2}}\right)dx
 \le C M_k.
\ea\ee

Set $z\triangleq\rho\phi^{k-1}$. By Young's inequality for $(B_k,B_k^*)$,
\be\nonumber\ba
 k\int z|F|\,dx
 \le C\|\nabla u\|_{L^2}^2
 \left[\int B_k^*(z)\,dx + \int B_k\left(\frac{|F|}{2C_1\|\nabla u\|_{L^2}^2}\right)dx\right].
\ea\ee
To bound $\int B_k^*(z)\,dx$, we fix a sufficiently large constant $C_2\ge r_1$ and split $\{\phi>0\}$ into two regions.
In the region $\{\rho\le C_2(1+\phi)L(1+\phi)^{\beta_D}\}$, we have $B_k^*(z)\le C\rho$ when $0<\phi<1$, whereas for $\phi\ge 1$,
\be\nonumber\ba
 B_k^*(z)
 \le C\rho\phi^k\left(\frac{\rho}{\phi}\right)^{1/k}L(z)^{m+\frac{m}{k}}
 \le C\rho\phi^k L(\phi)^{m+\frac{m+\beta_D}{k}}
 \le C\rho\phi^k L(\phi)^\xi.
\ea\ee
In the remaining active region $\{\rho>C_2(1+\phi)L(1+\phi)^{\beta_D},\ \phi>0\}$, choosing $C_2$ sufficiently large and using \eqref{extlow:struct} together with $\psi=\Lambda(\rho)-N-\phi$ yield $|\psi|\ge c\rho/L(\rho)^{\beta_D}$, so that
\be\nonumber\ba
 B_k^*(z)
 \le C\rho^{k+1}L(\rho)^{\frac{m(k+1)}{k}}
 \le C\left(\frac{\rho}{L(\rho)^{\beta_D}}\right)^{2k}
 \le C|\psi|^{2k}.
\ea\ee
Combining these estimates with \eqref{extlow:moment} and the elementary bound $C(\log(e+\phi))^\xi\le \delta\log\phi+C(1+\delta^{-\nu_*})$ (for $\phi\ge 1$, $\delta\in(0,1]$, and $\nu_*\triangleq\xi/(1-\xi)$), we arrive at
\be\la{extlow:comm}\ba
 k\int z|F|\,dx
 \le C(1+\delta^{-\nu_*})\|\nabla u\|_{L^2}^2 M_k
 +\delta\|\nabla u\|_{L^2}^2\int_{\{\phi\ge 1\}}\rho\phi^k\log\phi\,dx.
\ea\ee
Moreover, for any measurable function $f$ and $i\in\{1,2\}$, H\"older's inequality and \eqref{extlow:moment} give
\be\la{extlow:dualtest}\ba
 \int z|f|^i\,dx
 \le I_k^{1-\frac{1}{k}}\left(\int_{\{\phi>0\}}\rho|f|^{ik}\,dx\right)^{\frac{1}{k}}
 \le C M_k^{1-\frac{1}{k}+\frac{1}{k^2}}
 \left(\|f\|_{L^{ik}}+\|f\|_{L^{\frac{ik^2}{k-1}}}\right)^i .
\ea\ee

We now differentiate $I_k$ in time. Using $(\ref{fns})_1$ and \eqref{fp35}, we have
\be\nonumber\ba
 \frac{d}{dt}I_k = -k\int z(F+P+\ol{G})\,dx.
\ea\ee
Recall that $\ol{G}=\overline{\lambda\div u}-\bar P$, where $\bar P=\int\rho\theta\,dx\le E_*/c_v$. By Young's inequality, \eqref{extlow:struct}, and \eqref{extlow:moment},
\be\nonumber\ba
 |\overline{\lambda\div u}|
 \le \left(\int h(\rho)\,dx\right)^{1/2}\left(\int Q\,dx\right)^{1/2}
 \le \frac{1}{2}+\frac{1}{2}\left(\int h(\rho)\,dx\right)\int Q\,dx,
\ea\ee
and
\be\nonumber\ba
 \int h(\rho)\,dx
 \le C + C\int_{\{\rho>K_0\}}\rho\Lambda(\rho)^{1-\frac{1}{B_h}}\,dx
 \le C J_k^{\frac{1-1/B_h}{k}}
 \le C M_k^{\frac{1-1/B_h}{k}}.
\ea\ee
Since $\int z\,dx\le M_0^{1/k}I_k^{1-1/k}\le M_0^{1/k}M_k^{1-1/k}$, it follows that
\be\nonumber\ba
 -k\int z(P+\ol{G})\,dx
 \le C M_k\int Q\,dx + k\int z(C_*-P)\,dx,
\ea\ee
where $C_*\triangleq \frac{1}{2}+E_*/c_v$ is independent of $K_0$ and $N$.
We now fix $K_0\triangleq\max\{r_1,\,2C_*/h_0\}$ (and $N=2+\Lambda(K_0)$).
On $\{\rho>K_0,\ \theta>h_0/2\}$, we have $P=\rho\theta>C_*$, so $z(C_*-P)\le 0$.
On $\{\theta\le h_0/2\}$, if $\mathcal D_\theta\le 1$, H\"older's inequality, \eqref{extlow:moment}, and \eqref{extlow:cold} (with $a=h_0/2$ and $n\ge k^2/(k-1)$) yield
\be\nonumber\ba
 \int_{\{\theta\le h_0/2\}} z\,dx
 \le I_k^{1-\frac{1}{k}}\left(\|\rho\|_{L^k}|\{\theta\le h_0/2\}|^{1-\frac{1}{k}}\right)^{\frac{1}{k}}
 \le C M_k\,\mathcal D_\theta(t),
\ea\ee
while if $\mathcal D_\theta>1$ the same bound follows from $\int z\,dx\le C M_k$.
On $\{\rho\le K_0,\ \phi>0\}$, the choice $N=2+\Lambda(K_0)$ implies $1+\phi\le|\psi|$, so that $z=\rho\phi^{k-1}\le\rho\phi^{k-2}|\psi|^2$. Fixing $q\in(2,k_*)$, the Sobolev embedding, elliptic regularity, \eqref{ltp-low-u}, and \eqref{extlow:moment} give
\be\nonumber\ba
 \|\psi\|_{L^\infty}
 \le C\|\rho u\|_{L^q}
 \le C\|\rho\|_{L^k}\|u\|_{L^{\frac{kq}{k-q}}}
 \le C M_k^{\frac{1}{k}}\|\nabla u\|_{L^2},
\ea\ee
and hence
\be\nonumber\ba
 \int_{\{\rho\le K_0\}} z\,dx
 \le \|\psi\|_{L^\infty}^2\int\rho\phi^{k-2}\,dx
 \le C M_k^{1-\frac{2}{k}}\|\psi\|_{L^\infty}^2
 \le C M_k\|\nabla u\|_{L^2}^2.
\ea\ee
Combining the above estimates with \eqref{extlow:comm}, we obtain
\be\la{low:fixed-moment-rate}\ba
 \frac{d}{dt}I_k
 \le C(1+\delta^{-\nu_*})\,\mathfrak d\,M_k
 +\delta\|\nabla u\|_{L^2}^2\int_{\{\phi\ge 1\}}\rho\phi^k\log\phi\,dx.
\ea\ee

Finally, choose the decreasing exponent and parameter
\be\nonumber\ba
 k(t)\triangleq k_*+1-\frac{1}{1+U}\int_0^t\|\nabla u(s)\|_{L^2}^2\,ds \in [k_*,k_*+1],\quad
 \delta\triangleq \frac{1}{2(1+U)}.
\ea\ee
Differentiating $M_{k(t)}(t)=1+\int\rho\phi^{k(t)}\,dx$ and applying \eqref{low:fixed-moment-rate} give
\be\nonumber\ba
 \frac{d}{dt}M_{k(t)}(t)
 &\le C(1+U)^{\nu_*}\,\mathfrak d\,M_{k(t)}(t)
 -\frac{\|\nabla u\|_{L^2}^2}{2(1+U)}\int_{\{\phi\ge 1\}}\rho\phi^{k(t)}\log\phi\,dx \\
 &\quad -\frac{\|\nabla u\|_{L^2}^2}{1+U}\int_{\{0<\phi<1\}}\rho\phi^{k(t)}\log\phi\,dx \\
 & \le C(1+U)^{\nu_*}\,\mathfrak d\,M_{k(t)}(t) + C \| \na u \|_{L^2}^2 \\
 & \le C(1+U)^{\nu_*}\,\mathfrak d\,M_{k(t)}(t).
\ea\ee
Since $M_{k(0)}(0)\le C$, Gr\"onwall's inequality yields
\be\nonumber\ba
 \log M_{k(t)}(t)
 \le \log M_{k(0)}(0) + C(1+U)^{\nu_*}\int_0^t\mathfrak d(s)\,ds
 \le C(1+U)^{1+\nu_*}.
\ea\ee
As $k(t)\ge k_*>p$, \eqref{extlow:moment} and \eqref{extlow:struct} directly imply \eqref{extlow:NvsU}. This completes the proof of Lemma~\ref{low:moment-result}.
\end{proof}

\begin{lemma}\label{low:uniform}
Assume \textup{(A)} with \eqref{ext:iter} on the attained range,
\textup{(B)}, or \textup{(C)} with $\alpha>1-d$ and \eqref{new:gain} on the attained range, allowing the lower-bound constant $c_*/2$ in either data-class alternative. For $0<\alpha\le1$, a further
sufficient hypothesis is \eqref{ext:iter} with $D(r)\le C L(r)^\tau$ together with
\begin{equation}\tag{H1}\label{new:smallentropy}
 \sup_{0\le t\le T}\int_{\{\rho>K_0\}}\rho L(\rho)^q\,dx\le C_q,\quad q>\frac{1+\tau}{\alpha} .
\end{equation}
Then there exists a fixed constant $l\ge 2L_0+2$ such that
\begin{align}
 U+\mathcal P_l&\le C,\tag{H2}\label{new:smalldiss}\\
 \sup_{0\le t\le T}\|\rho(t)\|_{L^p}&\le C_p \quad(2\le p<\infty),\label{extlow:unif}
\end{align}
where
$$\mathcal P_l\triangleq\int_0^T\int_{\{\theta>l\}}\frac{P^2}{\mu+\lambda(\rho)}\,dx\,dt.$$
The constants are independent of $T$ and the provisional density ceiling whenever the bound in \eqref{new:smallentropy}, if invoked, is independent of both.
\end{lemma}

\begin{proof}
Let $0\le\nu<1$ satisfy $a\triangleq\alpha+\nu-1>0$, and let $l\ge 2L_0+2$. Define
\be\nonumber\ba
 f_l\triangleq\mathbf 1_{\{\theta>l\}}\frac{\rho^2\theta^\nu}{\mu+\lambda(\rho)},\quad
 F_l\triangleq 1+\sup_{0\le t\le T}\int f_l\log(e+f_l)\,dx,
\ea\ee
and consider the truncated function $z_l\triangleq\bigl(\theta^{1-\nu/2}-(l/2)^{1-\nu/2}\bigr)_+$.
Since $l/2>M_1$, we have $z_l=0$ on $\Gamma_t$, so Lemma~\ref{1pi} and \eqref{fp27} give
$$\|z_l\|_{H^1}\le C\|\nabla z_l\|_{L^2}.$$
On $\{\theta>l\}$, it holds that $\theta^{2-\nu}\le C z_l^2$, and hence $\frac{P^2}{\mu+\lambda}\le C f_l z_l^2$.
Applying the Trudinger--Moser inequality on $\T^2$ to $z_l/\|\nabla z_l\|_{L^2}$ together with Young's inequality $s_1 s_2\le s_1\log(e+s_1)+e^{s_2}$ yields
\be\nonumber\ba
 \int_{\{\theta>l\}}\frac{P^2}{\mu+\lambda(\rho)}\,dx
 \le C F_l\|\nabla z_l\|_{L^2}^2.
\ea\ee
Splitting the support of $\nabla z_l$ into $\{l/2<\theta\le 2l\}$ and $\{\theta>2l\}$, and noting that on $\{l/2<\theta\le 2l\}$,
\be\nonumber\ba
 |\nabla z_l|^2
 =\left(1-\frac{\nu}{2}\right)^2\theta^{-\nu}|\nabla\theta|^2
 \le \frac{C l^{2-\alpha-\nu}}{\kappa_0}\,\frac{\kappa(\theta)}{\theta^2}|\nabla\theta|^2
 = \frac{C l^{1-a}}{\kappa_0}\,\frac{\kappa(\theta)}{\theta^2}|\nabla\theta|^2,
\ea\ee
we integrate over $(0,T)$ and apply \eqref{fp01} to deduce
\be\nonumber\ba
 \mathcal P_l
 \le C F_l\left(\frac{l^{1-a}}{\kappa_0}
 +\int_0^T\int_{\{\theta>2l\}}\theta^{-\nu}|\nabla\theta|^2\,dx\,dt\right).
\ea\ee

Next, let $\chi_l\in W^{1,\infty}([0,\infty))$ rise linearly from $0$ at $l$ to $1$ at $2l$, with constant extensions $\chi_l=0$ on $[0,l]$ and $\chi_l=1$ on $[2l,\infty)$.
Multiplying $(\ref{fns})_3$ by $\chi_l(\theta)$, integrating over $\T^2\times(0,T)$, and using (\ref{fp01}) and Young's inequality, we derive
\be\nonumber\ba
 & \int_0^T\int\chi_l(\theta)Q\,dx\,dt \\
 & = c_v\int\rho\left(\int_0^\theta\chi_l(s)\,ds\right)dx\bigg|_0^T
 +\int_0^T\int\chi_l'(\theta)\kappa(\theta)|\nabla\theta|^2\,dx\,dt
+\int_0^T\int\chi_l(\theta)P\div u\,dx\,dt \\
& \le C + C l \int_0^T \int \frac{\kappa(\theta)}{\theta^2}|\nabla\theta|^2\,dx\,dt
+ \frac{1}{2} \int_0^T\int\chi_l(\theta)Q\,dx\,dt + \frac{1}{2} \int_0^T\int \mathbf 1_{\{\theta>l\}}\frac{P^2}{\mu+\lambda(\rho)} dx dt \\
& \le C + C l + \frac{1}{2} \int_0^T\int\chi_l(\theta)Q\,dx\,dt + \frac{1}{2} \mathcal P_l,
\ea\ee
which together with (\ref{fp01}) yields
\be\nonumber\ba
 \int_0^T\int Q\,dx\,dt\le C(1+l+\mathcal P_l).
\ea\ee
Similarly, testing $(\ref{fns})_3$ by $-\chi_l(\theta)(l/\theta)^a$ and integrating over $\T^2\times(0,T)$, we note that
$$0\le\int_0^\theta\chi_l(s) \left(\frac{l}{s}\right)^a\,ds\le\theta, \quad P\div u-Q\le \frac{P^2}{4(\mu+\lambda)},$$
that the positive part of $(\chi_l(\theta)(l/\theta)^a)'$ on $(l,2l)$ is at most $l^{-1}$, and that on $\{\theta>2l\}$ one has $-\kappa(\theta)\bigl((l/\theta)^a\bigr)'=\kappa_0 a l^a\theta^{-\nu}$. Thus,
\be\nonumber\ba
 \kappa_0 a l^a\int_0^T\int_{\{\theta>2l\}}\theta^{-\nu}|\nabla\theta|^2\,dx\,dt
 \le C\left(1 + l\int_0^T\mathcal D_\theta(t)\,dt + \mathcal P_l\right)
 \le C(1+l+\mathcal P_l).
\ea\ee
Collecting the three estimates above, we arrive at
\begin{equation}\label{extlow:band}
\begin{aligned}
& \mathcal P_l \le C F_l\left(\frac{l^{1-a}}{\kappa_0}
       +\int_0^T\int_{\{\theta>2l\}}\theta^{-\nu}|\nabla\theta|^2\,dx\,dt\right),\\
& \int_0^T\int_{\{\theta>2l\}}\theta^{-\nu}|\nabla\theta|^2\,dx\,dt
\le \frac{C}{\kappa_0 a l^a}(1+l+\mathcal P_l).
\end{aligned}
\end{equation}
The same estimates hold on every subinterval of $[0,T]$. Substituting the second inequality of \eqref{extlow:band} into the first gives
\be\nonumber\ba
 \mathcal P_l \le \frac{C(1+a^{-1})F_l}{\kappa_0 l^a}(1+l+\mathcal P_l).
\ea\ee
Consequently, using \eqref{fp13} and \eqref{fp01}, we conclude that
\be\nonumber\ba
\frac{C(1+a^{-1})F_l}{\kappa_0 l^a}\le \frac{1}{2},\quad
U+\mathcal P_l\le C(1+l).
\ea\ee

We now close the estimate in each regime.

\smallskip\noindent\textit{Case 1: $\alpha>1$.}
Set $\nu=0$ (so $a=\alpha-1>0$), fix $p>2$, and put $A\triangleq\max\{r_1,N_p^{2/(p-1)}\}$.
On $\{\theta>l,\ \rho\le A\}$, \eqref{extlow:struct} and $f_l\le\rho^2/\mu$ give $f_l\le CD(A)\rho$ and $\log(e+f_l)\le C L(\rho)$, whence
\be\nonumber\ba
 \int_{\{\theta>l,\ \rho\le A\}}f_l\log(e+f_l)\,dx
 \le C D(A)\int\rho L(\rho)\,dx
 \le C D(A).
\ea\ee
On $\{\theta>l,\ \rho>A\}$, since $f_l\le C\rho L(\rho)^{\beta_D}$ and $p>2$, we obtain
\be\nonumber\ba
 \int_{\{\theta>l,\ \rho>A\}}f_l\log(e+f_l)\,dx
 \le C\int_{\{\rho>A\}}\rho^{\frac{p+1}{2}}\,dx
 \le C A^{-\frac{p-1}{2}}\int_{\{\rho>A\}}\rho^p\,dx
 \le C.
\ea\ee
Hence $F_l\le C(1+D(A))$. Choosing $l=2L_0+2+C_*(1+D(A))^{1/(\alpha-1)}$ with $C_*$ large enough so that $C(1+a^{-1})F_l/(\kappa_0 l^{\alpha-1})\le \frac{1}{2}$, and applying \eqref{extlow:NvsU} on each initial subinterval $[0,t]\subset[0,T]$, we obtain for $U(t)\triangleq 1+\int_0^t\mathfrak d(s)\,ds$:
\be\nonumber\ba
 U(t)\le C_2\Bigl(1+D\bigl(\exp(C_1 U(t)^{c_p})\bigr)\Bigr)^{\frac{1}{\alpha-1}},\quad t\in[0,T].
\ea\ee
Note that for every $\gamma>0$, there exist arbitrarily large $Y>0$ such that $D(e^Y)\le Y^\gamma$: this is immediate for $D\equiv 1$ in \textup{(C)}, follows from \eqref{new:fixed} in \textup{(B)}, and holds in \textup{(A)} since $\int_1^\infty\frac{dy}{y D(e^y)}=\infty$ by \eqref{new:den}.
Choose $\gamma>0$ with $c_p\gamma<\alpha-1$, and select such a $Y$ satisfying
\be\nonumber\ba
 \left(\frac{Y}{C_1}\right)^{1/c_p}
 > \max\left\{1,\ C_2(1+Y^\gamma)^{\frac{1}{\alpha-1}}\right\}.
\ea\ee
Since $t\mapsto U(t)$ is continuous with $U(0)=1$, the value $(Y/C_1)^{1/c_p}$ can never be attained, for otherwise $U(t)\le C_2(1+D(e^Y))^{1/(\alpha-1)}\le C_2(1+Y^\gamma)^{1/(\alpha-1)}<(Y/C_1)^{1/c_p}$, a contradiction. Thus $U\le(Y/C_1)^{1/c_p}$ uniformly, which bounds $A$ and fixes a single cutoff $l$ yielding \eqref{new:smalldiss}.

\smallskip\noindent\textit{Case 2: $0<\alpha\le 1$ under \textup{(C)} with $\alpha>1-d$.}
Take $\nu=d\in[0,1)$, so that $a=\alpha+d-1>0$. By \eqref{new:gain} on the attained density range and Young's inequality, on $\{\theta>l\}$ we have $f_l\le C\rho\theta^d L(\rho)^{-d}$ and $\log(e+f_l)\le C(L(\rho)+\log(e+\theta))$, so that
\be\nonumber\ba
 f_l\log(e+f_l)
 &\le C\rho\Bigl(\theta^d L(\rho)^{1-d} + \theta^d L(\rho)^{-d}\log(e+\theta)\Bigr) \\
 &\le C\rho\bigl(1+\theta+L(\rho)\bigr).
\ea\ee
By \eqref{fp01}, this implies $F_l\le C$, and choosing a sufficiently large fixed $l$ gives \eqref{new:smalldiss}.

\smallskip\noindent\textit{Case 3: $0<\alpha\le 1$ under \eqref{ext:iter} and \eqref{new:smallentropy}.}
Since $q>(1+\tau)/\alpha$, we may choose $\nu\in\bigl(1-\alpha,\,1-(1+\tau)/q\bigr)$ and $\varepsilon>0$ small enough so that $\nu+\varepsilon<1$ and $\tau/(1-\nu-\varepsilon)<q$. Then $a=\alpha+\nu-1>0$, and on $\{\theta>l\}$,
\be\nonumber\ba
 f_l\log(e+f_l)
 &\le C\rho\Bigl(\theta^\nu L(\rho)^{1+\tau} + \theta^\nu L(\rho)^\tau\log(e+\theta)\Bigr) \\
 &\le C\rho\Bigl(\theta^\nu L(\rho)^{1+\tau} + \theta^{\nu+\varepsilon} L(\rho)^\tau\Bigr) \\
 &\le C\rho\Bigl(1+\theta+L(\rho)^{\frac{1+\tau}{1-\nu}}+L(\rho)^{\frac{\tau}{1-\nu-\varepsilon}}\Bigr)
 \le C\rho\bigl(1+\theta+L(\rho)^q\bigr),
\ea\ee
where we applied Young's inequality in the penultimate step. Integrating over $\T^2$ and using \eqref{fp01} together with \eqref{new:smallentropy} yield $F_l\le C$, so a sufficiently large fixed $l$ again gives \eqref{new:smalldiss}.

Finally, in all three cases, substituting \eqref{new:smalldiss} into \eqref{extlow:NvsU} bounds $N_p$ for every $p>2$, which together with the mass bound on $\{\rho\le r_1\}$ and H\"older's inequality yields \eqref{extlow:unif}. This completes the proof of Lemma~\ref{low:uniform}.
\end{proof}

We next derive finite-time dissipation bounds in the borderline regimes \textup{(D)} and \textup{(C)}, and then deduce high-temperature tail estimates and a lower temperature bound.

\begin{lemma}\label{lem:linear-moment-input}
Under \textup{(D)}, for each $T_0<\infty$, there exists a common cutoff $l\ge 2L_0+2$ such that
\be\nonumber\ba
 U+\mathcal P_l\le C_{T_0}
\ea\ee
on every $[0,T]\subset[0,T_0]$, independently of the initial density floor and the provisional density ceiling. Moreover, for any fixed $R_0>r_*$,
\be\nonumber\ba
 \sup_{0\le t\le T}\int_{\{\rho>2R_0\}}\frac{\rho L(\rho)^2}{D(\rho)}\,dx\le C_{T_0},
\ea\ee
and the profile $D=\Psi\circ J$ satisfies
\be\nonumber\ba
 \lim_{r\to\infty}\frac{\log D(r)}{\log L(r)}=0,\quad
 \int_{r_*}^\infty\frac{dr}{(e+r)L(r)D(r)}=\infty .
\ea\ee
\end{lemma}

\begin{proof}
Write $D(r)=\Psi(J(r))$. By \eqref{crit:profile}, integrating $\Psi'(s)/\Psi(s)\le p/s$ over $[x,tx]$ yields $\Psi(tx)\le t^p\Psi(x)$ for all $t,x\ge 1$. In particular, $D(r)\le C J(r)^p$ and, for $r\ge r_*$,
\be\nonumber\ba
 0\le \frac{rD'(r)}{D(r)}
 =\frac{r J'(r)\Psi'(J(r))}{\Psi(J(r))}
 \le \frac{C}{L(r)J(r)}.
\ea\ee
Set
\be\nonumber\ba
 e_u(t)\triangleq\int\frac{Q}{\theta}\,dx,\quad
 e_\theta(t)\triangleq\kappa_0\int\frac{|\nabla\theta|^2}{\theta}\,dx,\quad
 b(t)\triangleq\bigl[e_u(t)(1+e_\theta(t))\bigr]^{1/2}.
\ea\ee
Since $\alpha=1$, \eqref{fp01} and Young's inequality give
\be\nonumber\ba
 \int_0^T\bigl(e_u(t)+e_\theta(t)\bigr)\,dt\le C,\quad
 \int_0^T b(t)\,dt\le C_{T_0}.
\ea\ee

Choose a smooth nondecreasing cutoff $\chi\in C^\infty([0,\infty))$ satisfying $\chi=0$ on $[0,R_0]$ and $\chi=1$ on $[2R_0,\infty)$ for a fixed $R_0>r_*$. For $m\ge 0$, set $\sigma_m\triangleq 2^{-m}$, $w_m(r)\triangleq \chi(r)L(r)^{1+\sigma_m}/D(r)$, and
\be\nonumber\ba
 H_m(t)\triangleq 1+\int\rho w_m(\rho)\,dx.
\ea\ee
Note that $H_m(0)\le C$ uniformly in $m\ge 0$. We claim that
\be\la{crit:moment-rate}\ba
 \left|\frac{d}{dt}\sqrt{H_0(t)}\right|\le C\left(\int Q\,dx\right)^{1/2},\quad
 |H_m'(t)|\le C b(t)\sqrt{H_{m-1}(t)}\quad (m\ge 1).
\ea\ee
Indeed, using $(\ref{fns})_1$, we have
\be\nonumber\ba
 H_m'(t)=-\int\rho^2 w_m'(\rho)\div u\,dx,
\ea\ee
with $|r w_m'(r)|\le C L(r)^{\sigma_m}/D(r)$ for $r\ge 2R_0$ and $|r w_m'(r)|\le C$ for $R_0<r<2R_0$.
For $m=0$, H\"older's inequality, $Q\ge(\mu+\lambda)(\div u)^2$, and \eqref{crit:viscosity} give
\be\nonumber\ba
 |H_0'(t)|
 \le\left(\int\frac{\rho^4|w_0'(\rho)|^2}{\mu+\lambda(\rho)}\,dx\right)^{1/2}
 \left(\int Q\,dx\right)^{1/2}
 \le C\sqrt{H_0(t)}\left(\int Q\,dx\right)^{1/2},
\ea\ee
which proves the first bound in \eqref{crit:moment-rate}.
For $m\ge 1$, applying H\"older's inequality yields
\be\nonumber\ba
 |H_m'(t)|
 \le\left(\int\frac{\rho^4\theta|w_m'(\rho)|^2}{\mu+\lambda(\rho)}\,dx\right)^{1/2} e_u(t)^{1/2}
 \le C\left(1+\int_{\{\rho\ge 2R_0\}}\frac{\rho\theta L(\rho)^{2\sigma_m}}{D(\rho)}\,dx\right)^{1/2} e_u(t)^{1/2},
\ea\ee
where we used \eqref{fp01} on $\{R_0<\rho<2R_0\}$.
Fix a large constant $A>2L_0$. On $\{\rho\ge 2R_0,\ \theta\le A(1+e_\theta)L(\rho)\}$, since $2\sigma_m=\sigma_{m-1}$, we have
\be\nonumber\ba
 \int_{\{\rho\ge 2R_0,\ \theta\le A(1+e_\theta)L(\rho)\}}
 \frac{\rho\theta L(\rho)^{2\sigma_m}}{D(\rho)}\,dx
 \le A(1+e_\theta)\int\rho w_{m-1}(\rho)\,dx
 \le C(1+e_\theta)H_{m-1}.
\ea\ee
On the complementary set $\{\rho\ge 2R_0,\ \theta>A(1+e_\theta)L(\rho)\}$, we have $\theta>2L_0$ and $\rho\le\exp\bigl(\theta/[A(1+e_\theta)]\bigr)$, so that
\be\nonumber\ba
 \frac{\rho\theta L(\rho)^{2\sigma_m}}{D(\rho)}
 \le C(1+e_\theta)\exp\left(\frac{2\theta}{A(1+e_\theta)}\right).
\ea\ee
Since $(\sqrt\theta-\sqrt{L_0})_+$ vanishes on $\Gamma_t$ and satisfies $4\kappa_0\|\nabla\sqrt\theta\|_{L^2}^2=e_\theta$, Lemma~\ref{1pi} and the Trudinger--Moser inequality give a small constant $c>0$ such that
\be\nonumber\ba
 \int_{\{\theta>2L_0\}}\exp\left(\frac{c\theta}{1+e_\theta}\right)dx\le C.
\ea\ee
Choosing $A>\max\{2L_0,2/c\}$ bounds the second region by $C(1+e_\theta)$, which establishes \eqref{crit:moment-rate}.

Fix $M_{\rm in}>\max\{2,H_0(0)\}$ uniformly over the data, and set $M_j\triangleq 2^j M_{\rm in}$ for $j\ge 1$. Let $t_j$ be the first time at which $H_0$ attains $M_j$ (with $t_0\triangleq 0$), and put
\be\nonumber\ba
 I_j\triangleq [t_{j-1},t_j],\quad
 B_j\triangleq\int_{I_j}b(t)\,dt,\quad
 e_j\triangleq\int_{I_j}\bigl(e_u(t)+e_\theta(t)\bigr)\,dt.
\ea\ee
For each attained $j\ge 1$, select an integer $n_j\ge 1$ so that
\be\nonumber\ba
 \delta_j\triangleq 2^{-n_j}\in\left[\frac{1}{4\log(e+j)},\,\frac{1}{2\log(e+j)}\right],
\ea\ee
and set $K_j\triangleq\sup_{0\le t\le t_j}H_{n_j}(t)$. To compare $H_{n_j}$ with $H_0$, define
\be\nonumber\ba
 a_j\triangleq\frac{\delta_j}{2(1-\delta_j)},\quad
 M^*(t)\triangleq\max\left\{M_{\rm in},\,\sup_{0\le s\le t}H_0(s)\right\},\quad
 \mathcal V_j(t)\triangleq 1+\int_0^t b(s)M^*(s)^{a_j}\,ds.
\ea\ee
With $\gamma_k\triangleq 2(1-\sigma_k)$, we claim that
\be\la{crit:moment-induction}\ba
 H_k(t)\le C\,M^*(t)^{\sigma_k-a_j\gamma_k}\,\mathcal V_j(t)^{\gamma_k},
 \quad 0\le k\le n_j,\quad 0\le t\le t_j .
\ea\ee
For $k=0$, we have $\sigma_0=1$ and $\gamma_0=0$, so \eqref{crit:moment-induction} holds with constant $1$. Assuming \eqref{crit:moment-induction} holds for $k-1$ with constant $C_{k-1}$ ($1\le k\le n_j$), we substitute it into \eqref{crit:moment-rate} and use $\sigma_k=\sigma_{k-1}/2$, $\gamma_k=1+\gamma_{k-1}/2$, and $\mathcal V_j'(t)=b(t)M^*(t)^{a_j}$ to obtain
\be\nonumber\ba
 |H_k'(t)|
 &\le C\sqrt{C_{k-1}}\,b(t)\,M^*(t)^{\sigma_k-a_j(\gamma_k-1)}\,\mathcal V_j(t)^{\gamma_k-1} \\
 &= C\sqrt{C_{k-1}}\,M^*(t)^{\sigma_k-a_j\gamma_k}\,\mathcal V_j(t)^{\gamma_k-1}\mathcal V_j'(t).
\ea\ee
Since $k\le n_j$ implies $\sigma_k-a_j\gamma_k=(\sigma_k-\delta_j)/(1-\delta_j)\ge 0$ and $M^*$ is nondecreasing, integrating over $[0,t]$ yields
\be\nonumber\ba
 H_k(t)
 \le H_k(0)+\frac{C\sqrt{C_{k-1}}}{\gamma_k}M^*(t)^{\sigma_k-a_j\gamma_k}\mathcal V_j(t)^{\gamma_k}
 \le \bigl(C+C\sqrt{C_{k-1}}\bigr)M^*(t)^{\sigma_k-a_j\gamma_k}\mathcal V_j(t)^{\gamma_k}.
\ea\ee
Hence the constants are uniformly bounded in $k$ and $j$, proving \eqref{crit:moment-induction}.

At $k=n_j$, we have $\sigma_{n_j}-a_j\gamma_{n_j}=0$ and $a_j\gamma_{n_j}=\delta_j$. Since $M^*(s)\le M_i=2^{i-j}M_j$ on $I_i$ ($1\le i\le j$) and $a_j\ge\delta_j/2$, dividing $\mathcal V_j(t_j)$ by $M_j^{a_j}$ gives
\be\nonumber\ba
 \frac{\mathcal V_j(t_j)}{M_j^{a_j}}
 \le M_j^{-a_j}+\sum_{i=1}^j 2^{-a_j(j-i)}B_i
 \le M_j^{-a_j}+\sum_{i=1}^j e^{-c\delta_j(j-i)}B_i
 \triangleq h_j,
\ea\ee
and therefore \eqref{crit:moment-induction} yields
\be\la{crit:moment-record}\ba
\frac{K_j}{M_j^{\delta_j}} \le C h_j^{2(1-\delta_j)},\quad
 h_j=M_j^{-a_j}+\sum_{i=1}^j e^{-c\delta_j(j-i)}B_i,
\ea\ee
and
\be\la{crit:kernel-sum}\ba
\sum_{j=1}^N\delta_j h_j \le C\left(1+\int_0^{t_N}b(t)\,dt\right).
\ea\ee
Indeed, $M_j^{-a_j}\le e^{-c\delta_j j}$ is summable over $j\ge 1$. For each fixed $i\ge 1$, the kernel sum over $i\le j\le 2i+4$ satisfies $\sum_{j=i}^{2i+4}\delta_j e^{-c\delta_j(j-i)}\le C$ as $\delta_j\asymp\delta_i$, while for $j>2i+4$ we have $j-i>j/2$ and thus $\delta_j e^{-c\delta_j(j-i)}\le e^{-c' j/\log(e+j)}$, which is summable in $j$. Summing over $i$ yields \eqref{crit:kernel-sum}.

We now estimate the dissipation on $I_j$ via dyadic temperature bands. Fix a large threshold $\ell_0>8\max\{R_0,L_0,1\}$ and write $\delta=\delta_j$, $K=K_j$, and
\be\nonumber\ba
 \ell\triangleq\max\left\{\ell_0,\ M_j^{1/2},\ K^{1/\delta}\Psi(e+j)^{(2+\delta)/\delta}\right\}.
\ea\ee
Interpolating $L^{1+\delta}/D$ between $L/D$ and $L^2/D$ with respect to $\chi\rho\,dx$ and using \eqref{fp01} give $K\le C M_j^\delta$. Since $\log\Psi(e+j)\le C\log(e+j)$ and $\delta^{-1}\le 4\log(e+j)$, it follows that
\be\la{crit:cutoff-size}\ba
 c(1+j)\le\log(e+\ell)\le\log(e+\ell K)\le C(1+j).
\ea\ee
For $k\ge -2$, define $v_k\triangleq 2^k\ell$, $A_k\triangleq\{v_k<\theta\le 2v_k\}$, and
\be\nonumber\ba
 G_k\triangleq\int_{I_j}\int_{A_k}|\nabla\theta|^2\,dx\,dt,\quad
 P_k\triangleq\int_{I_j}\int_{A_k}\frac{P^2}{\mu+\lambda(\rho)}\,dx\,dt.
\ea\ee
For $k\ge -1$, set $f_k\triangleq\mathbf 1_{A_k}\rho^2/(\mu+\lambda(\rho))$ and $g_k\triangleq\sup_{t\in I_j}\|f_k(t)\|_{L\log L}$. We claim that
\be\la{low:band-pressure}\ba
 P_k\le C g_k G_{k-1}\quad (k\ge -1),\quad
 g_{-1}\le C,\quad
 \sum_{k\ge 0}g_k\le \frac{C}{\delta}.
\ea\ee
To verify the first inequality, let $w_k\triangleq(\min\{\theta,v_k\}-v_{k-1})_+$. Since $v_{k-1}\ge\ell/4>M_1$, $w_k$ vanishes on $\Gamma_t$; moreover, $w_k=v_k/2\ge\theta/4$ on $A_k$ and $\nabla w_k=\nabla\theta\,\mathbf 1_{A_{k-1}}$. Applying Orlicz--H\"older together with Lemma~\ref{1pi} and the Trudinger--Moser inequality on $\T^2$ gives
\be\nonumber\ba
 \int_{A_k}\frac{P^2}{\mu+\lambda(\rho)}\,dx
 \le 16\int f_k w_k^2\,dx
 \le C\|f_k\|_{L\log L}\|\nabla w_k\|_{L^2}^2
 \le C g_k\int_{A_{k-1}}|\nabla\theta|^2\,dx,
\ea\ee
and integrating over $I_j$ yields $P_k\le C g_k G_{k-1}$.

To bound $g_k$ for $k\ge -1$, set $r\triangleq 1/(1+\delta)$, $s\triangleq 1+r$, and
\be\nonumber\ba
 u_k\triangleq v_k^{-1/3}+K^r v_k^{r-1}\Psi\bigl(\log(e+v_k K)\bigr)^s.
\ea\ee
By \eqref{fp15}, $\int_{A_k}\rho\,dx\le C/v_k$. On $A_k\cap\{\rho\le v_k\}$, since $f_k\le C\rho D(v_k)$ and $u_k\ge v_k^{-1/3}$, we have
\be\nonumber\ba
 \int_{A_k\cap\{\rho\le v_k\}} f_k\log\left(e+\frac{f_k}{u_k}\right)dx
 \le C D(v_k)L(v_k)\int_{A_k}\rho\,dx
 \le C v_k^{-1}L(v_k)^2\le C u_k.
\ea\ee
On $A_k\cap\{\rho>v_k\}$, we have $f_k\log(e+f_k/u_k)\le C\rho D(\rho)L(\rho)$. Writing $Z\triangleq L(\rho)^{1+\delta}/D(\rho)$, we note that
$$D(\rho)L(\rho)=Z^r D(\rho)^s,\quad \int_{\{\rho>v_k\}}\rho Z\,dx\le K,\quad D(\rho)\le C\Psi(\log(e+Z)).$$
Splitting this region at $Z=v_k K$ gives
\be\nonumber\ba
 \int_{A_k\cap\{\rho>v_k,\ Z\le v_k K\}} \rho Z^r\Psi\bigl(\log(e+Z)\bigr)^s\,dx
 &\le (v_k K)^r\Psi\bigl(\log(e+v_k K)\bigr)^s\int_{A_k}\rho\,dx \\
 &\le C K^r v_k^{r-1}\Psi\bigl(\log(e+v_k K)\bigr)^s\le C u_k.
\ea\ee
On $\{Z>v_k K\}$, the logarithmic derivative of $z\mapsto z^{r-1}\Psi(\log(e+z))^s$ satisfies
\be\nonumber\ba
 \frac{d}{dz}\log\Bigl(z^{r-1}\Psi\bigl(\log(e+z)\bigr)^s\Bigr)
 \le \frac{1}{z}\left(r-1+\frac{sp}{\log(e+z)}\right)\le 0,
\ea\ee
provided $\ell_0$ is chosen sufficiently large, because $1-r\ge\delta/2\ge[8\log(e+j)]^{-1}$ while $\log(e+v_k K)\ge c(1+j)$ by \eqref{crit:cutoff-size}. Hence
\be\nonumber\ba
 \int_{A_k\cap\{\rho>v_k,\ Z>v_k K\}} \rho Z^r\Psi\bigl(\log(e+Z)\bigr)^s\,dx
 \le (v_k K)^{r-1}\Psi\bigl(\log(e+v_k K)\bigr)^s\int_{\{\rho>v_k\}}\rho Z\,dx
 \le u_k.
\ea\ee
Combining the above estimates shows that $g_k\le C u_k$. By \eqref{crit:cutoff-size}, the definition of $\ell$ (which ensures $K^r\ell^{r-1}\le\Psi(e+j)^{-s}$), and $\Psi(tx)\le t^p\Psi(x)$, we obtain
\be\nonumber\ba
 u_k
 \le v_k^{-1/3}+C 2^{-k(1-r)}\left(1+\frac{k}{e+j}\right)^{sp}
 \le v_k^{-1/3}+C 2^{-k\delta/2}\left(1+\frac{k}{e+j}\right)^{2p}\quad (k\ge 0),
\ea\ee
together with $u_{-1}\le C$. Summing over $k\ge 0$ proves \eqref{low:band-pressure}.

For each $k\ge 0$, let $b_k(\theta)$ be the continuous hat function supported on $[v_{k-1},v_{k+1}]$, piecewise linear in $\log\theta$, with $b_k(v_k)=1$. Then $\kappa(\theta)b_k'(\theta)=\kappa_0/\log 2$ on $A_{k-1}$ and $-\kappa_0/\log 2$ on $A_k$. Multiplying $(\ref{fns})_3$ by $-b_k(\theta)$, integrating over $\T^2\times I_j$, and using Young's inequality $b_k|P\div u|\le \frac{1}{2}b_k Q+C\mathbf 1_{A_{k-1}\cup A_k}P^2/(\mu+\lambda)$, we arrive at
\be\nonumber\ba
 \frac{\kappa_0}{\log 2}(G_k-G_{k-1})
 +\frac{1}{2}\int_{I_j}\int b_k(\theta)Q\,dx\,dt
 \le a_k+C(P_{k-1}+P_k),
\ea\ee
where $a_k\triangleq c_v\bigl|\int\rho\bigl(\int_0^\theta b_k(s)\,ds\bigr)\,dx\big|_{t_{j-1}}^{t_j}\bigr|$ satisfies $\sum_{k\ge 0}a_k\le C$ by $\sum_{k\ge 0}b_k\le 1$ and \eqref{fp15}.
Since $G_{-2}+G_{-1}\le C\ell e_j$ by \eqref{fp01}, substituting \eqref{low:band-pressure} shows that the running maximum $\ol G_k\triangleq\max_{-2\le i\le k}G_i$ satisfies
\be\nonumber\ba
 \ol G_k\le \bigl(1+C(g_{k-1}+g_k)\bigr)\ol G_{k-1}+C a_k\quad (k\ge 0).
\ea\ee
Iterating in $k$ and using \eqref{low:band-pressure} yield $\sup_{k\ge -1}G_k\le C\exp(C/\delta)(1+\ell e_j)$, and hence
\be\nonumber\ba
 \int_{I_j}\int_{\{\theta>\ell\}}\frac{P^2}{\mu+\lambda(\rho)}\,dx\,dt
 =\sum_{k\ge 0}P_k
 \le C\left(\sum_{k\ge 0}g_k\right)\sup_{k\ge -1}G_k
 \le C\exp(C/\delta)(1+\ell e_j).
\ea\ee
Testing $(\ref{fns})_3$ by $\chi_\ell(\theta)$ on $I_j$ as in the proof of Lemma~\ref{low:uniform} (with transition cost $C\ell e_j$ and $\int_{I_j}\int_{\{\theta\le 2\ell\}}Q\,dx\,dt\le 2\ell e_j$), we deduce
\be\la{crit:heat-window}\ba
 \int_{I_j}\int Q\,dx\,dt\le \exp(C/\delta)(1+\ell e_j).
\ea\ee

Finally, integrating the first inequality of \eqref{crit:moment-rate} over $I_j$ and using $H_0(t_j)=M_j=2H_0(t_{j-1})$ give
\be\nonumber\ba
 c M_j\le |I_j|\int_{I_j}\int Q\,dx\,dt
 \le |I_j|\exp(C/\delta_j)(1+\ell e_j).
\ea\ee
Because $\delta_j^{-1}\le 4\log(e+j)$ and $M_j=2^j M_{\rm in}$, the contributions of $1$ and $\max\{\ell_0,M_j^{1/2}\}e_j$ after dividing by $M_j$ are $O_{T_0}\bigl(e^{C/\delta_j}M_j^{-1/2}\bigr)\to 0$ as $j\to\infty$. Thus, for all sufficiently large attained $j$,
\be\nonumber\ba
 1\le C\exp(C/\delta_j)|I_j|e_j\left(\frac{K_j}{M_j^{\delta_j}}\right)^{\frac{1}{\delta_j}}\Psi(e+j)^{\frac{2+\delta_j}{\delta_j}}.
\ea\ee
Raising both sides to the power $\delta_j/2$ and substituting \eqref{crit:moment-record} yield
\be\nonumber\ba
 \Psi(e+j)^{-1-\frac{\delta_j}{2}}
 \le C |I_j|^{\frac{\delta_j}{2}} e_j^{\frac{\delta_j}{2}} h_j^{1-\delta_j}.
\ea\ee
Multiplying by $\delta_j$ and applying Young's inequality, we obtain
\be\nonumber\ba
 \delta_j\Psi(e+j)^{-1-\frac{\delta_j}{2}}
 \le C\bigl(|I_j|+e_j+\delta_j h_j\bigr).
\ea\ee
Summing over all large attained $j$ and using \eqref{crit:kernel-sum}, $\sum|I_j|\le T_0$, and $\sum e_j\le C$, we conclude that every attained partial sum of $\sum_j\delta_j\Psi(e+j)^{-1-\delta_j/2}$ is bounded by $C_{T_0}$. However, since $\Psi(e+j)^{\delta_j/2}\le C$, the general term is comparable to $[\Psi(e+j)\log(e+j)]^{-1}$, whose series diverges by \eqref{crit:profile}. Consequently, the number of attained levels $j$ is bounded in terms of $T_0$, which implies
\be\nonumber\ba
 \sup_{0\le t\le T}\int_{\{\rho>2R_0\}}\frac{\rho L(\rho)^2}{D(\rho)}\,dx
 \le \sup_{0\le t\le T}H_0(t)\le C_{T_0}.
\ea\ee
Since $D\le C J^p\le C L^{1/4}$, this yields
$$\sup_{0\le t\le T}\int_{\{\rho>2R_0\}}\rho L(\rho)^{3/2}\,dx\le C_{T_0}.$$
Applying Lemma~\ref{low:uniform} (\eqref{new:smallentropy}--\eqref{new:smalldiss}) with $q=3/2$ and $\tau=1/4$ gives $U+\mathcal P_l\le C_{T_0}$. Finally, $D\le C J^p=L^{o(1)}$ gives $\log D=o(\log L)$, and the change of variable $x=J(r)$ with $J'(r)\asymp[(e+r)L(r)]^{-1}$ yields $\int_{r_*}^\infty\frac{dr}{(e+r)L(r)D(r)}\asymp\int_{J(r_*)}^\infty\frac{dx}{\Psi(x)}=\infty$ by \eqref{crit:profile}. This completes the proof of Lemma~\ref{lem:linear-moment-input}.
\end{proof}

\begin{lemma}\label{lem:gain-moment-input}
Under the fixed-law alternative of \textup{(C)} with $\alpha\ge\frac{2(1-d)}{2+d}$, for each fixed $T_0<\infty$, there exists a common cutoff $l\ge 2L_0+2$ such that
\be\nonumber\ba
 U+\mathcal P_l\le C_{T_0}
\ea\ee
on every $[0,T]\subset[0,T_0]$, independently of the initial density floor and the provisional density ceiling.
\end{lemma}

\begin{proof}
Set $s\triangleq\frac{1-d}{2+d}\in(0,\frac{1}{2}]$, $\nu\triangleq 1-s$, and $a\triangleq\alpha+\nu-1=\alpha-s\ge s>0$. Fix $R_0>r_*$ and define
\be\nonumber\ba
 \varphi(r)\triangleq\bigl(L(r)^{1+\frac{d}{2}}-L(R_0)^{1+\frac{d}{2}}\bigr)_+,\quad
 \mathcal M(t)\triangleq\left(1+\int\rho\varphi(\rho)^2\,dx\right)^{1/2}.
\ea\ee
For any subinterval $I\subset[0,T]$, denote
\be\nonumber\ba
 \mathcal M_I\triangleq\sup_{t\in I}\mathcal M(t),\quad
 Q_I\triangleq\int_I\int Q\,dx\,dt,\quad
 \mathcal P_l(I)\triangleq\int_I\int_{\{\theta>l\}}\frac{P^2}{\mu+\lambda(\rho)}\,dx\,dt.
\ea\ee
By \eqref{new:gain}, we have $r|r\varphi'(r)|^2\le C r L(r)^d\le C(\mu+\lambda(r))$ for $r>R_0$. Differentiating $\mathcal M^2$ via $(\ref{fns})_1$ and applying H\"older's inequality yield
\be\nonumber\ba
 |\mathcal M'(t)|
& =\frac{1}{\mathcal M(t)}\left|\int\rho^2\varphi(\rho)\varphi'(\rho)\div u\,dx\right| \\
& \le \frac{1}{\mathcal M(t)}\left(\int\rho\varphi(\rho)^2\,dx\right)^{1/2}
 \left(\int\rho|\rho\varphi'(\rho)|^2(\div u)^2\,dx\right)^{1/2}
 \le C\left(\int Q\,dx\right)^{1/2}.
\ea\ee
Next, on $\{\theta>l\}$, let $f_l=\mathbf 1_{\{\theta>l\}}\rho^2\theta^\nu/(\mu+\lambda(\rho))$ be as in Lemma~\ref{low:uniform}. By \eqref{new:gain} and $\nu<1$, we have
\be\nonumber\ba
 f_l\log(e+f_l)
 \le C\rho\Bigl(\theta^\nu L(\rho)^{1-d}+\theta^\nu\log(e+\theta)\Bigr)
 \le C\rho\bigl(1+\theta+\theta^\nu L(\rho)^{1-d}\bigr).
\ea\ee
Applying H\"older's inequality and \eqref{fp01}, we obtain, for $t\in I$,
\be\nonumber\ba
 \int f_l\log(e+f_l)\,dx
 \le C+C\left(\int_{\{\theta>l\}}\rho\theta\,dx\right)^\nu
 \left(\int_{\{\theta>l\}}\rho L(\rho)^{2+d}\,dx\right)^s
 \le C\mathcal M_I^{2s},
\ea\ee
and hence $F_l\le C\mathcal M_I^{2s}$ on $I$. Taking $l=C_*\mathcal M_I^2$ in \eqref{extlow:band} on the interval $I$ with $C_*\ge 2L_0+2$ sufficiently large, since $a\ge s$ and $\mathcal M_I\ge 1$, we have
\be\nonumber\ba
 \frac{C(1+a^{-1})F_l}{\kappa_0 l^a}
 \le C\mathcal M_I^{2s}(C_*\mathcal M_I^2)^{-a}
 \le C C_*^{-a}\le \frac{1}{2},
\ea\ee
which together with the heat-test estimates in Lemma~\ref{low:uniform} yields
\be\la{low:gain-window}\ba
 |\mathcal M'(t)|\le C\left(\int Q\,dx\right)^{1/2},\quad
 Q_I+\mathcal P_{C_*\mathcal M_I^2}(I)\le C\mathcal M_I^2 .
\ea\ee

On any subinterval $I=[t_0,t_1]\subset[0,T]$, integrating the first inequality in \eqref{low:gain-window} and applying H\"older's inequality together with the second inequality give
\be\nonumber\ba
 \mathcal M_I
 \le \mathcal M(t_0)+C|I|^{1/2}Q_I^{1/2}
 \le \mathcal M(t_0)+C_0|I|^{1/2}\mathcal M_I.
\ea\ee
Taking $|I|\le (2C_0)^{-2}$ yields $\mathcal M_I\le 2\mathcal M(t_0)$. Covering $[0,T_0]$ with finitely many such subintervals bounds $\sup_{0\le t\le T}\mathcal M(t)\le C_{T_0}$. Applying \eqref{low:gain-window} on $[0,T]$ and enlarging the cutoff to $l=C_* C_{T_0}^2$, together with \eqref{fp01} and \eqref{fp13}, we obtain $U+\mathcal P_l\le C_{T_0}$. This completes the proof of Lemma~\ref{lem:gain-moment-input}.
\end{proof}

\begin{lemma}\label{low:temp-floor}
Assume that $U+\mathcal P_l\le C$ for a fixed $l\ge 2L_0+2$, and set $E\triangleq\{\theta>4l\}$. For any $q\ge 2$, $0<\epsilon\le (1+\alpha)/2$, and $0<\delta\le 1/2$, one has
\be\la{extlow:heattail}\ba
 \int_0^T\|\theta\mathbf 1_E\|_{L^q}^{1+\alpha-\epsilon}\,dt\le \frac{Cq}{\epsilon},\quad
 \int_0^T\|K(\theta)\mathbf 1_E\|_{L^q}^{1-\delta}\,dt\le \frac{Cq}{\delta} .
\ea\ee
With
$$H\triangleq c_v\int\rho\int_0^\theta K(v)\,dv\,dx,\quad H_E\triangleq c_v\int_E\rho\int_0^\theta K(v)\,dv\,dx,$$
for all sufficiently small $\zeta\in(0,\zeta_0]$ and every interval $I\subset[0,T]$ with $|I|\le 1$,
\be\la{extlow:potentialtail}\ba
 \int_0^T H_E^{1-\zeta}\,dt+\int_I(1+H)^{1-\zeta}\,dt\le C\zeta^{-2}.
\ea\ee
Moreover, for every $R\ge\max\bigl\{r_1,\sup_{0\le t\le T}\|\rho(t)\|_{L^\infty}\bigr\}$, it holds that
\be\la{extlow:templower}\ba
 \inf_{\T^2\times[0,T]}\theta\ge \frac{c}{1+D(R)}.
\ea\ee
The constants are independent of $T$ under \eqref{new:smalldiss} with a time-uniform bound, and may depend on $T_0$ under the finite-horizon bounds of Lemmas~\ref{lem:linear-moment-input} and~\ref{lem:gain-moment-input}.
\end{lemma}

\begin{proof}
Set $c_0\triangleq(1+\alpha-\epsilon)/2\ge(1+\alpha)/4$. Testing $(\ref{fns})_3$ by $-\chi_l(\theta)(l/\theta)^\epsilon$ as in the derivation of \eqref{extlow:band}, we obtain
\be\nonumber\ba
 \int_0^T\left\|\nabla\bigl(\theta^{c_0}-(2l)^{c_0}\bigr)_+\right\|_{L^2}^2\,dt
 = c_0^2\int_0^T\int_{\{\theta>2l\}}\theta^{\alpha-1-\epsilon}|\nabla\theta|^2\,dx\,dt
 \le \frac{C}{\epsilon}(1+l+\mathcal P_l)
 \le \frac{C}{\epsilon}.
\ea\ee
Since $(\theta^{c_0}-(2l)^{c_0})_+$ vanishes on $\Gamma_t$ and satisfies $(\theta^{c_0}-(2l)^{c_0})_+\ge(1-2^{-c_0})\theta^{c_0}$ on $E=\{\theta>4l\}$, applying Lemma~\ref{1pi} and \eqref{gn11} yields
\be\nonumber\ba
 \int_0^T\|\theta\mathbf 1_E\|_{L^q}^{1+\alpha-\epsilon}\,dt
 =\int_0^T\|\theta^{c_0}\mathbf 1_E\|_{L^{q/c_0}}^2\,dt
 \le C q\int_0^T\left\|\nabla\bigl(\theta^{c_0}-(2l)^{c_0}\bigr)_+\right\|_{L^2}^2\,dt
 \le \frac{Cq}{\epsilon},
\ea\ee
which gives the first estimate in \eqref{extlow:heattail}. Replacing $(q,\epsilon)$ by $\bigl((1+\alpha)q,(1+\alpha)\delta\bigr)$ and using $K(\theta)=\frac{\kappa_0}{\alpha+1}\theta^{\alpha+1}$ directly yield the second estimate in \eqref{extlow:heattail}.

Next, put $b\triangleq\zeta/[2(1-\zeta)]$. Since $\int_0^\theta K(v)\,dv=\frac{1}{\alpha+2}\theta K(\theta)$ and $\theta\le C K(\theta)$ on $E$, H\"older's inequality gives
\be\nonumber\ba
 H_E
 \le C\int_E\rho\theta^{1-b}K(\theta)^{1+b}\,dx
 \le C\left(\int_E\rho\theta\,dx\right)^{1-b}\|\rho\|_{L^2(E)}^b
 \|K(\theta)\mathbf 1_E\|_{L^{\frac{2(1+b)}{b}}}^{1+b}.
\ea\ee
Using \eqref{fp01} and \eqref{extlow:NvsU}, noting that $(1+b)(1-\zeta)=1-\zeta/2$ and $2(1+b)/b\le 4/\zeta$, and applying \eqref{extlow:heattail}, we deduce
\be\nonumber\ba
 \int_0^T H_E^{1-\zeta}\,dt
 \le C\int_0^T\|K(\theta)\mathbf 1_E\|_{L^{\frac{2(1+b)}{b}}}^{1-\frac{\zeta}{2}}\,dt
 \le C\zeta^{-2}.
\ea\ee
Since the contribution of $\T^2\setminus E=\{\theta\le 4l\}$ to $H$ is bounded by $C M_0$, \eqref{extlow:potentialtail} follows immediately.

We now prove \eqref{extlow:templower}. Fix $0<l_1<l_2<h_0$. For $k\in[l_1,l_2]$, define
\be\nonumber\ba
 w\triangleq\bigl(K(k)-K(\theta)\bigr)_+,\quad
 \Phi_k(s)\triangleq\int_s^{\max\{s,k\}}\bigl(K(k)-K(v)\bigr)\,dv,
\ea\ee
so that $\Phi_k(\theta)\asymp w^2\asymp(k-\theta)_+^2$ uniformly in $k\in[l_1,l_2]$. Multiplying $(\ref{fns})_3$ by $-w$, integrating by parts over $\T^2$, and using $(\ref{fns})_1$, $Q-P\div u\ge -P^2/[4(\mu+\lambda)]$, and \eqref{extlow:struct}, we derive
\be\nonumber\ba
 c_v\frac{d}{dt}\int\rho\Phi_k(\theta)\,dx+\int|\nabla w|^2\,dx
 \le \int\frac{w P^2}{4(\mu+\lambda(\rho))}\,dx
 \le C\int_{\{\theta<k\}}\rho D(\rho)\,dx.
\ea\ee
Fix $q>3$. For any measurable space-time subset $\mathcal S\subset\T^2\times[s,s+1]$ with mass $y\triangleq\iint_{\mathcal S}\rho\,dx\,dt$, H\"older's inequality and the density moment bounds from \eqref{extlow:NvsU} yield
\be\nonumber\ba
 \iint_{\mathcal S}\rho D(\rho)\,dx\,dt
 \le y^{1-\frac{1}{q}}\left(\int_s^{s+1}\int\rho D(\rho)^q\,dx\,dt\right)^{\frac{1}{q}}
 \le C_q y^{1-\frac{1}{q}}.
\ea\ee
Now let $[s,s+1]\subset[0,T]$. For $i\ge 0$, set
\be\nonumber\ba
 k_i\triangleq l_1+(l_2-l_1)2^{-i},\quad
 s_i\triangleq s+\frac{1-2^{-i}}{2},\quad
 y_i\triangleq\int_{s_i}^{s+1}\int_{\{\theta<k_i\}}\rho\,dx\,dt,\quad
 w_i\triangleq\bigl(K(k_i)-K(\theta)\bigr)_+.
\ea\ee
Multiplying the differential inequality for $w_i$ by a smooth time cutoff vanishing on $[s,s_i]$, equal to $1$ on $[s_{i+1},s+1]$, and with derivative at most $C 2^i$, we integrate in time to obtain
\be\nonumber\ba
 \sup_{t\in[s_{i+1},s+1]}\int\rho w_i^2\,dx
 +\int_{s_{i+1}}^{s+1}\|\nabla w_i\|_{L^2}^2\,dt
 \le C 2^i y_i^{1-\frac{1}{q}}.
\ea\ee
Since $w_i=0$ on $\{\theta\ge h_0\}$, Lemma~\ref{1pi}, Sobolev embedding, and $\|\rho\|_{L^2}\le C$ give
\be\nonumber\ba
 \left(\int\rho w_i^4\,dx\right)^{1/2}
 \le \|\rho\|_{L^2}^{1/2}\|w_i\|_{L^8}^2
 \le C\|\nabla w_i\|_{L^2}^2.
\ea\ee
Noting that $w_i\ge c 2^{-i}$ on $\{\theta<k_{i+1}\}$, we apply H\"older's inequality to deduce
\be\nonumber\ba
 y_{i+1}
 \le C 2^{3i}\int_{s_{i+1}}^{s+1}\int\rho w_i^3\,dx\,dt
 &\le C 2^{3i}\int_{s_{i+1}}^{s+1}
 \left(\int\rho w_i^2\,dx\right)^{1/2}\left(\int\rho w_i^4\,dx\right)^{1/2}dt \\
 &\le C C_1^i\, y_i^{\frac{3}{2}\left(1-\frac{1}{q}\right)}.
\ea\ee
Here $\frac{3}{2}(1-1/q)>1$, while H\"older's inequality and \eqref{extlow:cold} (with $n=2$) give
\be\nonumber\ba
 y_0
 =\int_s^{s+1}\int_{\{\theta<l_2\}}\rho\,dx\,dt
 \le \int_s^{s+1}\|\rho\|_{L^2}|\{\theta<l_2\}|^{1/2}\,dt
 \le C\int_s^{s+1}\mathcal D_\theta(t)\,dt.
\ea\ee
Hence there exists a small constant $\varepsilon_0>0$ such that whenever $\int_s^{s+1}\mathcal D_\theta(t)\,dt\le\varepsilon_0$, the iteration gives $y_i\to 0$ as $i\to\infty$, and consequently $\theta(\cdot,s+1)\ge l_1$.
By \eqref{fp01}, at most finitely many (uniformly bounded) disjoint unit windows $[m,m+1]\subset[0,T]$ can violate $\int_m^{m+1}\mathcal D_\theta\,dt\le\varepsilon_0$.
For any $t\in[0,T]$, starting from the last good endpoint $t_*\in[0,t]$ (where $\theta(\cdot,t_*)\ge l_1$) or from $t_*=0$ (where $\theta(\cdot,0)\ge\underline\theta_0$), the elapsed time $t-t_*$ is bounded uniformly in $T$. Applying the comparison bound \eqref{fp117} on $[t_*,t]$ with $a_0\le C(1+D(R))$ from \eqref{extlow:struct} yields \eqref{extlow:templower}. This completes the proof of Lemma~\ref{low:temp-floor}.
\end{proof}

Under \textup{(A)}, \textup{(B)}, or \textup{(C)} with $\alpha>1-d$, Lemma~\ref{low:uniform} provides \eqref{new:smalldiss}--\eqref{extlow:unif} uniformly in $T$. Under \textup{(D)} or in the remaining range of \textup{(C)}, Lemmas~\ref{lem:linear-moment-input}--\ref{lem:gain-moment-input} and \eqref{extlow:NvsU} yield $U+\mathcal P_l\le C_{T_0}$ and all fixed density moments on $[0,T]\subset[0,T_0]$ (in these two cases, constants below may depend on $T_0$).

Set
\be\nonumber\ba
 D_u\triangleq\int\bigl(h(\div u)^2+\mu\omega^2\bigr)\,dx,\quad
 Z\triangleq\|\nabla K(\theta)\|_{L^2},\quad
 B\triangleq\int\rho\kappa(\theta)|\dot\theta|^2\,dx.
\ea\ee
Let $s\triangleq\frac{2}{2+\alpha}$, and choose $\sigma\in\bigl(\frac{1+s}{2},1\bigr)$ close enough to $1$ so that $\frac{3}{2+\alpha}<\sigma$ if $\alpha>1$, or $r_c\triangleq 3-(2+\alpha)\sigma\in(0,1)$ if $0<\alpha\le 1$. Define
\be\la{exthigh:d}\ba
 S\triangleq 1+A_1^2+H,\quad S_0\triangleq 1+D_u+(1+H)^\sigma .
\ea\ee
Fix $\ell>\max\{6,1+\alpha,1+2/\alpha\}$ and put
\be\nonumber\ba
 a\triangleq\frac{\ell}{\ell-1},\quad
 q_\ell\triangleq\frac{\ell+1}{\ell-1},\quad
 \Theta\triangleq 1+\|\theta\|_{L^{2\ell}},\quad
 g\triangleq A_1^2+\Theta^{q_\ell},\quad
 g_*\triangleq 1+g+(1+H)^{(2s-1)_+}.
\ea\ee

\begin{lemma}\label{core:pressure-times}
For $i=2$ and, when $\alpha>1$, also for $i=3$, one has
\be\la{exthigh:pressure}\ba
 \int\frac{P^i}{h^{i-1}}\,dx\le C(1+H)^{\frac{i}{2+\alpha}} .
\ea\ee
When $0<\alpha\le 1$, the cubic integral satisfies
\be\nonumber\ba
 \int\frac{P^3}{h^2}\,dx\le C(1+H)^\sigma\bigl(1+\|\theta\mathbf 1_E\|_{L^{p_0}}^{r_c}\bigr),
\ea\ee
where $E\triangleq\{\theta>4l\}$ and $p_0\triangleq 2r_c/(1-\sigma)$. Moreover,
\be\la{exthigh:compare}\ba
 &\|\nabla u\|_{L^2}^2+A_1^2\le C S_0,\quad S_0\le C S,\quad S\le C S_0^{1/\sigma},\\
 &\|\nabla u\|_{L^2}^4\le C\bigl(A_1^4+(1+H)^{2s}\bigr).
\ea\ee
For sufficiently large fixed constants $S_*,g_0$, let $J\triangleq\{t\in[0,T]: S(t)>S_*\}\cup\{t\in[0,T]: g_*(t)>g_0\}$. Then for every interval $I\subset[0,T]$ with $|I|\le 1$ and every $0<\zeta\le\zeta_0$,
\be\la{exthigh:tails}\ba
 \int_I\bigl(g_*+\sqrt S\bigr)\,dt+\zeta^2\int_I S^{1-\zeta}\,dt&\le C,\\
 |J|+\int_J(g_*+S_0)\,dt+\zeta^2\int_J S^{1-\zeta}\,dt&\le C .
\ea\ee
\end{lemma}

\begin{proof}
On $E=\{\theta>4l\}$, \eqref{extlow:struct} gives $P^i/h^{i-1}\le C\rho D(\rho)^{i-1}\theta^i$, while $\int_E\rho\theta^{2+\alpha}\,dx\le C H$. For $i<2+\alpha$, H\"older's inequality and the boundedness of fixed density moments yield
\be\nonumber\ba
 \int_E\frac{P^i}{h^{i-1}}\,dx
 \le C\left(\int_E\rho\theta^{2+\alpha}\,dx\right)^{\frac{i}{2+\alpha}}
 \left(\int_E\rho D(\rho)^{\frac{(i-1)(2+\alpha)}{2+\alpha-i}}\,dx\right)^{1-\frac{i}{2+\alpha}}
 \le C H^{\frac{i}{2+\alpha}}.
\ea\ee
On $E^c$, we have $P\le 4l\rho$, so $\int_{E^c}P^i/h^{i-1}\,dx\le C$, which proves \eqref{exthigh:pressure}. When $i=3$ and $0<\alpha\le 1$, writing $\theta^3=\theta^{(2+\alpha)\sigma}\theta^{r_c}$ and applying H\"older's inequality give
\be\nonumber\ba
 \int_E\frac{P^3}{h^2}\,dx
 \le C H^\sigma\left(\int_E\rho^2 D(\rho)^{\frac{4}{1-\sigma}}\,dx\right)^{\frac{1-\sigma}{2}}
 \|\theta\mathbf 1_E\|_{L^{p_0}}^{r_c}
 \le C H^\sigma\|\theta\mathbf 1_E\|_{L^{p_0}}^{r_c}.
\ea\ee

Next, from $h\div u=G+P$, \eqref{dc1}, and \eqref{exthigh:pressure} ($i=2$), we deduce
\be\nonumber\ba
 \|\nabla u\|_{L^2}^2+D_u\le C\bigl(A_1^2+(1+H)^s\bigr),\quad
 A_1^2\le C\bigl(1+D_u+(1+H)^s\bigr),
\ea\ee
which together with $s<\sigma<1$ directly imply \eqref{exthigh:compare}.

Finally, $U+\mathcal P_l\le C$ and \eqref{fp13} yield
\be\la{fp02}\ba
 \int_0^T\left(\|\nabla u\|_{L^2}^2+D_u+\int_E\frac{P^2}{h}\,dx\right)dt\le C, \quad A_1^2\le C\bigl(1+D_u+\int_E\frac{P^2}{h}\,dx\bigr).
\ea\ee
On $E^c$, the contributions to $H$ and $\Theta$ are uniformly bounded. Since $\ell>1+2/\alpha$ implies $q_\ell<1+\alpha$, we may choose $\epsilon,\zeta_0>0$ so small that $q_\ell<1+\alpha-\epsilon$ and $\max\{\sigma,\frac{1}{2},(2s-1)_+\}<1-\zeta_0$. For $S_*,g_0$ sufficiently large, Chebyshev's inequality, \eqref{fp02}, and \eqref{extlow:heattail}--\eqref{extlow:potentialtail} give
\be\nonumber\ba
 |J|\le C\int_0^T\left(D_u+\int_E\frac{P^2}{h}\,dx
 +\|\theta\mathbf 1_E\|_{L^{2\ell}}^{1+\alpha-\epsilon}+H_E^{1-\zeta_0}\right)dt\le C.
\ea\ee
Combining this with \eqref{fp02} and \eqref{extlow:heattail}--\eqref{extlow:potentialtail} yields \eqref{exthigh:tails}. This completes the proof of Lemma~\ref{core:pressure-times}.
\end{proof}

Fix $r\ge\max\bigl\{r_*,\sup_{0\le t\le T}\|\rho(t)\|_{L^\infty}\bigr\}$ and set
\be\nonumber\ba
 R\triangleq 2+r,\quad D_R\triangleq 1+D(r),\quad
 L_R\triangleq R\log(e+R),\quad \mathcal L_T\triangleq R_T\log(e+R_T),
\ea\ee
together with
\be\nonumber\ba
 X\triangleq e+A_1^2+A_2^2+Z^2,\quad
 E_\theta\triangleq\frac{1}{2}Z^2-\int K(\theta)Q\,dx,
\ea\ee
and the thermal weight
\be\la{thermal:weight}\ba
 \Xi\triangleq\max\left\{1,\ \sup_{\T^2\times[0,T]}\frac{\rho}{h(\rho)\kappa(\theta)}\right\}.
\ea\ee
By \eqref{extlow:struct} and \eqref{extlow:templower}, we have
\be\la{exthigh:loss}\ba
 \frac{\rho}{h}\le C D_R,\quad \frac{1}{\kappa(\theta)}\le C D_R^\alpha,\quad
 \Xi\le C D_R^{1+\alpha},\quad \int\frac{\rho^2|\dot\theta|^2}{h}\,dx\le\Xi B .
\ea\ee

\begin{lemma}\label{core:flux-heat}
There exist constants $C,n>0$ independent of $r$ such that, for every $z\in(0,1]$,
\be\la{kvu:q4}\ba
 Q_4\triangleq\|\nabla u\|_{L^4}^4+\|\sqrt h\,\div u\|_{L^4}^4+\|P\div u\|_{L^2}^2+\|P\|_{L^4}^4
 \le C\mathcal L_T g X,
\ea\ee
\be\la{kvu:xprime}\ba
(A_1^2)'+A_2^2&\le z Z^2+Cz^{-n}\mathcal L_T g A_1^2+C R_T\|\nabla u\|_{L^2}^4 + C A_1\sqrt{\Xi B},
\ea\ee
\be\la{fp77}\ba
\int K(\theta)\bigl(Q+|P\div u|\bigr)\,dx \le z(A_2^2+Z^2)+Cz^{-n}A_1^4\log(e+A_1^2),
\ea\ee
and
\be\la{fp05}\ba
E_\theta'+\frac{c_v}{2}B \le z\left(\|\nabla\dot u\|_{L^2}^2+\int h|D_t(\div u)|^2\,dx\right)
+ C z^{-n}D_R^n g X(\log X+L_R).
\ea\ee
and
\end{lemma}

\begin{proof}
From \eqref{fp31}, we deduce that $\Delta G=\div(\rho\dot u)$ and $\mu\Delta\omega=\nabla^\perp\cdot(\rho\dot u)$. Standard periodic elliptic estimates give, for $k\ge 1$ and $1<p<\infty$,
\be\la{2fp32}\ba
 \|\nabla^k G\|_{L^p}+\|\nabla^k\omega\|_{L^p}
 \le C_{k,p}\|\nabla^{k-1}(\rho\dot u)\|_{L^p}.
\ea\ee
For fixed $p>2$, applying Sobolev embedding, \eqref{gn11}, and \eqref{2fp32} and $\|\sqrt\rho\|_{L^p}\le C_p$ bounds $\|G-\ol G\|_{L^p}\le C_p A_2$. For $m\ge 8$, the order-zero elliptic bound costs $Cm$, while $\int\rho\dot u\,dx=0$ gives $M_0|\ol{\dot u}|\le C\|\nabla\dot u\|_{L^2}$ and hence $\|\dot u\|_{L^m}\le C m^{1/2}\|\nabla\dot u\|_{L^2}$ by \eqref{gn11}. Consequently,
\be\la{core:periodic-flux-lp}\ba
 &\|G-\ol G\|_{L^p}\le C_p A_2,\quad
 \|\nabla G\|_{L^2}+\|\nabla\omega\|_{L^2}\le C R_T^{1/2}A_2,\\
 &\|\nabla G\|_{L^m}\le C m^{3/2}R_T\|\nabla\dot u\|_{L^2}.
\ea\ee
Since $|\ol G|\le\|G/\sqrt h\|_{L^2}\|\sqrt h\|_{L^2}\le C A_1$ and
$\ol\omega=0$, we have
$$\|G\|_{L^6}+\|\omega\|_{L^6}\le C(A_1+A_2).$$
Applying \eqref{dc1} together with $h\div u=G+P$, we obtain for $2\le p\le 4$,
\be\la{2fp003}\ba
 \|\nabla u\|_{L^p}+\|\sqrt h\,\div u\|_{L^p}+\|P\div u\|_{L^{p/2}}^{1/2}
 \le C\left(A_1+\|P\|_{L^p}+A_1^{\frac{6-p}{2p}}A_2^{\frac{3p-6}{2p}}\right).
\ea\ee
Interpolating $\|P\|_{L^1}\le C$ with $\|P\|_{L^\ell}\le\|\rho\|_{L^{2\ell}}\|\theta\|_{L^{2\ell}}\le C\Theta$, and noting that Lemma~\ref{1pi} on $\Gamma_t$ gives $\|K(\theta)\|_{H^1}\le C(1+Z)$, we deduce
\be\la{ltu-pressure}\ba
 \|P\|_{L^p}\le C\Theta^{a(1-1/p)}\quad (1\le p\le\ell),\quad
 \Theta\le C(1+Z)^{\frac{1}{1+\alpha}}.
\ea\ee
In particular, $\|\nabla u\|_{L^2}^2\le C(A_1^2+\|P\|_{L^2}^2)\le C g\le C X$, and since $3a-q_\ell<2a<2(1+\alpha)$,
\be\la{ltu-absorb}\ba
 \|P\|_{L^2}^2&\le C\Theta^a,&
 \|P\|_{L^4}^4&\le C\Theta^{3a}\le C\Theta^{q_\ell}(1+Z^2),\\
 \|P\|_{L^2}\|P\|_{L^4}^2&\le C\Theta^{2a}\le z Z^2+Cz^{-n},&
 \|P\|_{L^4}^2&\le C\Theta^{q_\ell/2}(1+Z).
\ea\ee

To prove \eqref{kvu:q4}, set $\tau\triangleq[4\log(e+R_T)]^{-1}$ and $p\triangleq 2(1+\tau)/\tau$. By \eqref{visc:envelope}, $\|G\|_{L^2}\le C R_T^{B_h/2}A_1$. Using H\"older's inequality, \eqref{gn11}, \eqref{core:periodic-flux-lp}, and $p^\tau+R_T^{B_h\tau/2}\le C$, we derive
\be\la{core:weighted-flux}\ba
 \left\|\frac{G^2}{h}\right\|_{L^2}
 \le C A_1^{1-\tau}\|G\|_{L^p}^{1+\tau}
 &\le C A_1^{1-\tau}\Bigl(|\ol G|^{1+\tau}
 +\bigl(C\sqrt p\bigr)^{1+\tau}\|G-\ol G\|_{L^2}^\tau\|\nabla G\|_{L^2}\Bigr)\\
 &\le C A_1^2+C\sqrt{\mathcal L_T}\,A_1 A_2.
\ea\ee
Similarly, $\|\omega\|_{L^4}^4\le C\|\omega\|_{L^2}^2\|\nabla\omega\|_{L^2}^2\le C R_T A_1^2 A_2^2$. Combining these bounds with \eqref{dc1}, $h\div u=G+P$, and \eqref{ltu-absorb} yields \eqref{kvu:q4}.

Next, to prove \eqref{kvu:xprime}, direct computation gives
\begin{align}
 \curl\dot u&=D_t\omega+\omega\div u,\label{fp44}\\
 \div\dot u&=D_t(\div u)+(\div u)^2-2\det\nabla u.\label{fp45}
\end{align}
Multiplying \eqref{fp31} by $2\dot u$, integrating by parts over $\T^2$, and using $D_t P+P\div u=\rho\dot\theta$ and $D_t h=-\rho\lambda'\div u$, we obtain
\be\la{fp43}\ba
 (A_1^2)'+2A_2^2
 &=-\mu\int\omega^2\div u\,dx+4\int G\det\nabla u\,dx
  -\int\frac{h+\rho\lambda'}{h^2}G^2\div u\,dx\\
 &\quad -2\int\frac{\rho\lambda'}{h^2}PG\div u\,dx
  -2\int\frac{\rho}{h}G\dot\theta\,dx
 \triangleq\sum_{i=1}^5 I_i.
\ea\ee
By H\"older's inequality, the $\mathcal H^1$--$\mathcal{BMO}$ pairing for $\det\nabla u$ \cite{CLMS}, and \eqref{core:periodic-flux-lp},
\be\nonumber\ba
 |I_1|+|I_2|
 \le C\|\omega\|_{L^4}^2\|\nabla u\|_{L^2}+C\|\nabla G\|_{L^2}\|\nabla u\|_{L^2}^2
 \le \frac{1}{4}A_2^2+C R_T\bigl(A_1^2\|\nabla u\|_{L^2}^2+\|\nabla u\|_{L^2}^4\bigr).
\ea\ee
Since $|\rho\lambda'|\le C h$ by \eqref{visc:regularity}, \eqref{core:weighted-flux}, \eqref{ltu-absorb}, and H\"older's inequality yield
\be\nonumber\ba
 |I_3|+|I_4|
 \le C\|\nabla u\|_{L^2}\left(\left\|\frac{G^2}{h}\right\|_{L^2}+\|P\|_{L^4}^2\right)
 \le \frac{1}{4}A_2^2+z Z^2+Cz^{-n}\mathcal L_T g A_1^2.
\ea\ee
For $I_5$, H\"older's inequality and \eqref{exthigh:loss} give
\be\la{fp411}\ba
 \frac{|I_5|}{2}
 \le\left\|\frac{G}{\sqrt h}\right\|_{L^2}\left\|\frac{\rho\dot\theta}{\sqrt h}\right\|_{L^2}
 \le A_1\sqrt{\Xi B}.
\ea\ee
Substituting the bounds for $I_1$--$I_5$ into \eqref{fp43} proves \eqref{kvu:xprime}.

We now prove \eqref{fp77}. For any $b\ge\ell$, let $p\triangleq 2b/(b-1)$. Then $\|K(\theta)\|_{L^{2b}}\le C\sqrt b(1+Z)$ by Lemma~\ref{1pi} and \eqref{gn11}, while \eqref{ltu-pressure} gives $\|P\|_{L^p}^2\le C\Theta^{q_\ell}$. Applying H\"older's inequality and \eqref{2fp003} yields
\be\la{fp74}\ba
 \int K(\theta)\bigl(Q+|P\div u|\bigr)\,dx
 \le C\sqrt b(1+Z)\left(g+A_1^{2-\frac{3}{b}}A_2^{\frac{3}{b}}\right).
\ea\ee
Replacing $K(\theta)$ by $K(\theta)^2$ gives the same bound with prefactor $C b(1+Z)^2$; choosing $b=\ell+2\log(e+A_2^2)$ thus yields
\be\la{fp517}\ba
 \int K(\theta)^2\bigl(|\nabla u|^2+\lambda(\div u)^2\bigr)\,dx
 \le C g X\log X .
\ea\ee
Choosing instead $b=\ell+\log(e+A_1^2)$ in \eqref{fp74} and applying Young's inequality together with \eqref{ltu-pressure}, we deduce
\be\nonumber\ba
 \int K(\theta)\bigl(Q+|P\div u|\bigr)\,dx
 &\le z(Z^2+A_2^2)+Cz^{-n}\left(b^{\frac{1}{1-q_\ell/(1+\alpha)}}+b A_1^4+b^{\frac{b}{b-3}}A_1^{4+\frac{6}{b-3}}\right)\\
 &\le z(Z^2+A_2^2)+Cz^{-n}A_1^4\log(e+A_1^2),
\ea\ee
which proves \eqref{fp77}. Multiplying $(\ref{fns})_3$ by $K(\theta)$, integrating over $\T^2$, and applying \eqref{fp77} give
\be\la{fp720}\ba
 H'+Z^2=\int K(\theta)(Q-P\div u)\,dx,\quad
 H'+\frac{3}{4}Z^2\le \frac{1}{4}A_2^2+C A_1^4\log(e+A_1^2).
\ea\ee
Finally, to prove \eqref{fp05}, we multiply $(\ref{fns})_3$ by $D_t K(\theta)=\kappa(\theta)\dot\theta$ and integrate over $\T^2$:
\be\nonumber\ba
 c_v B=\int\bigl(\Delta K(\theta)+Q-P\div u\bigr)D_t K(\theta)\,dx.
\ea\ee
For the diffusion term, periodic integration by parts and (\ref{gn11}) give
\be\nonumber\ba
 \int\Delta K(\theta)\,D_t K(\theta)\,dx
 \le -\frac{1}{2}\frac{d}{dt}Z^2+C\|\nabla u\|_{L^2}Z\|\Delta K(\theta)\|_{L^2}.
\ea\ee
Since $\Delta K(\theta)=c_v\rho\dot\theta-Q+P\div u$, \eqref{exthigh:loss} and \eqref{kvu:q4} yield
\be\nonumber\ba
 \|\Delta K(\theta)\|_{L^2}^2\le C R D_R^\alpha B+C L_R g X,
\ea\ee
whence
$$C\|\nabla u\|_{L^2}Z\|\Delta K(\theta)\|_{L^2}\le \frac{c_v}{4}B+C D_R^\alpha L_R g X.$$
For the production term, using $D_t\lambda=-\rho\lambda'\div u$, $D_t\nabla u=\nabla\dot u-(\nabla u)^2$, $|\rho\lambda'|\le C h$, and Young's inequality, we derive
\be\nonumber\ba
 \int Q\,D_t K(\theta)\,dx-\frac{d}{dt}\int K(\theta)Q\,dx
 &=-\int K(\theta)\bigl(D_t Q+Q\div u\bigr)\,dx\\
 &\le z\left(\|\nabla\dot u\|_{L^2}^2+\int h|D_t(\div u)|^2\,dx\right)\\
 &\quad +Cz^{-n}\left(Q_4+\int K(\theta)^2\bigl(|\nabla u|^2+h(\div u)^2\bigr)\,dx\right).
\ea\ee
For the pressure work term, since $\theta^2\kappa(\theta)\le C(1+K(\theta)^2)$ and $\rho/h\le C D_R$, Young's inequality gives
\be\la{fp57}\ba
 -\int P\div u\,D_t K(\theta)\,dx
 \le \frac{c_v}{4}B+C D_R\int\bigl(1+K(\theta)^2\bigr)h(\div u)^2\,dx.
\ea\ee
Collecting the above estimates and applying \eqref{kvu:q4} and \eqref{fp517} yield \eqref{fp05}. This completes the proof of Lemma~\ref{core:flux-heat}.
\end{proof}

\begin{lemma}\label{core:acceleration}
Let $v\triangleq 2\nabla\Delta^{-1}\det\nabla u$ and $w\triangleq\dot u+v$. Setting $E_w\triangleq\int\rho|w|^2\,dx$ and $D_w\triangleq\int\bigl(h|\div w|^2+\mu|\omega(w)|^2\bigr)\,dx$, we have
\be\la{kvu:corrected}\ba
 E_w'+c D_w\le C\Xi B+C\mathcal L_T g(X+E_w).
\ea\ee
\end{lemma}

\begin{proof}
Denote $\mathcal J(f,l)\triangleq\nabla f\cdot\nabla^\perp l$. Since $\det\nabla u=\mathcal J(u^1,u^2)$ has zero mean on $\T^2$, we have $\div v=2\det\nabla u$ and $\omega(v)=0$. Thus \eqref{fp44}--\eqref{fp45} give
\be\nonumber\ba
 \div w=D_t(\div u)+(\div u)^2,\quad
 \omega(w)=D_t\omega+\omega\div u,\quad
 \|\nabla w\|_{L^2}^2\le C D_w.
\ea\ee
By the Jacobian estimate $\|\nabla\Delta^{-1}\mathcal J(f,l)\|_{L^2}\le C\|\nabla f\|_{L^2}\|\nabla l\|_{L^2}$ \cite{CLMS} and $L^2$ elliptic regularity,
\be\la{kv:qsize}\ba
 \|v\|_{L^2}\le C\|\nabla u\|_{L^2}^2,\quad
 \|\nabla v\|_{L^2}\le C\|\nabla u\|_{L^4}^2.
\ea\ee
Differentiating $\det\nabla u$ gives $\partial_t(\det\nabla u)+\div(u\det\nabla u)=\mathcal J(\dot u^1,u^2)+\mathcal J(u^1,\dot u^2)$, so that
\be\nonumber\ba
 D_t v=2\nabla\Delta^{-1}\bigl(\mathcal J(\dot u^1,u^2)+\mathcal J(u^1,\dot u^2)\bigr)+\mathcal C,\quad
 \mathcal C\triangleq u\cdot\nabla v-\nabla\Delta^{-1}\div(u\div v).
\ea\ee
Writing $v=\nabla\varphi$, direct computation gives
\be\nonumber\ba
 \div\mathcal C=-\mathcal J(u^1,v^2)-\mathcal J(v^1,u^2),\quad
 \omega(\mathcal C)=-\sum_{j=1}^2\mathcal J(u^j,\partial_j\varphi),\quad
 \ol{\mathcal C}=-\int(\div u)v\,dx.
\ea\ee
By \eqref{dc1} and the Jacobian bound \cite{CLMS}, we have $\|\mathcal C\|_{L^2}\le C\|\nabla u\|_{L^2}(\|v\|_{L^2}+\|\nabla v\|_{L^2})$. Substituting $\dot u=w-v$ and using \eqref{kv:qsize} yield
\be\la{kv:qmaterial}\ba
 \|D_t v\|_{L^2}\le C\|\nabla u\|_{L^2}\bigl(\|\nabla w\|_{L^2}+\|\nabla u\|_{L^4}^2\bigr).
\ea\ee

Next, put $\mathbb J\triangleq\bigl(\begin{smallmatrix}0&1\\-1&0\end{smallmatrix}\bigr)$ and $\mathcal T\triangleq G\mathbb I_2+\mu\omega\mathbb J$. Differentiating $\rho\dot u=\div\mathcal T$ along the flow yields
\be\nonumber\ba
 \rho D_t\dot u=\div\bigl(D_t\mathcal T+(\div u)\mathcal T-\mathcal T(\nabla u)^\top\bigr).
\ea\ee
Using
$$D_t G+(\div u)G=h\div w-\rho\lambda'(\div u)^2-\rho\dot\theta,\quad D_t\omega+(\div u)\omega=\omega(w),$$
we multiply $\rho D_t w=\rho D_t\dot u+\rho D_t v$ by $w$ and integrate by parts over $\T^2$ to obtain
\be\nonumber\ba
 \frac{1}{2}E_w'+D_w
 =\int\bigl(\rho\lambda'(\div u)^2+\rho\dot\theta\bigr)\div w\,dx
 +\int\mathcal T(\nabla u)^\top:\nabla w\,dx
 +\int\rho D_t v\cdot w\,dx.
\ea\ee
Using $(\nabla u)^\top:\nabla w=(\div u)\div w-\mathcal J(u^1,w^2)-\mathcal J(w^1,u^2)$ in the $G$-part of the second integral, we pair $G$ with the two Jacobians via \cite{CLMS} and \eqref{core:periodic-flux-lp} to bound their contribution by $C\|\nabla G\|_{L^2}\|\nabla u\|_{L^2}\|\nabla w\|_{L^2}$, while the vorticity part is bounded by $C\|\nabla u\|_{L^4}^2\|\nabla w\|_{L^2}$. Applying Young's inequality with weight $h$ to the remaining divergence terms (using $|\rho\lambda'|\le C h$, \eqref{exthigh:loss}, and $G^2(\div u)^2/h\le 2h(\div u)^4+2P^2(\div u)^2/h$), and using \eqref{kv:qmaterial} for the last integral, we absorb the small multiples of $D_w$ to arrive at
\be\la{kv:corrected-raw}\ba
 E_w'+c D_w
 \le C\Xi B+C R_T\|\nabla u\|_{L^2}^2(E_w+A_2^2)
 +C\int\left(|\nabla u|^4+h(\div u)^4+\frac{P^2(\div u)^2}{h}\right)dx.
\ea\ee
Combining \eqref{kv:corrected-raw} with \eqref{kvu:q4} and $\|\nabla u\|_{L^2}^2\le C g$ yields \eqref{kvu:corrected}. This completes the proof of Lemma~\ref{core:acceleration}.
\end{proof}

\begin{lemma}\label{records:joint}
Fix $r\ge\max\bigl\{r_*,\sup_{0\le t\le T}\|\rho(t)\|_{L^\infty}\bigr\}$, $R\triangleq 2+r$, and set
\be\nonumber\ba
 \mathcal H\triangleq 1+H,\quad N\triangleq R+\log(e+A_1^2),\quad
 F_R\triangleq A_1^4 N+R\mathcal H^{2s},\quad
 Y_R\triangleq e+E_\theta+t_R E_w+C_R F_R,
\ea\ee
and $W_R\triangleq k_R\log Y_R+\log S$.
With $t_R\triangleq c_1 D_R^{-(1+\alpha)}$, one can choose $c_1>0$ small and $1\le C_R,k_R\le M_R\triangleq C D_R^C$ so that $Y_R\ge e$, $Y_R\asymp e+t_R A_2^2+Z^2+C_R F_R$, and
\be\la{exthigh:joint}\ba
 W_R'+M_R^{-1}\mathcal D_R\le M_R g_*(W_R+L_R),\quad
 \mathcal D_R\triangleq\frac{D_w+B}{Y_R}+\frac{Y_R+A_2^2+Z^2}{S}.
\ea\ee
For $\eta\triangleq 1/8$, a fixed sufficiently small $\zeta_*>0$, and $Q_R\triangleq L_R\exp(C D_R^C)$, one has
\be\la{exthigh:window}\ba
 &\sup_{0\le t\le T}W_R(t)+\sup_{|I|\le 1}\int_I\mathcal D_R\,dt\le Q_R,\\
 &\int_{\{S_0>R^\eta\}}\mathcal D_R\,dt\le Q_R
   \quad\text{if } D(r)\le(\log r)^{\zeta_*}\text{ and } r\gg 1,
\ea\ee
where all time sets lie in $[0,T]$. On the indicated large levels, $Q_R=R^{1+o(1)}$. Constants inherit the time dependence of the low-order bounds.
\end{lemma}

\begin{proof}
Since $\dot u=w-v$ and $D_t(\div u)=\div w-(\div u)^2$, \eqref{kv:qsize} and \eqref{kvu:q4} yield
\be\la{common:material-gradient}\ba
 \|\nabla\dot u\|_{L^2}^2+\int h|D_t(\div u)|^2\,dx
 \le C(D_w+Q_4)\le C(D_w+L_R g X).
\ea\ee
Substituting \eqref{common:material-gradient} into \eqref{fp05} gives, for $0<z\le 1$ and a fixed $n>0$,
\be\la{exthigh:tprime}\ba
 E_\theta'+c B\le z D_w+C z^{-n} D_R^n g X(\log X+L_R).
\ea\ee
Adding $t_R$ times \eqref{kvu:corrected} to \eqref{exthigh:tprime}, using $\Xi\le C D_R^{1+\alpha}$, and choosing $c_1>0$ small and then $z\ll t_R$ absorb $t_R C\Xi B$ and $z D_w$ on the right-hand side. Applying \eqref{fp77} with the same $z$ and noting
$$\int\rho|v|^2\,dx\le C R\|\nabla u\|_{L^2}^4\le C F_R,$$
we increase $C_R\le C D_R^C$ to obtain $Y_R\ge e$ and $Y_R\asymp e+t_R A_2^2+Z^2+C_R F_R$.

Next, we differentiate $F_R$ via $(A_1^4 N)'=\bigl(2A_1^2 N+A_1^4/(e+A_1^2)\bigr)(A_1^2)'$, \eqref{kvu:xprime}, and \eqref{fp720}. Since $C_R A_1^4 N\le C Y_R$ and $Z^2+A_2^2+R\|\nabla u\|_{L^2}^4\le C D_R^C Y_R$, applying Young's inequality to the $\sqrt{\Xi B}$ term yields
\be\nonumber\ba
 \frac{C_R(A_1^4 N)'}{Y_R}
 &\le \epsilon\frac{B}{Y_R}+C_\epsilon D_R^C g(L_R+\log Y_R),\\
 \frac{C_R R(\mathcal H^{2s})'}{Y_R}
 &\le C D_R^C\mathcal H^{(2s-1)_+}(R+\log Y_R).
\ea\ee
Choosing $\epsilon>0$ small absorbs $\epsilon B/Y_R$ into $c B/Y_R$ and controls $(\log Y_R)'$.
For $(\log S)'$, we add \eqref{fp720} and \eqref{kvu:xprime} and divide by $S\ge\max\{A_1^2,\mathcal H\}$. By the equivalence of $Y_R$ and Young's inequality,
\be\nonumber\ba
 & \frac{Y_R}{S}
 \le C D_R^C\left(\frac{A_2^2+Z^2}{S}+g_*(R+\log Y_R)\right),\\
 & \frac{D_R^{(1+\alpha)/2}A_1\sqrt B}{S}
 \le \epsilon\frac{Y_R}{S}+C_\epsilon D_R^{1+\alpha}\frac{B}{Y_R}.
\ea\ee
Taking $\epsilon=c D_R^{-m}$ ($m$ large, $c$ small) absorbs $\epsilon Y_R/S$ into $(A_2^2+Z^2)/S$, and adding $k_R(\log Y_R)'$ with $k_R\le C D_R^C$ sufficiently large absorbs $C_\epsilon D_R^{1+\alpha}B/Y_R$, proving \eqref{exthigh:joint}.

By \eqref{th13} and \eqref{ext:dataM}, $Y_R(0)\le C D_R^C(1+R)$, so Gr\"onwall's inequality bounds $W_R$ on $[0,1]\cap[0,T]$. On any later unit window $[b,b+1]\subset[0,T]$, set $A(t)\triangleq\exp\bigl(-M_R\int_b^t g_*(\tau)\,d\tau\bigr)$ and $V_*(t)\triangleq A(t)(W_R(t)+L_R)$. By \eqref{exthigh:tails}, $A\ge e^{-C M_R}$. Using $Y_R/S=e^{W_R/k_R}S^{-1-1/k_R}$ and $W_R-V_*\ge -L_R$, we deduce from \eqref{exthigh:joint} that
\be\nonumber\ba
 \frac{d}{dt}\bigl(e^{-V_*/k_R}\bigr)
 \ge \frac{e^{-C M_R-L_R/k_R}}{k_R M_R}S^{-1-\frac{1}{k_R}}.
\ea\ee
Integrating over $[b,b+1]$ and applying Jensen's inequality with $\int_b^{b+1}\sqrt S\,dt\le C$ from \eqref{exthigh:tails} yield $W_R(b+1)\le Q_R^0\triangleq L_R\exp(C_0 D_R^{C_0})$. Integrating \eqref{exthigh:joint} over any unit interval $I\subset[0,T]$ then gives the first line of \eqref{exthigh:window}.

For the second line, set $h_R\triangleq k_R\log(C M_R^2 Q_R^0)$. On $\{W_R>h_R\}\setminus J$, both $S$ and $g_*$ are bounded, so $Y_R/(M_R S)$ absorbs twice the right-hand side of \eqref{exthigh:joint}. On $\{S_0>R^\eta\}$, \eqref{exthigh:compare} gives $S\ge c R^\eta$ while $Y_R\ge c R/M_R$, so that
\be\nonumber\ba
 W_R-h_R\ge \eta\log R-C-k_R\bigl(C+\log\log(e+R)+3\log M_R+C_0 D_R^{C_0}\bigr)>0
\ea\ee
for $r\gg 1$ and $D(r)\le(\log r)^{\zeta_*}$ with $\zeta_*>0$ small enough. Integrating the differential inequality for $(W_R-h_R)_+$ over $[0,T]$ yields
\be\nonumber\ba
 \int_{\{S_0>R^\eta\}}\mathcal D_R\,dt
 \le C M_R\bigl(W_R(0)-h_R\bigr)_++C M_R^2 Q_R^0\int_J g_*\,dt\le Q_R.
\ea\ee
Decreasing $\zeta_*$ once more ensures $Q_R=R^{1+o(1)}$. This completes the proof of Lemma~\ref{records:joint}.
\end{proof}

Set $r_j\triangleq 2^j r_0$ with $r_0>2\|\rho_0\|_{L^\infty}$ sufficiently large, $R_j\triangleq 2+r_j$, and $L_j\triangleq R_j\log(e+R_j)$.
Whenever $\|\rho(t)\|_{L^\infty}$ reaches $r_j$ on $[0,T]$, let $t_j$ be the first such time.
For each such $j$ satisfying $D(r_j)\le(\log r_j)^{\zeta_*}$, let $j^-$ be the largest index $i<j$ satisfying the same inequality.
If no such $i$ exists, set $t_{j^-}=0$. Define
\be\nonumber\ba
 I_j^0\triangleq [t_{j^-},t_j],\quad
 E_j\triangleq I_j^0\cap\{S_0>R_j^\eta\},\quad
 e_j\triangleq\int_{E_j}S_0\,dt,\quad
 d_j\triangleq L_j^{-1}\int_{E_j}\frac{A_2^2}{S_0}\,dt .
\ea\ee

\begin{lemma}\label{records:sum}
Suppose that $\|\rho(t)\|_{L^\infty}$ reaches $r_N$ by time $T$. Then
\be\la{exthigh:sum}\ba
 \sum_{\substack{j\le N\\D(r_j)\le(\log r_j)^{\zeta_*}}}(d_j+e_j)\le C,
\ea\ee
where $C$ is independent of $N$.
\end{lemma}

\begin{proof}
Multiplying \eqref{fp31} by $\dot u$, integrating by parts over $\T^2$, and using $\div\dot u=\div w-2\det\nabla u$ give
\be\la{per:mechanical-identity}\ba
 \frac{1}{2}D_u'+A_2^2
 &=\int P\div w\,dx+2\int G\det\nabla u\,dx\\
 &\quad -\frac{1}{2}\int(h+\rho\lambda')(\div u)^3\,dx
  -\frac{\mu}{2}\int\omega^2\div u\,dx .
\ea\ee
By the Jacobian bound \cite{CLMS}, \eqref{2fp32}, and the $L^4$ vorticity estimate, the second and fourth integrals on the right-hand side are bounded by $C R^{1/2}A_2\|\nabla u\|_{L^2}^2$. Using $h\div u=G+P$, $|\rho\lambda'|\le C h$, \eqref{core:weighted-flux}, and Lemma~\ref{core:pressure-times}, we estimate the cubic terms by
\be\nonumber\ba
 \int\frac{|G|^3}{h^2}\,dx
 \le C A_1\left\|\frac{G^2}{h}\right\|_{L^2}
 \le \epsilon A_2^2+C_\epsilon L_R S_0^2,\quad
 \int\frac{|P|^3}{h^2}\,dx\le C S_0(1+F_\alpha),
\ea\ee
where $F_\alpha\triangleq 0$ for $\alpha>1$ and $F_\alpha\triangleq\|\theta\mathbf 1_E\|_{L^{p_0}}^{r_c}$ for $0<\alpha\le 1$. Set $b\triangleq 1+S_0+F_\alpha$, so $\int_J b\,dt\le C$ by \eqref{exthigh:tails} and \eqref{extlow:heattail}. Adding $\sigma(1+H)^{\sigma-1}$ times \eqref{fp720} to twice \eqref{per:mechanical-identity}, dividing by $S_0$, and using $\|\nabla u\|_{L^2}^2+A_1^2\le C S_0$ from \eqref{exthigh:compare} yield
\be\la{exthigh:mech}\ba
 (\log S_0)'+c\frac{A_2^2}{S_0}
 \le \frac{2\left|\int P\div w\,dx\right|}{S_0}+C b(L_R+\log S_0).
\ea\ee
By H\"older's inequality, \eqref{exthigh:pressure}, and $1+H\le S\le C S_0^{1/\sigma}$, with $\delta\triangleq 1-(1+s)/(2\sigma)>0$, we have
\be\la{exthigh:work}\ba
 \frac{\left|\int P\div w\,dx\right|}{S_0}
 \le C\frac{(1+H)^{s/2}\sqrt S}{S_0}\left(\frac{D_w}{Y_R}\frac{Y_R}{S}\right)^{1/2}
 \le C S_0^{-\delta}\mathcal D_R .
\ea\ee

Fix $j\le N$ satisfying $D(r_j)\le(\log r_j)^{\zeta_*}$. Set $r=r_j$ on $[0,t_j]$, and enlarge $r_0$ so that $E_j\subset J$. By \eqref{exthigh:window} and $Q_{R_j}=R_j^{1+o(1)}$,
\be\nonumber\ba
 L_j^{-1}\int_{E_j}\frac{\left|\int P\div w\,dx\right|}{S_0}\,dt
 \le C R_j^{-\eta\delta}\frac{Q_{R_j}}{L_j}
 \le C R_j^{-\eta\delta/2}.
\ea\ee
Define $U_j\triangleq(\log S_0-\eta\log R_j)_+$ and $B_*(t)\triangleq\exp\bigl(-C\int_0^t\mathbf 1_J b\,d\tau\bigr)\ge c$. On $I_j^0$, almost everywhere,
\be\nonumber\ba
 (B_* U_j)'=B_*\mathbf 1_{E_j}(\log S_0)'-C\mathbf 1_J b B_* U_j .
\ea\ee
Since $E_j\subset J$, the last term cancels the contribution of $b U_j$ from \eqref{exthigh:mech}. Multiplying \eqref{exthigh:mech} by $B_*\mathbf 1_{E_j}/L_j$ and integrating over $I_j^0$ give
\be\la{per:layer-endpoints}\ba
 c d_j+\left[\frac{B_* U_j}{L_j}\right]_{t_{j^-}}^{t_j}
 \le C\int_{E_j}b\,dt+C R_j^{-\eta\delta/2}.
\ea\ee
Consider the indices $j\le N$ satisfying $D(r_j)\le(\log r_j)^{\zeta_*}$. For two consecutive such indices $j<j^+$, we have $U_{j^+}/L_{j^+}\le U_j/L_j$ at the common endpoint $t_j$.
Summing \eqref{per:layer-endpoints} over all such indices and using the disjoint interiors of the intervals $I_j^0$, we obtain
\be\nonumber\ba
 c\sum_{\substack{j\le N\\D(r_j)\le(\log r_j)^{\zeta_*}}} d_j
 \le C\left(1+\int_J b\,dt+\sum_{j\le N} R_j^{-\eta\delta/2}\right)\le C,\quad
 \sum_{\substack{j\le N\\D(r_j)\le(\log r_j)^{\zeta_*}}} e_j\le\int_J S_0\,dt\le C,
\ea\ee
which proves \eqref{exthigh:sum}.
\end{proof}

\begin{lemma}\label{records:flux}
Set $f\triangleq G-\ol G$. In this lemma, $j$ ranges over the indices satisfying $D(r_j)\le(\log r_j)^{\zeta_*}$ for which $\|\rho(t)\|_{L^\infty}$ reaches $r_j$ by time $T$. There exists $\xi\in(0,1)$ such that, for every interval $I\subset[0,t_j]$,
\be\nonumber\ba
 \int_I|\ol G|\,dt\le C(1+|I|),\quad
 \int_{I\cap\{S_0\le R_j^\eta\}}\|f\|_{L^\infty}\,dt\le C R_j^\xi(1+|I|).
\ea\ee
Moreover, on the exceptional sets $E_j$, there exist weights $x_j\ge 0$ and $w_j\in(0,1)$ such that
\be\la{per:exceptional-flux}\ba
 \int_{E_j}\|f\|_{L^\infty}\,dt\le C L_j x_j^{w_j},\quad
 \sum_j x_j\le C,\quad
 1-w_j=O(1/j).
\ea\ee
\end{lemma}

\begin{proof}
Set $T_2\triangleq 1+A_2+\|\nabla u\|_{L^2}$ and $T_3\triangleq 1+A_1+T_2+\|\nabla\dot u\|_{L^2}+Z^{1/(1+\alpha)}$. By \eqref{core:periodic-flux-lp} and \eqref{visc:envelope},
\be\la{per:cross-elliptic}\ba
 \|f\|_{L^p}&\le C_p T_2\quad (p>2\text{ fixed}),&
 \|f\|_{L^2}&\le C R^{B_h/2}A_1,\\
 \|\nabla f\|_{L^n}&\le C n^8 R D_R T_3\quad (n\ge 8),&
 \|\nabla f\|_{L^2}&\le C R^{1/2}T_2 .
\ea\ee
Since $\int(h+P)\,dx\le C$, we have $|\ol G|\le C(1+\sqrt{D_u})$, which together with \eqref{fp02} gives $\int_I|\ol G|\,dt\le C(1+|I|)$.
For any $p,n\ge 8$, balancing the average $C l^{-2/p}\|f\|_{L^p}$ on a ball of radius $l$ against the Morrey oscillation $C l^{1-2/n}\|\nabla f\|_{L^n}$ yields, uniformly in $p,n\ge 8$,
\be\la{records:ball}\ba
 \|f\|_{L^\infty}\le C\|f\|_{L^p}^{1-v}\|\nabla f\|_{L^n}^v,\quad
 v\triangleq\frac{2}{p(1-2/n)+2}.
\ea\ee

On $E_j$ and on every unit window, \eqref{exthigh:tails} with $\zeta=1/Q_{R_j}$ and $\log S\le Q_{R_j}$ imply $\int S\,dt\le C Q_{R_j}^2$. By \eqref{common:material-gradient}, $X/Y_R\le C D_R^C$, and \eqref{exthigh:window}, $\int T_3^2/Y_R\,dt$ is bounded by a fixed power of $R_j$. Applying H\"older's inequality to $(T_3^2/Y_R)^{1/3}(Y_R/S)^{1/3}S^{1/3}$ gives $\int T_3^{2/3}\,dt\le R_j^{c_3}$ for a fixed $c_3>0$.

On the subset of any unit window where $S_0\le R_j^\eta$, \eqref{exthigh:compare} and \eqref{exthigh:window} yield $\int T_2^2\,dt\le C R_j^{\eta/\sigma}Q_{R_j}$. Taking $n=8$ and $p$ fixed and sufficiently large in \eqref{records:ball} and integrating via H\"older's inequality produce a power of $R_j$ tending to $(1+\eta/\sigma)/2<1$ as $p\to\infty$. Covering $I$ by unit windows gives the second bound.

On $E_j$, take $p\triangleq 16+10\log(e+R_j)$, $n\triangleq e^p$, $a_j\triangleq(1-2/n)v$, and $w_j\triangleq 1-(2-1/n)v$. Applying \eqref{gn11} to $\|f\|_{L^p}$ in \eqref{records:ball} and using $2(1-v)/p=a_j$, $(n^8 R_j^{1/2}D_{R_j})^v R_j^{(B_h-1)a_j/2}\le C$, and $T_2^{v/n}\le C$, we deduce
\be\nonumber\ba
 \|f\|_{L^\infty}\le C\sqrt{p R_j}\,A_1^{a_j}T_2^{w_j}T_3^v .
\ea\ee
Set $\widetilde d_j\triangleq L_j^{-1}\int_{E_j}T_2^2/S_0\,dt\le C(d_j+|E_j|/L_j)$ and $x_j\triangleq(e_j\widetilde d_j)^{1/2}$. By \eqref{exthigh:sum}, $\sum_j x_j\le C$ and $\int_{E_j}T_2\,dt\le\sqrt{L_j}\,x_j$. Since $a_j/2+w_j+3v/2=1$ and $\int_{E_j}A_1^2\,dt\le C\int_J S_0\,dt\le C$, H\"older's inequality yields
\be\nonumber\ba
 \int_{E_j}\|f\|_{L^\infty}\,dt
 \le C\sqrt{p R_j}\bigl(\sqrt{L_j}\,x_j\bigr)^{w_j}R_j^{\frac{3 c_3 v}{2}}
 \le C L_j x_j^{w_j}.
\ea\ee
As $1-w_j=O(v)=O(1/\log R_j)=O(1/j)$, this proves \eqref{per:exceptional-flux}.
\end{proof}

\begin{proposition}[Periodic density bound and continuation]\label{order:periodic-cap}
Let $(\rho,u,\theta)$ be a smooth positive-density periodic solution with the data and coefficients of Theorem~\ref{th1} on $[0,T]\subset[0,T_0]$, $T_0<\infty$. Then there exist positive constants $c,C$, independent of the initial density floor, such that
\be\la{construction:common-cap}\ba
 &\inf_{\T^2\times[0,T]}\theta\ge c,\\
 &\sup_{0\le t\le T}\bigl(\|\rho\|_{L^\infty}+A_1+A_2+\|u\|_{H^1}
   +\|K(\theta)\|_{H^1}+\|\theta\|_{H^1}\bigr)\\
 &\quad +\sup_{\substack{I\subset[0,T]\\ |I|\le 1}}
   \int_I\bigl(\|\nabla\dot u\|_{L^2}^2+B\bigr)\,dt\le C .
\ea\ee
The constants are independent of $T_0$ in \textup{(A)}, \textup{(B)}, and \textup{(C)} with $\alpha>1-d$. In \eqref{ext:powerfamily}, for fixed $D,c_*,r_*,H_0,C_\lambda$ and physical parameters, the density ceiling is common to all data satisfying \eqref{ext:dataM} and all indicated viscosities with $0\le\varepsilon\le\varepsilon_*(M)$.
The same uniformity holds in \eqref{ext:gainfamily} for fixed $d,\alpha,c_*,r_*,H_0,C_\lambda$ and physical parameters, with $0\le\varepsilon\le\varepsilon_{*,C}(M)$.

More generally, on a fixed horizon $[0,T_0]$, \eqref{construction:common-cap} holds under \eqref{ext:heatlaw} and \eqref{visc:regularity} whenever $\mu+\lambda(r)\ge cr/D(r)$ for $r\ge r_*$, $D\ge 1$ is nondecreasing,
\be\la{per:finite-input}\ba
 \lim_{r\to\infty}\frac{\log D(r)}{\log L(r)}=0,\quad
 \int_2^\infty\frac{dr}{(e+r)L(r)D(r)}=\infty ,
\ea\ee
and $U+\mathcal P_l\le C_{T_0}$ holds on every smooth prefix with a common cutoff $l$ and bound independent of the provisional ceiling and initial density floor.
\end{proposition}

\begin{proof}
Fix an index $j$ satisfying $D(r_j)\le(\log r_j)^{\zeta_*}$ for which $r_j$ is reached by time $T$.
At $t_j$, choose a point where $\rho(\cdot,t_j)=r_j$ and follow its trajectory backwards to the last time the density equals $r_j/2$.
The resulting crossing interval $I_j$ lies in $I_j^0$, since $r_{j^-}\le r_j/2$ when $j^-$ exists and $\|\rho_0\|_{L^\infty}<r_j/2$ otherwise. In the time-uniform range, \eqref{extlow:templower} gives $P\ge c r_j/D_{R_j}$ on $I_j$. Under \eqref{per:finite-input} (which includes the fixed-law alternative of \textup{(C)} with $D\equiv 1$ and \textup{(D)} via Lemma~\ref{lem:linear-moment-input}), \eqref{extlow:NvsU} and Lemma~\ref{low:temp-floor} yield the fixed density moments and \eqref{exthigh:pressure}--\eqref{exthigh:tails}, while \eqref{fp117} with $a_0\le C D_{R_j}$ gives $P\ge c_{T_0}r_j/D_{R_j}$ on $I_j$. Under \eqref{per:finite-input}, $D=L^{o(1)}$ implies that all sufficiently large indices satisfy the stated bound on $D(r_j)$, and Gr\"onwall's inequality in \eqref{exthigh:joint} with $\int_0^T g_*\,dt\le C_{T_0}$ gives
\be\nonumber\ba
 \sup_{0\le t\le T}W_R(t)+\int_0^T\mathcal D_R\,dt
 \le C_{T_0}L_R\exp\bigl(C_{T_0}D_R^{C_{T_0}}\bigr)=R^{1+o(1)}.
\ea\ee

Integrating $D_t\Lambda(\rho)=-P-G$ along the trajectory over $I_j$ and applying Lemma~\ref{records:flux}, we derive
\be\nonumber\ba
 \Lambda(r_j)-\Lambda(r_j/2)+\frac{c r_j}{D_{R_j}}|I_j|
 \le C R_j^\xi(1+|I_j|)+C L_j x_j^{w_j}.
\ea\ee
Since $D_{R_j}=1+D(r_j)\le 1+(\log r_j)^{\zeta_*}$ and $\xi<1$, we have $C R_j^\xi\le c r_j/D_{R_j}$ for all sufficiently large $j$ under consideration. Thus
\be\la{order:crossing}\ba
 \frac{\Lambda(r_j)-\Lambda(r_j/2)}{L_j}
 \le \frac{C R_j^\xi}{L_j}+C x_j^{w_j}.
\ea\ee
Here $\sum_j R_j^\xi/L_j<\infty$, while $1-w_j=O(1/j)$ implies $x_j^{w_j}\le C x_j$ if $x_j\ge e^{-j}$ and $x_j^{w_j}\le e^{-j/2}$ if $x_j<e^{-j}$ (for large $j$). By \eqref{per:exceptional-flux}, whenever $r_N$ is reached by time $T$, summing the left-hand side of \eqref{order:crossing} over $j\le N$ with $D(r_j)\le(\log r_j)^{\zeta_*}$ gives a bound independent of $N$.

On the other hand, under \textup{(A)}, \textup{(C)}, or \eqref{per:finite-input}, \eqref{extlow:struct} bounds each normalized increment below by $c/[j D(r_j)]$. By monotonicity and \eqref{new:den} or \eqref{per:finite-input}, the series $\sum_{j\ge1}1/[jD(r_j)]$ diverges. Its terms with $D(r_j)>(\log r_j)^{\zeta_*}$ have finite sum, bounded by $C\sum_{j\ge1}j^{-1-\zeta_*}$.
Thus the subseries over indices $j$ satisfying $D(r_j)\le(\log r_j)^{\zeta_*}$ still diverges.
Under \textup{(B)}, all sufficiently large indices satisfy this last inequality, and each increment is comparable to $\int_{r_j/2}^{r_j}\frac{h(z)}{z^2 L(z)}\,dz$, whose sum diverges by \eqref{new:fixed} and Tonelli's theorem (since $\int_z^\infty\frac{dr}{r^2 L(r)}\asymp\frac{1}{z L(z)}$). Hence some finite level is excluded in every case.

For the data-class family \eqref{ext:powerfamily}, first choose an excluded level $R_*=R_*(M)>\max\{r_*,4M+4\}$ by the preceding argument with lower-bound constant $c_*/2$ and data bound $4M+4$, keeping $D,c_*,r_*,H_0,C_\lambda$ and the physical parameters fixed. This choice is independent of $\varepsilon$, $T$, and the particular viscosity law, since the preceding estimates use the lower bound only on the attained density range. Setting
\be\la{per:family-level}\ba
 \varepsilon_*\triangleq\min\left\{1,\ \frac{\log 2}{\log(1+R_*)}\right\}
\ea\ee
ensures, for $0\le\varepsilon\le\varepsilon_*$, that
\be\nonumber\ba
 &(1+r)^\varepsilon\le(1+R_*)^{\varepsilon_*}\le2,\\
 &h(r)\ge\frac{c_*r}{(1+r)^\varepsilon D(r)}
       \ge\frac{c_*r}{2D(r)}\qquad(r_*\le r\le R_*).
\ea\ee
Suppose that $\|\rho(t)\|_{L^\infty}$ first reaches $R_*$ at some finite time $t_*$. On $[0,t_*]$, all the preceding estimates apply with lower-bound constant $c_*/2$, and summing \eqref{order:crossing} up to $R_*$ contradicts the choice of this excluded level. Thus $R_*$ is again excluded uniformly, without imposing \eqref{ext:iter} above $R_*$.

For \eqref{ext:gainfamily}, fix $d$ and $\alpha>1-d$, and take $D\equiv1$. Apply the preceding time-uniform argument on attained ranges where \eqref{new:gain} holds with $c_*/2$. By \eqref{visc:envelope}, Lemma~\ref{low:uniform}, and the initial bounds entering \eqref{exthigh:window}, its constants are controlled by the fixed quantities in \textup{(C)} and the data bound. Thus, using data bound $4M+4$, first choose a corresponding excluded level $R_{*,C}=R_{*,C}(M)>\max\{r_*,4M+4\}$, independently of $\varepsilon$, $T$, and the particular viscosity law. Set
\begin{equation}\label{per:gain-family-level}
 \varepsilon_{*,C}(M)\triangleq
 \min\left\{1,\ \frac{\log2}{\log(1+R_{*,C})}\right\}.
\end{equation}
Then \eqref{ext:gainfamily} gives
\[
 \mu+\lambda(r)\ge\frac{c_*rL(r)^d}{(1+r)^\varepsilon}
 \ge\frac{c_*}{2}rL(r)^d
 \qquad(r_*\le r\le R_{*,C},\quad0\le\varepsilon\le\varepsilon_{*,C}(M)).
\]
If $\|\rho(t)\|_\infty$ first reached $R_{*,C}$ at some finite time, all preceding estimates, including \eqref{new:smalldiss} with $\alpha>1-d$, would hold up to that time with lower-bound constant $c_*/2$. Summing \eqref{order:crossing} with $D\equiv1$ would contradict the choice of $R_{*,C}$. Hence this level is excluded uniformly, without imposing \eqref{new:gain} above $R_{*,C}$.

Finally, with the density ceiling $R$ fixed, \eqref{exthigh:window} bounds $Y_R$, $S$, and $\int_I(D_w+B)\,dt\le (\sup_I Y_R)\int_I\mathcal D_R\,dt\le C$, and \eqref{common:material-gradient} then bounds $\int_I\|\nabla\dot u\|_{L^2}^2\,dt$. Combining \eqref{extlow:templower}, \eqref{ltp-low-u}, Lemma~\ref{1pi} on $\Gamma_t$, and $\|\nabla\theta\|_{L^2}\le Z/\kappa(c)$ gives \eqref{construction:common-cap}. This completes the proof of Proposition~\ref{order:periodic-cap}.
\end{proof}

\subsection{Higher-order estimates}

Assume the periodic hypotheses of Theorem~\ref{th1} on $[0,T]\subset[0,T_0]$. All estimates below are independent of the lower bound for the initial density in the smooth approximations.

\begin{lemma}\label{fpl9}
For any $0<\tau<T$,
\be\la{fp09}\ba
 \sup_{\tau\le t\le T}\bigl(\|\theta\|_{H^2}+\|\sqrt\rho\,\theta_t\|_{L^2}\bigr)
 +\sup_{\substack{I\subset[\tau,T]\\ |I|\le 1}}\int_I
   \bigl(\|\dot\theta\|_{H^1}^2+\|\dot u\|_{H^1}^2\bigr)\,dt\le C_{\tau,T_0}.
\ea\ee
Under \textup{(A)}, \textup{(B)}, or \textup{(C)} with $\alpha>1-d$, the constant depends on $\min\{1,\tau\}$ and is independent of $T_0$.
\end{lemma}

\begin{proof}
By \eqref{dc1}, \eqref{2fp32}, and \eqref{construction:common-cap}, we have for any fixed $2\le p<\infty$,
\be\la{fp103}\ba
 \|\nabla u\|_{L^p}\le C_p\bigl(\|G\|_{H^1}+\|\omega\|_{H^1}+\|\rho\theta\|_{L^p}\bigr)
 \le C_p(1+A_2).
\ea\ee
Since $\|\rho\|_{L^\infty}\le C$ and $M_0>0$, the mass-weighted Poincar\'e inequality holds:
\be\la{ltp-wpoin}\ba
 \left\|f-M_0^{-1}\int\rho f\,dx\right\|_{L^p}\le C_p\|\nabla f\|_{L^2}\quad (2\le p<\infty).
\ea\ee
Together with $\int\rho u\,dx=\int\rho\dot u\,dx=0$ and \eqref{construction:common-cap}, \eqref{fp103}--\eqref{ltp-wpoin} bound $\|u\|_{W^{1,p}}$ uniformly in time and $\int_I\|\dot u\|_{H^1}^2\,dt$ for $|I|\le 1$.

Set $k\triangleq K(\theta)$, $v\triangleq D_t k=\kappa(\theta)\dot\theta$, and $m_\kappa\triangleq\rho/\kappa(\theta)$, so that $B=\int m_\kappa v^2\,dx$. By \eqref{construction:common-cap}, $\|m_\kappa\|_{L^\infty}\le C$, and for a fixed large $L>0$, mass and energy conservation give
\be\la{fp93}\ba
 \int m_\kappa\,dx
 \ge \frac{M_0-L^{-1}\int\rho\theta\,dx}{\max_{[0,L]}\kappa}
 \ge \frac{M_0}{2\max_{[0,L]}\kappa}.
\ea\ee
Subtracting the spatial mean of $v$ and applying \eqref{gn11} yield
\be\la{fp96}\ba
 \|v\|_{L^p}\le C\bigl(\sqrt B+p^{1/2}\|\nabla v\|_{L^2}\bigr),\quad p\ge 2.
\ea\ee

Applying $\partial_t+\div(u\,\cdot)$ to $(\ref{fns})_3$ gives
\be\la{fp97}\ba
 c_v\rho D_t\dot\theta-\Delta v
 =\partial_i C_i+D_t Q+Q\div u-\rho\dot\theta\div u-P D_t(\div u),
\ea\ee
where $C_i\triangleq(\div u)\partial_i k-(\partial_i u^j+\partial_j u^i)\partial_j k$. Multiplying \eqref{fp97} by $v=\kappa(\theta)\dot\theta$ and integrating by parts over $\T^2$ produce $\frac{c_v}{2}B'+\|\nabla v\|_{L^2}^2$ on the left-hand side and the cubic term $\frac{c_v}{2}\int\rho\kappa'(\theta)\dot\theta^3\,dx$ on the right-hand side. Since $|\kappa'(\theta)|/\kappa(\theta)^2\le C$ above the temperature floor, \eqref{fp96}, and Young's inequality yield, for $p\ge 5$,
\be\la{fp910}\ba
 \left|\int\rho\kappa'(\theta)\dot\theta^3\,dx\right|
 &\le C B^{\frac{p-3}{p-2}}\|v\|_{L^p}^{\frac{p}{p-2}}\\
 &\le \varepsilon\|\nabla v\|_{L^2}^2
 +C_\varepsilon\left(1+B^2+p^{\frac{p}{p-4}}B^{\frac{2(p-3)}{p-4}}\right).
\ea\ee
Choosing $p=4+\log(e+B)$ bounds the last parenthesis by $C\bigl(1+B^2\log(e+B)\bigr)$.
Moreover, $(\ref{fns})_3$ and \eqref{fp103} imply
\be\la{fp912}\ba
 \|k\|_{H^2}\le C\bigl(1+\|\Delta k\|_{L^2}\bigr)\le C\bigl(1+\sqrt B\bigr),
\ea\ee
so
$$\left|\int C_i\partial_i v\,dx\right|\le C\|\nabla u\|_{L^4}\|\nabla k\|_{L^4}\|\nabla v\|_{L^2}\le\varepsilon\|\nabla v\|_{L^2}^2+C_\varepsilon(1+B).$$
Using $D_t\nabla u=\nabla\dot u-(\nabla u)^2$, $D_t\lambda=-\rho\lambda'\div u$, \eqref{fp103}, and \eqref{fp96}, we bound the remaining tested source terms in \eqref{fp97} by
\be\la{fp914}\ba
 C\bigl(1+\|\nabla\dot u\|_{L^2}+\sqrt B\bigr)\|v\|_{L^4}
 \le \varepsilon\|\nabla v\|_{L^2}^2+C_\varepsilon\bigl(1+B+\|\nabla\dot u\|_{L^2}^2\bigr).
\ea\ee
Combining the above bounds yields
\be\la{fp915}\ba
 \frac{c_v}{2}B'+\frac{1}{2}\|\nabla v\|_{L^2}^2
 \le C\bigl(1+B^2\log(e+B)+\|\nabla\dot u\|_{L^2}^2\bigr).
\ea\ee

For each $t\in[\tau,T]$, set $\delta_t\triangleq\min\{1,t\}$. By \eqref{construction:common-cap}, there exists $s\in[t-\delta_t,t-\delta_t/2]$ with $B(s)\le 2C/\delta_t$. Dividing \eqref{fp915} by $(e+B)\log(e+B)$, integrating over $[s,t]$, and using \eqref{construction:common-cap} bound $\log\log(e+B(t))$ uniformly for $t\ge\tau$. Integrating \eqref{fp915} once more gives
\be\la{fp918}\ba
 \sup_{\tau\le t\le T}B(t)
 +\sup_{\substack{I\subset[\tau,T]\\ |I|\le 1}}\int_I\|\nabla v\|_{L^2}^2\,dt
 \le C_{\tau,T_0}.
\ea\ee
By \eqref{fp912}, Sobolev embedding, and the regularity of $K^{-1}$ above the temperature floor,
\be\la{fp920}\ba
 \|\theta\|_{L^\infty}+\|\theta\|_{H^2}\le C_{\tau,T_0}\quad (t\ge\tau).
\ea\ee
From $\theta_t=v/\kappa(\theta)-u\cdot\nabla\theta$ we obtain $\|\sqrt\rho\,\theta_t\|_{L^2}\le C_{\tau,T_0}$, and from $\|\nabla\dot\theta\|_{L^2}\le C(\|\nabla v\|_{L^2}+\|v\|_{L^4}\|\nabla\theta\|_{L^4})$, \eqref{fp96}, and \eqref{fp918} we derive
\be\la{fp925}\ba
 \sup_{\substack{I\subset[\tau,T]\\ |I|\le 1}}\int_I\|\dot\theta\|_{H^1}^2\,dt
 \le C_{\tau,T_0},
\ea\ee
which gives \eqref{fp09} and completes the proof of Lemma \ref{fpl9}.
\end{proof}

\begin{lemma}\label{init:common-lemma}
For the exponent $q>2$ in \eqref{th12} and every finite $T>0$,
\be\la{init:prelayer}\ba
 &\sup_{0\le t\le T}\bigl(\|\rho\|_{W^{1,q}}+\|\rho_t\|_{L^q}+\|u\|_{H^2}
   +\|\sqrt\rho\,u_t\|_{L^2}+\|\theta\|_{H^1}+\|\theta\|_{L^\infty}\bigr)\\
 &\quad +\int_0^T\bigl(\|u\|_{W^{2,q}}^2+\|u_t\|_{H^1}^2
   +\|\theta\|_{H^2}^2+\|\sqrt\rho\,\theta_t\|_{L^2}^2\bigr)\,dt\le C_T .
\ea\ee
Moreover, for every $0<\tau<T$,
\be\la{fp010}\ba
\int_\tau^T\bigl(\|\theta\|_{H^3}^2+\|\theta_t\|_{H^1}^2\bigr)\,dt\le C_{\tau,T}.
\ea\ee
\end{lemma}

\begin{proof}
By \eqref{fp912} and \eqref{construction:common-cap}, we have $\int_0^T\|k\|_{H^2}^2\,dt\le C_T$ and $\sup_{[0,T]}\|k\|_{H^1}\le C$. Using $\|\nabla k\|_{L^4}^4\le C\|\nabla k\|_{L^2}^2\|k\|_{H^2}^2$ and composing with $K^{-1}$ give $\theta\in L^2(0,T;H^2)$, while $\theta_t=\dot\theta-u\cdot\nabla\theta$ bounds $\sqrt\rho\,\theta_t$ in $L^2(0,T;L^2)$.

Set $F\triangleq 1+\|\dot u\|_{H^1}+\|\theta\|_{H^2}$ and $Y\triangleq\|\nabla\rho\|_{L^q}$. By \eqref{2fp32}, $\|G\|_{W^{1,q}}+\|\omega\|_{W^{1,q}}\le C F$. Differentiating $h\div u=G+\rho\theta$ yields
\be\nonumber\ba
 \nabla\div u=\frac{1}{h}\bigl(\nabla G+\rho\nabla\theta+(\theta-\lambda'\div u)\nabla\rho\bigr),
\ea\ee
whence $\|\div u\|_{L^\infty}+\|\omega\|_{L^\infty}\le C F$ and $\|\nabla\div u\|_{L^q}+\|\nabla^2 u\|_{L^q}\le C F(1+Y)$. Applying Lemma~\ref{bkm} gives
\be\la{fp107}\ba
 \|\nabla u\|_{L^\infty}\le C F^2\log(e+Y).
\ea\ee
Differentiating $(\ref{fns})_1$ gives $Y'\le C\|\nabla u\|_{L^\infty}Y+C\|\nabla\div u\|_{L^q}$, so that $\bigl(\log(e+Y)\bigr)'\le C F^2\log(e+Y)$.
By Gr\"onwall's inequality,
\be\la{fp1010}\ba
 \sup_{0\le t\le T}Y(t)+\int_0^T\bigl(\|\nabla^2 u\|_{L^q}^2+\|\nabla u\|_{L^\infty}\bigr)\,dt\le C_T.
\ea\ee
Applying H\"older's inequality and \eqref{fp103} in the identity for $\nabla\div u$ bounds $\|(\theta-\lambda'\div u)\nabla\rho\|_{L^2}\le C_T$, and \eqref{2fp32} then yields $\sup_{[0,T]}\|u\|_{H^2}\le C_T$. The bounds for $\rho_t$, $\sqrt\rho\,u_t$, and $u_t$ follow directly from $(\ref{fns})_1$ and $u_t=\dot u-u\cdot\nabla u$.

Next, to bound $\|\theta\|_{L^\infty}$ on $[0,\min\{1,T\}]$, let $\psi$ be the periodic solution to $(-\Delta+\rho)\psi=Q$. By \eqref{ltp-wpoin} and $Q\ge 0$, the operator is coercive and $\psi\ge 0$. Since $\|Q\|_{L^4}\le C\|\nabla u\|_{L^8}^2\le C_T$, elliptic regularity gives $\sup_{[0,T]}\|\psi\|_{W^{2,4}}\le C_T$. Differentiating in time yields $(-\Delta+\rho)\psi_t=Q_t-\rho_t\psi$, with
\be\nonumber\ba
 \|Q_t-\rho_t\psi\|_{L^{4/3}}
 \le C\|\nabla u\|_{L^4}\|\nabla u_t\|_{L^2}
 +C\|\rho_t\|_{L^2}\bigl(\|\nabla u\|_{L^8}^2+\|\psi\|_{L^4}\bigr)
 \in L^2(0,T),
\ea\ee
so $\psi_t\in L^2(0,T;L^\infty)$. Define $\bar\theta\triangleq K^{-1}(\xi+\psi)$, where
\be\nonumber\ba
 \xi(0)\triangleq 1+\|K(\theta_0)\|_{L^\infty},\quad
 \xi'=C_*(1+\xi)\bigl(1+\|\psi_t\|_{L^\infty}+\|\div u\|_{L^\infty}\bigr).
\ea\ee
By \eqref{fp1010}, $\xi$ stays bounded. Since $\kappa(\bar\theta)(1+\bar\theta)\le C(1+\xi)$ and $-\Delta K(\bar\theta)=Q-\rho\psi$, choosing $C_*$ sufficiently large gives
\be\nonumber\ba
 c_v\rho D_t\bar\theta-\Delta K(\bar\theta)+\rho\bar\theta\div u-Q
 =\rho\left(\frac{c_v(\xi'+D_t\psi)}{\kappa(\bar\theta)}-\psi+\bar\theta\div u\right)\ge 0.
\ea\ee
Since $\{\theta>\bar\theta\}=\{K(\theta)>K(\bar\theta)\}$, subtracting $(\ref{fns})_3$ and applying the positive-part inequality for the Laplacian together with $(\ref{fns})_1$ yield
\be\nonumber\ba
 c_v\frac{d}{dt}\int\rho(\theta-\bar\theta)_+\,dx
 \le \|\div u\|_{L^\infty}\int\rho(\theta-\bar\theta)_+\,dx.
\ea\ee
As $(\theta_0-\bar\theta(0))_+=0$, we obtain $\theta\le\bar\theta$ on $[0,\min\{1,T\}]$. Together with \eqref{fp920} for $t\ge 1$, this establishes \eqref{init:prelayer}.

Finally, differentiating $\Delta k=c_v\rho\dot\theta+P\div u-Q$ in space and applying \eqref{fp09}, \eqref{fp1010}, and H\"older's inequality for $t\ge\tau$ give
\be\la{fp1015}\ba
 \|\nabla^3 k\|_{L^2}\le C\|\nabla\Delta k\|_{L^2}
 \le C_{\tau,T}\bigl(1+\|\dot\theta\|_{H^1}+\|\nabla^2 u\|_{L^q}\bigr).
\ea\ee
Since the first three derivatives of $K^{-1}$ are bounded on the attained temperature range and $k$ is bounded in $H^2$, we have $\|\theta\|_{H^3}\le C_{\tau,T}(1+\|k\|_{H^3})$. Differentiating $\theta_t=\dot\theta-u\cdot\nabla\theta$ and using \eqref{fp925} yield
\be\la{fp1018}\ba
 \int_\tau^T\bigl(\|\theta\|_{H^3}^2+\|\theta_t\|_{H^1}^2\bigr)\,dt\le C_{\tau,T},
\ea\ee
which yields \eqref{fp010} and finishes the proof of Lemma \ref{init:common-lemma}.
\end{proof}

\begin{lemma}\label{reg:tor-decay-norms}
Assume \textup{(A)}, \textup{(B)}, or \textup{(C)} with $\alpha>1-d$. Then for each $1\le p<\infty$, there exist constants $C_p,\eta_p>0$ such that, for all $t\ge 1$,
\be\nonumber\ba
 \|\rho(t)-\ol\rho_0\|_{L^p}+\|u(t)-u_\infty\|_{W^{1,p}}
 +\|\theta(t)-\theta_\infty\|_{W^{1,p}}\le C_p e^{-\eta_p t}.
\ea\ee
\end{lemma}

\begin{proof}
Working in the zero-momentum frame ($u_\infty=0$), set $r\triangleq\rho-\ol\rho_0$ and $E_k\triangleq\frac{1}{2}\int\rho|u|^2\,dx$. Energy conservation gives $c_v M_0(\theta_\infty-m)=E_k$. Define the relative free energy
\be\nonumber\ba
 \mathcal F\triangleq E_k
 +c_v\int\rho\left(\theta-\theta_\infty-\theta_\infty\log\frac{\theta}{\theta_\infty}\right)dx
 +\theta_\infty\int\left(\rho\log\frac{\rho}{\ol\rho_0}-\rho+\ol\rho_0\right)dx.
\ea\ee
By the time-uniform bounds \eqref{construction:common-cap} and \eqref{fp920} (at $\tau=1$), we have
\be\nonumber\ba
 \mathcal F\asymp E_k+\int\rho|\theta-\theta_\infty|^2\,dx+\|r\|_{L^2}^2,
\ea\ee
and the entropy identity \eqref{fp118} yields
\be\nonumber\ba
 -\mathcal F'
 =\theta_\infty\int\left(\frac{Q}{\theta}+\frac{\kappa(\theta)|\nabla\theta|^2}{\theta^2}\right)dx
 \ge c\bigl(\|\nabla u\|_{L^2}^2+\|\nabla\theta\|_{L^2}^2\bigr).
\ea\ee

To recover density dissipation, set $w\triangleq\nabla\Delta^{-1}r$ and $I\triangleq\int\rho u\cdot w\,dx$. Since $\div w=r$ and $r_t=-\div(\rho u)$, we have $\|w\|_{H^1}\le C\|r\|_{L^2}$ and $\|w_t\|_{L^2}\le C\|\rho u\|_{L^2}$. Multiplying $(\ref{fns})_2$ by $w$ and integrating by parts over $\T^2$ give
\be\nonumber\ba
 I'=\int\bigl(P r+\rho u\otimes u:\nabla w-\mathbb S:\nabla w+\rho u\cdot w_t\bigr)\,dx .
\ea\ee
Using Poincar\'e's inequality, \eqref{ltp-wpoin}, and Young's inequality, we arrive at
\be\nonumber\ba
\int P r\,dx
=\int\theta r^2\,dx+\ol\rho_0\int(\theta-\ol\theta)r\,dx
\ge c\|r\|_{L^2}^2-C\|\nabla\theta\|_{L^2}^2,
\ea\ee
and
\be\nonumber\ba
 \left|\int\bigl(\rho u\otimes u:\nabla w-\mathbb S:\nabla w+\rho u\cdot w_t\bigr)\,dx\right|
\le \varepsilon\|r\|_{L^2}^2+C_\varepsilon\|\nabla u\|_{L^2}^2.
\ea\ee
Since $|I|\le C\sqrt{E_k}\|r\|_{L^2}$, choosing $\varepsilon,\delta>0$ sufficiently small yields
\be\nonumber\ba
 \mathcal F-\delta I\asymp\mathcal F,\quad
 (\mathcal F-\delta I)'+c\bigl(\|r\|_{L^2}^2+\|\nabla u\|_{L^2}^2+\|\nabla\theta\|_{L^2}^2\bigr)\le 0.
\ea\ee
Because $E_k\le C\|\nabla u\|_{L^2}^2$ and
\be\nonumber\ba
 \int\rho|\theta-\theta_\infty|^2\,dx
 =\int\rho|\theta-m|^2\,dx+M_0(m-\theta_\infty)^2
 \le C\|\nabla\theta\|_{L^2}^2+C E_k^2,
\ea\ee
with $E_k$ bounded, Gr\"onwall's inequality yields $\mathcal F(t)\le C e^{-\gamma t}$ for $t\ge 1$.

Finally, for $f=u$ or $\theta-\theta_\infty$, we have
\be\nonumber\ba
 \ol\rho_0\|f\|_{L^2}^2\le\int\rho|f|^2\,dx+\|r\|_{L^2}\|f\|_{L^4}^2,
\ea\ee
where $\|f\|_{L^4}$ is uniformly bounded by \eqref{fp103}, \eqref{ltp-wpoin}, and \eqref{fp09}. Hence $u$, $\theta-\theta_\infty$, and $P-\ol\rho_0\theta_\infty$ decay exponentially in $L^2$. Testing $(\ref{fns})_2$ by $u$ gives
\be\nonumber\ba
 \mu\|\nabla u\|_{L^2}^2+\int(\mu+\lambda)(\div u)^2\,dx
 =\int(P-\ol\rho_0\theta_\infty)\div u\,dx-\int\rho\dot u\cdot u\,dx,
\ea\ee
while $\|\nabla\theta\|_{L^2}^2=-\int(\theta-\theta_\infty)\Delta\theta\,dx$. Since $\|\nabla u\|_{L^2}+A_2+\|\theta\|_{H^2}\le C$ for $t\ge 1$, both $\|\nabla u\|_{L^2}$ and $\|\nabla\theta\|_{L^2}$ decay exponentially. Interpolating with the uniform density ceiling and the $W^{1,s}$ bounds for $u,\theta$ gives the decay for $p\ge 2$, and H\"older's inequality covers $1\le p<2$. Galilean invariance restores the general frame. This completes the proof of Lemma~\ref{reg:tor-decay-norms}.
\end{proof}

\section{A priori estimates on a bounded domain}\label{tri:bd}

Let $\Omega$ and the initial data satisfy the assumptions of Theorem~\ref{tri:main}, and let $(\rho,u,\theta)$ be a smooth positive-density strong solution to \eqref{fns}, \eqref{cztj}, and \eqref{tri:bd-boundary} on $\Omega\times(0,T]$.

Set $u_\tau\triangleq u\cdot\tau$, $h\triangleq 2\mu+\lambda$, and $G\triangleq h\div u-P$, and let a bar denote the spatial average on $\Omega$. With $Q$, $\mathcal D_\theta$, and $R_T$ as in \eqref{fp52} and \eqref{mdsj}, we define
\be\nonumber\ba
 D_\partial&\triangleq\mu\int_{\partial\Omega}A u_\tau^2\,dS,&
 A_2^2&\triangleq\int_\Omega\rho|\dot u|^2\,dx,\\
 A_1^2&\triangleq 1+\int_\Omega\left(\frac{G^2}{h}+\mu\omega^2\right)dx+D_\partial,&
 D_u&\triangleq\int_\Omega\bigl(h(\div u)^2+\mu\omega^2\bigr)\,dx+D_\partial .
\ea\ee
Recall that on the simply connected $C^{3,1}$ domain $\Omega$, any vector field $v\in H^1(\Omega;\rr)$ with $v\cdot n=0$ on $\partial\Omega$ satisfies the tangent div--curl and Poincar\'e estimates $\|v\|_{H^1}\le C\|\nabla v\|_{L^2}\le C(\|\div v\|_{L^2}+\|\omega(v)\|_{L^2})$.

\subsection{Low-order estimates}

\begin{lemma}\label{bdthree:wall}
The total mass is conserved, $\int_\Omega\rho\,dx=M_0$, and the total energy $E(t)\triangleq\int_\Omega\bigl(\frac{1}{2}\rho|u|^2+c_v\rho\theta\bigr)\,dx$ satisfies
\be\la{wallthree:sign}\ba
 E(t)+\mu\int_0^t\int_{\partial\Omega}(A+2\chi)u_\tau^2\,dS\,ds=E(0).
\ea\ee
Moreover, there exists a constant $C>0$, independent of $T$ and the initial density floor, such that
\be\la{bd41:entropy}\ba
 \sup_{0\le t\le T}\int_\Omega
   \rho\bigl(|u|^2+\theta+(-\log\theta)_++\log(1+\rho)\bigr)\,dx
 +\int_0^T\int_\Omega
   \left(\frac{Q}{\theta}+\frac{\kappa(\theta)|\nabla\theta|^2}{\theta^2}\right)dx\,dt
 \le C.
\ea\ee
If $\sup_{\Omega\times[0,T]}\rho\le r$, then
\be\la{yjywdxj}\ba
 \theta(x,t)\ge(\underline\theta_0^{-1}+a_r t)^{-1},\quad
 a_r\triangleq\frac{1}{4c_v}\sup_{0\le z\le r}\frac{z}{\mu+\lambda(z)} .
\ea\ee
\end{lemma}

\begin{proof}
Since $u\cdot n=0$ and $\partial_\tau n=\chi\tau$ on $\partial\Omega$, the Navier-slip condition \eqref{tri:bd-boundary} gives $2D(u)n\cdot\tau=\omega+2(\partial_\tau u)\cdot n=-(A+2\chi)u_\tau$. Multiplying $(\ref{fns})_2$ by $u$ and integrating by parts yield
\be\nonumber\ba
 D_u=\int_\Omega Q\,dx+\mu\int_{\partial\Omega}(A+2\chi)u_\tau^2\,dS,\quad
 \frac{1}{2}\frac{d}{dt}\int_\Omega\rho|u|^2\,dx+D_u=\int_\Omega P\div u\,dx .
\ea\ee
Since $\partial_n K(\theta)=0$ on $\partial\Omega$, adding the integrated heat equation $(\ref{fns})_3$ gives \eqref{wallthree:sign}.
Similarly, testing $(\ref{fns})_3$ by $-\theta^{-1}$ and using $(\ref{fns})_1$ show that $\mathscr S\triangleq\int_\Omega(c_v\rho\log\theta-\rho\log\rho)\,dx$ satisfies
\be\nonumber\ba
 \mathscr S'(t)=\int_\Omega\frac{Q}{\theta}\,dx+\mathcal D_\theta(t).
\ea\ee
Combining this with \eqref{wallthree:sign} and $A+2\chi\ge 0$ yields \eqref{bd41:entropy} as in Lemma~\ref{fpl1}.

Finally, since $u\cdot n=\partial_n K(\theta)=0$ on $\partial\Omega$, the temperature comparison argument in Lemma~\ref{fpl1} applies verbatim and yields \eqref{yjywdxj}.
\end{proof}

To handle the density transport on $\Omega$, for any vector field $v$ let $\mathcal N(v)$ denote the unique zero-mean solution of the Neumann problem
\be\nonumber\ba
 \begin{cases}
  \Delta\mathcal N(v)=\div v & \text{in }\Omega,\\
  \partial_n\mathcal N(v)=v\cdot n & \text{on }\partial\Omega,
 \end{cases}\quad
 \int_\Omega\mathcal N(v)\,dx=0.
\ea\ee
Define the auxiliary potentials and commutator
\be\la{bd42:potentials}\ba
 f_2\triangleq\mathcal N(\rho u),\quad
 f_1\triangleq\mathcal N\bigl(\div(\rho u\otimes u)\bigr),\quad
 F\triangleq f_1-u\cdot\nabla f_2 .
\ea\ee

\begin{lemma}\label{bdthree:mom}
With $H_A\triangleq\mathcal N(\mu\nabla^\perp\omega)$, the mass and momentum equations yield
\be\la{bd42:transport}\ba
 D_t\bigl(\Lambda(\rho)+f_2\bigr)
 =-F-\ol{\lambda\div u}+H_A+\ol P-P .
\ea\ee
Moreover, for every $q>2$,
\be\la{bd42:commutator}\ba
 \|F\|_{L^q}\le C_q\|\nabla u\|_{L^2}\|\rho u\|_{L^q},\quad
 \|H_A\|_{L^q}\le C_q\|\nabla u\|_{L^2}.
\ea\ee
\end{lemma}

\begin{proof}
Since $u\cdot n=0$ on $\partial\Omega$, we have $\partial_n f_2=0$ and $\partial_n f_1=\rho u_\tau(\partial_\tau u)\cdot n=-\chi\rho u_\tau^2$ on $\partial\Omega$. Applying $\mathcal N$ to the conservative momentum equation $\partial_t(\rho u)+\div(\rho u\otimes u)=\nabla G+\mu\nabla^\perp\omega$ and using $n\cdot\nabla^\perp\omega=-\partial_\tau\omega=\partial_\tau(A u_\tau)$ on $\partial\Omega$, we obtain
\be\la{bd42:flux}\ba
 \partial_t f_2+f_1=G-\ol G+H_A ,
\ea\ee
where $H_A=\mathcal N(\mu\nabla^\perp\omega)$ solves the harmonic Neumann problem
\be\la{bd42:harmonic}\ba
 \begin{cases}
  \Delta H_A=0 & \text{in }\Omega,\\
  \partial_n H_A=\mu\partial_\tau(A u_\tau) & \text{on }\partial\Omega,
 \end{cases}\quad
 \int_\Omega H_A\,dx=0.
\ea\ee
By standard Neumann regularity, the trace theorem, and the tangent Poincar\'e inequality,
\be\la{bd42:hbound}\ba
 \|\nabla H_A\|_{L^2}
 &\le C\|\partial_\tau(A u_\tau)\|_{H^{-1/2}(\partial\Omega)}
 \le C\|u\|_{H^1}\le C\|\nabla u\|_{L^2},\\
 \|H_A\|_{L^q}
 &\le C_q\|\nabla H_A\|_{L^2}\le C_q\|\nabla u\|_{L^2}\quad (1\le q<\infty).
\ea\ee

Next, for $q>2$, set $s\triangleq 2q/(q+2)$. By the Navier-slip commutator estimate of Fan and Huang \cite[Lemma 4.1]{FanHuang}, we have
\be\nonumber\ba
 \|\nabla F\|_{L^s}\le C_q\|\nabla u\|_{L^2}\|\rho u\|_{L^q}.
\ea\ee
Since $\ol{f_1}=0$ and $u\cdot n=0$ on $\partial\Omega$, integrating by parts gives
\be\nonumber\ba
 |\ol F|
 =\frac{1}{|\Omega|}\left|\int_\Omega f_2\div u\,dx\right|
 \le C\|\nabla u\|_{L^2}\|f_2\|_{L^2}
 \le C_q\|\nabla u\|_{L^2}\|\rho u\|_{L^q}.
\ea\ee
Combining the last two bounds with the Sobolev embedding proves the first inequality in \eqref{bd42:commutator}.

Finally, $(\ref{fns})_1$ gives $D_t\Lambda(\rho)=-(G+P)$, while $u\cdot n=0$ implies $\int_\Omega\div u\,dx=0$ and hence $\ol G=\ol{\lambda\div u}-\ol P$. Adding $u\cdot\nabla f_2$ to \eqref{bd42:flux} and eliminating $G$ yield \eqref{bd42:transport}.
\end{proof}

Choose $D$ as in Section~\ref{per-apriori} so that \eqref{extlow:struct} holds on $\Omega$ (in the data-class alternatives \eqref{ext:powerfamily} and \eqref{ext:gainfamily}, restrict to attained ranges where \eqref{ext:iter} and \eqref{new:gain}, respectively, hold with $c_*/2$). Let $\mathfrak d, U, N_p, \mathcal P_l$ denote the same quantities as in Section~\ref{per-apriori} with integrals over $\Omega$.

\begin{lemma}\label{bd:density-moment}
For each $p>2$, there exist positive constants $C_p,c_p$, independent of $T$, the initial density floor, and the provisional density ceiling, such that
\be\nonumber\ba
 \log N_p\le C_p(1+U)^{c_p}.
\ea\ee
\end{lemma}

\begin{proof}
From \eqref{bd41:entropy}, Chebyshev's inequality as in Lemma~\ref{fpl1} yields fixed constants $0<h_0<L_0$ and $M_2>0$ such that $|\Gamma_t|\ge M_0/(2M_2)$ and $|\{\theta\ge h_0\}|\ge M_0/(2M_2)$, whence the cold-set bound \eqref{extlow:cold} holds on $\Omega$.
With $\psi\triangleq -f_2$, Neumann regularity and Sobolev embedding give $\|\psi\|_{L^{2k}}\le C_k\|\rho u\|_{L^{2k/(k+1)}}\le C_k\|\rho\|_{L^k}^{1/2}$, so \eqref{extlow:moment} holds as in Lemma~\ref{low:moment-result}.
For any tangent field $g$ ($g\cdot n=0$ on $\partial\Omega$), the Neumann commutator $T_u(g)\triangleq\mathcal N(\div(u\otimes g))-u\cdot\nabla\mathcal N(g)$ obeys the $L^{q_\pm}$ bounds $\|T_u(g)\|_{L^{q_\pm}}\le C\|\nabla u\|_{L^2}\|g\|_{L^{q_\pm}}$ by \cite[Lemma~2.4]{FJLL2} (cf.\ \eqref{bd42:commutator}). Since level-set truncations of $\rho u$ preserve tangency and $\|u\|_{H^1}\le C\|\nabla u\|_{L^2}$, the same distribution-function interpolation and exponential velocity bound give \eqref{extlow:comm}.

Compared with \eqref{fp35}, \eqref{bd42:transport} contains the additional harmonic term $H_A$. Setting $\phi\triangleq(\Lambda(\rho)+f_2-N)_+$, $z\triangleq\rho\phi^{k-1}$, and $M_k\triangleq 1+\int_\Omega\rho\phi^k\,dx$, we apply \eqref{bd42:hbound} and \eqref{extlow:dualtest} to obtain
\be\nonumber\ba
 \int_\Omega z|H_A|\,dx
 \le \int_\Omega z\,dx + C M_k^{1-\frac{1}{k}+\frac{1}{k^2}}\|\nabla u\|_{L^2}^2.
\ea\ee
The first term merely shifts the constant $C_*$ in the pressure drift of Lemma~\ref{low:moment-result}, while the second is absorbed into $C M_k\|\nabla u\|_{L^2}^2$. Since $u\cdot n=0$ eliminates all boundary fluxes in the decreasing-exponent test, \eqref{low:fixed-moment-rate} and Gr\"onwall's inequality conclude the proof exactly as in Lemma~\ref{low:moment-result}.
\end{proof}

\begin{lemma}\label{bd:low-transfer}
There exists a common cutoff $l\ge 2L_0+2$ such that $U+\mathcal P_l\le C$ in \textup{(A)}, \textup{(B)}, and \textup{(C)} with $\alpha>1-d$ (on the attained ranges specified above), with $C$ independent of $T$. In \textup{(D)} and in the remaining range of \textup{(C)}, the same bound holds on $[0,T]\subset[0,T_0]$ with $l,C$ depending on $T_0$. In every case, \eqref{extlow:heattail}--\eqref{extlow:templower} hold on $\Omega$ with the corresponding time dependence.
\end{lemma}

\begin{proof}
Because $\partial_n\theta=0$ and $u\cdot n=0$ on $\partial\Omega$, all boundary fluxes vanish in the heat-equation cutoff tests of Lemmas~\ref{low:uniform}--\ref{low:temp-floor}, while the endpoint costs are bounded by $\sup_t\int_\Omega\rho\theta\,dx\le C$ via \eqref{wallthree:sign}.
Moreover, tangent div--curl and Lemma~\ref{bdthree:wall} give
\be\nonumber\ba
 \|\nabla u\|_{L^2}^2\le C D_u
 =C\int_\Omega Q\,dx+C\mu\int_{\partial\Omega}(A+2\chi)u_\tau^2\,dS,
\ea\ee
whose boundary term has uniformly bounded time integral by \eqref{wallthree:sign}. Consequently, \eqref{extlow:band} holds on $\Omega$, so the closure in Lemma~\ref{low:uniform} (including under \eqref{new:smallentropy}) carries over verbatim.
Likewise, zero normal flux yields \eqref{low:gain-window} under \textup{(C)}, and \eqref{crit:moment-rate} and \eqref{low:band-pressure} under \textup{(D)}, so Lemmas~\ref{lem:gain-moment-input} and~\ref{lem:linear-moment-input} give $U+\mathcal P_l\le C_{T_0}$. Finally, Lemma~\ref{bd:density-moment} controls every fixed density moment, and the hot and cold truncations of Lemma~\ref{low:temp-floor} yield \eqref{extlow:heattail}--\eqref{extlow:templower} on $\Omega$.
\end{proof}

Define
$$H\triangleq c_v\int_\Omega\rho\mathcal K(\theta)\,dx,\quad
Z\triangleq\|\nabla K(\theta)\|_{L^2},\quad
B\triangleq\int_\Omega\rho\kappa(\theta)|\dot\theta|^2\,dx,\quad
\Theta\triangleq 1+\|\theta\|_{L^{2\ell}}.$$
With the exponents of Lemma~\ref{core:pressure-times}, set
$$ S\triangleq 1+A_1^2+H,\quad S_0\triangleq 1+D_u+(1+H)^\sigma,\quad X \triangleq e+A_1^2+A_2^2+Z^2,$$
and define $g,g_*,J$ as in Section~\ref{per-apriori}.
For $r\ge\max\{r_*,\sup_{\Omega\times[0,T]}\rho\}$, put
$$ R\triangleq 2+r,\quad L_R\triangleq R\log(e+R),\quad D_R\triangleq 1+D(r),$$
so that Lemma~\ref{bd:low-transfer} and \eqref{extlow:struct} give
$$ \Xi\triangleq\max\left\{1,\ \sup_{\Omega\times[0,T]}\frac{\rho}{h(\rho)\kappa(\theta)}\right\} \le C D_R^{1+\alpha}.$$

\begin{lemma}\label{est:wall-diss}
The solution satisfies
\be\la{bdthree:mech}\ba
 \int_0^T\left[\|\nabla u\|_{L^2}^2
   +\int_\Omega\bigl(h(\div u)^2+Q\bigr)\,dx
   +\int_{\partial\Omega}A u_\tau^2\,dS\right]dt\le C.
\ea\ee
For every interval $I\subset[0,T]$ with $|I|\le 1$,
\be\la{bd41:pressure}\ba
 \int_I\int_\Omega\left(\frac{P^2}{\mu+\lambda}+P|\div u|\right)dx\,dt\le C,
\ea\ee
and for any $1\le p<\infty$ and $0<\varepsilon\le\min\{1,\alpha/2\}$,
\be\la{bd41:temp}\ba
 \int_I\|\theta\|_{L^p(\Omega)}^{1+\alpha-\varepsilon}\,dt\le \frac{Cp}{\varepsilon}.
\ea\ee
\end{lemma}

\begin{proof}
Estimate \eqref{bdthree:mech} follows from Lemma~\ref{bd:low-transfer}, tangent div--curl, and \eqref{wallthree:sign}. Splitting $P^2/(\mu+\lambda)$ at $\theta=l$ and applying Young's inequality with \eqref{bdthree:mech} yield \eqref{bd41:pressure}, while Lemma~\ref{bd:low-transfer} (\eqref{extlow:heattail}) and $|\Omega|<\infty$ yield \eqref{bd41:temp}.
\end{proof}

\begin{lemma}\label{bdnew:spatial}
There exists a constant $C>0$ such that
\be\la{bd2fp03}\ba
 &\|\nabla u\|_{L^4}+\|\sqrt h\,\div u\|_{L^4}+\|P\div u\|_{L^2}^{1/2}\\
 &\le C R_T^{\frac{B_h}{2}}A_1+C R_T^{\frac{B_h}{2}}A_1^{\frac{1}{2}}\bigl(A_2+\|\nabla u\|_{L^2}\bigr)^{\frac{1}{2}}+C\|P\|_{L^4}.
\ea\ee
For any $2\le p\le 4$,
\be\la{bd2fp003}\ba
 &\|\nabla u\|_{L^p}+\|\sqrt h\,\div u\|_{L^p}+\|P\div u\|_{L^{p/2}}^{1/2}\\
 &\le C A_1+C\|P\|_{L^p}+C A_1^{\frac{6-p}{2p}}\bigl(A_2+\|\nabla u\|_{L^2}\bigr)^{\frac{3p-6}{2p}}.
\ea\ee
Moreover,
\be\la{bdltu-pressure}\ba
 \|P\|_{L^p}\le C\Theta^{a(1-1/p)}\quad (1\le p\le\ell),\quad
 \Theta\le C(1+Z)^{\frac{1}{\alpha+1}} .
\ea\ee
\end{lemma}

\begin{proof}
Extend $A\tau$ smoothly to $\ol\Omega$. By \eqref{fp31} and \eqref{tri:bd-boundary}, $\omega$ satisfies
\be\nonumber\ba
 \mu\Delta\omega=\nabla^\perp\cdot(\rho\dot u)\quad\text{in }\Omega,\quad
 \omega+(A\tau)\cdot u=0\quad\text{on }\partial\Omega .
\ea\ee
Dirichlet elliptic regularity for $\omega+(A\tau)\cdot u$, \eqref{fp31}, and tangent Poincar\'e give, for $q\in\{3/2,2\}$,
\be\nonumber\ba
 \|\nabla G\|_{L^q}+\|\nabla\omega\|_{L^q}
 \le C\bigl(\|\rho\dot u\|_{L^q}+\|u\|_{W^{1,q}}\bigr)
 \le C\bigl(\|\rho\dot u\|_{L^q}+\|\nabla u\|_{L^2}\bigr).
\ea\ee
Since $\|\rho\|_{L^3}+\int_\Omega h\,dx\le C$ (by Lemma~\ref{bd:density-moment} and \eqref{visc:envelope}), we have $|\ol G|+\|\omega\|_{L^2}\le C A_1$, and Sobolev embedding yields
\be\nonumber\ba
 \|G\|_{L^6}+\|\omega\|_{L^6}\le C\bigl(A_1+A_2+\|\nabla u\|_{L^2}\bigr),\quad
 \|\nabla G\|_{L^2}+\|\nabla\omega\|_{L^2}\le C\bigl(R_T^{1/2}A_2+\|\nabla u\|_{L^2}\bigr).
\ea\ee
Applying tangent div--curl to $h\div u=G+P$ and interpolating as in \eqref{2fp003} give \eqref{bd2fp03}--\eqref{bd2fp003}, while Lemma~\ref{1pi} on $\Gamma_t$ yields \eqref{bdltu-pressure} as in \eqref{ltu-pressure}.
\end{proof}

\begin{lemma}\label{bd17:neumann}
Let $v\triangleq\nabla\varphi$ and $w\triangleq\dot u+v$, where $\varphi$ is the unique solution of
\be\nonumber\ba
 \begin{cases}
  \Delta\varphi=2\det\nabla u & \text{in }\Omega,\\
  \partial_n\varphi=\chi u_\tau^2 & \text{on }\partial\Omega,
 \end{cases}\quad
 \int_\Omega\varphi\,dx=0 .
\ea\ee
Then $w\cdot n=0$ on $\partial\Omega$, $\div w=D_t(\div u)+(\div u)^2$, $\omega(w)=D_t\omega+(\div u)\omega$, and
\be\la{bd17:liftbounds}\ba
 \|v\|_{L^2}&\le C\|\nabla u\|_{L^2}^2,\quad
 \|v\|_{H^1}\le C\|\nabla u\|_{L^4}^2,\\
 \|D_t v\|_{L^2}&\le C\|\nabla u\|_{L^2}\bigl(\|w\|_{H^1}+\|\nabla u\|_{L^4}^2\bigr).
\ea\ee
\end{lemma}

\begin{proof}
Since $u\cdot n=0$ and $\partial_\tau n=\chi\tau$ on $\partial\Omega$, we have
\be\nonumber\ba
 \dot u\cdot n=-\chi u_\tau^2,\quad
 2\int_\Omega\det\nabla u\,dx
 =\int_{\partial\Omega}\bigl(u\div u-u\cdot\nabla u\bigr)\cdot n\,dS
 =\int_{\partial\Omega}\chi u_\tau^2\,dS ,
\ea\ee
which verifies Neumann compatibility and $w\cdot n=0$; the identities for $\div w$ and $\omega(w)$ follow from \eqref{fp44}--\eqref{fp45}.

For $B_0(a,b)\triangleq(\div a)(\div b)-\operatorname{tr}(\nabla a\,\nabla b)$ with $a\cdot n=b\cdot n=0$ on $\partial\Omega$, we claim that
\be\la{bd17:jac}\ba
 \left|\int_\Omega f B_0(a,b)\,dx\right|
 \le C\|f\|_{H^1}\|a\|_{H^1}\|b\|_{H^1},\quad f\in H^1(\Omega).
\ea\ee
Indeed, in a local boundary frame $(\tau,n)$, the terms with two differentiated components are $\mathcal J(a_\tau,b_n)+\mathcal J(b_\tau,a_n)$ (with $\mathcal J$ as in Lemma~\ref{core:acceleration}); extending $a_n,b_n$ by zero across $\partial\Omega$ and the remaining factors by $H^1$ extension, \cite{CLMS} bounds both interior and boundary terms, while frame-derivative terms follow from $H^1(\Omega)\hookrightarrow L^4(\Omega)$.

Testing the Neumann problem by $\varphi$ and applying \eqref{bd17:jac} (with $a=b=u$), trace, and Poincar\'e inequalities give
\be\nonumber\ba
 \|v\|_{L^2}^2
 =\int_{\partial\Omega}\chi u_\tau^2\varphi\,dS-\int_\Omega\varphi B_0(u,u)\,dx
 \le C\|\nabla u\|_{L^2}^2\|v\|_{L^2}.
\ea\ee
Using smooth extensions of $\chi,\tau,n$ to $\Omega$, the algebra property of $W^{1,4}(\Omega)$ gives $\|\chi u_\tau^2 n\|_{W^{1,4}}\le C\|\nabla u\|_{L^4}^2$, so Neumann $H^2$ regularity yields
\be\nonumber\ba
 \|v\|_{H^1}\le C\bigl(\|\det\nabla u\|_{L^2}+\|\chi u_\tau^2\|_{H^{1/2}(\partial\Omega)}\bigr)
 \le C\|\nabla u\|_{L^4}^2 .
\ea\ee

Finally, differentiating $\div v=2\det\nabla u$, $\omega(v)=0$, and $v\cdot n=\chi u_\tau^2$ along the flow yields
\be\nonumber\ba
 \div(D_t v)&=2B_0(u,w)-3B_0(u,v),\quad
 \omega(D_t v)=-\sum_{j=1}^2\mathcal J(u^j,v^j),\\
 (D_t v)\cdot n&=(\partial_\tau\chi)u_\tau^3+2\chi u_\tau w_\tau-3\chi u_\tau v_\tau \quad\text{on }\partial\Omega .
\ea\ee
Applying \eqref{bd17:jac} (with $b=w$ and $b=v-\chi u_\tau^2 n$), the zero-extended $H^1_0(\Omega)$ Jacobian bound as in \eqref{kv:qmaterial}, and the trace inequality gives
\be\nonumber\ba
 \|\div(D_t v)\|_{(H^1(\Omega))'}+\|\omega(D_t v)\|_{H^{-1}(\Omega)}
 +\|(D_t v)\cdot n\|_{L^2(\partial\Omega)}
 \le C\|\nabla u\|_{L^2}\bigl(\|w\|_{H^1}+\|\nabla u\|_{L^4}^2\bigr).
\ea\ee
The div--curl decomposition on $\Omega$ then yields the last bound in \eqref{bd17:liftbounds}.
\end{proof}

\begin{lemma}\label{bd17:acc}
Let $w$ be as in Lemma~\ref{bd17:neumann} and $\mathcal L_T\triangleq R_T\log(e+R_T)$. For every $\eta\in(0,1)$,
\be\la{bd07a}\ba
 \|\nabla u\|_{L^4}^4+\|\sqrt h\,\div u\|_{L^4}^4
 +\|P\div u\|_{L^2}^2+\|P\|_{L^4}^4
 \le C\mathcal L_T g X,
\ea\ee
\be\la{bd07b}\ba
 (A_1^2)'+A_2^2
 \le \eta Z^2+C_\eta\mathcal L_T g A_1^2
 +C R_T\|\nabla u\|_{L^2}^4+C A_1\sqrt{\Xi B},
\ea\ee
and
\be\la{bd07c}\ba
 &\frac{d}{dt}\int_\Omega\rho|w|^2\,dx
 +c\int_\Omega\bigl(h|\div w|^2+\mu|\omega(w)|^2\bigr)\,dx\\
 &\le C\Xi B+C\mathcal L_T g\left(X+\int_\Omega\rho|w|^2\,dx\right),
\ea\ee
with constants independent of $R_T$.
\end{lemma}

\begin{proof}
By \eqref{visc:envelope} and Lemma~\ref{bdnew:spatial}, $\|G\|_{L^2}\le C R_T^{B_h/2}A_1$ and $\|\nabla G\|_{L^2}+\|\nabla\omega\|_{L^2}\le C(R_T^{1/2}A_2+\|\nabla u\|_{L^2})$. Applying the weighted interpolation \eqref{core:weighted-flux} yields
\be\nonumber\ba
 \left\|\frac{G^2}{h}\right\|_{L^2}
 \le C A_1^2+C\sqrt{\log(e+R_T)}\,A_1\bigl(R_T^{1/2}A_2+\|\nabla u\|_{L^2}\bigr),
\ea\ee
which together with tangent div--curl and \eqref{ltu-absorb} gives \eqref{bd07a}.

Next, multiplying \eqref{fp31} by $2\dot u$ and using $\dot u\cdot n=-\chi u_\tau^2$ and $\dot u\cdot\tau=\partial_t u_\tau+u_\tau\partial_\tau u_\tau$ on $\partial\Omega$, we obtain
\be\nonumber\ba
 (A_1^2)'+2A_2^2
 =\sum_{i=1}^5 I_i-2\int_{\partial\Omega}\chi G u_\tau^2\,dS
 +\frac{2\mu}{3}\int_{\partial\Omega}(\partial_\tau A)u_\tau^3\,dS ,
\ea\ee
where $I_i$ are the volume integrals in \eqref{fp43} over $\Omega$. By \eqref{bd17:jac}, the trace theorem, and Lemma~\ref{bdnew:spatial},
\be\nonumber\ba
 \left|\int_\Omega G\det\nabla u\,dx\right|
 +\int_{\partial\Omega}\bigl(|G|u_\tau^2+|u_\tau|^3\bigr)\,dS
 &\le C\bigl(A_1+R_T^{1/2}A_2+\|\nabla u\|_{L^2}\bigr)\|\nabla u\|_{L^2}^2\\
 &\le \varepsilon A_2^2+C_\varepsilon\bigl(R_T\|\nabla u\|_{L^2}^4+g A_1^2\bigr).
\ea\ee
Estimating $I_1,I_3,I_4,I_5$ exactly as in the proof of \eqref{kvu:xprime} yields \eqref{bd07b}.

To prove \eqref{bd07c}, test the differentiated momentum equation by $w$ as in Lemma~\ref{core:acceleration}. With $\mathcal T=G\mathbb I_2+\mu\omega\mathbb J$, \eqref{tri:bd-boundary} and $(\nabla u)^\top n=(\div u-\partial_\tau u_\tau)n-\chi u_\tau\tau$ give on $\partial\Omega$:
\be\nonumber\ba
 &\bigl[\bigl(D_t\mathcal T+(\div u)\mathcal T-\mathcal T(\nabla u)^\top\bigr)n\bigr]\cdot\tau\\
 &=-\mu A w_\tau+\mu A v_\tau-\mu(\partial_\tau A)u_\tau^2
   -\frac{\mu}{2}A\,\partial_\tau(u_\tau^2)+\chi G u_\tau .
\ea\ee
Since $w\cdot n=0$, pairing with $w_\tau$ produces $-\mu\int_{\partial\Omega}A w_\tau^2\,dS\le 0$, while the remaining four boundary terms are bounded via \eqref{bd17:liftbounds} and trace duality ($\|u_\tau^2\|_{H^{1/2}(\partial\Omega)}\le C\|\nabla u\|_{L^4}^2$) by
\be\nonumber\ba
 C\bigl(\|\nabla u\|_{L^4}^2+\|G\|_{H^1}\|\nabla u\|_{L^2}\bigr)\|w\|_{H^1}.
\ea\ee
For the Jacobian and $D_t v$ volume integrals, \eqref{bd17:jac} and \eqref{bd17:liftbounds} yield
\be\nonumber\ba
 \left|\int_\Omega G\bigl((\div u)(\div w)-(\nabla u)^\top:\nabla w\bigr)\,dx\right|
\le C\|G\|_{H^1}\|\nabla u\|_{L^2}\|w\|_{H^1},
\ea\ee
and
\be\nonumber\ba
 \left|\int_\Omega\rho D_t v\cdot w\,dx\right|
\le C R_T^{1/2}\|\sqrt\rho\,w\|_{L^2}\|\nabla u\|_{L^2}
   \bigl(\|w\|_{H^1}+\|\nabla u\|_{L^4}^2\bigr).
\ea\ee
Handling the remaining volume terms as in \eqref{kv:corrected-raw} and using $\|w\|_{H^1}^2\le C\int_\Omega(h|\div w|^2+\mu|\omega(w)|^2)\,dx$ yield \eqref{bd07c}.
This completes the proof of Lemma \ref{bd17:acc}.
\end{proof}

\begin{lemma}\label{est:wall-heat}
On an interval with $\rho\le r$ and $\theta\ge c/D_R$, one has
\be\nonumber\ba
 H'+\frac{3}{4}Z^2 &\le \frac{1}{4}A_2^2+C A_1^4\log(e+A_1^2),\\
 \|\Delta K(\theta)\|_{L^2}^2 &\le C R D_R^\alpha B+C L_R g X,\\
 \|\nabla\dot u\|_{L^2}^2+\int_\Omega h|D_t(\div u)|^2\,dx &\le C D_w+C L_R g X .
\ea\ee
\end{lemma}

\begin{proof}
Testing $(\ref{fns})_3$ by $K(\theta)$ gives
$$H'+Z^2=\int_\Omega K(\theta)(Q-P\div u)\,dx.$$
The extra term $A_1^{2-3/b}\|\nabla u\|_{L^2}^{3/b}$ from \eqref{bd2fp003} is bounded by $A_1^2+\|\nabla u\|_{L^2}^2\le C g$ for $b\ge\ell$, so \eqref{fp77}, \eqref{fp720}, and \eqref{fp517} apply on $\Omega$. By $\partial_n K(\theta)=0$, Neumann regularity gives $\|\nabla K\|_{L^4}^2\le C Z\|\Delta K\|_{L^2}$, and $(\ref{fns})_3$, $\kappa^{-1}\le C D_R^\alpha$, and \eqref{bd07a} yield the second bound. The last bound follows from $\dot u=w-v$, $D_t(\div u)=\div w-(\div u)^2$, and \eqref{bd17:liftbounds}.
\end{proof}

\begin{lemma}\label{tri:bdcap}
Let $r\ge\sup_{0\le t\le T}\|\rho(t)\|_{L^\infty}$, $R=2+r$, and $D_R=1+D(r)$. For any fixed $p>2$ and every $m\ge 8$,
\be\la{bdfp08}\ba
 \|G-\ol G\|_{L^p}\le C_p\bigl(A_2+\|\nabla u\|_{L^2}\bigr),\quad
 \|\nabla G\|_{L^2}\le C\bigl(R^{1/2}A_2+\|\nabla u\|_{L^2}\bigr),
\ea\ee
and
\be\la{bdfp008}\ba
 \|\nabla G\|_{L^m}\le C m^{15/2}R D_R
 \bigl(1+A_1+A_2+\|\nabla u\|_{L^2}+\|\nabla\dot u\|_{L^2}+Z^{\frac{1}{1+\alpha}}\bigr),
\ea\ee
with constants independent of $r$ and $m$.
\end{lemma}

\begin{proof}
For $q=2p/(p+2)$, Sobolev embedding and Lemma~\ref{bdnew:spatial} give
\be\nonumber\ba
 \|G-\ol G\|_{L^p}\le C_p\|\nabla G\|_{L^q}
 \le C_p\bigl(\|\sqrt\rho\|_{L^p}A_2+\|\nabla u\|_{L^2}\bigr),
\ea\ee
which together with Lemma~\ref{bdnew:spatial} proves \eqref{bdfp08}. For $m\ge 8$, the Dirichlet and tangent div--curl estimates \cite[Lemmas~2.2--2.3]{FJLL2} and $P/h\le C D_R\theta$ give
\be\nonumber\ba
 \|\nabla G\|_{L^m}&\le C m\bigl(\|\rho\dot u\|_{L^m}+\|u\|_{W^{1,m}}\bigr),\\
 \|\nabla u\|_{L^m}&\le C m^6\bigl(\|G\|_{L^m}+\|\omega\|_{L^m}+D_R\|\theta\|_{L^m}\bigr).
\ea\ee
By \eqref{gn11}, Lemma~\ref{bdnew:spatial}, and mass-weighted Poincar\'e,
\be\nonumber\ba
& \|\dot u\|_{L^m} \le C m^{1/2}\bigl(A_2+\|\nabla\dot u\|_{L^2}\bigr),\quad
 \|\theta\|_{L^m}\le C m^{\frac{1}{2(1+\alpha)}}(1+Z)^{\frac{1}{1+\alpha}},\\
& \|G\|_{L^m}+\|\omega\|_{L^m}+\|u\|_{L^m} \le C m^{1/2}\bigl(A_1+R^{1/2}A_2+\|\nabla u\|_{L^2}\bigr).
\ea\ee
Substituting these estimates yields \eqref{bdfp008}.
\end{proof}

\begin{lemma}\label{bd:record-transfer}
Under the hypotheses of Lemma~\ref{bd:low-transfer}, the differential inequalities \eqref{exthigh:joint}, \eqref{exthigh:work}, and \eqref{exthigh:mech} hold on $\Omega$ (with $D_\partial$ included in $D_u$ and $A_1^2$, and with the Neumann correction of Lemma~\ref{bd17:neumann}). Consequently, for the compatible data of Theorem~\ref{tri:main}, estimates \eqref{exthigh:window} and \eqref{exthigh:sum}, as well as Lemma~\ref{records:flux} for $f\triangleq G-\ol G$, hold on $\Omega$ with the same time dependence as in Lemma~\ref{bd:low-transfer}.
\end{lemma}

\begin{proof}
Lemma~\ref{bd:low-transfer} yields the pressure bounds of Lemma~\ref{core:pressure-times} and
\be\nonumber\ba
 \int_0^T\left(D_u+\int_{\{\theta>4l\}}\frac{P^2}{h}\,dx\right)dt\le C .
\ea\ee
Since the same nonnegative term $D_\partial$ enters both $A_1^2$ and $D_u$, tangent div--curl gives \eqref{exthigh:compare}, and \eqref{exthigh:tails} follows as on the torus. By Lemmas~\ref{bd17:acc} and~\ref{est:wall-heat}, \eqref{common:material-gradient} and \eqref{exthigh:tprime} hold on $\Omega$, so Lemma~\ref{records:joint} yields \eqref{exthigh:joint} and \eqref{exthigh:window}.

In the mechanical identity \eqref{per:mechanical-identity}, $\dot u\cdot n=-\chi u_\tau^2$ and tangential integration by parts add the boundary term
\be\nonumber\ba
 -\int_{\partial\Omega}\chi G u_\tau^2\,dS
 +\frac{\mu}{3}\int_{\partial\Omega}(\partial_\tau A)u_\tau^3\,dS
\ea\ee
to the right-hand side. By the trace and Jacobian estimates in Lemma~\ref{bd17:acc}, this boundary term together with $2\int_\Omega G\det\nabla u\,dx$ is bounded by $C\bigl(A_1+R^{1/2}A_2+\|\nabla u\|_{L^2}\bigr)\|\nabla u\|_{L^2}^2$. Since $D_w$ controls $\int_\Omega h|\div w|^2\,dx$, both \eqref{exthigh:work} and \eqref{exthigh:mech} hold (before using $A_2(0)$). Compatibility bounds $A_2(0)$ and the trace of $u_0$ bounds $D_\partial(0)$.
If $r_N$ is reached by time $T$, summing \eqref{per:layer-endpoints} over $j\le N$ with $D(r_j)\le(\log r_j)^{\zeta_*}$ gives \eqref{exthigh:sum}, as in Lemma~\ref{records:sum}.

Finally, $\int_\Omega(h+P)\,dx\le C$ and \eqref{bdthree:mech} give $\int_I|\ol G|\,dt\le C(1+|I|)$. For $f=G-\ol G$, Lemma~\ref{tri:bdcap} supplies the elliptic inputs \eqref{per:cross-elliptic}, and the ball interpolation \eqref{records:ball} holds in boundary charts.
The time-integration argument of Lemma~\ref{records:flux} then applies verbatim.
\end{proof}

\begin{lemma}\label{order:wall-cap}
For any smooth positive-density solution with the data of Theorem~\ref{tri:main},
\be\nonumber\ba
 \sup_{0\le t\le T}\|\rho(t)\|_{L^\infty}\le C_{T_0}\quad (T\le T_0),
\ea\ee
independently of the initial density floor. Under \textup{(A)}, \textup{(B)}, or \textup{(C)} with $\alpha>1-d$, the bound is independent of $T_0$.
\end{lemma}

\begin{proof}
Consider the indices $j$ satisfying $D(r_j)\le(\log r_j)^{\zeta_*}$. For each such index whose level $r_j$ is reached by time $T$, choose $y_j\in\ol\Omega$ with $\rho(y_j,t_j)=r_j$.
Since $u\cdot n=0$, the backward trajectory from $y_j$ stays in $\ol\Omega$. Its density is at most $r_j/2$ at $t_{j^-}$, so its last crossing of $r_j/2$ lies in $[t_{j^-},t_j]$. Along this crossing, $D_t\Lambda(\rho)+P=-G$ with $P\ge c r_j/D_{R_j}$ (with $c=c_{T_0}$ in the finite-horizon regimes). Lemma~\ref{bd:record-transfer} therefore yields \eqref{order:crossing}, and the divergence argument of Proposition~\ref{order:periodic-cap} (together with \eqref{per:family-level} for \eqref{ext:powerfamily}) shows that at least one level satisfying the coefficient condition cannot be reached.
For \eqref{ext:gainfamily}, take $D\equiv1$ and first choose $R_{*,C}$ by this bounded-domain argument with lower-bound constant $c_*/2$ in \eqref{new:gain} and data bound $4M+4$. Since $\alpha>1-d$, this choice is independent of time and uniform over the stated coefficients and data. Define $\varepsilon_{*,C}(M)$ by \eqref{per:gain-family-level} with this level. The comparison in Proposition~\ref{order:periodic-cap} preserves \eqref{new:gain} up to the first attainment of $R_{*,C}$, so the same crossing argument excludes this level.
Fixing the resulting ceiling in \eqref{exthigh:window} and applying Lemma~\ref{est:wall-heat}, tangent Poincar\'e for $u$, and Lemma~\ref{1pi} on $\Gamma_t$ for $K(\theta)$ also yield \eqref{construction:common-cap} on $\Omega$.
\end{proof}

\begin{lemma}\label{wallthree:uniformlow}
Under \textup{(A)}, \textup{(B)}, or \textup{(C)} with $\alpha>1-d$, for every $1\le p<\infty$ there exist constants $C_p,\theta_*>0$, independent of $T$ and the minimum initial density, such that
\be\la{wallthree:lower}\ba
 &\sup_{0\le t\le T}\left(E(t)+\frac{M_0}{\int_\Omega\rho\theta\,dx}
   +\|\rho(t)\|_{L^p}\right)\\
 &\quad +\int_0^T\int_\Omega
   \left(|\nabla u|^2+Q+\frac{Q}{\theta}+\frac{\kappa(\theta)|\nabla\theta|^2}{\theta^2}\right)dx\,dt
   +\int_0^T\int_{\partial\Omega}u_\tau^2\,dS\,dt\le C_p,
\ea\ee
and
\be\nonumber\ba
 \theta(x,t)\ge\theta_*,\quad (x,t)\in\ol\Omega\times[0,T].
\ea\ee
\end{lemma}

\begin{proof}
By Lemmas~\ref{bd:low-transfer} and~\ref{order:wall-cap}, $\theta\ge\theta_*>0$ and $\int_\Omega\rho\theta\,dx\ge M_0\theta_*$. Combining \eqref{bd41:entropy}, \eqref{bdthree:mech}, and the trace estimate $\int_{\partial\Omega}u_\tau^2\,dS\le C\|\nabla u\|_{L^2}^2$ gives \eqref{wallthree:lower}.
\end{proof}

\subsection{Higher-order estimates}

\begin{lemma}\label{reg:wall-initial}\label{reg:wall-barrier}
For the exponent $q>2$ in \eqref{tri:bd-data} and every finite $T>0$, the bound \eqref{init:prelayer} holds on $\Omega$, independently of the initial density floor.
\end{lemma}

\begin{proof}
At the fixed density ceiling, Lemma~\ref{bdnew:spatial}, $h\div u=G+P$, and tangent div--curl bound $\|u\|_{W^{1,p}}$ for every finite $p\ge 2$, while Neumann $H^2$ regularity for $K(\theta)$ gives $\theta\in L^2(0,T;H^2)$ as in Lemma~\ref{init:common-lemma}. Applying the Dirichlet $W^{1,q}$ bound to $\omega+(A\tau)\cdot u$ and tangent $W^{2,q}$ div--curl to $u$ yields
\be\nonumber\ba
 \|\div u\|_{L^\infty}+\|\omega\|_{L^\infty}\le C F,\quad
 \|\nabla^2 u\|_{L^q}\le C F\bigl(1+\|\nabla\rho\|_{L^q}\bigr),
\ea\ee
where $F\triangleq 1+\|\dot u\|_{H^1}+\|\theta\|_{H^2}\in L^2(0,T)$. Since $u\cdot n=0$ eliminates boundary fluxes in the differentiated continuity equation, the logarithmic Gr\"onwall argument of Lemma~\ref{init:common-lemma} yields all velocity and density-gradient bounds in \eqref{init:prelayer}.

For the initial temperature ceiling on $[0,\min\{1,T\}]$, replace the periodic auxiliary problem in Lemma~\ref{init:common-lemma} by the Neumann problem
\be\nonumber\ba
 (-\Delta+\rho)\psi=Q\quad\text{in }\Omega,\quad
 \partial_n\psi=0\quad\text{on }\partial\Omega .
\ea\ee
Mass-weighted Poincar\'e gives coercivity and $\psi\ge 0$, while Neumann $W^{2,4}$ and $W^{2,4/3}$ regularity bound $\psi$ and $\psi_t$ exactly as in Lemma~\ref{init:common-lemma}. Since $\partial_n K(\theta)=\partial_n K(\bar\theta)=0$ on $\partial\Omega$, the positive-part comparison applies verbatim and bounds $\|\theta\|_{L^\infty}$ on $[0,\min\{1,T\}]$.
\end{proof}

\begin{lemma}\la{tri:bdhigher}
For every $0<\tau<T\le T_0$, there exists a constant $C=C(\min\{1,\tau\},T_0)$ (independent of $T_0$ under \textup{(A)}, \textup{(B)}, or \textup{(C)} with $\alpha>1-d$) such that
\be\la{yjyho1}\ba
 \sup_{\tau\le t\le T}\bigl(\|\theta\|_{H^2}+\|\sqrt\rho\,\theta_t\|_{L^2}\bigr)
 +\sup_{\substack{I\subset[\tau,T]\\ |I|\le 1}}
   \int_I\bigl(\|\dot\theta\|_{H^1}^2+\|\dot u\|_{H^1}^2\bigr)\,dt\le C.
\ea\ee
\end{lemma}

\begin{proof}
With $k\triangleq K(\theta)$ and $v\triangleq D_t k$, testing \eqref{fp97} by $v$ produces the boundary integral $\int_{\partial\Omega}v\,n\cdot(\nabla v+C)\,dS$, where $C_i=(\div u)\partial_i k-(\partial_i u^j+\partial_j u^i)\partial_j k$. Differentiating $n\cdot\nabla k=0$ and $u\cdot n=0$ tangentially along $\partial\Omega$ shows that $n_i u^j\partial_{ij}k=n_i(\partial_j u^i)\partial_j k=-\chi u_\tau\partial_\tau k$, and hence
\be\nonumber\ba
 n\cdot(\nabla v+C)
 =n_i u^j\partial_{ij}k-n_i(\partial_j u^i)\partial_j k=0\quad\text{on }\partial\Omega .
\ea\ee
Moreover, \eqref{fp93} and \eqref{fp96} follow from mass and energy conservation, \eqref{fp103} from Lemma~\ref{bdnew:spatial} and tangent div--curl, and \eqref{fp912} from the Neumann bound $\|k\|_{H^2}\le C(\|\Delta k\|_{L^2}+\|k\|_{L^2})$. Thus \eqref{fp910}--\eqref{fp918} carry over directly from Lemma~\ref{fpl9}, giving \eqref{yjyho1}.
\end{proof}

\begin{lemma}\label{reg:wall-higher}
For every $0<\tau<T<\infty$ and the exponent $q>2$ in \eqref{tri:bd-data}, there exists a constant $C(\tau,T)>0$ such that
\be\la{yjyho2}\ba
 &\sup_{\tau\le t\le T}\bigl(\|\rho\|_{W^{1,q}}+\|\rho_t\|_{L^q}+\|u\|_{H^2}+\|\sqrt\rho\,u_t\|_{L^2}\bigr)\\
 &\quad +\int_\tau^T\bigl(\|u\|_{W^{2,q}}^2+\|\theta\|_{H^3}^2+\|u_t\|_{H^1}^2+\|\theta_t\|_{H^1}^2\bigr)\,dt
 \le C(\tau,T).
\ea\ee
\end{lemma}

\begin{proof}
Combining Lemma~\ref{reg:wall-initial}, Lemma~\ref{tri:bdhigher}, and the spatial-differentiation step of Lemma~\ref{init:common-lemma} (replacing periodic $H^3$ regularity with the Neumann estimate $\|K(\theta)\|_{H^3}\le C(\|\nabla\Delta K(\theta)\|_{L^2}+\|K(\theta)\|_{H^2})$) directly yields \eqref{yjyho2}.
\end{proof}

\begin{lemma}\label{reg:wall-norms}
Assume \textup{(A)}, \textup{(B)}, or \textup{(C)} with $\alpha>1-d$. Set $\rho_*\triangleq M_0/|\Omega|$ and $E_\infty\triangleq\lim_{t\to\infty}E(t)$. Then $\theta_{\rm eq}=E_\infty/(c_v M_0)\ge\theta_*>0$ satisfies \eqref{tri:bd-equilibrium}, and for each $1\le p<\infty$, there exist constants $C_p,\eta_p>0$ such that, for all $t\ge 1$,
\be\nonumber\ba
 \|\rho(t)-\rho_*\|_{L^p}+\|u(t)\|_{W^{1,p}}
 +\|\theta(t)-\theta_{\rm eq}\|_{W^{1,p}}\le C_p e^{-\eta_p t}.
\ea\ee
\end{lemma}

\begin{proof}
Set $r\triangleq\rho-\rho_*$, $E_k\triangleq\frac{1}{2}\int_\Omega\rho|u|^2\,dx$, $m\triangleq M_0^{-1}\int_\Omega\rho\theta\,dx$, $\theta_E\triangleq E(t)/(c_v M_0)$, and $B_\partial\triangleq\mu\int_{\partial\Omega}(A+2\chi)u_\tau^2\,dS$. By Lemmas~\ref{wallthree:uniformlow}, \ref{reg:wall-barrier}, and~\ref{tri:bdhigher} (at $\tau=1$), $0<c\le m\le\theta_E\le C$ for $t\ge 1$. Define the entropy deficit
\be\nonumber\ba
 \mathcal F\triangleq c_v M_0\log\theta_E-M_0\log\rho_*-\mathscr S
 =c_v M_0\log\frac{\theta_E}{m}
 +c_v\int_\Omega\rho\log\frac{m}{\theta}\,dx
 +\int_\Omega\left(\rho\log\frac{\rho}{\rho_*}-\rho+\rho_*\right)dx .
\ea\ee
By Taylor's formula on the compact density and temperature ranges, $\mathcal F\asymp E_k+\int_\Omega\rho|\theta-m|^2\,dx+\|r\|_{L^2}^2$. Differentiating $\mathcal F$ via \eqref{wallthree:sign} ($E'=-B_\partial\le 0$) and Lemma~\ref{bdthree:wall}, and noting $\int_\Omega Q\,dx+B_\partial=D_u\ge c\|\nabla u\|_{L^2}^2$, we obtain
\be\la{extslip:diss}\ba
 -\mathcal F'
 =\int_\Omega\left(\frac{Q}{\theta}+\frac{\kappa(\theta)|\nabla\theta|^2}{\theta^2}\right)dx
 +\frac{B_\partial}{\theta_E}
 \ge c\bigl(\|\nabla u\|_{L^2}^2+\|\nabla\theta\|_{L^2}^2\bigr).
\ea\ee

Let $\phi$ be the zero-mean Neumann solution of $\Delta\phi=r$ in $\Omega$, $\partial_n\phi=0$ on $\partial\Omega$, and put $v\triangleq\nabla\phi$ and $I\triangleq\int_\Omega\rho u\cdot v\,dx$. Testing $(\ref{fns})_2$ by $v$ produces the same volume terms as in Lemma~\ref{reg:tor-decay-norms} plus the boundary term $-\mu\int_{\partial\Omega}(A+2\chi)u_\tau v_\tau\,dS$, which is bounded by $C\|\nabla u\|_{L^2}\|r\|_{L^2}$ via the trace theorem. Estimating the volume terms as in Lemma~\ref{reg:tor-decay-norms} (using tangent Poincar\'e and $\|v_t\|_{L^2}\le\|\rho u\|_{L^2}$) gives
\be\nonumber\ba
 I'\ge c\|r\|_{L^2}^2-C\bigl(\|\nabla u\|_{L^2}^2+\|\nabla\theta\|_{L^2}^2\bigr),\quad
 |I|\le C\sqrt{E_k}\|r\|_{L^2}.
\ea\ee
Choosing $\delta>0$ small so that $\mathcal F-\delta I\asymp\mathcal F$ and applying Gr\"onwall's inequality via \eqref{extslip:diss} yield $\mathcal F(t)\le C e^{-ct}$ for $t\ge 1$. Integrating the equality in \eqref{extslip:diss} over $[t,\infty)$ then gives
\be\nonumber\ba
 0\le E(t)-E_\infty=\int_t^\infty B_\partial\,ds
 \le\left(\sup_{s\ge 1}\theta_E(s)\right)\mathcal F(t),\quad
 \theta_{\rm eq}=\frac{E_\infty}{c_v M_0}\ge\theta_*>0,
\ea\ee
which proves \eqref{tri:bd-equilibrium} and the exponential decay of $|m-\theta_{\rm eq}|+\|r\|_{L^2}+\|\sqrt\rho\,u\|_{L^2}+\|\sqrt\rho(\theta-\theta_{\rm eq})\|_{L^2}$. Upgrading to unweighted $L^p\times W^{1,p}\times W^{1,p}$ norms follows by the same argument as in Lemma~\ref{reg:tor-decay-norms}, using Lemmas~\ref{bdnew:spatial} and~\ref{tri:bdhigher}.
\end{proof}

\section{A priori estimates on the whole plane}\label{wh3:section}

Consider a smooth positive-density solution of \eqref{fns} with the data of
Theorem~\ref{tri:whole} and the far-field state \eqref{tri:farfield}.
Spatial norms are taken over $\mathbb R^2$ unless a subset is specified.
Write
\[
 \begin{gathered}
 P_\infty=\rho_\infty\theta_\infty,\qquad
 A_1^2=1+\int_{\mathbb R^2}\left(\frac{(G+P_\infty)^2}{h}+\mu\omega^2\right)dx,\\
 H=c_v\int_{\mathbb R^2}\rho\int_{\theta_\infty}^{\theta}(K(v)-K(\theta_\infty))\,dv\,dx.
 \end{gathered}
\]
All estimates below are independent of the minimum initial density.

\subsection{Low-order estimates}
\begin{lemma}\label{wh3:low}
 There exists a positive constant $C$, independent of $T$, such that
\begin{equation}\label{wh3:primitives}
\begin{aligned}
 &\sup_{0\le t\le T}\int_{\mathbb R^2}
 \bigg[\frac{1}{2}\rho|u|^2
    +c_v\rho\left(\theta-\theta_\infty
                 -\theta_\infty\log\frac{\theta}{\theta_\infty}\right)+\theta_\infty
   \left(\rho\log\frac{\rho}{\rho_\infty}-\rho+\rho_\infty\right)\bigg]\,dx\\
 &+\int_0^T\!\int_{\mathbb R^2}
    \left(\frac{2\mu|D(u)|^2+\lambda(\div u)^2}{\theta}
          +\frac{\kappa(\theta)|\nabla\theta|^2}{\theta^2}\right)dx\,dt\le C.
\end{aligned}
\end{equation}
\end{lemma}
\begin{proof}
Let $y(t)=(\underline\theta_0^{-1}+a_0 t)^{-1}$ be as in \eqref{fp117}, which satisfies $y'=-a_0 y^2$. Together with \eqref{fp19}, this gives
\begin{equation}\label{wh3:temperature-comparison}
 c_v\rho(\dot\theta-y')-\Delta K(\theta)
 \ge -c_v a_0\rho(\theta^2-y^2)
 =-c_v a_0\rho(\theta+y)(\theta-y).
\end{equation}
Put $w=(y-\theta)_+$. Since $y\le\theta_\infty$, we have $0\le w\le|\theta-\theta_\infty|$ and $w\in H^1$. Multiplying \eqref{wh3:temperature-comparison} by $-w$ and integrating over $\mathbb{R}^2$ gives
\[
 \frac{c_v}{2}\frac{d}{dt}\int_{\mathbb R^2}\rho w^2\,dx+\int_{\mathbb R^2}\kappa(\theta)|\nabla w|^2\,dx
 \le -c_v a_0\int_{\mathbb R^2}\rho(\theta+y)w^2\,dx\le 0,
\]
which together with $w(0)=0$ yields $\theta\ge y(t)>0$.

Multiplying $\eqref{fns}_2$ by $u$ and integrating by parts over $\mathbb{R}^2$ gives
\begin{equation}\label{wh3:relative-kinetic}
 \frac{1}{2}\frac{d}{dt}\int_{\mathbb R^2}\rho|u|^2\,dx+\int_{\mathbb R^2}\bigl(2\mu|D(u)|^2+\lambda(\div u)^2\bigr)\,dx
  =\int_{\mathbb R^2}(P-P_\infty)\div u\,dx.
\end{equation}
Multiplying $\eqref{fns}_3$ by $1-\theta_\infty/\theta$ and integrating over $\mathbb{R}^2$, we deduce from $\eqref{fns}_1$ that
\begin{equation}\label{wh3:relative-thermal}
\begin{aligned}
 &c_v\frac{d}{dt}\int_{\mathbb R^2}\rho
    \left(\theta-\theta_\infty-\theta_\infty
                   \log\frac{\theta}{\theta_\infty}\right)dx
   +\theta_\infty\int_{\mathbb R^2}
             \frac{\kappa(\theta)|\nabla\theta|^2}{\theta^2}\,dx\\
 =&\int_{\mathbb R^2}
        \left(1-\frac{\theta_\infty}{\theta}\right)\bigl(2\mu|D(u)|^2+\lambda(\div u)^2-P\div u\bigr)\,dx.
\end{aligned}
\end{equation}
Multiplying $\eqref{fns}_1$ by $\theta_\infty\log(\rho/\rho_\infty)$ and integrating over $\mathbb{R}^2$ gives
\begin{equation}\label{wh3:relative-density}
 \theta_\infty\frac{d}{dt}\int_{\mathbb R^2}
      \left(\rho\log\frac{\rho}{\rho_\infty}-\rho+\rho_\infty\right)dx
    =-\theta_\infty\int_{\mathbb R^2}(\rho-\rho_\infty)\div u\,dx.
\end{equation}
Notice that
\[
 \left[P-P_\infty
       -\left(1-\frac{\theta_\infty}{\theta}\right)P
       -\theta_\infty(\rho-\rho_\infty)\right]\div u
 =\theta_\infty\left(\frac{P}{\theta}-\rho\right)\div u=0.
\]
Combining \eqref{wh3:relative-kinetic} and \eqref{wh3:relative-thermal} with \eqref{wh3:relative-density}, integrating over $[0,t]$ and using \eqref{tri:whole-data}, we obtain
\[
\begin{aligned}
 &\int_{\mathbb R^2}\bigg[\frac{1}{2}\rho|u|^2
    +c_v\rho\left(\theta-\theta_\infty-\theta_\infty\log\frac{\theta}{\theta_\infty}\right)
    +\theta_\infty\left(\rho\log\frac{\rho}{\rho_\infty}-\rho+\rho_\infty\right)\bigg]\,dx\\
 &+\theta_\infty\int_0^t\!\int_{\mathbb R^2}
    \left(\frac{2\mu|D(u)|^2+\lambda(\div u)^2}{\theta}
          +\frac{\kappa(\theta)|\nabla\theta|^2}{\theta^2}\right)dx\,ds\\
 =&\int_{\mathbb R^2}\bigg[\frac{1}{2}\rho_0|u_0|^2
    +c_v\rho_0\left(\theta_0-\theta_\infty-\theta_\infty\log\frac{\theta_0}{\theta_\infty}\right)
    +\theta_\infty\left(\rho_0\log\frac{\rho_0}{\rho_\infty}-\rho_0+\rho_\infty\right)\bigg]\,dx\\
 \le& C\bigl(\|u_0\|_2^2+\|\theta_0-\theta_\infty\|_2^2+\|\rho_0-\rho_\infty\|_2^2\bigr)\\
 \le& C,
\end{aligned}
\]
which yields \eqref{wh3:primitives}.
\end{proof}

To extract temperature bounds from \eqref{wh3:primitives}, split at $\theta=2\theta_\infty$ and $\rho=\rho_\infty/2$.

\begin{lemma}\label{plane:thermal-sets}
There is a constant $C$ independent of $T$ such that
\begin{equation}\label{wh3:thermal-sets-bound}
 \sup_{0\le t\le T}\left(
 \int_{\{\theta>2\theta_\infty\}}\rho\theta\,dx+
 \int_{\{\theta\le2\theta_\infty\}}|\theta-\theta_\infty|^2\,dx
 \right)\le C .
\end{equation}
\end{lemma}
\begin{proof}
By the density entropy in \eqref{wh3:primitives},
\begin{equation}\label{wh3:low-density-measure}
 |\{\rho<\rho_\infty/2\}|
 \le C\|\sqrt\rho-\sqrt{\rho_\infty}\|_2^2
 \le C\int_{\mathbb R^2}\left(\rho\log\frac{\rho}{\rho_\infty}-\rho+\rho_\infty\right)dx\le C.
\end{equation}
On $\{\theta>2\theta_\infty\}$, the thermal entropy in \eqref{wh3:primitives} gives
\begin{equation}\label{wh3:theta-hot-energy}
 \int_{\{\theta>2\theta_\infty\}}\rho\theta\,dx
 \le\frac{2}{1-\log 2}\int_{\{\theta>2\theta_\infty\}}\rho
    \left(\theta-\theta_\infty-\theta_\infty\log\frac{\theta}{\theta_\infty}\right)dx
 \le C.
\end{equation}
On $\{\theta\le 2\theta_\infty\}$, Taylor's formula bounds $|\theta-\theta_\infty|^2$ by the thermal entropy density. Splitting at $\rho=\rho_\infty/2$ gives
\begin{equation}\label{wh3:theta-cool-ltwo}
 \int_{\{\theta\le2\theta_\infty\}}|\theta-\theta_\infty|^2\,dx
 \le\frac{2}{\rho_\infty}\int_{\{\theta\le2\theta_\infty\}}\rho|\theta-\theta_\infty|^2\,dx
     +\theta_\infty^2|\{\rho<\rho_\infty/2\}|\le C.
\end{equation}
Both temperature bounds are uniform in $t$.
\end{proof}

Throughout the following estimates, we take $D$ from \textup{(A)}, $D_h$ from \textup{(B)}, $D=1$ in \textup{(C)}, and $D=\Psi\circ J$ in \textup{(D)}, and fix $r_1>2\rho_\infty$ sufficiently large so that \eqref{extlow:struct} holds for $\rho\ge r_1$ on the attained density range.
In the data-class alternatives \eqref{ext:powerfamily} and \eqref{ext:gainfamily}, we first restrict to time intervals on which \eqref{ext:iter} and \eqref{new:gain}, respectively, hold with $c_*/2$ on this range, as in Section~\ref{per-apriori}; the restrictions are removed in Lemma~\ref{wh3:permanentlow}. Set
\[
 A_2=\|\nabla u\|_2^2+
 \int_{\mathbb R^2}\bigl(2\mu|D(u)|^2+\lambda(\div u)^2+
              \kappa(\theta)|\nabla\log\theta|^2\bigr)\,dx .
\]
We next derive the density-moment estimate in terms of $A_2$.

\begin{lemma}\label{plane:drift}
For every $p>2$, there exist constants $C_p,c_p>0$, independent of $T$ and $\sup_{0\le t\le T}\|\rho(t)\|_\infty$, such that
\begin{equation}\label{wh3:drift-bound}
 \log\left(2+\sup_{0\le t\le T}\int_{\{\rho>r_1\}}\rho^p\,dx\right)\le C_p\left(2+\int_0^T A_2\,dt\right)^{c_p}.
\end{equation}
\end{lemma}
\begin{proof}
Take $k_*>\max\{8,p,B_h\}$ as in Lemma~\ref{low:moment-result} and $k\in[k_*,k_*+1]$.
For $K_0\ge r_1$ and a fixed $N>1+\Lambda(K_0)$, set
\[
 \phi=\Lambda(\rho)-(1-\Delta)^{-1}\div(\rho u)-N.
\]
Let $G_B$ be the convolution kernel of $(1-\Delta)^{-1}$. By \cite[(4,2)--(4,6), pp.~416--417]{AronszajnSmith1961}, we have $G_B>0$, $\int_{\mathbb R^2} G_B\,dx=1$, and
\begin{equation}\label{wh3:bessel-kernel}
\begin{gathered}
 G_B(x)\le C\log\!\left(e+\frac{1}{|x|}\right)e^{-|x|/2},\qquad
 |\nabla G_B(x)|\le\frac{e^{-|x|/\sqrt2}}{\pi|x|}\quad(x\ne0),\\
 \nabla G_B\in L^q(\mathbb R^2)\quad(1\le q<2).
\end{gathered}
\end{equation}
Observe that $\rho u=(\sqrt\rho-\sqrt{\rho_\infty})\sqrt\rho\,u+\sqrt{\rho_\infty}\sqrt\rho\,u$. Hence we obtain from \eqref{wh3:primitives}, $|\sqrt\rho-\sqrt{\rho_\infty}|^2\le|\rho-\rho_\infty|$ and H\"older's inequality that
\[
 \|(\sqrt\rho-\sqrt{\rho_\infty})\sqrt\rho\,u\|_1\le C,\qquad
 \|(\sqrt\rho-\sqrt{\rho_\infty})\sqrt\rho\,u\|_{2k/(k+1)}
 \le C\|\rho-\rho_\infty\|_k^{1/2}.
\]
By \eqref{wh3:bessel-kernel} and the Hardy--Littlewood--Sobolev and Young inequalities, we have
\begin{equation}\label{s5:bessel-l2k}
 \|(1-\Delta)^{-1}\div(\rho u)\|_{2k}\le C(1+\|\rho-\rho_\infty\|_k^{1/2}),
\end{equation}
and
\begin{equation}\label{s5:bessel-l1l2}
 \|(1-\Delta)^{-1}\div(\rho u)\|_{L^1+L^2}\le C.
\end{equation}
Since $|(1-\Delta)^{-1}\div(\rho u)|>1$ on $\{\rho\le K_0,\ \phi>0\}$, it follows from \eqref{wh3:primitives} and \eqref{s5:bessel-l1l2} that
\begin{equation}\label{s5:truncated-mass}
 \int_{\{\phi>0\}}\rho\,dx
 \le K_0|\{|(1-\Delta)^{-1}\div(\rho u)|>1\}|+\int_{\{\rho>K_0\}}\rho\,dx
 \le C.
\end{equation}
By \eqref{extlow:struct} on $\{\rho>K_0\}$ and \eqref{wh3:primitives} on $\{\rho\le K_0\}$, we have
\[
 \|\rho-\rho_\infty\|_k^k\le C\left(1+\int_{\{\rho>K_0\}}\rho\Lambda(\rho)^k\,dx\right).
\]
On $\{\rho>K_0\}$, since
\[
 \Lambda(\rho)\le N+\phi_+
                      +|(1-\Delta)^{-1}\div(\rho u)|,
\]
by \eqref{wh3:primitives}, \eqref{s5:bessel-l2k} and H\"older's inequality, we obtain
\[
\begin{aligned}
 1+\int_{\{\rho>K_0\}}\rho\Lambda(\rho)^k\,dx
\le C\left(1+\int_{\mathbb R^2}\rho\phi_+^k\,dx\right)
+C\left(1+\int_{\{\rho>K_0\}}\rho\Lambda(\rho)^k\,dx\right)^{\frac{1}{2}+\frac{1}{2(k-1)}}.
\end{aligned}
\]
Since $\frac{1}{2}+\frac{1}{2(k-1)}<1$, by \eqref{s5:bessel-l2k} and Young's inequality, we get
\begin{equation}\label{s5:moment-comparison}
\begin{aligned}
 1+\int_{\{\rho>K_0\}}\rho\Lambda(\rho)^k\,dx
 &\le C\left(1+\int_{\mathbb R^2}\rho\phi_+^k\,dx\right),\\
 \|(1-\Delta)^{-1}\div(\rho u)\|_{2k}^{2k}+\|\rho-\rho_\infty\|_k^k
 &\le C\left(1+\int_{\{\rho>K_0\}}\rho\Lambda(\rho)^k\,dx\right).
\end{aligned}
\end{equation}

By $\eqref{fns}_1$, we have
\begin{equation}\label{wh3:lambda-transport}
 \dot{\Lambda(\rho)}=-(2\mu+\lambda)\div u,\qquad
 \Lambda(r)-\Lambda(\rho_\infty)=\int_{\rho_\infty}^r
                                      \frac{2\mu+\lambda(z)}{z}\,dz.
\end{equation}
Applying $(1-\Delta)^{-1}\div$ to $\rho\dot u^i=\partial_t(\rho u^i)+\partial_l(\rho u^i u^l)$ and using $\Delta G=\div(\rho\dot u)$, we obtain
\begin{equation}\label{wh3:psi-transport}
 \begin{aligned}
 \partial_t\bigl((1-\Delta)^{-1}\div(\rho u)\bigr)
+(1-\Delta)^{-1}\partial_i\partial_l(\rho u^i u^l)
=-\bigl(1-(1-\Delta)^{-1}\bigr)(G+P_\infty).
 \end{aligned}
\end{equation}
Combining \eqref{wh3:lambda-transport} with \eqref{wh3:psi-transport} and $G=(2\mu+\lambda)\div u-P$ yields
\begin{equation}\label{wh3:bessel}
 \dot\phi+P-P_\infty
 =-[u^l,\partial_i\partial_l(1-\Delta)^{-1}](\rho u^i)
            -(1-\Delta)^{-1}(G+P_\infty).
\end{equation}

Set
\[
 z=\rho\phi_+^{k-1},\qquad
 T_u g=[u^l,\partial_i\partial_l(1-\Delta)^{-1}]g^i.
\]
Multiplying \eqref{wh3:bessel} by $kz$, integrating over $\mathbb{R}^2$ and using $G=h\div u-P$ and $(1-\Delta)^{-1}P_\infty=P_\infty$, we obtain
\begin{equation}\label{s5:Ik-id}
 \frac{d}{dt}\int_{\mathbb R^2}\rho\phi_+^k\,dx
 =-k\int_{\mathbb R^2} z\bigl[T_u(\rho u)+P+(1-\Delta)^{-1}(h\div u-P)\bigr]\,dx.
\end{equation}

For $f\in H^1$, by \eqref{gn11}, \eqref{wh3:low-density-measure} and H\"older's inequality, we have
\[
 \|f\|_2^2\le C\|\sqrt\rho\,f\|_2^2+C\|f\|_2\|\nabla f\|_2,
\]
and Young's inequality yields
\begin{equation}\label{wh3:anchored}
 \|f\|_2\le C\bigl(\|\sqrt\rho\,f\|_2+\|\nabla f\|_2\bigr).
\end{equation}
By \eqref{gn11}, we have $\|u\|_{2j}^{2j}\le (Cj)^j\|u\|_2^2\|\nabla u\|_2^{2j-2}$ for all $j\ge 1$. Using \eqref{wh3:primitives} and \eqref{wh3:anchored} with $f=u$, we obtain, for $c>0$ sufficiently small and $\|\nabla u\|_2>0$,
\begin{equation}\label{wh3:velocity-exponential}
\begin{aligned}
 \int_{\mathbb R^2}\left(e^{c|u|^2/\|\nabla u\|_2^2}-1\right)dx
= \sum_{j=1}^{\infty}\frac{c^j}{j!\,\|\nabla u\|_2^{2j}}\|u\|_{2j}^{2j}
\le \frac{\|u\|_2^2}{\|\nabla u\|_2^2}\sum_{j=1}^{\infty}\frac{(Cc j)^j}{j!}
 \le C\bigl(1+\|\nabla u\|_2^{-2}\bigr).
\end{aligned}
\end{equation}

For $\|\nabla u\|_2>0$, we decompose $\rho u=g_c+g_s$, where
\[
 \begin{aligned}
 g_c =(\rho-\rho_\infty)u\,{\bf1}_{\{|u|\le C_0\|\nabla u\|_2\sqrt{\log(e+|\rho-\rho_\infty|)}\}},\quad
 g_s = \rho u-g_c.
 \end{aligned}
\]
With $\frac{1}{2}+\beta_D<m<1$ in the Young function $B_k$ from \eqref{extlow:comm}, the commutator estimate \cite{CRW} gives, for fixed $q>2$,
\[
 \|T_u g\|_q\le C_q \|\nabla u\|_2\|g\|_q,\quad
 c\,\frac{s^{k+1}}{[\log(e+s)]^{m(k+1)}}\le B_k(s)\le
 C\,\frac{s^{k+1}}{[\log(e+s)]^{m(k+1)}}.
\]
Since $m>\frac{1}{2}+\beta_D$, for large $|\rho-\rho_\infty|$ we have
\[
 B_k\bigl(|\rho-\rho_\infty|\sqrt{\log(e+|\rho-\rho_\infty|)}\bigr)+|\rho-\rho_\infty|^k [\log(e+|\rho-\rho_\infty|)]^{k/2}
 \le C\frac{\rho^{k+1}}{D(\rho)^k}.
\]
Using \eqref{wh3:primitives} for bounded $|\rho-\rho_\infty|$ and \eqref{s5:moment-comparison} for large $|\rho-\rho_\infty|$, we arrive at
\[
 \int_{\mathbb R^2} B_k\!\left(\frac{|g_c|}{C\|\nabla u\|_2}\right)dx
 \le C\left(1+\int_{\mathbb R^2}\rho\phi_+^k\,dx\right),
\]
and
\[
 \|g_c\|_k
 \le C\|\nabla u\|_2\left(1+\int_{\mathbb R^2}\rho\phi_+^k\,dx\right)^{1/k}.
\]
Similar to \eqref{extlow:comm}, with $m+(m+\beta_D)/k_*<\xi<1$, we obtain, for $0<\delta\le 1$,
\begin{equation}\label{s5:comm-gc}
\begin{aligned}
 k\int_{\mathbb R^2} z|T_u g_c|\,dx
 &\le C(1+\delta^{-\xi/(1-\xi)})\|\nabla u\|_2^2
       \left(1+\int_{\mathbb R^2}\rho\phi_+^k\,dx\right)\\
 &\quad+\delta\|\nabla u\|_2^2\frac{\partial}{\partial k}
 \left(\int_{\mathbb R^2}\rho\bigl[\max\{1,\phi\}^k-1\bigr]\,dx\right).
\end{aligned}
\end{equation}
On $\{|u|>C_0\|\nabla u\|_2\sqrt{\log(e+|\rho-\rho_\infty|)}\}$, we have $|\rho-\rho_\infty|\le\exp\bigl(|u|^2/(C_0\|\nabla u\|_2)^2\bigr)$. Taking $C_0$ sufficiently large, we deduce from \eqref{gn11}, \eqref{wh3:primitives}, \eqref{wh3:anchored}, \eqref{wh3:velocity-exponential} and the commutator estimate \cite{CRW} that
\[
 \|g_s\|_q\le C\bigl(\|\nabla u\|_2+\|\nabla u\|_2^{1-2/q}\bigr),\qquad
 \|T_u g_s\|_q\le C\bigl(\|\nabla u\|_2^2+\|\nabla u\|_2^{2-2/q}\bigr).
\]
For $i=1,2$, by \eqref{s5:moment-comparison} and H\"older's inequality, we have
\begin{equation}\label{s5:dual-fi}
\begin{aligned}
 \int_{\mathbb R^2} z|f|^i\,dx
 &\le\left(\int_{\mathbb R^2}\rho\phi_+^k\,dx\right)^{1-1/k}
           \left(\int_{\{\phi>0\}}\rho|f|^{ik}\,dx\right)^{1/k}\\
 &\le C\left(1+\int_{\mathbb R^2}\rho\phi_+^k\,dx\right)^{1-1/k+1/k^2}
           \bigl(\|f\|_{ik}+\|f\|_{ik^2/(k-1)}\bigr)^i.
\end{aligned}
\end{equation}
Taking $i=1$ in \eqref{s5:dual-fi} for $\|\nabla u\|_2\ge 1$, and using $|T_u g_s|\le 1+|T_u g_s|^2$ with $i=2$ for $\|\nabla u\|_2<1$, we obtain
\begin{equation}\label{s5:comm-gs}
 \int_{\mathbb R^2} z|T_u g_s|\,dx\le C \left(1+\int_{\mathbb R^2}\rho\phi_+^k\,dx\right) \|\nabla u\|_2^2+\int_{\mathbb R^2} z\,dx.
\end{equation}

By \eqref{visc:envelope}, \eqref{wh3:primitives} and H\"older's inequality, we have
\[
 \begin{aligned}
 \int_{\mathbb R^2}\frac{|h-h(\rho_\infty)|^2}{h}\,dx
\le C+C\int_{\{\rho>K_0\}}h\,dx
\le C\left(1+\int_{\{\rho>K_0\}}\rho\Lambda(\rho)^k\,dx\right)^{(1-1/B_h)/k}.
 \end{aligned}
\]
Writing $h\div u=h(\rho_\infty)\div u+(h-h(\rho_\infty))\div u$, by Young's convolution inequality and Cauchy--Schwarz, we obtain, for the exponents in \eqref{s5:dual-fi},
\[
\begin{aligned}
 \|(1-\Delta)^{-1}(h\div u)\|_q
 &\le C_q\left(1+\int_{\{\rho>K_0\}}\rho\Lambda(\rho)^k\,dx\right)^{(1-1/B_h)/(2k)}\\
 &{}\times\left(\int_{\mathbb R^2}\bigl(2\mu|D(u)|^2+\lambda(\div u)^2\bigr)dx\right)^{1/2}.
\end{aligned}
\]
Substituting this into \eqref{s5:dual-fi} with $i=2$ and using \eqref{s5:moment-comparison} and $k\ge B_h$, we obtain
\begin{equation}\label{s5:B-hdiv}
\begin{aligned}
 \int_{\mathbb R^2} z|(1-\Delta)^{-1}(h\div u)|^2\,dx
\le C\left(1+\int_{\mathbb R^2}\rho\phi_+^k\,dx\right)
\int_{\mathbb R^2}\bigl(2\mu|D(u)|^2+\lambda(\div u)^2\bigr)\,dx.
\end{aligned}
\end{equation}

Since $e^{cG_B}\in L^1(B_1)$ and $G_B\in L^1\cap L^\infty(\mathbb{R}^2\setminus B_1)$, it follows from \eqref{wh3:primitives} and Young's inequality that
\[
 \|(1-\Delta)^{-1}(P{\bf1}_{\{\theta\le2\theta_\infty\}})\|_\infty
 \le 2\theta_\infty\|(1-\Delta)^{-1}\rho\|_\infty\le C.
\]
For $0<h_0<\theta_\infty<L_0$, since $\theta-\theta_\infty-\theta_\infty\log(\theta/\theta_\infty)\ge c>0$ on $\{\theta<h_0\}\cup\{\theta>L_0\}$, it follows from \eqref{wh3:primitives} and \eqref{wh3:low-density-measure} that
\begin{equation}\label{wh3:cutoff-set-measure}
 |\{\theta<h_0\}|+|\{\theta>L_0\}|
 \le|\{\rho<\rho_\infty/2\}|+
       C\int_{\mathbb R^2}\rho\left(\theta-\theta_\infty-\theta_\infty\log\frac{\theta}{\theta_\infty}\right)\,dx\le C.
\end{equation}
For $f\in H^1$ with $|\operatorname{supp} f|<\infty$, by Sobolev's and H\"older's inequalities, we have
\begin{equation}\label{wh3:sobolev-support}
 \|f\|_2\le C\|\nabla f\|_1
      \le C|\operatorname{supp} f|^{1/2}\|\nabla f\|_2.
\end{equation}
Applying \eqref{gn11} and \eqref{wh3:sobolev-support} to $(\theta^{\alpha/2}-(3\theta_\infty/2)^{\alpha/2})_+$ and using \eqref{wh3:cutoff-set-measure} with $j\ge\max\{2,2q/\alpha,2nq\}$, we obtain, for $q,n\ge1$,
\begin{equation}\label{s5:small-hot-tail}
 \|\theta{\bf1}_{\{\theta>2\theta_\infty\}}\|_q
 \le C_{q,n}\|\sqrt{\kappa(\theta)}\,\nabla\log\theta\|_2^{2n},
\end{equation}
provided that $\|\sqrt{\kappa(\theta)}\,\nabla\log\theta\|_2\le1$.

By \eqref{wh3:theta-hot-energy}, we have $\|P{\bf1}_{\{\theta>2\theta_\infty\}}\|_1\le C$. When $\|\sqrt{\kappa(\theta)}\,\nabla\log\theta\|_2\le 1$, combining this with \eqref{s5:moment-comparison}, \eqref{s5:small-hot-tail} (with $q=k/(k-1)$ and $n=2$) and H\"older's inequality yields
\begin{equation}\label{s5:Ph-L1}
 \|P{\bf1}_{\{\theta>2\theta_\infty\}}\|_1
 \le C\left(1+\int_{\mathbb R^2}\rho\phi_+^k\,dx\right)^{1/(2k)}
      \|\sqrt{\kappa(\theta)}\,\nabla\log\theta\|_2^2 .
\end{equation}
When $\|\sqrt{\kappa(\theta)}\,\nabla\log\theta\|_2>1$, \eqref{s5:Ph-L1} follows directly from $\|P{\bf1}_{\{\theta>2\theta_\infty\}}\|_1\le C$. By \eqref{s5:dual-fi} with $i=1$, \eqref{s5:Ph-L1} and Young's convolution inequality, we obtain
\begin{equation}\label{s5:B-Ph}
 \int_{\mathbb R^2} z(1-\Delta)^{-1}(P{\bf1}_{\{\theta>2\theta_\infty\}})\,dx
 \le C\left(1+\int_{\mathbb R^2}\rho\phi_+^k\,dx\right)
      \|\sqrt{\kappa(\theta)}\,\nabla\log\theta\|_2^2 .
\end{equation}

For $0<a<h_0<\theta_\infty$, taking $w=(h_0-\max\{a,\theta\})_+$ in \eqref{gn11} and \eqref{wh3:sobolev-support} gives, for $n\ge1$,
\begin{equation}\label{s5:cold-set}
 |\{\theta\le a\}|\le\frac{\|w\|_{2n}^{2n}}{(h_0-a)^{2n}}
 \le C_{a,n}\|\nabla w\|_2^{2n}
 \le C_{a,n}\|\sqrt{\kappa(\theta)}\,\nabla\log\theta\|_2^{2n}.
\end{equation}

Combining \eqref{s5:Ik-id} with \eqref{s5:comm-gc}, \eqref{s5:comm-gs}, \eqref{s5:B-hdiv}, \eqref{s5:B-Ph} and Young's inequality, we obtain
\begin{equation}\label{s5:Ik-reduced}
\begin{aligned}
 \frac{d}{dt}\int_{\mathbb R^2}\rho\phi_+^k\,dx
 &\le C(1+\delta^{-\xi/(1-\xi)})A_2\left(1+\int_{\mathbb R^2}\rho\phi_+^k\,dx\right)\\
 &{}+\delta\|\nabla u\|_2^2\frac{\partial}{\partial k}
 \left(\int_{\mathbb R^2}\rho\bigl[\max\{1,\phi\}^k-1\bigr]\,dx\right)
 +k\int_{\mathbb R^2} z(C_*-P)\,dx,
\end{aligned}
\end{equation}
where $C_*>0$ is a constant independent of $K_0$ and $N$. Choose $K_0\ge r_1$ with $aK_0>C_*$ and $N>1+\Lambda(K_0)$. Since $P=\rho\theta$ and $z=0$ on $\{\phi\le 0\}$, we have $C_*-P\le C_*-aK_0<0$ on $\{\rho>K_0,\ \theta>a\}$, and hence
\begin{equation}\label{s5:drift-split}
 k\int_{\mathbb R^2} z(C_*-P)\,dx
 \le kC_*\left(\int_{\{\theta\le a\}}z\,dx
                 +\int_{\{\phi>0,\ \rho\le K_0\}}z\,dx\right).
\end{equation}
By \eqref{s5:truncated-mass}, \eqref{s5:moment-comparison}, \eqref{s5:cold-set} (with $n$ large when $\|\sqrt{\kappa(\theta)}\,\nabla\log\theta\|_2\le 1$) and H\"older's inequality, we obtain
\begin{equation}\label{s5:z-cold}
 \int_{\{\theta\le a\}}z\,dx
 \le C\left(1+\int_{\mathbb R^2}\rho\phi_+^k\,dx\right)
      \|\sqrt{\kappa(\theta)}\,\nabla\log\theta\|_2^2 .
\end{equation}
On $\{\phi>0,\ \rho\le K_0\}$, we have $1+\phi_+\le|(1-\Delta)^{-1}\div(\rho u)|$. By \eqref{wh3:bessel-kernel} and the $L^k$ estimates for $g_c$ and $g_s$, we get
\[
\begin{aligned}
 \|(1-\Delta)^{-1}\div(\rho u)\|_\infty
 &\le C\|g_c+g_s\|_k\\
 &\le C\|\nabla u\|_2\left(1+\int_{\mathbb R^2}\rho\phi_+^k\,dx\right)^{1/k}
+ C \bigl(\|\nabla u\|_2+\|\nabla u\|_2^{1-2/k}\bigr).
\end{aligned}
\]
Taking $i=2$ for $\|\nabla u\|_2\ge 1$ and $i=4$ for $\|\nabla u\|_2<1$ yields
\begin{equation}\label{s5:z-lowrho}
\begin{aligned}
 \int_{\{\phi>0,\ \rho\le K_0\}}z\,dx
 &\le C\left(1+\int_{\mathbb R^2}\rho\phi_+^k\,dx\right)^{1-i/k}\|(1-\Delta)^{-1}\div(\rho u)\|_\infty^i\\
 &\le C\|\nabla u\|_2^2\left(1+\int_{\mathbb R^2}\rho\phi_+^k\,dx\right).
\end{aligned}
\end{equation}
Substituting \eqref{s5:drift-split}--\eqref{s5:z-lowrho} into \eqref{s5:Ik-reduced}, we obtain, for $0<\delta\le 1$ and $k\in[k_*,k_*+1]$,
\[
\begin{aligned}
 \frac{d}{dt}\int_{\mathbb R^2}\rho\phi_+^k\,dx
 \le& C(1+\delta^{-\xi/(1-\xi)})A_2\left(1+\int_{\mathbb R^2}\rho\phi_+^k\,dx\right)\\
 &+\delta\|\nabla u\|_2^2\frac{\partial}{\partial k}
 \left(\int_{\mathbb R^2}\rho\bigl[\max\{1,\phi\}^k-1\bigr]\,dx\right).
\end{aligned}
\]

Choose
\[
 \begin{aligned}
 k(t)=k_*+1-\frac{\int_0^t\|\nabla u\|_2^2\,ds}{2+\int_0^T A_2\,ds},\quad
 \delta=\frac{1}{2\bigl(2+\int_0^T A_2\,ds\bigr)},
 \end{aligned}
\]
so that $k(t)\in[k_*,k_*+1]$ and $k'+\delta\|\nabla u\|_2^2\le 0$. By \eqref{s5:truncated-mass}, we arrive at
\[
 \frac{d}{dt}\log\left(1+\int_{\mathbb R^2}\rho\phi_+^{k(t)}\,dx\right)
 \le C(1+\delta^{-\xi/(1-\xi)})A_2+C|k'|.
\]
Integrating over $[0,t]$ and using \eqref{extlow:struct}, \eqref{s5:moment-comparison} and $k(t)>p$ yields \eqref{wh3:drift-bound}.
\end{proof}

We next derive the estimates for $\int_0^T A_2\,dt$ and the $L^p$ norms of $\rho-\rho_\infty$.

\begin{lemma}\label{plane:low-transfer}
\label{wh3:moment}
In case \textup{(A)}, assume that
\[
 h(r)\ge\frac{c_*r}{2D(r)}
 \qquad\left(r_*\le r\le\sup_{0\le t\le T}\|\rho(t)\|_\infty\right).
\]
In the data-class alternative \eqref{ext:gainfamily}, assume likewise that
\[
 \mu+\lambda(r)\ge\frac{c_*}{2}rL(r)^d
 \qquad\left(r_*\le r\le\sup_{0\le t\le T}\|\rho(t)\|_\infty\right).
\]
Then there exist constants $l>4\max\{1,\theta_\infty\}$, $C$, and
$C_p$, independent of $\sup_{0\le t\le T}\|\rho(t)\|_\infty$, such that
\begin{equation}\label{wh3:rho-moment-1}
 \int_0^T A_2\,dt
 +\int_0^T\!\int_{\{\theta>l\}}
       \frac{P^2}{\mu+\lambda}\,dx\,dt
 \le C,
\end{equation}
and
\begin{equation}\label{wh3:rho-moment-2}
 \sup_{0\le t\le T}\|\rho(t)-\rho_\infty\|_p
 \le C_p
 \qquad(2\le p<\infty).
\end{equation}
Here $l$, $C$, and $C_p$ are independent of $T$ under \textup{(A)},
\textup{(B)}, or \textup{(C)} with $\alpha>1-d$, whereas under
\textup{(D)} or \textup{(C)} with
$2(1-d)/(2+d)\le\alpha\le 1-d$ they depend on $T_0$ for
$0<T\le T_0<\infty$.
\end{lemma}

\begin{proof}
By \eqref{wh3:low-density-measure} and \eqref{wh3:theta-hot-energy},
we get
\begin{equation}\label{wh3:theta-hot-measure}
 |\{\theta>2\theta_\infty\}|
 \le |\{\rho<\rho_\infty/2\}|
      +\frac{1}{P_\infty}
       \int_{\{\theta>2\theta_\infty\}}\rho\theta\,dx
 \le C.
\end{equation}
For $l>4\theta_\infty$ and $0\le\nu<2$, it follows from
\eqref{wh3:theta-hot-energy} and \eqref{wh3:theta-hot-measure} that
\[
 \sup_{0\le t\le T}\int_{\{\theta>l\}}\rho\theta\,dx\le C,
 \qquad
 \bigl|\operatorname{supp}
   \bigl(\theta^{1-\nu/2}-(l/2)^{1-\nu/2}\bigr)_+\bigr|\le C.
\]
Moreover, for any $L_0>\theta_\infty$ and $z\in H^1$ supported in
$\{\theta>L_0\}$ with $\|\nabla z\|_2>0$, we have
$|\{\theta>L_0\}|\le C$ by \eqref{wh3:cutoff-set-measure}. Thus the
P\'olya--Szeg\H{o} inequality
\cite[(6.22)]{BianchiGardnerGronchiKiderlen2024} and the
Trudinger--Moser inequality \cite[Theorem~1]{Moser1971} yield
\begin{equation}\label{wh3:pressure-hot-bound}
 \int_{\{z\ne0\}}
 \exp\!\left(\frac{4\pi z^2}{\|\nabla z\|_2^2}\right)dx
 \le C|\{z\ne0\}|\le C.
\end{equation}
Below we apply \eqref{wh3:pressure-hot-bound} with $L_0=l/2$ and
$z=\bigl(\theta^{1-\nu/2}-(l/2)^{1-\nu/2}\bigr)_+$.

Since
$2|D(u)|^2=|\nabla u|^2+
\sum_{i,j=1}^2\partial_i u^j\,\partial_j u^i$,
integrating by parts over $\mathbb R^2$ gives
\[
 \sum_{i,j=1}^2\int_{\mathbb R^2}
     \partial_i u^j\,\partial_j u^i\,dx
 =-\int_{\mathbb R^2}u\cdot\nabla(\div u)\,dx
 =\int_{\mathbb R^2}(\div u)^2\,dx,
\]
and hence
\[
 \int_{\mathbb R^2}
   \bigl(2\mu|D(u)|^2+\lambda(\div u)^2\bigr)\,dx
 =\int_{\mathbb R^2}
   \bigl(\mu|\nabla u|^2+(\mu+\lambda)(\div u)^2\bigr)\,dx
 \ge\mu\|\nabla u\|_2^2.
\]

Let $l\ge4(1+\theta_\infty)$. Since $r_1>2\rho_\infty$,
the relative entropy in \eqref{wh3:primitives} controls
$\int_{\{\rho>r_1\}}\rho L(\rho)\,dx$. Therefore, by
\eqref{wh3:theta-hot-measure},
\[
 \int_{\{\theta>l\}}\rho L(\rho)\,dx
 \le r_1L(r_1)|\{\theta>l\}|
     +\int_{\{\rho>r_1\}}\rho L(\rho)\,dx
 \le C
\]
uniformly in $l$ and $t$.

For $0\le\nu<1$ with $a=\alpha+\nu-1>0$, set
\[
 f_l=\mathbf 1_{\{\theta>l\}}
       \frac{\rho^2\theta^\nu}{\mu+\lambda(\rho)}.
\]
On any interval $I\subset[0,T]$, write
\[
 F_l(I)=1+\sup_{s\in I}\int_{\mathbb R^2}f_l\log(e+f_l)\,dx,
 \qquad
 \mathcal P_l(I)=\int_I\!\int_{\{\theta>l\}}
                   \frac{P^2}{\mu+\lambda(\rho)}\,dx\,ds.
\]
The two temperature tests used to derive \eqref{extlow:band} apply
on $\mathbb R^2$. Indeed, the absolute values of their primitives
are bounded by $\theta\mathbf 1_{\{\theta>l\}}$, so
\eqref{wh3:theta-hot-energy} controls the endpoint terms.
Using \eqref{wh3:pressure-hot-bound} for the truncated temperature,
\eqref{wh3:primitives}, and the preceding bound for
$\|\nabla u\|_2^2$, we obtain
\[
 \begin{aligned}
 \int_I A_2(s)\,ds
 &\le C\bigl(1+l+\mathcal P_l(I)\bigr),\\
 \mathcal P_l(I)
 &\le \frac{C(1+a^{-1})F_l(I)}{\kappa_0 l^a}
          \bigl(1+l+\mathcal P_l(I)\bigr).
 \end{aligned}
\]
The constants are independent of $I$ and $l$.

Suppose first that $\alpha>1$. Take $\nu=0$, so $a=\alpha-1$, and set
\[
 U(t)=2+\int_0^t A_2(s)\,ds.
\]
By \eqref{wh3:drift-bound} with $p=3$ on $[0,t]$, there exist
$C_0,c_3>0$, independent of $T$ and the density upper bound, such that
\[
 \sup_{0\le s\le t}\int_{\{\rho>r_1\}}\rho^3\,dx
 \le \exp\!\bigl(C_0U(t)^{c_3}\bigr).
\]
Fix $R\ge r_1$. On $\{\theta>l,\ \rho\le R\}$,
\eqref{extlow:struct}, the monotonicity of $D$, and
$\mu+\lambda\ge\mu$ give
\[
 f_l\le C D(R)\rho,
 \qquad
 \log(e+f_l)\le C L(\rho).
\]
On $\{\rho>R\}$, using $D(\rho)\le C L(\rho)^{\beta_D}$, we have
\[
 f_l\log(e+f_l)
 \le C\rho L(\rho)^{1+\beta_D}
 \le C\rho^2
 \le \frac{C}{R}\rho^3.
\]
Splitting the integral at $\rho=R$ thus gives
\[
 \begin{aligned}
 F_l([0,t])
 &\le 1+C D(R)\sup_{0\le s\le t}
               \int_{\{\theta>l\}}\rho L(\rho)\,dx
       +\frac{C}{R}\sup_{0\le s\le t}
               \int_{\{\rho>R\}}\rho^3\,dx\\
 &\le C\bigl(1+D(R)\bigr)
       +\frac{C}{R}\exp\!\bigl(C_0U(t)^{c_3}\bigr).
 \end{aligned}
\]

Choose $0<\gamma<a/c_3$. There are arbitrarily large $R$ such that
$D(R)\le(\log R)^\gamma$. In case \textup{(A)}, the opposite
inequality for all sufficiently large $R$ would make the integral
in \eqref{new:den} finite. In case \textup{(B)}, this follows from
\eqref{new:fixed} with $D=D_h$, and in case \textup{(C)} we have
$D=1$. For such an $R$, set
\[
 M_R=\left(\frac{\log R}{C_0}\right)^{1/c_3},
 \qquad
 l=4(1+\theta_\infty)+K\bigl(1+D(R)\bigr)^{1/a}.
\]
Whenever $U(t)\le M_R$, we have
$R^{-1}\exp(C_0U(t)^{c_3})\le1$, and hence
$F_l([0,t])\le C(1+D(R))$. We can therefore fix $K$, independently
of $R$ and $T$, so large that
\[
 \frac{C(1+a^{-1})F_l([0,t])}{\kappa_0l^a}\le\frac12
 \qquad\text{whenever }U(t)\le M_R.
\]
Absorbing the term containing $\mathcal P_l([0,t])$ gives
\[
 U(t)+\mathcal P_l([0,t])\le C_1(1+l)
 \qquad\text{whenever }U(t)\le M_R.
\]
Since $\gamma/a<1/c_3$ and $D(R)\le(\log R)^\gamma$ along the
chosen sequence, we may now fix $R$ so large that
\[
 2<M_R,\qquad C_1(1+l)<M_R.
\]
If the continuous function $U$ first attained $M_R$ at some
$t_*\le T$, the preceding estimate would give
$M_R=U(t_*)\le C_1(1+l)<M_R$. Since $U(0)=2$, it follows that
$U(t)<M_R$ throughout $[0,T]$. Taking $t=T$ with this fixed $l$
proves \eqref{wh3:rho-moment-1}.

For \textup{(C)} with $1-d<\alpha\le1$, take $\nu=d$, so that
$a=\alpha+d-1>0$. By \eqref{new:gain} on the attained density range and Young's inequality, as in
the proof of Lemma~\ref{low:uniform},
\[
 f_l\log(e+f_l)
 \le C\mathbf 1_{\{\theta>l\}}\rho\bigl(\theta+L(\rho)\bigr).
\]
The bounds for $\int_{\{\theta>l\}}\rho\theta\,dx$ and
$\int_{\{\theta>l\}}\rho L(\rho)\,dx$ give $F_l([0,T])\le C$.
Choosing $l$ sufficiently large and absorbing as above yields
\eqref{wh3:rho-moment-1}.

Under \textup{(D)}, take $R_0>\max\{r_*,2\rho_\infty\}$ in
\eqref{crit:moment-rate} and fix $l_0>2\theta_\infty$. With $b$
as in the proof of Lemma~\ref{lem:linear-moment-input},
\eqref{wh3:primitives} and \eqref{wh3:theta-hot-energy} give
\[
 \int_{\{\theta>l_0\}}\rho\theta\,dx\le C,
 \qquad
 \int_0^T b(t)\,dt\le C_{T_0}.
\]
Applying \eqref{wh3:pressure-hot-bound} to
$(\sqrt\theta-\sqrt{l_0})_+$ gives the exponential estimate used
in \eqref{crit:moment-rate}. For $v_k>2l_0$, we also have
\[
 \int_{\{v_k<\theta\le2v_k\}}\rho\,dx\le\frac{C}{v_k},
 \qquad
 \bigl|\operatorname{supp}
   \bigl(\min\{\theta,v_k\}-v_k/2\bigr)_+\bigr|\le C.
\]
Thus \eqref{wh3:pressure-hot-bound} applies to the temperature
truncation on each band. Using H\"older's inequality for the
weights $L(\rho)/D(\rho)$ and $L(\rho)^2/D(\rho)$ as in the proof
of Lemma~\ref{lem:linear-moment-input} gives
\eqref{low:band-pressure}. Combining \eqref{crit:moment-record}
with \eqref{crit:heat-window} as in that lemma yields
\[
 \sup_{0\le t\le T}
 \int_{\{\rho>2R_0\}}\rho L(\rho)^{3/2}\,dx\le C_{T_0}.
\]
Since $D(r)\le C L(r)^{1/4}$ and $3/2>1+1/4$, the condition in
\eqref{new:smallentropy} holds with $q=3/2$, $\tau=1/4$, and
$\alpha=1$. The estimate for $f_l\log(e+f_l)$ in the proof of
Lemma~\ref{low:uniform}, integrated over $\{\theta>l\}$, therefore
gives $F_l([0,T])\le C_{T_0}$ for the choice of $\nu$ made there.
Taking $l=l(T_0)$ sufficiently large yields
\eqref{wh3:rho-moment-1} for $0<T\le T_0$.

For \textup{(C)} with $2(1-d)/(2+d)\le\alpha\le1-d$, take
$R_0>\max\{r_*,2\rho_\infty\}$ and let
\[
 \mathcal M(t)=
 \left(1+\int_{\mathbb R^2}\rho
  \bigl(L(\rho)^{1+d/2}-L(R_0)^{1+d/2}\bigr)_+^2\,dx\right)^{1/2}
\]
as in the proof of Lemma~\ref{lem:gain-moment-input}. By
\eqref{wh3:theta-hot-energy}, \eqref{wh3:theta-hot-measure}, and the
definition of $\mathcal M$, we have
\[
 \int_{\{\theta>l\}}\rho\theta\,dx\le C,
 \qquad
 \int_{\{\theta>l\}}\rho L(\rho)^{2+d}\,dx
 \le C\mathcal M(t)^2.
\]
Set $s=(1-d)/(2+d)$ and $\nu=1-s$. On each subinterval
$I\subset[0,T]$, write $\mathcal M_I=\sup_{t\in I}\mathcal M(t)$.
Since $a=\alpha-s\ge s$, taking $l=C_*\mathcal M_I^2$ with $C_*$
sufficiently large gives, as in \eqref{low:gain-window},
\[
 F_l(I)\le C\mathcal M_I^{2s},
 \qquad
 \frac{C(1+a^{-1})F_l(I)}{\kappa_0l^a}\le\frac12.
\]
Consequently, the heat estimates on $I$ and the differential
inequality for $\mathcal M$ in \eqref{low:gain-window} yield
\[
 \mathcal M_I
 \le\mathcal M(\inf I)+C|I|^{1/2}\mathcal M_I.
\]
Iterating over subintervals with $C|I|^{1/2}\le1/2$ gives
$\sup_{0\le t\le T}\mathcal M(t)\le C_{T_0}$. Applying the heat
estimates on $[0,T]$ with $l=C_*C_{T_0}^2$ proves
\eqref{wh3:rho-moment-1} for $0<T\le T_0$.

Finally, for $p>2$, using \eqref{wh3:primitives} on
$\{\rho\le r_1\}$ and \eqref{wh3:rho-moment-1} together with
Lemma~\ref{plane:drift} on $\{\rho>r_1\}$, we obtain
\[
 \begin{aligned}
 \|\rho-\rho_\infty\|_p^p
 &\le C_p\int_{\{\rho\le r_1\}}
       \left(\rho\log\frac{\rho}{\rho_\infty}-\rho+\rho_\infty\right)dx
       +C_p\int_{\{\rho>r_1\}}\rho^p\,dx\\
 &\le C_p.
 \end{aligned}
\]
For $p=2$, since $\rho^2\le\rho^3/r_1$ on $\{\rho>r_1\}$, we also
have $\|\rho-\rho_\infty\|_2\le C$, which completes the proof of
\eqref{wh3:rho-moment-2}.
\end{proof}

Using Lemma~\ref{plane:low-transfer}, we next derive the following estimates for the dissipation, the relative pressure, and the temperature, with the same time dependence of the constants as in Lemma~\ref{plane:low-transfer}.

\begin{lemma}\label{wh3:thermal}
There exists a positive constant $C$ such that
\begin{equation}\label{wh3:thermal-dissipation}
 \int_0^T\!\int_{\mathbb R^2}
 \bigl(|\nabla u|^2+2\mu|D(u)|^2+\lambda(\div u)^2
               +\kappa(\theta)|\nabla\log\theta|^2\bigr)\,dx\,dt\le C.
\end{equation}
Moreover, for $0\le t_1<t_2\le T$ with $t_2-t_1\le 1$, we have
\begin{equation}\label{wh3:pressure-interval}
 \int_{t_1}^{t_2}\!\int_{\mathbb R^2}
       \frac{(P-P_\infty)^2}{2\mu+\lambda}\,dx\,dt\le C,
\end{equation}
and, for $2\le p<\infty$ and $0<\varepsilon\le\min\{1,\alpha\}$,
\begin{equation}\label{wh3:theta-interval}
 \int_{t_1}^{t_2}\|\theta-\theta_\infty\|_p^{1+\alpha-\varepsilon}\,dt
      \le \frac{Cp}{\varepsilon}.
\end{equation}
\end{lemma}
\begin{proof}
The estimate \eqref{wh3:thermal-dissipation} follows directly from \eqref{wh3:rho-moment-1}. Let $I=[t_1,t_2]$ and $L\ge l$. By \eqref{wh3:theta-cool-ltwo}, \eqref{wh3:rho-moment-2}, \eqref{wh3:theta-hot-measure} and $P-P_\infty=\theta(\rho-\rho_\infty)+\rho_\infty(\theta-\theta_\infty)$, we have
\[
 \|(P-P_\infty){\bf1}_{\{\theta\le L\}}\|_2
 \le L\|\rho-\rho_\infty\|_2+\rho_\infty\|(\theta-\theta_\infty){\bf1}_{\{\theta\le L\}}\|_2
 \le C_L,
\]
and since $h\ge 2\mu$,
\begin{equation}\label{wh3:pressure-cool}
 \int_{\{\theta\le L\}}\frac{|P-P_\infty|^2}{2\mu+\lambda}\,dx\le C_L.
\end{equation}
Combining \eqref{wh3:rho-moment-1} and \eqref{wh3:theta-hot-measure} with \eqref{wh3:pressure-cool}, we obtain
\[
\begin{aligned}
\int_I\!\int_{\mathbb R^2}\frac{(P-P_\infty)^2}{h}\,dx\,dt
\le C_L|I|+2\int_0^T\!\int_{\{\theta>l\}}\frac{P^2}{\mu+\lambda}\,dx\,dt+C P_\infty^2\int_I|\{\theta>L\}|\,dt
 \le C,
\end{aligned}
\]
which gives \eqref{wh3:pressure-interval}.

Set $c=(1+\alpha-\varepsilon)/2$ and $w=(\theta^c-(2L)^c)_+$. Multiplying $\eqref{fns}_3$ by $\chi_L(\theta)(L/\theta)^\varepsilon$ as in \eqref{extlow:heattail} and using \eqref{gn11} and \eqref{wh3:sobolev-support} with $j=\max\{2,p/c\}$, we obtain
\[
 \int_0^T\|\nabla w\|_2^2\,dt\le\frac{C}{\varepsilon},\qquad
 \|\theta{\bf1}_{\{\theta>4L\}}\|_p^{2c}\le C p\|\nabla w\|_2^2,
\]
and hence
\begin{equation}\label{s5:theta-hot-tail}
 \int_0^T\|\theta{\bf1}_{\{\theta>4L\}}\|_p^{1+\alpha-\varepsilon}\,dt\le\frac{Cp}{\varepsilon}.
\end{equation}
On $\{\theta\le 4L\}$, by \eqref{wh3:theta-cool-ltwo}, \eqref{wh3:theta-hot-measure} and H\"older's inequality, we have
\begin{equation}\label{s5:theta-cool-lp}
 \|\theta-\theta_\infty\|_{L^p(\{\theta\le4L\})}
 \le\|\theta-\theta_\infty\|_{L^2(\{\theta\le4L\})}^{2/p}\|\theta-\theta_\infty\|_{L^\infty(\{\theta\le4L\})}^{1-2/p}\le C_L.
\end{equation}
Combining \eqref{s5:theta-hot-tail} with \eqref{s5:theta-cool-lp} yields \eqref{wh3:theta-interval}.

Replacing $p$ by $(1+\alpha)p$ and $\varepsilon$ by $(1+\alpha)\delta$ in \eqref{s5:theta-hot-tail} gives
\begin{equation}\label{s5:KE-tail}
 \int_0^T\|K(\theta){\bf1}_{\{\theta>4L\}}\|_p^{1-\delta}\,dt
 \le C\int_0^T\|\theta{\bf1}_{\{\theta>4L\}}\|_{(1+\alpha)p}^{(1+\alpha)(1-\delta)}\,dt
 \le\frac{Cp}{\delta}.
\end{equation}
For $s>2\theta_\infty$, we have
\[
 \int_{\theta_\infty}^s\bigl(K(v)-K(\theta_\infty)\bigr)\,dv
 =\int_{\theta_\infty}^s(s-v)\kappa(v)\,dv,
\]
which yields
\begin{equation}\label{wh3:potential-hot}
 c s^{2+\alpha}\le\int_{\theta_\infty}^s
             \bigl(K(v)-K(\theta_\infty)\bigr)\,dv
 \le sK(s)\qquad(s>2\theta_\infty).
\end{equation}
For $0<s\le 2\theta_\infty$, we have
\[
 \kappa(s)\le C\frac{\theta_\infty}{s^2}
 =C\frac{d^2}{ds^2}\left(s-\theta_\infty-\theta_\infty\log\frac{s}{\theta_\infty}\right).
\]
Since both $\int_{\theta_\infty}^s\bigl(K(v)-K(\theta_\infty)\bigr)\,dv$ and $s-\theta_\infty-\theta_\infty\log(s/\theta_\infty)$, together with their first derivatives, vanish at $s=\theta_\infty$, integrating twice from $\theta_\infty$ yields
\begin{equation}\label{wh3:potential-cool}
 \int_{\theta_\infty}^s\bigl(K(v)-K(\theta_\infty)\bigr)\,dv
 \le C\left(s-\theta_\infty-\theta_\infty\log\frac{s}{\theta_\infty}\right) \quad(0<s\le 2\theta_\infty).
\end{equation}
For small $\zeta>0$, set $b=\zeta/[2(1-\zeta)]$. Since $\theta\le C K(\theta)$ on $\{\theta>4L\}$, H\"older's inequality yields
\[
\begin{aligned}
 &c_v\int_{\{\theta>4L\}}\rho\int_{\theta_\infty}^{\theta}\bigl(K(v)-K(\theta_\infty)\bigr)\,dv\,dx\\
 & \le \, C\int_{\{\theta>4L\}}(\rho\theta)^{1-b}\rho^b K(\theta)^{1+b}\,dx\\
 & \le \, C\left(\int_{\{\theta>4L\}}\rho\theta\,dx\right)^{1-b}
    \|\rho\|_{L^2(\{\theta>4L\})}^b
    \|K(\theta){\bf1}_{\{\theta>4L\}}\|_{2(1+b)/b}^{1+b}.
\end{aligned}
\]
By \eqref{wh3:theta-hot-energy}, \eqref{wh3:rho-moment-2}, \eqref{wh3:theta-hot-measure} and \eqref{s5:KE-tail} with $\delta=\zeta/2$, since $(1+b)(1-\zeta)=1-\zeta/2$, we obtain
\begin{equation}\label{s5:H-hot-int}
 \int_0^T\left(c_v\int_{\{\theta>4L\}}\rho\int_{\theta_\infty}^{\theta}\bigl(K(v)-K(\theta_\infty)\bigr)\,dv\,dx\right)^{1-\zeta}dt\le\frac{C}{\zeta^2}.
\end{equation}
By \eqref{wh3:theta-hot-energy} on $\{2\theta_\infty<\theta\le 4L\}$ and \eqref{wh3:potential-cool} on $\{\theta\le 2\theta_\infty\}$, we have $c_v\int_{\{\theta\le 4L\}}\rho\int_{\theta_\infty}^\theta\bigl(K(v)-K(\theta_\infty)\bigr)\,dv\,dx\le C_L$. Together with \eqref{s5:H-hot-int}, this yields
\[
 \int_I(1+H)^{1-\zeta}\,dt\le C_L|I|+\frac{C}{\zeta^2}\le\frac{C_L}{\zeta^2}.
\]

Next, let $r\ge\max\{r_1,\sup_{0\le t\le T}\|\rho(t)\|_\infty\}$ and $0<l_1<l_2<h_0<\theta_\infty$. For $l_1\le k\le l_2$, set $w=(K(k)-K(\theta))_+$. By \eqref{wh3:cutoff-set-measure}, we have $|\{\theta<k\}|\le|\{\theta<h_0\}|\le C$, and it follows from \eqref{gn11}, \eqref{wh3:sobolev-support} and \eqref{wh3:rho-moment-2} that
\[
 \|\rho\|_{L^2(\{\theta<k\})}\le\|\rho-\rho_\infty\|_2+\rho_\infty|\{\theta<k\}|^{1/2}\le C,
\]
and
\[
 \left(\int_{\mathbb R^2}\rho w^4\,dx\right)^{1/2}\le C\|w\|_8^2\le C\|\nabla w\|_2^2.
\]
On any $[b,b+1]\subset[0,T]$, use the levels and time cutoffs from Lemma~\ref{low:temp-floor}.
\[
 \begin{gathered}
 k_i=l_1+2^{-i}(l_2-l_1),\quad s_i=b+(1-2^{-i})/2,\quad
 y_i=\int_{s_i}^{b+1}\!\int_{\{\theta<k_i\}}\rho\,dx\,d\tau .
 \end{gathered}
\]
The cutoff estimates from that proof with $q=6$ give
\be\nonumber\ba
y_{i+1} \le C C_1^i y_i^{5/4},\quad C_1\ge2,\quad
y_0 \le C\int_b^{b+1}\|\sqrt{\kappa(\theta)}\,\nabla\log\theta\|_2^2\,d\tau,
\ea\ee
and
\be\nonumber\ba
\int_{s_{i+1}}^{b+1}\|\nabla(K(k_i)-K(\theta))_+\|_2^2\,d\tau
\le C 2^i y_i^{5/6}.
\ea\ee
Choose $\epsilon>0$ with $C\epsilon\le (C C_1^4)^{-4}$, using one $C$ for both bounds on $y_i$. Then
\[
 \int_b^{b+1}\|\sqrt{\kappa(\theta)}\,\nabla\log\theta\|_2^2\,d\tau
 \le\epsilon
 \ \Longrightarrow\
 y_i\le C\epsilon C_1^{-4i}\quad(i\ge0).
\]
By \eqref{wh3:sobolev-support} and $|\{\theta<k_i\}|\le C$,
\[
 \int_{b+1/2}^{b+1}\|(K(l_1)-K(\theta))_+\|_2^2\,d\tau
 \le C2^i y_i^{5/6}\longrightarrow0.
\]
By continuity, $\theta(x,b+1)\ge l_1$ for all $x\in\mathbb R^2$.

For the time-uniform bound, fix $t\in[0,T]$ and set
\[
 s=\max\!\left(\{0\}\cup
 \left\{j\in\{1,\ldots,\lfloor t\rfloor\}:
 \int_{j-1}^{j}\|\sqrt{\kappa(\theta)}\,\nabla\log\theta\|_2^2\,d\tau\le\epsilon
 \right\}\right).
\]
Each unit interval $[j-1,j]$ with $s<j\le\lfloor t\rfloor$ has dissipation greater than $\epsilon$. By \eqref{wh3:primitives},
\[
 \begin{aligned}
 \epsilon(\lfloor t\rfloor-s)
\le \sum_{j=s+1}^{\lfloor t\rfloor}\int_{j-1}^{j}
       \|\sqrt{\kappa(\theta)}\,\nabla\log\theta\|_2^2\,d\tau
= \int_s^{\lfloor t\rfloor}\|\sqrt{\kappa(\theta)}\,\nabla\log\theta\|_2^2\,d\tau
 \le C.
 \end{aligned}
\]
Thus $t-s\le1+C/\epsilon$. At $s$ the lower bound is
\[
 \theta(x,s)\ge
 \begin{cases}
  l_1,&s\ge1\quad\text{by the cutoff estimate on }[s-1,s],\\
  \underline\theta_0,&s=0.
 \end{cases}
\]
Apply \eqref{wh3:temperature-comparison} on $[s,t]$ with initial comparison value
$\min\{l_1,\underline\theta_0\}$. Using $a_0\le C(1+D(r))$ from \eqref{fp117} and \eqref{extlow:struct}, we obtain
\[
 \theta(x,t)^{-1}
 \le\max\{l_1^{-1},\underline\theta_0^{-1}\}+C(1+D(r))(t-s)
 \le C(1+D(r)).
\]
Under \textup{(D)}, or \textup{(C)} with $2(1-d)/(2+d)\le\alpha\le 1-d$, applying \eqref{wh3:temperature-comparison} on $[0,T]$ with $0<T\le T_0$ and $\sup_{0\le z\le r}z/(\mu+\lambda(z))\le C(1+D(r))$ gives
\begin{equation}\label{wh3:provisional-temperature-floor}
 \inf_{\mathbb R^2\times[0,T]}\theta\ge
 \begin{cases}
 c/(1+D(r)), & \text{for \textup{(A)}, \textup{(B)}, or \textup{(C)} with $\alpha>1-d$},\\
 c_{T_0}/(1+D(r)), & \text{for \textup{(D)}, or \textup{(C)} with $2(1-d)/(2+d)\le\alpha\le 1-d$},
 \end{cases}
\end{equation}
where both $c$ and $c_{T_0}$ are independent of $r$.
\end{proof}

Fix $L_0\ge 4l$, $\ell>\max\{6,1+\alpha,1+2/\alpha\}$, and $\sigma\in\bigl(\frac{1}{2}+\frac{1}{2+\alpha},1\bigr)$ sufficiently close to $1$ so that $3/(2+\alpha)<\sigma$ when $\alpha>1$ and $0<3-(2+\alpha)\sigma<1$ when $0<\alpha\le 1$. Set
\[
 S_0=1+(1+H)^\sigma+\int_{\mathbb R^2}\bigl(h(\div u)^2+\mu\omega^2\bigr)\,dx.
\]
We next estimate the relative pressure in terms of $H$.

\begin{lemma}\label{plane:relative-pressure}
There exists a positive constant $C$ such that, for $i=2$, and also for $i=3$ when $\alpha>1$,
\begin{equation}\label{wh3:relative-pressure-i}
 \int_{\mathbb R^2}\frac{|P-P_\infty|^i}{h^{i-1}}\,dx\le C(1+H)^{i/(2+\alpha)}.
\end{equation}
Moreover, when $0<\alpha\le 1$, we have
\begin{equation}\label{wh3:relative-pressure-sublinear}
 \int_{\mathbb R^2}\frac{|P-P_\infty|^3}{h^2}\,dx
 \le C(1+H)^\sigma
 \left(1+\|\theta{\bf1}_{\{\theta>L_0\}}\|_{\frac{2[3-(2+\alpha)\sigma]}{1-\sigma}}^{3-(2+\alpha)\sigma}\right).
\end{equation}
\end{lemma}
\begin{proof}
On $\{\theta\le L_0\}$, by \eqref{wh3:theta-cool-ltwo}, \eqref{wh3:rho-moment-2} and \eqref{wh3:theta-hot-measure}, we have, for $i\ge 2$,
\begin{equation}\label{wh3:relative-pressure-cool-lp}
 \|P-P_\infty\|_{L^i(\{\theta\le L_0\})}
 \le L_0\|\rho-\rho_\infty\|_i+\rho_\infty\|\theta-\theta_\infty\|_{L^i(\{\theta\le L_0\})}
 \le C_i.
\end{equation}
On $\{\theta>L_0\}$, by \eqref{wh3:cutoff-set-measure} and $h\ge 2\mu$, we have $\int_{\{\theta>L_0\}}P_\infty^i/h^{i-1}\,dx\le C|\{\theta>L_0\}|\le C$.
For $i<2+\alpha$, by \eqref{extlow:struct} and H\"older's inequality, we obtain
\begin{equation}\label{s5:P-hot-i}
\begin{aligned}
 \int_{\{\theta>L_0\}}\frac{P^i}{h^{i-1}}\,dx
 &\le C\left(\int_{\{\theta>L_0\}}\rho\theta^{2+\alpha}\,dx\right)^{i/(2+\alpha)}\\
 &\quad\times\left(\int_{\{\theta>L_0\}}\rho D(\rho)^{(i-1)(2+\alpha)/(2+\alpha-i)}\,dx\right)^{1-i/(2+\alpha)}\\
 &\le C H^{i/(2+\alpha)},
\end{aligned}
\end{equation}
where the last inequality uses \eqref{wh3:rho-moment-2}, \eqref{wh3:potential-hot} and $|\{\theta>L_0\}|\le C$. Combining \eqref{wh3:relative-pressure-cool-lp} with \eqref{s5:P-hot-i} yields \eqref{wh3:relative-pressure-i}.

For $0<\alpha\le 1$, by \eqref{extlow:struct} and H\"older's inequality on $\{\theta>L_0\}$, we obtain
\begin{equation}\label{s5:P-hot-sub}
\begin{aligned}
 \int_{\{\theta>L_0\}}\frac{P^3}{h^2}\,dx
 &\le C\left(\int_{\{\theta>L_0\}}\rho\theta^{2+\alpha}\,dx\right)^\sigma
     \|\rho D(\rho)^{2/(1-\sigma)}\|_{L^2(\{\theta>L_0\})}^{1-\sigma}
\|\theta{\bf1}_{\{\theta>L_0\}}\|_{\frac{2[3-(2+\alpha)\sigma]}{1-\sigma}}^{3-(2+\alpha)\sigma}\\
 &\le C H^\sigma \|\theta{\bf1}_{\{\theta>L_0\}}\|_{\frac{2[3-(2+\alpha)\sigma]}{1-\sigma}}^{3-(2+\alpha)\sigma},
\end{aligned}
\end{equation}
where the last inequality uses \eqref{wh3:rho-moment-2} and \eqref{wh3:potential-hot}. Combining \eqref{wh3:relative-pressure-cool-lp} with \eqref{s5:P-hot-sub} and $|\{\theta>L_0\}|\le C$ gives \eqref{wh3:relative-pressure-sublinear}.

By \eqref{wh3:relative-pressure-i} with $i=2$, $h\div u=G+P$ and the div--curl estimate, we obtain
\[
 \|\nabla u\|_2^2+\int_{\mathbb R^2}\bigl(h(\div u)^2+\mu\omega^2\bigr)\,dx\le C\bigl[A_1^2+(1+H)^{2/(2+\alpha)}\bigr],
\]
and
\[
 A_1^2\le C\left(1+\int_{\mathbb R^2}\bigl(h(\div u)^2+\mu\omega^2\bigr)\,dx
            +\int_{\{\theta>L_0\}}\frac{P^2}{h}\,dx\right),
\]
so \eqref{exthigh:compare} holds with $1+A_1^2+H$ in place of $S$. By \eqref{wh3:rho-moment-1} and \eqref{wh3:thermal-dissipation}, we have
\[
 \int_0^T\left(\int_{\mathbb R^2}\bigl(h(\div u)^2+\mu\omega^2\bigr)\,dx
              +\int_{\{\theta>L_0\}}\frac{P^2}{h}\,dx\right)dt\le C.
\]
Set
\[
 g_*=1+A_1^2+(1+\|\theta-\theta_\infty\|_{2\ell})^{\frac{\ell+1}{\ell-1}}
           +(1+H)^{(4/(2+\alpha)-1)_+}.
\]
Take $L=L_0/4$ in \eqref{s5:theta-hot-tail} and \eqref{s5:H-hot-int}. Choose $\varepsilon$ and $\zeta_0$ with
\[
 \begin{gathered}
 0<\varepsilon\le\min\{1,\alpha\},\quad
 \frac{\ell+1}{\ell-1}\le 1+\alpha-\varepsilon,\quad
 \max\!\left\{\frac12,\sigma,\left(\frac{4}{2+\alpha}-1\right)_+\right\}<1-\zeta_0<1.
 \end{gathered}
\]
For $0<\zeta\le\zeta_0$, the two estimates give
\[
 \int_0^T\|\theta{\bf1}_{\{\theta>L_0\}}\|_{2\ell}^{1+\alpha-\varepsilon}\,dt\le C
\]
and
\[
 \int_0^T\left(c_v\int_{\{\theta>L_0\}}\rho
    \int_{\theta_\infty}^{\theta}\bigl(K(v)-K(\theta_\infty)\bigr)\,dv\,dx
 \right)^{1-\zeta}dt\le C\zeta^{-2}.
\]
On $\{\theta\le L_0\}$, \eqref{s5:theta-cool-lp}, \eqref{wh3:potential-cool} and \eqref{wh3:theta-hot-energy} give
\[
 \|\theta-\theta_\infty\|_{L^{2\ell}(\{\theta\le L_0\})}
 +c_v\int_{\{\theta\le L_0\}}\rho
    \int_{\theta_\infty}^{\theta}\bigl(K(v)-K(\theta_\infty)\bigr)\,dv\,dx
 \le C_{L_0}.
\]
The preceding mechanical estimates give, for any measurable $I\subset[0,T]$,
\[
 \int_I\left[A_1^2+\int_{\mathbb R^2}
 \bigl(h(\div u)^2+\mu\omega^2\bigr)\,dx\right]dt\le C(1+|I|).
\]
Splitting the temperature terms at $L_0$ and applying the two integral bounds above yields
\[
 \int_I\bigl[g_*+S_0+(1+A_1^2+H)^{1-\zeta}\bigr]\,dt
 \le C\bigl(|I|+\zeta^{-2}\bigr).
\]
For $J=\{1+A_1^2+H>S_*\}\cup\{g_*>g_0\}$, the estimate with $\zeta=\zeta_0$ gives
\[
 \begin{aligned}
 \min\{S_*^{1-\zeta_0},g_0\}|J|
 &\le\int_J\bigl[g_*+(1+A_1^2+H)^{1-\zeta_0}\bigr]\,dt\\
 &\le C\bigl(|J|+\zeta_0^{-2}\bigr).
 \end{aligned}
\]
Choose $S_*,g_0$ so the minimum exceeds $2C$. Then $|J|\le C$, and
\[
 \int_I\bigl[g_*+\sqrt{1+A_1^2+H}\,\bigr]\,dt
 +\zeta^2\int_I(1+A_1^2+H)^{1-\zeta}\,dt\le C
 \quad(|I|\le1)
\]
and
\[
 |J|+\int_J(g_*+S_0)\,dt
 +\zeta^2\int_J(1+A_1^2+H)^{1-\zeta}\,dt\le C,
\]
which are the estimates in \eqref{exthigh:tails}.
\end{proof}

On any subinterval $[0,t]\subset[0,T]$ where $\rho\le r$ with $r\ge\max\{r_1,r_*\}$, set $R=2+r$,
$D_R=1+D(r)$, $L_R=R\log(e+R)$ and $R_T=2+\sup_{0\le s\le T}\|\rho(s)\|_\infty$.
By \eqref{wh3:provisional-temperature-floor} and the viscosity lower bound,
\[
 \Xi=\max\{1,\sup\rho/(h\kappa(\theta))\}\le C D_R^{1+\alpha},\qquad
 B=\int_{\mathbb R^2}\rho\kappa(\theta)|\dot\theta|^2\,dx .
\]

We next establish the differential estimate for $A_1^2$.

\begin{lemma}\label{wh3:quartic}
For every sufficiently small $\eta>0$, there exist positive constants $C$ and $C_\eta$ such that
\begin{equation}\label{wh3:flux-differential}
\begin{aligned}
 (A_1^2)'+\|\sqrt\rho\,\dot u\|_2^2
 &\le\eta \|\nabla K(\theta)\|_2^2+C R_T\|\nabla u\|_2^4+C A_1\sqrt{\Xi B}\\
 &\quad+C_\eta R_T\log(e+R_T)
   \bigl[A_1^2+(1+\|\theta-\theta_\infty\|_{2\ell})^{\frac{\ell+1}{\ell-1}}\bigr]A_1^2 .
\end{aligned}
\end{equation}
\end{lemma}
\begin{proof}
By \eqref{wh3:theta-hot-energy}, \eqref{wh3:theta-hot-measure}, \eqref{wh3:relative-pressure-cool-lp} and H\"older's inequality, we have
\[
 \|P-P_\infty\|_{L^p(\{\theta\le2\theta_\infty\})}+\|P-P_\infty\|_{L^1(\{\theta>2\theta_\infty\})}\le C,
\]
and
\[
\begin{aligned}
 \|P-P_\infty\|_{L^\ell(\{\theta>2\theta_\infty\})}
 &\le\|\rho\|_{L^{2\ell}(\{\theta>2\theta_\infty\})}\|\theta\|_{L^{2\ell}(\{\theta>2\theta_\infty\})}\\
 &{}+P_\infty|\{\theta>2\theta_\infty\}|^{1/\ell}\\
 &\le C(1+\|\theta-\theta_\infty\|_{2\ell}).
\end{aligned}
\]
By H\"older's inequality on $\{\theta>2\theta_\infty\}$ with $1/p=1-b+b/\ell$, we obtain
\begin{equation}\label{wh3:pressure-lp}
 \|P-P_\infty\|_p\le C(1+\|\theta-\theta_\infty\|_{2\ell})^{\frac{\ell(1-1/p)}{\ell-1}},\qquad 2\le p\le\ell.
\end{equation}
Applying \eqref{wh3:sobolev-support} to $(K(\theta)-K(2\theta_\infty))_+$ and \eqref{wh3:theta-cool-ltwo} on $\{\theta\le2\theta_\infty\}$, and using $K(\theta)=\frac{\kappa_0}{1+\alpha}\theta^{1+\alpha}$, we obtain
\begin{equation}\label{wh3:potential-hone}
 \|K(\theta)-K(\theta_\infty)\|_{H^1}\le C(1+\|\nabla K(\theta)\|_2),\qquad 1+\|\theta-\theta_\infty\|_{2\ell}\le C(1+\|\nabla K(\theta)\|_2)^{\frac{1}{1+\alpha}}.
\end{equation}

Writing
\begin{equation}\label{wh3:momentum-dot-split}
 \rho\dot u=(\sqrt\rho-\sqrt{\rho_\infty})\sqrt\rho\,\dot u+
                  \sqrt{\rho_\infty}\,\sqrt\rho\,\dot u,
\end{equation}
we deduce from \eqref{wh3:rho-moment-2} and H\"older's inequality that
\[
 \begin{aligned}
 \|(\sqrt\rho-\sqrt{\rho_\infty})\sqrt\rho\,\dot u\|_{3/2}
 &\le\|\sqrt\rho-\sqrt{\rho_\infty}\|_6\|\sqrt\rho\,\dot u\|_2
 \le C\|\sqrt\rho\,\dot u\|_2,\\
 \|\sqrt{\rho_\infty}\,\sqrt\rho\,\dot u\|_2
 &=\sqrt{\rho_\infty}\,\|\sqrt\rho\,\dot u\|_2.
 \end{aligned}
\]
For $f=G+P_\infty,\omega$, the elliptic estimates for $\Delta G=\div(\rho\dot u)$ and $\mu\Delta\omega=\nabla^\perp\cdot(\rho\dot u)$ together with the Sobolev inequality for $f^3$ give
\[
 \|f\|_6^3\le C\|f^2\nabla f\|_1
 \le C \|\sqrt\rho\,\dot u\|_2\bigl(\|f\|_6^2+\|f\|_4^2\bigr).
\]
By H\"older's inequality on $\{\rho\le 2\rho_\infty\}$ and $\{\rho>2\rho_\infty\}$, we get
\[
 \|f\|_4\le C(A_1+\|f\|_6),\qquad \|f\|_6\le C(A_1+\|\sqrt\rho\,\dot u\|_2).
\]
Applying H\"older's inequality to $(G+P_\infty)/\sqrt h$ and $\omega$ gives \eqref{2fp003} with $P-P_\infty$ in place of $P$. By \eqref{visc:envelope} and the $L^2$ elliptic estimate, we also have
\begin{equation}\label{wh3:flux-gradient-lfour}
 \|\nabla G\|_2+\|\nabla\omega\|_2\le C R_T^{1/2}\|\sqrt\rho\,\dot u\|_2,\qquad
 \|G+P_\infty\|_2\le C R_T^{B_h/2}A_1.
\end{equation}
Taking $p=2+8\log(e+R_T)$ in \eqref{core:weighted-flux} and using \eqref{gn11} and \eqref{wh3:flux-gradient-lfour} yields
\begin{equation}\label{wh3:flux-square-ltwo}
 \|(G+P_\infty)^2/h\|_2\le C\sqrt{R_T\log(e+R_T)}\,\|(G+P_\infty)/\sqrt h\|_2 \|\sqrt\rho\,\dot u\|_2.
\end{equation}
Applying the div--curl estimate to $h\div u=G+P$ and using \eqref{wh3:pressure-lp} and \eqref{wh3:flux-square-ltwo} as in \eqref{kvu:q4}, we obtain
\begin{equation}\label{wh3:quartic-bound}
 \begin{aligned}
 &\|\nabla u\|_4^4+\|\sqrt h\,\div u\|_4^4+
 \|(P-P_\infty)\div u\|_2^2+\|P-P_\infty\|_4^4\\
 &\le C R_T\log(e+R_T)
   \bigl[A_1^2+(1+\|\theta-\theta_\infty\|_{2\ell})^{\frac{\ell+1}{\ell-1}}\bigr]
   (e+A_1^2+\|\sqrt\rho\,\dot u\|_2^2+\|\nabla K(\theta)\|_2^2).
 \end{aligned}
\end{equation}

Multiplying $\rho\dot u=\nabla G+\mu\nabla^\perp\omega$ by $\dot u$ and integrating by parts over $\mathbb{R}^2$ gives
\begin{equation}\label{wh3:momentum-testing}
 \|\sqrt\rho\,\dot u\|_2^2=-\int_{\mathbb R^2}\bigl((G+P_\infty)\div\dot u+\mu\omega\curl\dot u\bigr)\,dx.
\end{equation}
The continuity equation and $P=\rho\theta$ give $\dot P=\rho\dot\theta-P\div u$, hence
\begin{equation}\label{wh3:pressure-transport}
 \dot P+(P-P_\infty)\div u=\rho\dot\theta-P_\infty\div u.
\end{equation}
Substituting $h\div u=G+P$ and \eqref{wh3:pressure-transport} into the derivation of \eqref{fp43}, the pressure term contributes
\[
 \begin{aligned}
 &-2\int_{\mathbb R^2}\frac{G+P_\infty}{h}\bigl(\dot P+(P-P_\infty)\div u\bigr)\,dx\\
 &=-2\int_{\mathbb R^2}\frac{G+P_\infty}{h}\rho\dot\theta\,dx
       +2P_\infty\int_{\mathbb R^2}\frac{G+P_\infty}{h}\div u\,dx.
 \end{aligned}
\]
and hence
\begin{equation}\label{wh3:flux-equality}
 (A_1^2)'+2\|\sqrt\rho\,\dot u\|_2^2=\sum_{i=1}^5 I_i(G+P_\infty,P-P_\infty)
                   +2P_\infty\int_{\mathbb R^2}\frac{G+P_\infty}{h}\div u\,dx,
\end{equation}
where $I_i(G+P_\infty,P-P_\infty)$ denotes the $i$th integral in \eqref{fp43} with $G,P$ replaced by $G+P_\infty,P-P_\infty$. By \eqref{wh3:pressure-lp}, Cauchy--Schwarz and the div--curl estimate, we have
\begin{equation}\label{s5:flux-bg}
 2P_\infty\left|\int_{\mathbb R^2}\frac{G+P_\infty}{h}\div u\,dx\right|
 \le C A_1 \|\nabla u\|_2\le C\bigl(A_1^2+\|\nabla u\|_2^2\bigr)
\end{equation}
and
\begin{equation}\label{s5:flux-gradu}
 \|\nabla u\|_2^2\le C\bigl(A_1^2+\|P-P_\infty\|_2^2\bigr)
 \le C\bigl[A_1^2+(1+\|\theta-\theta_\infty\|_{2\ell})^{\frac{\ell}{\ell-1}}\bigr].
\end{equation}
Estimating $\sum_{i=1}^5 I_i(G+P_\infty,P-P_\infty)$ as in \eqref{kvu:xprime} and substituting \eqref{s5:flux-bg}--\eqref{s5:flux-gradu} into \eqref{wh3:flux-equality} gives \eqref{wh3:flux-differential}.
\end{proof}

Set $v=2\nabla\Delta^{-1}\det\nabla u$ and $w=\dot u+v$. Using \eqref{wh3:quartic-bound}, we derive the following estimate for $w$.

\begin{lemma}\label{wh3:acceleration}
There exist positive constants $c$ and $C$ such that
\begin{equation}\label{wh3:acceleration-differential}
\begin{aligned}
 &\frac{d}{dt}\int_{\mathbb R^2}\rho|w|^2\,dx+c\int_{\mathbb R^2}\bigl(h|\div w|^2+\mu|\curl w|^2\bigr)\,dx\\
 &\quad\le C\Xi B+C R_T\log(e+R_T)
     \bigl[A_1^2+(1+\|\theta-\theta_\infty\|_{2\ell})^{\frac{\ell+1}{\ell-1}}\bigr]\\
 &\qquad\times\left(e+A_1^2+\|\sqrt\rho\,\dot u\|_2^2+\|\nabla K(\theta)\|_2^2+\int_{\mathbb R^2}\rho|w|^2\,dx\right).
\end{aligned}
\end{equation}
\end{lemma}
\begin{proof}
Since $\div v=2\det\nabla u$ and $\curl v=0$, the Hardy-space Jacobian estimate \cite{CLMS} and the $L^2$ elliptic estimate give
\begin{equation}\label{wh3:corrected-velocity-estimate}
 \|v\|_2\le C \|\nabla u\|_2^2,\qquad \|\nabla v\|_2\le C\|\nabla u\|_4^2.
\end{equation}
Because every $L^2$ harmonic field on $\mathbb{R}^2$ vanishes, \eqref{kv:qmaterial} yields
\[
 \|\dot v\|_2\le C \|\nabla u\|_2\bigl(\|\nabla w\|_2+\|\nabla u\|_4^2\bigr).
\]
Using \eqref{wh3:pressure-transport} and $\div w=\dot{(\div u)}+(\div u)^2$, we obtain
\[
 \dot G+(\div u)(G+P_\infty)
 =h\div w-\rho\lambda'(\div u)^2-\rho\dot\theta+P_\infty\div u.
\]
By Cauchy--Schwarz, Young's inequality and $\int_{\mathbb R^2}\frac{\rho^2|\dot\theta|^2}{h}\,dx\le\Xi B$, we obtain
\begin{equation}\label{wh3:acceleration-constant-pressure}
 P_\infty\left|\int_{\mathbb R^2}(\div u)\div w\,dx\right|\le\eta \int_{\mathbb R^2}\bigl(h|\div w|^2+\mu|\curl w|^2\bigr)\,dx+C_\eta \|\nabla u\|_2^2
\end{equation}
and
\begin{equation}\label{wh3:acceleration-temperature-coupling}
 \left|\int_{\mathbb R^2}\rho\dot\theta\div w\,dx\right|\le\eta \int_{\mathbb R^2}\bigl(h|\div w|^2+\mu|\curl w|^2\bigr)\,dx+C_\eta\Xi B.
\end{equation}
Replacing $G,P$ by $G+P_\infty,P-P_\infty$ in \eqref{kv:corrected-raw}, applying \eqref{wh3:quartic-bound} and \eqref{wh3:acceleration-constant-pressure}--\eqref{wh3:acceleration-temperature-coupling}, and taking $\eta$ small gives \eqref{wh3:acceleration-differential}.
\end{proof}

We next derive the differential estimate for $H$.

\begin{lemma}\label{est:relative-heat}
There exists a positive constant $C$ such that
\begin{equation}\label{wh3:heat-potential-test}
 H'+\frac{3}{4}\|\nabla K(\theta)\|_2^2\le\frac{1}{4}\|\sqrt\rho\,\dot u\|_2^2+C A_1^4\log(e+A_1^2).
\end{equation}
\end{lemma}
\begin{proof}
By \eqref{gn11} and \eqref{wh3:potential-hone}, we have
\[
 \|K(\theta)-K(\theta_\infty)\|_j\le C\sqrt{j}\,(1+\|\nabla K(\theta)\|_2),\qquad 2\le j<\infty.
\]
Taking $b=\ell+\log(e+A_1^2)$ and $p=\frac{2b}{b-1}$, we deduce from \eqref{2fp003}, \eqref{wh3:pressure-lp} and H\"older's and Young's inequalities that, for $0<\eta\le1$,
\begin{equation}\label{wh3:heat-Q}
 \begin{aligned}
 &\int_{\mathbb R^2}|K(\theta)-K(\theta_\infty)|\bigl(2\mu|D(u)|^2+\lambda(\div u)^2+|(P-P_\infty)\div u|\bigr)\,dx\\
 &\le C\sqrt b\,(1+\|\nabla K(\theta)\|_2)
 \bigl[A_1^2+(1+\|\theta-\theta_\infty\|_{2\ell})^{\frac{\ell+1}{\ell-1}}
       +A_1^{2-3/b}\|\sqrt\rho\,\dot u\|_2^{3/b}\bigr]\\
 &\le\eta\bigl(\|\sqrt\rho\,\dot u\|_2^2+\|\nabla K(\theta)\|_2^2\bigr)+C_\eta A_1^4\log(e+A_1^2).
 \end{aligned}
\end{equation}

Multiplying $\eqref{fns}_3$ by $K(\theta)-K(\theta_\infty)$ and integrating over $\mathbb{R}^2$ gives
\begin{equation}\label{s5:H-id}
 H'+\|\nabla K(\theta)\|_2^2=\int_{\mathbb R^2}\bigl(K(\theta)-K(\theta_\infty)\bigr)
          \bigl(2\mu|D(u)|^2+\lambda(\div u)^2-P\div u\bigr)\,dx.
\end{equation}
Integrating by parts on the background pressure term and using Young's inequality, we get
\begin{equation}\label{s5:H-bg}
 -P_\infty\int_{\mathbb R^2}\bigl(K(\theta)-K(\theta_\infty)\bigr)\div u\,dx
 =P_\infty\int_{\mathbb R^2} u\cdot\nabla K(\theta)\,dx
 \le\eta \|\nabla K(\theta)\|_2^2+C_\eta\|u\|_2^2.
\end{equation}
By \eqref{wh3:primitives}, \eqref{wh3:anchored}, \eqref{wh3:pressure-lp}--\eqref{wh3:potential-hone} and the div--curl estimate, since $\ell/[(\ell-1)(1+\alpha)]<2$, we have
\begin{equation}\label{s5:H-u2}
 \|u\|_2^2\le C\bigl(1+\|\nabla u\|_2^2\bigr)
 \le C A_1^2+C(1+\|\nabla K(\theta)\|_2)^{\frac{\ell}{(\ell-1)(1+\alpha)}}
 \le C A_1^2+\epsilon \|\nabla K(\theta)\|_2^2+C_\epsilon.
\end{equation}
Substituting \eqref{wh3:heat-Q} and \eqref{s5:H-bg}--\eqref{s5:H-u2} into \eqref{s5:H-id} and choosing $\eta,\epsilon$ sufficiently small yields \eqref{wh3:heat-potential-test}.

Multiplying $\eqref{fns}_3$ by $\kappa(\theta)\dot\theta$ and integrating over $\mathbb{R}^2$ also yields
\begin{equation}\label{wh3:heat-material-test}
 \begin{aligned}
 c_v B-\int_{\mathbb R^2}\Delta K(\theta)\,\kappa(\theta)\dot\theta\,dx
 &=\int_{\mathbb R^2}\bigl(2\mu|D(u)|^2+\lambda(\div u)^2-P\div u\bigr)\kappa(\theta)\dot\theta\,dx .
 \end{aligned}
\end{equation}
Subtracting the far-field pressure gives
\[
 \|P\div u\|_2\le\|(P-P_\infty)\div u\|_2+P_\infty\|\nabla u\|_2.
\]
The first norm on the right is controlled by \eqref{wh3:quartic-bound}. Since $\rho/h\le C D_R$ and $\kappa(\theta)=\kappa_0\theta^\alpha$, we have
\[
 \rho\theta^2\kappa(\theta)
 \le C D_R h\bigl(1+|K(\theta)-K(\theta_\infty)|^2\bigr).
\]
Young's inequality for $P(\div u)\kappa(\theta)\dot\theta$ yields
\begin{equation}\label{wh3:heat-pressure}
 |P(\div u)\kappa(\theta)\dot\theta|
 \le\eta\rho\kappa(\theta)|\dot\theta|^2+C_\eta D_R\bigl(1+|K(\theta)-K(\theta_\infty)|^2\bigr)h|\div u|^2.
\end{equation}
\end{proof}

We next establish the upper bound for the density and the positive lower bound for the temperature.

\begin{lemma}\label{wh3:permanentlow}
For any $T_0\in(0,\infty)$, there exist positive constants $C_{T_0}$ and $c_{T_0}$ such that, for all $0<T\le T_0$,
\[
 \sup_{0\le t\le T}\|\rho(t)\|_\infty\le C_{T_0},\qquad
 \inf_{\mathbb R^2\times[0,T]}\theta\ge c_{T_0}.
\]
Moreover, under \textup{(A)}, \textup{(B)}, or \textup{(C)} with $\alpha>1-d$, there exist positive constants $C$ and $\theta_*$, independent of $T$, such that
\begin{equation}\label{wh3:theta-positive}
 \sup_{0\le t\le T}\|\rho(t)\|_\infty\le C,\qquad
 \inf_{\mathbb R^2\times[0,T]}\theta\ge \theta_*.
\end{equation}
\end{lemma}
\begin{proof}
Take $r_0>2\max\{\rho_\infty,\|\rho_0\|_\infty,r_1,r_*\}$ and set $r_j=2^jr_0$. As in Lemma~\ref{records:sum}, we use the indices $j$ for which $D(r_j)\le(\log r_j)^{\zeta_*}$. The constants have the time dependence specified in Lemma~\ref{plane:low-transfer}.
For $w=\dot u+2\nabla\Delta^{-1}\det\nabla u$, it follows from \eqref{wh3:acceleration-constant-pressure}--\eqref{wh3:acceleration-temperature-coupling} that
\[
 \begin{aligned}
 P_\infty\left|\int_{\mathbb R^2}(\div u)\div w\,dx\right|
 &\le\varepsilon \int_{\mathbb R^2}\bigl(h|\div w|^2+\mu|\curl w|^2\bigr)\,dx+C_\varepsilon\|\nabla u\|_2^2,\\
 \left|\int_{\mathbb R^2}\rho\dot\theta\div w\,dx\right|
 &\le\varepsilon \int_{\mathbb R^2}\bigl(h|\div w|^2+\mu|\curl w|^2\bigr)\,dx+C_\varepsilon\Xi B.
 \end{aligned}
\]
Using $\div w=\dot{(\div u)}+(\div u)^2$ together with \eqref{wh3:quartic-bound}, \eqref{wh3:corrected-velocity-estimate} and \eqref{wh3:heat-material-test}--\eqref{wh3:heat-pressure}, we obtain
\[
\begin{aligned}
 &\|\Delta K(\theta)\|_2^2+\|\nabla\dot u\|_2^2+\int_{\mathbb R^2} h|\dot{(\div u)}|^2\,dx\\
 &\le C\left(R D_R^\alpha B+\int_{\mathbb R^2}\bigl(h|\div w|^2+\mu|\curl w|^2\bigr)\,dx\right)\\
 &{}+C L_R\bigl[A_1^2+(1+\|\theta-\theta_\infty\|_{2\ell})^{\frac{\ell+1}{\ell-1}}\bigr]
(e+A_1^2+\|\sqrt\rho\,\dot u\|_2^2+\|\nabla K(\theta)\|_2^2).
\end{aligned}
\]
Estimating the diffusion term in \eqref{wh3:heat-material-test} by $\|\Delta K(\theta)\|_2$ and bounding $\int_{\mathbb R^2}\bigl(K(\theta)-K(\theta_\infty)\bigr)^2(|\nabla u|^2+\lambda(\div u)^2)\,dx$ as in \eqref{fp517} gives, for $0<z\le 1$,
\[
 \begin{aligned}
 &\frac{d}{dt}\biggl[\frac{1}{2}\|\nabla K(\theta)\|_2^2-\int_{\mathbb R^2}\bigl(K(\theta)-K(\theta_\infty)\bigr)
                  \bigl(2\mu|D(u)|^2+\lambda(\div u)^2\bigr)\,dx\biggr]+cB\\
 &\le z\int_{\mathbb R^2}\bigl(h|\div w|^2+\mu|\curl w|^2\bigr)\,dx\\
 &{}+C z^{-n}D_R^n
    \bigl[A_1^2+(1+\|\theta-\theta_\infty\|_{2\ell})^{\frac{\ell+1}{\ell-1}}\bigr]
(e+A_1^2+\|\sqrt\rho\,\dot u\|_2^2+\|\nabla K(\theta)\|_2^2)\\
 & \quad \times\bigl[\log(e+A_1^2+\|\sqrt\rho\,\dot u\|_2^2+\|\nabla K(\theta)\|_2^2)+L_R\bigr].
 \end{aligned}
\]
By \eqref{wh3:acceleration-differential} and $\Xi\le C D_R^{1+\alpha}$,
the viscous term with coefficient $z$ and the thermal term $t_R\Xi B$ are absorbed
for $t_R=c_1D_R^{-(1+\alpha)}$, with $c_1>0$ small and $z\ll t_R$.

With $W_R=k_R\log Y_R+\log(1+A_1^2+H)$ as in Lemma~\ref{records:joint}, differentiating gives
\[
 W_R'=k_R\frac{Y_R'}{Y_R}+\frac{(A_1^2)'+H'}{1+A_1^2+H}.
\]
Using \eqref{wh3:flux-differential}, \eqref{wh3:heat-potential-test}, $A_1^2\le 1+A_1^2+H$ and Young's inequality, we have
\[
 \frac{A_1\sqrt{\Xi B}}{1+A_1^2+H}
 \le\epsilon\frac{Y_R}{1+A_1^2+H}+C_\epsilon\Xi\frac{B}{Y_R},
\]
and choosing $k_R$ as in Lemma~\ref{records:joint} yields
\[
 W_R'+M_R^{-1}\mathcal D_R\le M_R g_*(W_R+L_R).
\]
At $t=0$, we have $\|\sqrt{\rho_0}\,\dot u(0)\|_2\le\|g_1\|_2$, $\|\nabla K(\theta_0)\|_2\le\kappa(\|\theta_0\|_\infty)\|\nabla\theta_0\|_2$, and
\[
 H(0)\le \frac{c_v}{2}\|\rho_0\|_\infty
 \left(\max_{0\le s\le\max\{\|\theta_0\|_\infty,\theta_\infty\}}\kappa(s)\right)
 \|\theta_0-\theta_\infty\|_2^2\le C.
\]
Hence \eqref{exthigh:tails} and Lemma~\ref{records:joint} give \eqref{exthigh:window}.

Using \eqref{wh3:momentum-testing} and $\div w=\div\dot u+2\det\nabla u$, we obtain
\[
\begin{aligned}
 &\frac{1}{2}\frac{d}{dt}\int_{\mathbb R^2}\bigl(h(\div u)^2+\mu\omega^2\bigr)\,dx+\|\sqrt\rho\,\dot u\|_2^2\\
 &=\int_{\mathbb R^2}(P-P_\infty)\div w\,dx+2\int_{\mathbb R^2}(G+P_\infty)\det\nabla u\,dx\\
 & \quad -\frac{1}{2}\int_{\mathbb R^2}(h+\rho\lambda')(\div u)^3\,dx
       -\frac{\mu}{2}\int_{\mathbb R^2}\omega^2\div u\,dx.
\end{aligned}
\]
Since $\int_{\mathbb R^2}(P-P_\infty)^2/h\,dx\le C(1+H)^{2/(2+\alpha)}$ by Lemma~\ref{plane:relative-pressure} and $1+H\le 1+A_1^2+H\le C S_0^{1/\sigma}$, Cauchy--Schwarz yields
\[
\begin{aligned}
 \frac{\bigl|\int_{\mathbb R^2}(P-P_\infty)\div w\,dx\bigr|}{S_0}
 &\le \frac{C(1+H)^{1/(2+\alpha)}\sqrt{1+A_1^2+H}}{S_0}\\
 &{}\times\left(\frac{\int_{\mathbb R^2}(h|\div w|^2+\mu|\curl w|^2)\,dx}{Y_R}\right)^{1/2}
       \left(\frac{Y_R}{1+A_1^2+H}\right)^{1/2}\\
 &\le C S_0^{-\delta}\mathcal D_R,\qquad
 \delta=1-\frac{1+2/(2+\alpha)}{2\sigma}>0.
\end{aligned}
\]
Differentiate the expression for $S_0$ and divide by $S_0$. Use the preceding mechanical identity, the Jacobian estimate, \eqref{wh3:flux-square-ltwo} and Lemma~\ref{plane:relative-pressure}.
\begin{equation}\label{wh3:record-mechanical}
 (\log S_0)'+\frac{c \|\sqrt\rho\,\dot u\|_2^2}{S_0}
 \le\frac{2}{S_0}\left|\int_{\mathbb R^2}(P-P_\infty)\div w\,dx\right|+C b(L_R+\log S_0).
\end{equation}
Let $N$ be any index for which $\|\rho(t)\|_\infty$ reaches $r_N$ by time $T$. Summing \eqref{per:layer-endpoints} over $j\le N$ with $D(r_j)\le(\log r_j)^{\zeta_*}$, as in Lemma~\ref{records:sum}, gives
\[
 \begin{aligned}
 \sum_{\substack{j\le N\\D(r_j)\le(\log r_j)^{\zeta_*}}} d_j &\le C\left(1+\int_J b\,dt+\sum_{j\le N}R_j^{-\eta\delta/2}\right)\le C,\\
 \sum_{\substack{j\le N\\D(r_j)\le(\log r_j)^{\zeta_*}}} e_j &\le \int_J S_0\,dt\le C.
 \end{aligned}
\]

Applying $(1-\Delta)^{-1}$ to $\Delta G=\div(\rho\dot u)$ and using $(1-\Delta)^{-1}\Delta=(1-\Delta)^{-1}-I$, we have
\begin{equation}\label{wh3:flux-elliptic}
 \bigl(1-(1-\Delta)^{-1}\bigr)(G+P_\infty)=-(1-\Delta)^{-1}\div(\rho\dot u).
\end{equation}
In \eqref{wh3:momentum-dot-split}, H\"older's inequality gives
\[
 \begin{aligned}
 \|(\sqrt\rho-\sqrt{\rho_\infty})\sqrt\rho\,\dot u\|_{\frac{2p}{p+2}}
 &\le\|\sqrt\rho-\sqrt{\rho_\infty}\|_p\|\sqrt\rho\,\dot u\|_2,\\
 \|\sqrt{\rho_\infty}\sqrt\rho\,\dot u\|_2
 &=\sqrt{\rho_\infty}\|\sqrt\rho\,\dot u\|_2.
 \end{aligned}
\]
By \eqref{wh3:bessel-kernel}, \eqref{wh3:flux-elliptic} and the Hardy--Littlewood--Sobolev and Young inequalities, we obtain
\begin{equation}\label{wh3:flux-high-interim}
 \|\bigl(1-(1-\Delta)^{-1}\bigr)(G+P_\infty)\|_p
 \le C_p\bigl(1+\|\sqrt\rho-\sqrt{\rho_\infty}\|_p\bigr)\|\sqrt\rho\,\dot u\|_2.
\end{equation}
For $p>2$, \eqref{wh3:primitives} and \eqref{wh3:rho-moment-2} give
\[
 \int_{\mathbb R^2}|\sqrt\rho-\sqrt{\rho_\infty}|^p\,dx
 \le C_p\int_{\{\rho\le 2\rho_\infty\}}|\sqrt\rho-\sqrt{\rho_\infty}|^2\,dx
     +C_p\int_{\{\rho>2\rho_\infty\}}\rho^{p/2}\,dx\le C_p,
\]
Combining this with \eqref{wh3:flux-high-interim} gives
\begin{equation}\label{wh3:high-frequency}
 \|\bigl(1-(1-\Delta)^{-1}\bigr)(G+P_\infty)\|_p\le C_p\|\sqrt\rho\,\dot u\|_2\quad(p>2).
\end{equation}

For the low-frequency part, we write
\begin{equation}\label{wh3:flux-low-split}
 G+P_\infty
 =\left(\sqrt{2\mu+\lambda}-\sqrt{2\mu+\lambda(\rho_\infty)}\right)
   \frac{G+P_\infty}{\sqrt{2\mu+\lambda}}
 +\sqrt{2\mu+\lambda(\rho_\infty)}\,\frac{G+P_\infty}{\sqrt{2\mu+\lambda}}.
\end{equation}
Since $\|\sqrt{h(\rho)}-\sqrt{h(\rho_\infty)}\|_4\le C$ by \eqref{visc:regularity}, \eqref{visc:envelope} and \eqref{wh3:rho-moment-2}, Young's convolution inequality and H\"older's inequality give
\begin{equation}\label{wh3:low-frequency}
\begin{aligned}
 \|(1-\Delta)^{-1}(G+P_\infty)\|_\infty
 &\le\Bigl(\|G_B\|_4\|\sqrt{h(\rho)}-\sqrt{h(\rho_\infty)}\|_4
       +\sqrt{h(\rho_\infty)}\,\|G_B\|_2\Bigr)
\left\|\frac{G+P_\infty}{\sqrt h}\right\|_2\\
 &\le C A_1.
\end{aligned}
\end{equation}

By \eqref{gn11} and \eqref{wh3:anchored}, we have $\|\dot u\|_m\le C\sqrt{m}\,(\|\sqrt\rho\,\dot u\|_2+\|\nabla\dot u\|_2)$ for $m\ge 8$. Combining this with \eqref{visc:envelope} and the $L^p$ elliptic estimates for $\Delta G=\div(\rho\dot u)$, we obtain
\begin{equation}\label{wh3:flux-gradient-interpolation}
\begin{aligned}
 \|G+P_\infty\|_2&\le C R_T^{B_h/2} A_1,\qquad
 \|\nabla G\|_2\le C R_T^{1/2} \|\sqrt\rho\,\dot u\|_2,\\
 \|\nabla G\|_m&\le C m\|\rho\dot u\|_m
 \le C m^{3/2}R_T\bigl(\|\sqrt\rho\,\dot u\|_2+\|\nabla\dot u\|_2\bigr),\qquad m\ge 8.
\end{aligned}
\end{equation}

For $f=(1-(1-\Delta)^{-1})(G+P_\infty)$, $\|G_B\|_1=1$ gives
\[
 \|\nabla f\|_n\le 2\|\nabla G\|_n\qquad(n=2\ \text{or }n\ge 8).
\]
By \eqref{exthigh:tails}, \eqref{records:ball} (as in Lemma~\ref{records:flux}), \eqref{wh3:high-frequency}, \eqref{wh3:low-frequency} and \eqref{wh3:flux-gradient-interpolation}, for each $j\le N$ satisfying $D(r_j)\le(\log r_j)^{\zeta_*}$ and each interval $I\subset[0,t_j]$ we obtain
\begin{equation}\label{wh3:record-flux}
 \begin{aligned}
 \int_I\|(1-\Delta)^{-1}(G+P_\infty)\|_\infty\,dt
 &\le C\int_I A_1\,dt\le C(1+|I|),\\
 \int_{I\cap\{S_0\le R_j^\eta\}}\|f\|_\infty\,dt
 &\le C R_j^\xi(1+|I|), \qquad \xi<1,\\
 \int_{E_j}\|f\|_\infty\,dt
 &\le C L_j x_j^{w_j},\qquad
 \sum_{\substack{j\le N\\D(r_j)\le(\log r_j)^{\zeta_*}}} x_j\le C,\qquad 1-w_j=O(1/j).
 \end{aligned}
\end{equation}
Here $x_j\ge0$ and $0<w_j<1$ are the weights in \eqref{per:exceptional-flux}.

For each index $j$ in the sum in \eqref{wh3:record-flux}, the first time the density reaches $r_j$ is
\[
 t_j=\min\{t\in[0,T]:\|\rho(t)\|_\infty=r_j\}.
\]
By \eqref{tri:farfield}, $\rho(t_j,x)\to\rho_\infty<r_j/2$ as $|x|\to\infty$. Choose $y_j\in\mathbb R^2$ with $\rho(t_j,y_j)=r_j$ and solve
\[
 \dot X(t)=u(t,X(t)),\quad X(t_j)=y_j.
\]
If $j$ has a predecessor $j^-$ in this sum, then at the left endpoint of $I_j^0$ we have
\[
 \rho(t_{j^-},X(t_{j^-}))
 \le \|\rho(t_{j^-})\|_\infty
 =r_{j^-}\le r_j/2.
\]
If no such index exists, then $t_{j^-}=0$ and $\rho(0,X(0))\le\|\rho_0\|_\infty<r_0/2\le r_j/2$. By continuity,
\[
 \begin{gathered}
 s_j=\max\{t\in I_j^0:\rho(t,X(t))=r_j/2\},\\
 I_j=[s_j,t_j]\subset I_j^0,\quad
 \rho(t,X(t))\ge r_j/2\quad(t\in I_j).
 \end{gathered}
\]
Applying \eqref{wh3:provisional-temperature-floor} on $[0,t_j]$ with $r=r_j$, we obtain, for $t\in I_j$,
\[
 P(t,X(t))=\rho(t,X(t))\theta(t,X(t))\ge \frac{c r_j}{1+D(r_j)}.
\]
Integrating $\dot{\Lambda(\rho)}+P-P_\infty=-(G+P_\infty)$ along $X$ over $I_j=[s_j,t_j]$ and using \eqref{wh3:record-flux}, we obtain
\[
 \Lambda(r_j)-\Lambda(r_j/2)
 +\left(\frac{c r_j}{D_{R_j}}-P_\infty\right)|I_j|
 \le C(1+R_j^\xi)(1+|I_j|)+C L_j x_j^{w_j}.
\]
Since $\xi<1$ and $D_{R_j}\le 1+(\log r_j)^{\zeta_*}$ for these indices, $C(1+R_j^\xi)+P_\infty\le c r_j/D_{R_j}$ for all sufficiently large $j$ under consideration. Summing as in \eqref{order:crossing} gives
\[
 \sum_{\substack{j\le N\\D(r_j)\le(\log r_j)^{\zeta_*}}}\frac{\Lambda(r_j)-\Lambda(r_j/2)}{L_j}
 \le C\sum_{\substack{j\le N\\D(r_j)\le(\log r_j)^{\zeta_*}}}\left(\frac{R_j^\xi}{L_j}+x_j^{w_j}\right)\le C.
\]
By the divergence argument in Proposition~\ref{order:periodic-cap}, we may choose an index $j_*$ satisfying $D(r_{j_*})\le(\log r_{j_*})^{\zeta_*}$ such that
\[
 \sum_{\substack{j\le j_*\\D(r_j)\le(\log r_j)^{\zeta_*}}}
      \frac{\Lambda(r_j)-\Lambda(r_j/2)}{L_j}>C,
\]
where $C$ is the preceding bound, independent of $N$. If $r_{j_*}$ were reached by time $T$, that bound would hold with $N=j_*$. This contradiction gives
\[
 \sup_{0\le t\le T}\|\rho(t)\|_\infty<r_{j_*}.
\]
For \eqref{ext:powerfamily}, first choose the excluded level $R_*$ by the preceding whole-space argument with lower-bound constant $c_*/2$ and data bound $4M+4$, and then define $\varepsilon_*$ by \eqref{per:family-level}. We have
\[
 \begin{gathered}
 (1+r)^\varepsilon\le2,\quad h(r)\ge\frac{c_*r}{2D(r)},\\
 r_*\le r\le R_*,\quad 0\le\varepsilon\le\varepsilon_*.
 \end{gathered}
\]
Thus \eqref{order:crossing} up to the first attainment of $R_*$ gives $\sup_{0\le t\le T}\|\rho(t)\|_\infty<R_*$.
For \eqref{ext:gainfamily}, take $D\equiv1$ and use the preceding whole-space argument with lower-bound constant $c_*/2$ in \eqref{new:gain} and data bound $4M+4$ to choose $R_{*,C}$. The constants are uniform in $T$ because $\alpha>1-d$, and depend only on the fixed quantities and the data bound. Define $\varepsilon_{*,C}(M)$ by \eqref{per:gain-family-level} with this whole-space level. Then
\[
 \mu+\lambda(r)\ge\frac{c_*rL(r)^d}{(1+r)^\varepsilon}
 \ge\frac{c_*}{2}rL(r)^d
 \qquad(r_*\le r\le R_{*,C},\quad0\le\varepsilon\le\varepsilon_{*,C}(M)).
\]
Thus Lemma~\ref{plane:low-transfer} and the preceding flux estimates apply up to a first attainment of $R_{*,C}$. Summing \eqref{order:crossing} would contradict its choice, and hence $\sup_{0\le t\le T}\|\rho(t)\|_\infty<R_{*,C}$ uniformly over the stated coefficients and data.

Under \textup{(D)}, or \textup{(C)} with $2(1-d)/(2+d)\le\alpha\le 1-d$, taking
\[
 D(r)=
 \begin{cases}
 \Psi(J(r)),&\text{under \textup{(D)}},\\
 1,&\text{under \textup{(C)} with }\dfrac{2(1-d)}{2+d}\le\alpha\le 1-d,
 \end{cases}
\]
we have from \eqref{crit:profile} and Lemma~\ref{plane:low-transfer} that
\[
 \int^\infty\frac{dr}{(e+r)\log(e+r)D(r)}=\infty,\qquad
 \int_0^T\!\int_{\{\theta>l\}}\frac{P^2}{\mu+\lambda}\,dx\,dt\le C_{T_0}.
\]
Applying \eqref{order:crossing} on $[0,T_0]$ with $P\ge c_{T_0}r_j/D_{R_j}$ along $I_j$ gives $\sup_{0\le t\le T}\|\rho(t)\|_\infty\le C_{T_0}$.

Once $\sup_{0\le t\le T}\|\rho(t)\|_\infty$ is bounded, \eqref{exthigh:window} and \eqref{wh3:provisional-temperature-floor} give
\[
 \inf\theta\ge\frac{c}{D_R},\qquad
 \sup_t\bigl(A_1+\|\sqrt\rho\,\dot u\|_2+\|\nabla K(\theta)\|_2\bigr)+\sup_{|I|\le 1}\int_I\bigl(B+\|\nabla\dot u\|_2^2\bigr)\,dt\le C,
\]
and \eqref{wh3:anchored}, \eqref{wh3:potential-hone} and the monotonicity of $K$ yield
\[
 \begin{aligned}
 \|u\|_2&\le C\bigl(\|\sqrt\rho\,u\|_2+\|\nabla u\|_2\bigr),\\
 \|\theta-\theta_\infty\|_{H^1}
 &\le C\|K(\theta)-K(\theta_\infty)\|_{H^1}\le C(1+\|\nabla K(\theta)\|_2),
 \end{aligned}
\]
which gives \eqref{construction:common-cap} and \eqref{wh3:theta-positive}.
\end{proof}

\subsection{Higher-order estimates}

Based on Lemma~\ref{wh3:permanentlow}, we now derive the higher-order estimates up to the initial time.

\begin{lemma}\label{reg:plane-initial}
Let $q>2$ be as in Theorem~\ref{tri:whole}. For any finite $T>0$, there exists a positive constant $C_{T,q}$, independent of $\inf_{\mathbb{R}^2}\rho_0$, such that
\[
\begin{aligned}
 &\sup_{0\le t\le T}\bigl(\|\nabla\rho\|_2+\|\nabla\rho\|_q+
 \|\rho_t\|_2+\|\rho_t\|_q+\|u\|_{H^2}+\|\sqrt\rho\,u_t\|_2+\|\theta-\theta_\infty\|_{H^1}\bigr)\\
 &\quad+\int_0^T\bigl(\|\nabla^2u\|_q^2+\|u_t\|_{H^1}^2+
 \|\dot u\|_{H^1}^2+\|K(\theta)-K(\theta_\infty)\|_{H^2}^2+\|\theta-\theta_\infty\|_{H^2}^2\bigr)\,dt\le C_{T,q}.
\end{aligned}
\]
\end{lemma}
\begin{proof}
By Lemma~\ref{wh3:permanentlow}, \eqref{wh3:anchored}, \eqref{wh3:flux-gradient-lfour} and the div--curl estimate, we have, for $2\le p<\infty$,
\[
 \int_0^T\|\dot u\|_{H^1}^2\,dt
 +\sup_{0\le t\le T}\bigl(\|K(\theta)-K(\theta_\infty)\|_{H^1}+\|\theta-\theta_\infty\|_{H^1}+\|\nabla u\|_p\bigr)\le C_{T,p}.
\]
From $\eqref{fns}_3$, \eqref{wh3:quartic-bound} and Lemma~\ref{wh3:permanentlow}, we get
\[
 \begin{aligned}
 \|\Delta K(\theta)\|_2
 &\le c_v\|\rho\dot\theta\|_2+\|(P-P_\infty)\div u\|_2+P_\infty\|\nabla u\|_2\\
 &{}+\|2\mu|D(u)|^2+\lambda(\div u)^2\|_2
 \le C_T\bigl(1+\sqrt B\bigr),
 \end{aligned}
\]
which, together with the $H^2$ elliptic estimate and the lower bound for $\theta$, yields
\[
 \int_0^T\bigl(\|K(\theta)-K(\theta_\infty)\|_{H^2}^2+\|\theta-\theta_\infty\|_{H^2}^2\bigr)\,dt\le C_T.
\]

By elliptic estimates for $G+P_\infty$ and $\omega$ and $h\div u=G+P$, we have
\[
 \begin{aligned}
 \|\div u\|_\infty+\|\omega\|_\infty&\le C_T\bigl(1+\|\dot u\|_{H^1}+\|\theta-\theta_\infty\|_{H^2}\bigr),\\
 \|\nabla^2 u\|_q&\le C_T\bigl(1+\|\dot u\|_{H^1}+\|\theta-\theta_\infty\|_{H^2}\bigr)\bigl(1+\|\nabla\rho\|_q\bigr).
 \end{aligned}
\]
Differentiating $\eqref{fns}_1$ and applying Lemma~\ref{bkm} as in \eqref{fp107} gives
\[
 \begin{aligned}
 \frac{d}{dt}\|\nabla\rho\|_q
 &\le C\|\nabla u\|_\infty\|\nabla\rho\|_q+C\|\nabla^2 u\|_q,\\
 \frac{d}{dt}\log\bigl(e+\|\nabla\rho\|_q\bigr)
 &\le C_T\bigl(1+\|\dot u\|_{H^1}^2+\|\theta-\theta_\infty\|_{H^2}^2\bigr)\log\bigl(e+\|\nabla\rho\|_q\bigr).
 \end{aligned}
\]
Applying Gronwall's inequality yields
\[
 \sup_{0\le t\le T}\|\nabla\rho\|_q+\int_0^T\bigl(\|\nabla^2 u\|_q^2+\|\nabla u\|_\infty\bigr)\,dt\le C_{T,q}.
\]

By \eqref{th13}, we have
\[
\begin{aligned}
 \bigl(\theta_0-\lambda'(\rho_0)\div u_0\bigr)\nabla\rho_0&=g,\\
 g&=\sqrt{\rho_0}\,g_1+\mu\Delta u_0
    +(\mu+\lambda(\rho_0))\nabla\div u_0-\rho_0\nabla\theta_0\in L^2.
\end{aligned}
\]
Dividing by $\theta_0-\lambda'(\rho_0)\div u_0\ge\underline\theta_0/2$ on $\{|\lambda'(\rho_0)\div u_0|\le\underline\theta_0/2\}$ and applying Chebyshev's and H\"older's inequalities on $\{|\lambda'(\rho_0)\div u_0|>\underline\theta_0/2\}$, we obtain
\begin{equation}\label{s5:grad-rho0}
 \begin{aligned}
 |\{|\lambda'(\rho_0)\div u_0|>\underline\theta_0/2\}|&\le C\|\div u_0\|_2^2,\\
 \|\nabla\rho_0\|_2
 &\le 2\underline\theta_0^{-1}\|g\|_2+
 C\|\div u_0\|_2^{1-2/q}\|\nabla\rho_0\|_q.
 \end{aligned}
\end{equation}
For the approximating density $r_\delta=(1-\delta)\rho_0+\delta\rho_\infty$, \eqref{s5:grad-rho0} gives
\[
 \|\nabla r_\delta\|_2=(1-\delta)\|\nabla\rho_0\|_2\le C .
\]
Differentiate $h\div u=G+P$. The div--curl estimate and H\"older's inequality with $p=2q/(q-2)$ yield
\[
\begin{aligned}
 \|\nabla^2 u\|_2
 &\le C\bigl(\|\nabla G\|_2+\|\nabla\omega\|_2+\|\nabla\theta\|_2
        +\theta_\infty\|\nabla\rho\|_2\bigr)
+ C \bigl(\|\theta-\theta_\infty\|_p+\|\div u\|_p\bigr)\|\nabla\rho\|_q\\
 &\le C_T\bigl(1+\|\nabla\rho\|_2\bigr).
\end{aligned}
\]
By the differentiated continuity equation, we get
\[
 \frac{d}{dt}\log\bigl(1+\|\nabla\rho\|_2\bigr)\le C_T\bigl(1+\|\nabla u\|_\infty\bigr).
\]
Integrating over $[0,t]$ and using the estimate for $\|\nabla^2u\|_2$ gives
\begin{equation}\label{wh3:density-gradient-ltwo}
 \sup_{0\le t\le T}\|\nabla\rho(t)\|_2+\sup_{0\le t\le T}\|u(t)\|_{H^2}\le C_T.
\end{equation}
From $\rho_t=-u\cdot\nabla\rho-\rho\div u$ and $u_t=\dot u-u\cdot\nabla u$, we then obtain
\[
 \sup_{0\le t\le T}\bigl(\|\rho_t\|_{L^2\cap L^q}+\|\sqrt\rho\,u_t\|_2\bigr)
 +\|u_t\|_{L^2(0,T;H^1)}\le C_{T,q}.
\]
This completes the proof of Lemma \ref{reg:plane-initial}.
\end{proof}

Using Lemma~\ref{reg:plane-initial}, we next establish a short-time $L^\infty$ estimate for $\theta$.

\begin{lemma}\label{reg:plane-barrier}
There exists a positive constant $C$, independent of $T$ and $\inf_{\mathbb{R}^2}\rho_0$, such that
\[
 \sup_{0\le t\le\min\{T,1\}}\|\theta(t)\|_\infty\le C.
\]
\end{lemma}
\begin{proof}
On $[0,\min\{1,T\}]$, let $\psi$ solve
\[
 (-\Delta+\rho)\psi=2\mu|D(u)|^2+\lambda(\div u)^2.
\]
By \eqref{wh3:anchored}, we have
\[
 \|\psi\|_{H^1}^2\le C\bigl(\|\nabla\psi\|_2^2+\|\sqrt\rho\,\psi\|_2^2\bigr).
\]
Testing with $\psi_-=(-\psi)_+$ and integrating in $x$ gives
\[
 \int_{\mathbb R^2}\bigl(|\nabla\psi_-|^2+\rho\psi_-^2\bigr)\,dx
 =-\int_{\mathbb R^2}\bigl(2\mu|D(u)|^2+\lambda(\div u)^2\bigr)\psi_-\,dx
 \le0.
\]
Hence $\psi\ge0$ by \eqref{wh3:anchored}. For the $L^4$ elliptic estimate, write
\[
 (-\Delta+\rho_\infty)\psi
 =2\mu|D(u)|^2+\lambda(\div u)^2+(\rho_\infty-\rho)\psi.
\]
Together with Lemma~\ref{reg:plane-initial}, this gives
\[
 \|\psi\|_{W^{2,4}}\le C\bigl(\|\nabla u\|_8^2+
                         \|(\rho_\infty-\rho)\psi\|_4\bigr)\le C.
\]
Differentiating in time, we see from Lemma~\ref{reg:plane-initial} that $b=\partial_t\psi$ satisfies
\[
 (-\Delta+\rho)b=f,\qquad
 f=\partial_t\bigl(2\mu|D(u)|^2+\lambda(\div u)^2\bigr)-\rho_t\psi
      \in L^2_t L^{4/3}_x.
\]
By coercivity and the $L^{4/3}$ elliptic estimate for $-\Delta+\rho_\infty$, we get
\[
 \|b\|_{H^1}\le C\|f\|_{4/3},\qquad
 \|b\|_{W^{2,4/3}}\le C\bigl(\|f\|_{4/3}+\|\rho-\rho_\infty\|_2\|b\|_4\bigr)\le C\|f\|_{4/3},
\]
so $\|b\|_\infty\in L^2_t$ by Sobolev embedding.

As in Lemma~\ref{init:common-lemma}, we set
\[
\begin{aligned}
 \xi(0)&=1+\max\{\|K(\theta_0)\|_\infty,K(\theta_\infty)\},\\
 \xi'&=C(1+\xi)(1+\|b\|_\infty+\|\div u\|_\infty),\qquad
 \bar\theta=K^{-1}(\xi+\psi).
\end{aligned}
\]
Since $\|b\|_\infty+\|\div u\|_\infty\in L^1_t$, we have $\bar\theta\le C$. For $C$ sufficiently large, $\bar\theta$ satisfies
\begin{equation}\label{s5:theta-bar}
 c_v\rho \dot{\bar\theta}-\Delta K(\bar\theta)+\rho\bar\theta\div u-\bigl(2\mu|D(u)|^2+\lambda(\div u)^2\bigr)\ge 0.
\end{equation}
Since $\bar\theta\ge\theta_\infty+a$ for some $a>0$, we have $(\theta-\bar\theta)_+\le a^{-1}|\theta-\theta_\infty|^2$. Combining $\eqref{fns}_3$ with \eqref{s5:theta-bar}, multiplying by $\mathbf{1}_{\{\theta>\bar\theta\}}$ and integrating over $\mathbb{R}^2$, we obtain
\[
 c_v\frac{d}{dt}\int_{\mathbb R^2}\rho(\theta-\bar\theta)_+\,dx
 \le\|\div u\|_\infty\int_{\mathbb R^2}\rho(\theta-\bar\theta)_+\,dx.
\]
Since $(\theta-\bar\theta)_+(0)=0$, Gronwall's inequality yields $\theta\le\bar\theta\le C$.
\end{proof}

We next derive the higher-order estimates for $\theta$ and $u$ at positive times.

\begin{lemma}\label{wh3:temperaturehigh}
For any $0<\tau<T$, $2\le p<\infty$, and every interval $I\subset[\tau,T]$ with $|I|\le 1$, there exists a positive constant $C_{p,\tau,T}$ such that
\begin{equation}\label{wh3:temperature-regularity}
\begin{aligned}
 &\sup_{\tau\le t\le T}
 \left(\|\sqrt{\rho\kappa(\theta)}\,\dot\theta\|_2^2+\|\theta-\theta_\infty\|_{H^2}
            +\|\sqrt\rho\,\theta_t\|_2+\|u\|_{W^{1,p}}\right)\\
 &\quad+\int_I\left(\|\nabla(\kappa(\theta)\dot\theta)\|_2^2+
                \|\dot\theta\|_{H^1}^2+\|\dot u\|_{H^1}^2\right)dt\le C_{p,\tau,T} .
\end{aligned}
\end{equation}
Moreover, under \textup{(A)}, \textup{(B)}, or \textup{(C)} with $\alpha>1-d$, the constant depends only on $p$ and $\tau$ and is independent of $T$.
\end{lemma}
\begin{proof}
Write $v=\kappa(\theta)\dot\theta$. Since $\rho/\kappa(\theta)\ge c$ on $\{\rho\ge\rho_\infty/2,\ \theta\le 2\theta_\infty\}$, by \eqref{gn11}, \eqref{wh3:low-density-measure} and \eqref{wh3:theta-hot-measure}, we have
\[
\begin{aligned}
 \|v\|_2^2&\le CB+|\{\rho<\rho_\infty/2\}\cup\{\theta>2\theta_\infty\}|^{1/2}\|v\|_4^2\le CB+C\|v\|_2\|\nabla v\|_2,\\
 \|v\|_p&\le C\bigl(\sqrt B+\sqrt p\,\|\nabla v\|_2\bigr)\qquad (p\ge 5).
\end{aligned}
\]
Replacing $\theta$ by $\theta-\theta_\infty$ in the $L^4$ products of \eqref{fp914}, we estimate the constant-temperature terms by Young's inequality as
\[
 C\theta_\infty\bigl(\|\nabla\dot u\|_2+\|\nabla u\|_4^2\bigr)\|v\|_2
 \le\epsilon\|\nabla v\|_2^2+C_\epsilon\bigl(1+B+\|\nabla\dot u\|_2^2\bigr).
\]
In the elliptic estimate for $\Delta K(\theta)$, the background pressure contributes $P_\infty\|\div u\|_2\le C$, so \eqref{fp910}--\eqref{fp915} yield
\begin{equation}\label{s5:B-diff}
 \frac{c_v}{2}B'+\frac{1}{2}\|\nabla v\|_2^2
 \le C\bigl(1+B^2\log(e+B)+\|\nabla\dot u\|_2^2\bigr).
\end{equation}
For $t\ge\tau$, let $h=\min\{1,\tau\}$. By Lemma~\ref{wh3:permanentlow}, there exists $s\in[t-h,t-h/2]$ with $B(s)\le C/h$. Dividing \eqref{s5:B-diff} by $(e+B)\log(e+B)$ and integrating over $[s,t]$, we obtain
\[
 \log\log\bigl(e+B(t)\bigr)
 \le\log\log\bigl(e+B(s)\bigr)
   +C\int_s^t\bigl(1+B+\|\nabla\dot u\|_2^2\bigr)\,d\xi
 \le C_\tau.
\]
Integrating \eqref{s5:B-diff} in time over $I=[t_1,t_2]$ then yields
\[
 \int_I\|\nabla v\|_2^2\,dt
 \le C B(t_1)+C\int_I\bigl(1+B^2\log(e+B)+\|\nabla\dot u\|_2^2\bigr)\,dt
 \le C_\tau.
\]
Applying the $H^2$ elliptic estimate to $\Delta K(\theta)=c_v\rho\dot\theta+P\div u-\bigl(2\mu|D(u)|^2+\lambda(\div u)^2\bigr)$ and using the lower bound for $\theta$ together with $\theta_t=v/\kappa(\theta)-u\cdot\nabla\theta$ gives \eqref{wh3:temperature-regularity}.
\end{proof}

With \eqref{wh3:temperature-regularity}, the differentiated temperature equation gives $\nabla^3 K(\theta)\in L^2(\tau,T;L^2)$.

\begin{lemma}\label{reg:plane-higher}
Let $q>2$ be as in Theorem~\ref{tri:whole}. For every $0<\tau<T$, there exists a positive constant $C_{\tau,T,q}$ such that
\begin{equation}\label{wh3:temperature-regularity2}
\begin{aligned}
 &\sup_{\tau\le t\le T}\bigl(\|\rho-\rho_\infty\|_{H^1\cap W^{1,q}}
   +\|u\|_{H^2}+\|\theta-\theta_\infty\|_{H^2}
   +\|\sqrt\rho\,\theta_t\|_2\bigr)\\
 &\quad+\int_\tau^T\bigl(\|\nabla u\|_\infty+
     \|\nabla^2u\|_q^2+\|\nabla\theta_t\|_2^2\bigr)\,dt\le C_{\tau,T,q}.
\end{aligned}
\end{equation}
\end{lemma}
\begin{proof}
On $[\tau,T]$, Lemmas~\ref{reg:plane-initial} and~\ref{wh3:temperaturehigh} give
\[
 \nabla P=\theta\nabla\rho+\rho\nabla\theta,\qquad
 \|\nabla P\|_{L^2\cap L^q}\le C_{\tau,T,q}.
\]
Differentiating $\Delta K(\theta)=c_v\rho\dot\theta+P\div u-\bigl(2\mu|D(u)|^2+\lambda(\div u)^2\bigr)$ and using \eqref{fp1015} with $1/2=1/q+1/p$, we obtain
\[
 \|\nabla^3 K(\theta)\|_2\le C\|\nabla\Delta K(\theta)\|_2
 \le C_{\tau,T,q}\bigl(1+\|\dot\theta\|_{H^1}+\|\nabla^2 u\|_q\bigr),
\]
so $\theta-\theta_\infty\in L^2(\tau,T;H^3)$. Using $\nabla\theta_t=\nabla\dot\theta-(\nabla u)^{\mathsf T}\nabla\theta-u\cdot\nabla^2\theta$, we have
\[
 \|\nabla\theta_t\|_2
 \le\|\nabla\dot\theta\|_2+\|\nabla u\|_4\|\nabla\theta\|_4+\|u\|_\infty\|\nabla^2\theta\|_2.
\]
Integrating in $L^2(\tau,T)$ and using Lemmas~\ref{reg:plane-initial} and~\ref{wh3:temperaturehigh} yields \eqref{wh3:temperature-regularity2}.
\end{proof}

In the remainder of this section, we assume \textup{(A)}, \textup{(B)}, or \textup{(C)} with $\alpha>1-d$, and set
\[
 D_\rho=\bigl\langle \rho-\rho_\infty,\,
     (1-(1-\Delta)^{-1})(\rho-\rho_\infty)\bigr\rangle.
\]
\begin{lemma}\label{reg:plane-energy}
Let
\[
 \mathcal F=\int_{\mathbb R^2}\bigg[\frac{1}{2}\rho|u|^2
   +c_v\rho\left(\theta-\theta_\infty-\theta_\infty\log\frac{\theta}{\theta_\infty}\right)
   +\theta_\infty\left(\rho\log\frac{\rho}{\rho_\infty}-\rho+\rho_\infty\right)\bigg]\,dx.
\]
There exist positive constants $a$, $c$, and $C$ such that
$\mathcal E=\mathcal F+a\int_{\mathbb R^2}\rho u\cdot\nabla(1-\Delta)^{-1}(\rho-\rho_\infty)\,dx$ satisfies, for all $t\ge 1$,
\begin{equation}\label{extplane:compensated}
\begin{aligned}
& c \|(\rho-\rho_\infty,u,\theta-\theta_\infty)\|_2^2
 \le \mathcal E\le C\|(\rho-\rho_\infty,u,\theta-\theta_\infty)\|_2^2, \\
& \mathcal E'+c\bigl(D_\rho+\|\nabla u\|_2^2+\|\nabla\theta\|_2^2\bigr)
 \le 0.
\end{aligned}
\end{equation}
\end{lemma}
\begin{proof}
For $t\ge 1$, Lemmas~\ref{wh3:permanentlow} and~\ref{wh3:temperaturehigh} give
\[
\begin{aligned}
 &c\int_{\mathbb R^2}\bigl((\rho-\rho_\infty)^2+\rho|u|^2+\rho(\theta-\theta_\infty)^2\bigr)\,dx\\
 &\le\mathcal F
 \le C\int_{\mathbb R^2}\bigl((\rho-\rho_\infty)^2+\rho|u|^2+\rho(\theta-\theta_\infty)^2\bigr)\,dx.
\end{aligned}
\]
Since $|u|^2+(\theta-\theta_\infty)^2\le C(\rho-\rho_\infty)^2$ on $\{\rho<\rho_\infty/2\}$, \eqref{wh3:primitives} and Korn's inequality yield
\begin{equation}\label{s5:F-equiv}
\begin{aligned}
 c\|(\rho-\rho_\infty,u,\theta-\theta_\infty)\|_2^2
 &\le\mathcal F\le C\|(\rho-\rho_\infty,u,\theta-\theta_\infty)\|_2^2,\\
 \mathcal F'+c\bigl(\|\nabla u\|_2^2+\|\nabla\theta\|_2^2\bigr)
 &\le 0.
\end{aligned}
\end{equation}

Writing $\delta f=f(x)-f(y)$ and $d\mu=G_B(x-y)\,dx\,dy$, we have
\begin{equation}\label{s5:increment}
 \langle f,(1-(1-\Delta)^{-1})g\rangle=\frac{1}{2}\iint_{\mathbb R^2\times\mathbb R^2}\delta f\,\delta g\,d\mu,\qquad
 \delta P=\theta(x)\delta\rho+\rho(y)\delta\theta.
\end{equation}
By $\theta\ge\theta_*>0$, Young's inequality and $\langle\theta-\theta_\infty,(1-(1-\Delta)^{-1})(\theta-\theta_\infty)\rangle\le\|\nabla\theta\|_2^2$, we obtain
\begin{equation}\label{s5:P-coercive}
 \bigl\langle P-P_\infty,(1-(1-\Delta)^{-1})(\rho-\rho_\infty)\bigr\rangle
 \ge c D_\rho-C\|\nabla\theta\|_2^2.
\end{equation}
Let $I=\int_{\mathbb R^2}\rho u\cdot\nabla(1-\Delta)^{-1}(\rho-\rho_\infty)\,dx$ and
$\mathbb S=2\mu D(u)+\lambda(\div u)\mathbb{I}_2$.
By $\eqref{fns}_1$ and $\eqref{fns}_2$, we have
\[
 \rho_t=-\div(\rho u),\qquad
 (\rho u)_t=\div(\mathbb S-\rho u\otimes u)-\nabla P.
\]
Using these two equations in $I'$ and integrating by parts, we obtain
\[
 \begin{aligned}
 I'={}&-\bigl\langle P-P_\infty,(1-(1-\Delta)^{-1})(\rho-\rho_\infty)\bigr\rangle
      +\bigl\langle\div(\rho u),(1-\Delta)^{-1}\div(\rho u)\bigr\rangle \\
      &+\int_{\mathbb R^2}(\rho u\otimes u-\mathbb S):\nabla^2(1-\Delta)^{-1}(\rho-\rho_\infty)\,dx.
 \end{aligned}
\]
By the Fourier multiplier estimate and the estimates on $\rho$ and $u$,
\[
 \begin{aligned}
& \|\nabla^2(1-\Delta)^{-1}(\rho-\rho_\infty)\|_2 \le D_\rho^{1/2},\\
& \|\rho u\otimes u\|_2+\|\mathbb S\|_2 \le C(1+\|u\|_2)\|\nabla u\|_2.
 \end{aligned}
\]
By Young's inequality, the stress and transport integral in $I'$ is estimated by
$\delta D_\rho+C_\delta\|\nabla u\|_2^2$. By \eqref{s5:increment} and the Fourier multiplier estimate, we have
\begin{equation}\label{s5:I-conv}
 \begin{aligned}
 \bigl\langle\div(\rho u),(1-\Delta)^{-1}\div(\rho u)\bigr\rangle
 &\le\bigl\langle \rho u,(1-(1-\Delta)^{-1})(\rho u)\bigr\rangle\\
 &\le C\|\nabla u\|_2^2+C\|u\|_\infty^2 D_\rho.
 \end{aligned}
\end{equation}
Since $\|u\|_\infty\le C\|u\|_4^{1/2}\|\nabla u\|_4^{1/2}\le C\|\nabla u\|_2^{1/4}$ and $D_\rho\le C$, Young's inequality yields
\begin{equation}\label{s5:I-uinf}
 \|u\|_\infty^2 D_\rho\le\delta D_\rho+C_\delta\|\nabla u\|_2^2.
\end{equation}
Combining \eqref{s5:P-coercive} with \eqref{s5:I-conv}--\eqref{s5:I-uinf} and choosing $\delta>0$ sufficiently small, we arrive at
\[
 \begin{aligned}
 I'+c D_\rho&\le C\bigl(\|\nabla u\|_2^2+\|\nabla\theta\|_2^2\bigr),\\
 |I|&\le\|\rho u\|_2\|\nabla(1-\Delta)^{-1}(\rho-\rho_\infty)\|_2\le C\mathcal F.
 \end{aligned}
\]
By \eqref{s5:F-equiv} and the estimates for $I$ and $I'$, \eqref{extplane:compensated} follows for sufficiently small $a>0$.
\end{proof}

Using Lemma~\ref{reg:plane-energy}, we establish the $L^2$ decay of $(\rho-\rho_\infty,u,\theta-\theta_\infty)$.

\begin{lemma}\label{reg:plane-decay-ltwo}
Let $(\rho,u,\theta)$ be a global solution under \textup{(A)}, \textup{(B)}, or \textup{(C)} with $\alpha>1-d$ satisfying the estimates of Lemmas~\ref{wh3:permanentlow} and~\ref{wh3:temperaturehigh}. Then
\begin{equation}\label{wh3:decay-ltwo-limit}
 \lim_{t\to\infty}\bigl(\|\rho(t)-\rho_\infty\|_2+\|u(t)\|_2+
                         \|\theta(t)-\theta_\infty\|_2\bigr)=0 .
\end{equation}
\end{lemma}
\begin{proof}
Set
\[
 U=(\rho-\rho_\infty,v,y),\qquad
 v=\frac{\rho u}{\rho_\infty},\qquad
 y=\frac{\rho\bigl(c_v(\theta-\theta_\infty)+|u|^2/2\bigr)}{c_v\rho_\infty}.
\]
From $\eqref{fns}_1$, $\eqref{fns}_2$ and $\eqref{ffns}_3$, we have
\[
 \begin{pmatrix}
 \rho_t+\rho_\infty\div v\\
 v_t+(\theta_\infty/\rho_\infty)\nabla\rho+\nabla y\\
 y_t+(\theta_\infty/c_v)\div v
 \end{pmatrix}=\div(F_1+F_2),
\]
where
\[
\begin{aligned}
 |F_1|&\le C\rho\bigl(|u|^2+|\theta-\theta_\infty|\,|u|+|u|^3\bigr),\\
 |F_2|&\le C\bigl((1+|u|)|\nabla u|+|\nabla K(\theta)|\bigr).
\end{aligned}
\]
By \eqref{wh3:primitives} and Lemmas~\ref{wh3:permanentlow} and~\ref{wh3:temperaturehigh}, we have
\[
 \sup_{t\ge 1}\|F_1(t)\|_1\le C,\qquad
 \int_1^\infty\|F_2(t)\|_2^2\,dt\le C.
\]
For the homogeneous acoustic system, its Fourier equations give
\[
 \frac{d}{dt}\left(
  \frac{\theta_\infty}{\rho_\infty^2}|\widehat{\rho-\rho_\infty}|^2
  +|\widehat v|^2+\frac{c_v}{\theta_\infty}|\widehat y|^2
 \right)=0 .
\]
Thus the Fourier propagator is uniformly bounded. For $0<k\le1$, apply Duhamel's formula and Plancherel's identity with
\[
 P_k={\bf1}_{\{|\xi|\le k\}}(D),\quad \||\xi|\|_{L^2(B_k)}=Ck^2 .
\]
Cauchy--Schwarz in time then gives
\[
\begin{aligned}
 \|P_k U(t)\|_2
 &\le C\|P_k U(1)\|_2
     +C\int_1^t\Bigl(\||\xi|\|_{L^2(B_k)}\|F_1(s)\|_1+k\|F_2(s)\|_2\Bigr)\,ds\\
 &\le C\|P_k U(1)\|_2+C t k^2+C\sqrt t\,k.
\end{aligned}
\]
The differences from $u$ and $\theta-\theta_\infty$ satisfy
\[
 \begin{aligned}
 v-u&=\frac{(\rho-\rho_\infty)u}{\rho_\infty},\\
 y-(\theta-\theta_\infty)
 &=\frac{(\rho-\rho_\infty)(\theta-\theta_\infty)}{\rho_\infty}
   +\frac{\rho|u|^2}{2c_v\rho_\infty}.
 \end{aligned}
\]
By Plancherel's identity,
\[
 \|P_k(v-u,y-(\theta-\theta_\infty))\|_2\le C k\|(v-u,y-(\theta-\theta_\infty))\|_1\le C k.
\]
Choosing $k(t)=[(e+t)\log(e+t)]^{-1/2}$, we have $tk(t)^2+\sqrt t\,k(t)\to 0$ and $\|P_{k(t)}U(1)\|_2\to 0$, so
\[
 \lim_{t\to\infty}\|P_{k(t)}(\rho-\rho_\infty,u,\theta-\theta_\infty)\|_2=0.
\]
Since $|\xi|^2/(1+|\xi|^2)\ge k(t)^2/2$ and $|\xi|^2\ge k(t)^2/2$ for $|\xi|>k(t)$, it follows from \eqref{extplane:compensated} that for each $\epsilon>0$ there exists $t_\epsilon\ge 1$ such that, for all $t\ge t_\epsilon$,
\[
 \mathcal E'+c k(t)^2\mathcal E\le\epsilon k(t)^2,
\]
and hence
\[
 \mathcal E(t)\le\mathcal E(t_\epsilon)\exp\!\left(-c\int_{t_\epsilon}^t k(s)^2\,ds\right)+\frac{\epsilon}{c}.
\]
Since $\int_1^\infty k(s)^2\,ds=\infty$, letting $t\to\infty$ and then $\epsilon\downarrow 0$ yields $\mathcal E(t)\to 0$, which together with \eqref{extplane:compensated} gives \eqref{wh3:decay-ltwo-limit}.
\end{proof}

We next extend the $L^2$ decay to $L^p$ and $W^{1,p}$ norms.

\begin{lemma}\label{reg:plane-decay-norms}
Let $(\rho,u,\theta)$ be as in Lemma~\ref{reg:plane-decay-ltwo}. Then, for every $2\le p<\infty$,
\[
 \lim_{t\to\infty}\bigl(\|\rho(t)-\rho_\infty\|_p+\|u(t)\|_{W^{1,p}}+
                         \|\theta(t)-\theta_\infty\|_{W^{1,p}}\bigr)=0 .
\]
\end{lemma}
\begin{proof}
Multiplying $\eqref{fns}_2$ by $u$ and integrating by parts over $\mathbb{R}^2$ gives
\[
 \mu\|\nabla u\|_2^2+\int_{\mathbb R^2}(\mu+\lambda)(\div u)^2\,dx
 =\int_{\mathbb R^2}(P-P_\infty)\div u\,dx-\int_{\mathbb R^2}\rho\dot u\cdot u\,dx.
\]
By Lemma~\ref{reg:plane-decay-ltwo} and Lemmas~\ref{wh3:permanentlow} and~\ref{wh3:temperaturehigh}, we have
\[
\begin{aligned}
 \|P-P_\infty\|_2&\le\|\theta\|_\infty\|\rho-\rho_\infty\|_2+\rho_\infty\|\theta-\theta_\infty\|_2\longrightarrow 0,\\
 \mu\|\nabla u\|_2^2&\le\|P-P_\infty\|_2\|\div u\|_2+\|\sqrt\rho\,\dot u\|_2\|\sqrt\rho\,u\|_2\longrightarrow 0.
\end{aligned}
\]
By \eqref{gn11}, Lemma~\ref{wh3:temperaturehigh} and H\"older's inequality, we obtain, for $2\le p<\infty$,
\[
\begin{aligned}
 \|\nabla u\|_p&\le\|\nabla u\|_2^{1/(p-1)}\|\nabla u\|_{2p}^{(p-2)/(p-1)}\longrightarrow 0,\\
 \|\nabla\theta\|_p&\le C_p\|\theta-\theta_\infty\|_2^{1/p}\|\nabla^2\theta\|_2^{1-1/p}\longrightarrow 0,
\end{aligned}
\]
and for $f\in\{\rho-\rho_\infty,u,\theta-\theta_\infty\}$,
\[
 \|f\|_p\le\|f\|_2^{2/p}\|f\|_\infty^{1-2/p}\longrightarrow 0.
\]
This completes the proof.
\end{proof}

Together with the density and temperature estimates of Lemmas~\ref{wh3:permanentlow},
\ref{reg:plane-barrier}, and~\ref{wh3:temperaturehigh}, this proves the large-time
assertions of Theorem~\ref{tri:whole}.

\section{Proofs of Theorems \ref{th1}--\ref{tri:whole}}\label{init:construction}

With the a priori estimates established in Sections 3, 4 and 5 at hand, we are now in a position to prove the main results.

Throughout this section, set
\be\nonumber\ba
 (\rho_b,\theta_b)\triangleq
 \begin{cases}
  (0,0), & \Omega=\T^2\text{ or bounded slip domain},\\
  (\rho_\infty,\theta_\infty), & \Omega=\rr,
 \end{cases}\quad
 \rho_{\rm ref}\triangleq
 \begin{cases}
  \ol\rho_0, & |\Omega|<\infty,\\
  \rho_\infty, & \Omega=\rr.
 \end{cases}
\ea\ee

\begin{lemma}\label{approx:data}
Under the data assumptions of Theorems~\ref{th1}--\ref{tri:whole}, for $0<\varepsilon<1/2$, there exist smooth approximations $(r_\varepsilon,u_{0,\varepsilon},\theta_{0,\varepsilon})$ preserving the corresponding boundary or far-field conditions such that $r_\varepsilon\ge\varepsilon\rho_{\rm ref}>0$, $\underline\theta_0\le\theta_{0,\varepsilon}\le\|\theta_0\|_{L^\infty}$,
\be\nonumber\ba
 (r_\varepsilon-\rho_0,\,u_{0,\varepsilon}-u_0,\,\theta_{0,\varepsilon}-\theta_0)
 \longrightarrow(0,0,0)\quad\text{in } W^{1,q}\times H^2\times H^1 ,
\ea\ee
and $\sup_{\varepsilon}\|g_{1,\varepsilon}\|_{L^2}\le C$, where $g_{1,\varepsilon}$ denotes the compatibility function in \eqref{th13}. Moreover, $\int_\Omega r_\varepsilon\,dx=M_0$ when $|\Omega|<\infty$, and $r_\varepsilon-\rho_0\to 0$ in $H^1$ when $\Omega=\rr$.
\end{lemma}

\begin{proof}
Set $r_\varepsilon\triangleq(1-\varepsilon)\rho_0+\varepsilon\rho_{\rm ref}$.
Then $r_\varepsilon\ge\varepsilon\rho_{\rm ref}>0$, and $r_\varepsilon-\rho_0\to 0$ in $W^{1,q}$ (and in $H^1$ on $\rr$). Smoothing $\theta_0-\theta_b$ by the periodic, Neumann, or planar heat semigroup at a sufficiently small time $s_\varepsilon\in(0,\varepsilon]$ gives $\theta_{0,\varepsilon}$ satisfying
\be\nonumber\ba
 \|\theta_{0,\varepsilon}-\theta_0\|_{H^1}\le\varepsilon,\quad
 \underline\theta_0\le\theta_{0,\varepsilon}\le\|\theta_0\|_{L^\infty}.
\ea\ee
Write $\mathcal L_\rho v\triangleq -\mu\Delta v-\nabla\bigl((\mu+\lambda(\rho))\div v\bigr)$, and set
\be\nonumber\ba
 d_\varepsilon\triangleq\bigl(\lambda(r_\varepsilon)-\lambda(\rho_0)\bigr)\div u_0
 -(r_\varepsilon-\rho_0)\theta_0,\quad
 u_{0,\varepsilon}\triangleq u_0+w_\varepsilon ,
\ea\ee
where $w_\varepsilon$ solves $(I+\mathcal L_{r_\varepsilon})w_\varepsilon=\nabla d_\varepsilon$ in weak form,
\be\nonumber\ba
 \int_\Omega w_\varepsilon\cdot v\,dx
 +\int_\Omega\bigl(h(r_\varepsilon)\div w_\varepsilon\,\div v
 +\mu\omega(w_\varepsilon)\omega(v)\bigr)\,dx
 +\mu\int_{\partial\Omega}A w_{\varepsilon,\tau}v_\tau\,dS
 =-\int_\Omega d_\varepsilon\div v\,dx ,
\ea\ee
for admissible $v\in H^1$ (with the boundary integral present only on the slip domain). Coercivity and taking the curl yield
\be\nonumber\ba
 \|w_\varepsilon\|_{H^1}\le C\|d_\varepsilon\|_{L^2}\le C\varepsilon,\quad
 (1-\mu\Delta)\omega(w_\varepsilon)=0 .
\ea\ee
Hence $\omega(w_\varepsilon)=0$ on $\T^2$ and $\rr$, whereas on the slip domain the Dirichlet condition $\omega(w_\varepsilon)=-A w_{\varepsilon,\tau}$ on $\partial\Omega$ gives $\|\omega(w_\varepsilon)\|_{H^1}\le C\|w_\varepsilon\|_{H^1}\le C\varepsilon$. Since $\|\nabla d_\varepsilon\|_{L^2}\to 0$, the resolvent equation implies
\be\nonumber\ba
 \|\nabla(h(r_\varepsilon)\div w_\varepsilon)\|_{L^2}
 \le\|w_\varepsilon\|_{L^2}+\|\nabla d_\varepsilon\|_{L^2}
 +\mu\|\nabla\omega(w_\varepsilon)\|_{L^2}\longrightarrow 0 .
\ea\ee
Set $a_\varepsilon=1/h(r_\varepsilon)$ and
$Z_\varepsilon=h(r_\varepsilon)\div w_\varepsilon$.
Since $r_\varepsilon$ is uniformly bounded in $L^\infty$,
$h\ge 2\mu>0$, and
$\nabla r_\varepsilon=(1-\varepsilon)\nabla\rho_0$, the chain rule gives
\[
 \|a_\varepsilon\|_{L^\infty}\le\frac{1}{2\mu},
 \qquad
 \|\nabla a_\varepsilon\|_{L^q}
 =\left\|\frac{h'(r_\varepsilon)}{h(r_\varepsilon)^2}
                 \nabla r_\varepsilon\right\|_{L^q}
 \le C\|\nabla\rho_0\|_{L^q}.
\]
Moreover, $\|Z_\varepsilon\|_{L^2}\le C\|w_\varepsilon\|_{H^1}
\le C\varepsilon$, and the preceding estimate gives
$\|\nabla Z_\varepsilon\|_{L^2}\to0$. Thus $Z_\varepsilon\to0$ in $H^1$.
Taking $p=2q/(q-2)$ and using $H^1(\Omega)\hookrightarrow L^p(\Omega)$,
we obtain
\[
 \begin{aligned}
 \|\div w_\varepsilon\|_{H^1}
 &=\|a_\varepsilon Z_\varepsilon\|_{H^1}\\
 &\le \|a_\varepsilon\|_{L^\infty}\|Z_\varepsilon\|_{H^1}
      +\|\nabla a_\varepsilon\|_{L^q}\|Z_\varepsilon\|_{L^p}
 \le C\|Z_\varepsilon\|_{H^1}\longrightarrow0.
 \end{aligned}
\]
The div--curl estimate, with $w_\varepsilon\cdot n=0$ on the slip
boundary, now yields
\[
 \|w_\varepsilon\|_{H^2}
 \le C\bigl(\|w_\varepsilon\|_{L^2}
           +\|\div w_\varepsilon\|_{H^1}
           +\|\omega(w_\varepsilon)\|_{H^1}\bigr)
 \longrightarrow 0.
\]
Using $(I+\mathcal L_{r_\varepsilon})w_\varepsilon=\nabla d_\varepsilon$ and \eqref{th13}, the approximate compatibility function reads
\be\nonumber\ba
 g_{1,\varepsilon}
 =\frac{\mathcal L_{r_\varepsilon}u_{0,\varepsilon}+\nabla(r_\varepsilon\theta_{0,\varepsilon})}{\sqrt{r_\varepsilon}}
 =\sqrt{\frac{\rho_0}{r_\varepsilon}}\,g_1
 -\frac{w_\varepsilon}{\sqrt{r_\varepsilon}}
 +\frac{\nabla\bigl(r_\varepsilon(\theta_{0,\varepsilon}-\theta_0)\bigr)}{\sqrt{r_\varepsilon}}.
\ea\ee
Since $\rho_0/r_\varepsilon\le 2$, $r_\varepsilon\ge\varepsilon\rho_{\rm ref}$, $\|w_\varepsilon\|_{L^2}\le C\varepsilon$, and $H^1\hookrightarrow L^{2q/(q-2)}$, we obtain
\be\nonumber\ba
 \|g_{1,\varepsilon}\|_{L^2}
 &\le \sqrt 2\,\|g_1\|_{L^2}+C\sqrt\varepsilon
 +C\|\nabla(\theta_{0,\varepsilon}-\theta_0)\|_{L^2}
 +C\varepsilon^{-1/2}\|\nabla\rho_0\|_{L^q}\|\theta_{0,\varepsilon}-\theta_0\|_{L^{\frac{2q}{q-2}}}\\
 &\le C.
\ea\ee
Finally, at each fixed $\varepsilon>0$, further smoothing $r_\varepsilon$ by the heat semigroup and $u_{0,\varepsilon}$ by the constant-coefficient Lam\'e semigroup (or convolution on $\rr$) yields smooth data preserving the boundary or far-field conditions while perturbing $\|g_{1,\varepsilon}\|_{L^2}$ by at most $\varepsilon$.
\end{proof}

\begin{lemma}\label{approx:limit}
The positive-data solutions $(\rho_\varepsilon,u_\varepsilon,\theta_\varepsilon)$ corresponding to Lemma~\ref{approx:data} exist globally on $\Omega\times(0,\infty)$, satisfy the $\varepsilon$-uniform a priori bounds of Sections~\ref{per-apriori}--\ref{wh3:section} together with
\be\nonumber\ba
 \inf_{\Omega\times[0,T]}\rho_\varepsilon\ge\varepsilon\rho_{\rm ref}e^{-C_T},\quad
 \inf_{\Omega\times[0,T]}\theta_\varepsilon\ge c_T>0 ,
\ea\ee
and converge along a subsequence (locally in space when $\Omega=\rr$) to a global strong solution $(\rho,u,\theta)$ satisfying \eqref{th15} (or \eqref{th35} on $\rr$), \eqref{main:finite-bounds}, and
\be\la{init:heat-trace}\ba
 \int_\Omega\rho(t)|\theta(t)-\theta_0|^2\,dx\le C(t+t^2)\quad (0<t\le 1),\quad
 \|\rho(t)\theta(t)-\rho_0\theta_0\|_{L^2}\longrightarrow 0\quad (t\downarrow 0).
\ea\ee
\end{lemma}

\begin{proof}
For \eqref{ext:gainfamily}, the viscosity law and its exponent remain fixed throughout the approximation. For sufficiently small approximation parameters, Lemma~\ref{approx:data} gives \eqref{ext:dataM} with $4M+4$ in place of $M$: the data converge in the stated spaces, the temperature bounds are preserved, and $\|g_{1,\varepsilon}\|_2\le\sqrt2\|g_1\|_2+o(1)$. Since the levels $R_{*,C}$ in Proposition~\ref{order:periodic-cap} and Lemmas~\ref{order:wall-cap} and~\ref{wh3:permanentlow} were chosen with this enlarged bound, the same thresholds $\varepsilon_{*,C}(M)$ and density ceilings apply to these approximations.
For each fixed $\varepsilon>0$, Lemma~\ref{lct} gives a local strong solution. Integrating $(\ref{fns})_1$ along characteristics and applying the $\varepsilon$-uniform $L^1(0,T;L^\infty)$ bound on $\div u_\varepsilon$ from Sections~\ref{per-apriori}--\ref{wh3:section} yield
\be\nonumber\ba
 \rho_\varepsilon\bigl(t,X_\varepsilon(t,x)\bigr)
 =r_\varepsilon(x)\exp\left(-\int_0^t\div u_\varepsilon\bigl(s,X_\varepsilon(s,x)\bigr)\,ds\right)
 \ge\varepsilon\rho_{\rm ref}e^{-C_T},
\ea\ee
while Proposition~\ref{order:periodic-cap} and Lemmas~\ref{order:wall-cap} and~\ref{wh3:permanentlow} give $\theta_\varepsilon\ge c_T>0$. Combined with the $\varepsilon$-uniform initial- and positive-time estimates (Lemmas~\ref{fpl9}, \ref{init:common-lemma}, \ref{reg:wall-initial}--\ref{reg:wall-higher}, and~\ref{reg:plane-initial}--\ref{reg:plane-higher}), these positive floors allow Lemma~\ref{lct} to extend $(\rho_\varepsilon,u_\varepsilon,\theta_\varepsilon)$ to $[0,\infty)$.

By compact spatial embeddings, equicontinuity, and diagonal extraction (on balls by exhaustion when $\Omega=\rr$), there is a subsequence such that, for every $0<\tau<T$,
\be\nonumber\ba
 \rho_\varepsilon\to\rho\text{ in }C([0,T];C_{\rm loc}),\quad
 u_\varepsilon\to u\text{ in }C([0,T];H^1_{\rm loc}),\quad
 \theta_\varepsilon\to\theta\text{ in }C([\tau,T];H^1_{\rm loc})
\ea\ee
(globally in space when $|\Omega|<\infty$). In particular, $\rho_\varepsilon\theta_\varepsilon\to\rho\theta$ in $L^2_{\rm loc}$, $Q_\varepsilon\to Q$ in $L^1_{\rm loc}$, and $K(\theta_\varepsilon)\to K(\theta)$ in $L^2_{\rm loc}$ on $\Omega\times[\tau,T]$, which justifies passing to the limit in \eqref{fns}. The remaining bounds pass by weak lower semicontinuity: the transport formula for $\nabla\rho$ yields $\rho-\rho_b\in C([0,T];W^{1,q})$ (and in $C([0,T];L^2)$ on $\rr$ via \eqref{wh3:density-gradient-ltwo}), $u_t\in L^2(0,T;H^1)$ yields $u\in C([0,T];H^1)$, and $\theta-\theta_b\in L^2(\tau,T;H^3)$ with $\theta_t\in L^2(\tau,T;H^1)$ yields $\theta-\theta_b\in C([\tau,T];H^2)$.

To identify the initial thermal trace, applying $(\ref{fns})_1$ to $\rho_\varepsilon|\theta_\varepsilon-\theta_{0,\varepsilon}|^2$ and using $\kappa(\theta_\varepsilon)\ge\kappa_0 c_1^\alpha>0$ on $[0,1]$ give, for $0<t\le 1$,
\be\nonumber\ba
 \|\sqrt{\rho_\varepsilon(t)}\bigl(\theta_\varepsilon(t)-\theta_{0,\varepsilon}\bigr)\|_{L^2}
 \le\int_0^t\bigl(\|\sqrt{\rho_\varepsilon}\,\dot\theta_\varepsilon\|_{L^2}
 +\|\sqrt{\rho_\varepsilon}\,u_\varepsilon\cdot\nabla\theta_{0,\varepsilon}\|_{L^2}\bigr)\,ds
 \le C(\sqrt t+t).
\ea\ee
Since $\|\sqrt{\rho_\varepsilon(t)}(\theta_{0,\varepsilon}-\theta_0)\|_{L^2}\le C\varepsilon$, letting $\varepsilon\to 0$ (by exhaustion on $\rr$) yields the first inequality in \eqref{init:heat-trace}, and hence
\be\nonumber\ba
 \|\rho(t)\theta(t)-\rho_0\theta_0\|_{L^2}
 \le \|\rho(t)\|_{L^\infty}^{1/2}\|\sqrt{\rho(t)}(\theta(t)-\theta_0)\|_{L^2}
 +\|\theta_0\|_{L^\infty}\|\rho(t)-\rho_0\|_{L^2}\longrightarrow 0\quad (t\downarrow 0),
\ea\ee
which establishes \eqref{cztj} and $\rho\theta-\rho_b\theta_b\in C([0,T];L^2)$.
\end{proof}

\begin{proof}[Proof of Theorems~\ref{th1}--\ref{tri:whole}]
Lemmas~\ref{approx:data} and~\ref{approx:limit} provide the global strong solutions, regularity \eqref{th15} and \eqref{th35}, and finite-time bounds \eqref{main:finite-bounds} in all three domains. Under \textup{(A)}, \textup{(B)}, or \textup{(C)} with $\alpha>1-d$, the time-uniform density ceiling and upper/lower temperature bounds (combining the $[0,1]$ barrier with the positive-time $H^2$ bound at $\tau=1$ on $[1,\infty)$) pass to the limit and yield \eqref{ltp-main-uniform}, \eqref{tri:bd-uniform}, and \eqref{tri:whole-uniform}. The common density ceilings in the data-class alternatives \eqref{ext:powerfamily} and \eqref{ext:gainfamily} pass to the limit as well.

Finally, on $\T^2$, the approximate equilibria converge to \eqref{ltp-equilibrium} and the constants in Lemma~\ref{reg:tor-decay-norms} are uniform in $\varepsilon$, so passing to the limit at each $t\ge 1$ gives \eqref{ltp-main-decay}. On the slip domain $\Omega$ and on $\rr$, the time-uniform bounds on $[1,\infty)$ justify applying Lemma~\ref{reg:wall-norms} and Lemmas~\ref{reg:plane-decay-ltwo}--\ref{reg:plane-decay-norms} directly to the limiting strong solution (on $\Omega$, the energy balance from $t=0$ follows from the initial traces and $u\in C([0,T];H^1)$, while on $\rr$ the low-frequency limit at $t=1$ applies to the limiting $L^2$ fields), yielding \eqref{tri:bd-decay}--\eqref{tri:bd-equilibrium} and \eqref{tri:whole-decay}.
Moreover, the uniqueness of $(\rho,u,\theta)$ follows from arguments similar to those in the proof of \cite[Theorem 1]{ChoKim2006}; we omit the details.
This completes the proofs of Theorems~\ref{th1}--\ref{tri:whole}.
\end{proof}

\bigskip

\noindent\textbf{Data availability.} No data was used for the research described in the article.

\bigskip

\noindent\textbf{Conflict of interest.} The authors declare that they have no conflict of interest.

\section*{Acknowledgments}
The authors acknowledge the use of generative AI tools in developing and checking mathematical arguments and in drafting and revising this manuscript.

\end{document}